\documentclass[11pt]{article}

\usepackage[a4paper,margin=0.8in]{geometry}
\usepackage{authblk}
\usepackage{hrefhide}
\usepackage{mathrsfs}
\usepackage{dsfont}
\usepackage{amsmath,amssymb,amsthm,mathtools}

\usepackage{tikz}
\newcommand{\MMN}{\mathcal M_{\NN_0}}

\newcommand{\ZZJ}{{Z_J^c\times[0,T]}}
	\DeclareMathOperator{\dist}{dist}
	\DeclareMathOperator{\supp}{supp}
	\newcommand{\Lt}{\widetilde{L}}

\newcommand{\seque}{{h_r}}

	\renewcommand{\restriction}{|}
\newcommand{\rest}{|}

	\newcommand{\jt}{{j}}
	\newcommand{\diff}{ \Delta }%
\newcommand{\jty}{J}
								\newcommand{\myp}{{u}}

								\newcommand{\N}{{\mathbb{N}}}
								\newcommand{\up}{{\bf p}}
								\newcommand{\uq}{{\bf q}}

								\newcommand{\Leb}{\mbox{Leb}}
								\newcommand{\zaehler}{\kappa}
								\newcommand{\shiftpro}{\pro^{\mbox{s}}}
								\newcommand{\INT}{[0,T]}
								\newcommand{\MA}{\mathfrak{A}}
								\newcommand{\lk}{\left}
								\newcommand{\lqq}{\lefteqn}
								\newcommand{\rk}{\right}
								\newcommand{\la}{{\langle}}
								\newcommand{\ra}{{\rangle}}
								\newcommand{\ep} {\epsilon }

								\usepackage{amsfonts}
								\usepackage[english]{babel}
								\usepackage{latexsym}
								\usepackage{enumitem}

								\newcommand{\CM}{{{ \mathcal M }}}
								\newcommand{\bNN}{{\bar{\mathbb{N}}}}
								\newcommand{\PP}{\mathbb{P}}
								\newcommand{\MscrZ}{{\mathscr{Z}}}

								\newcommand{\BF}{{{ \mathbb{F} }}}
								\newcommand{\CF}{{{ \mathcal F }}}

								\newcommand{\Pro}{\mbox{Proj}}

								\usepackage{xcolor}

								\newcommand{\bcase}{\begin{cases}}
									\newcommand{\ecase}{\end{cases}}

								\newcommand{\eps}{\epsilon}

								\newcommand{\RR}{{\mathbb{R}}}
								\newcommand{\DD}{\mathbb{D}}
								\newcommand{\NN}{\mathbb{N}}
								\newcommand{\EE}{\mathbb{E}}

								\newcommand{\DEQS}{\begin{eqnarray*}}
									\newcommand{\EEQS}{\end{eqnarray*}}
								\newcommand{\DEQSZ}{\begin{eqnarray}}
									\newcommand{\EEQSZ}{\end{eqnarray}}

									\newcounter{lil7711a}
	\newenvironment{steps_a1}
	{\begin{list} {  Step (A.\arabic{lil7711a})} {  \usecounter{lil7711a}
				\setlength{\leftmargin}{0.0cm}
				\setlength{\topsep}{0.2cm} \setlength{\itemsep}{0.0cm}
				\setlength{\parsep}{0.1cm} \setlength{\itemindent}{0.8cm}
				\setlength{\parskip}{0.0cm}}} {\end{list}}

								\newcounter{lil77}
								\newenvironment{steps}
								{\begin{list} { \bf Step (\Roman{lil77})} {  \usecounter{lil77}
											\setlength{\leftmargin}{0.0cm}
											\setlength{\topsep}{0.2cm} \setlength{\itemsep}{0.0cm}
											\setlength{\parsep}{0.1cm} \setlength{\itemindent}{0.8cm}
											\setlength{\parskip}{0.0cm}}} {\end{list}}

								\newcounter{lil22}
								\newenvironment{steps-2}
								{\em \begin{list} {  (\roman{lil22})} {\usecounter{lil22} \em
											\setlength{\leftmargin}{1.3cm}
											\setlength{\topsep}{0.2cm} \setlength{\itemsep}{0.0cm}
											\setlength{\parsep}{0.1cm} \setlength{\itemindent}{0.4cm}
											\setlength{\parskip}{0.0cm}}} {\end{list}}

								\newcommand{\n}{\,\hbox{\vrule width 1pt height 8ptdepth 2pt}\,}

		\newcommand{\Law}{{{ \mathcal L }aw}}
\newcommand{\NNtilde}{{\tilde \NN}}
		\newcommand{\cB}{\mathcal{B}}
		\newcommand{\cF}{\mathcal{F}}
		\newcommand{\xhat}{\hat{x}}
		\newcommand{\zhat}{\hat{z}}
		\newcommand{\lbar}{\bar{\ell}}
		\newcommand{\ind}{\mathbf{1}}
		\newcommand{\mr}{{m_r}}
		\newcommand{\Jr}{{r,J}}
				\newcommand{\nr}{{n_r}}

		\numberwithin{equation}{section}
		\theoremstyle{plain}
		\newtheorem{notation}{Notation}[section]

		\newtheorem{theorem}{Theorem}[section]
		\newtheorem{lemma}{Lemma}[theorem]%
		\newtheorem{corollary}[theorem]{Corollary}%
		\newtheorem{proposition}[theorem]{Proposition}%

		\theoremstyle{definition}
		\newtheorem{hypo}[theorem]{Hypothesis}%
		\newtheorem{definition}[theorem]{Definition}%
		\newtheorem{remark}[theorem]{Remark}%

\title{Reconstruction of Poisson Random Measures in the Variational Framework with Applications to Jump-Driven SPDEs}

\author[1]{Erika Hausenblas}
\author[2]{Akash Ashirbad Panda}
\author[1]{Ankit Kumar\thanks{Corresponding author.}}

\affil[1]{Department Mathematik und Informationstechnologie, Montanuniversit\"at Leoben, Leoben, Austria, (erika.hausenblas@unileoben.ac.at, ankit.kumar@unileoben.ac.at)}

\affil[2]{Department of Mathematics,
Indian Institute of Technology Indore, Madhya Pradesh, India, (akashpanda@iitbbs.ac.in)}

\date{}

\begin{document}
%

\maketitle

\vspace{-0.8em}

\begin{abstract}
We establish a reconstruction principle for Poisson random measures arising in compactness arguments for jump-driven stochastic partial differential equations.
A basic difficulty in such arguments is that, after a Skorokhod-type representation, almost sure convergence of the solution processes does not by itself identify the limiting stochastic integral, since integration with respect to a compensated Poisson random measure is not a pathwise continuous functional of the solution and the noise.

We address this problem in the variational framework of a Gelfand triple $V\hookrightarrow H\equiv H^\ast\hookrightarrow V^\ast$ for equations with locally monotone drift. Starting from a tight sequence of approximating solutions in $D([0,T];H)\cap L^\alpha(0,T;V),$ we construct, on the represented probability space, a time-homogeneous Poisson random measure with the original L\'evy intensity and compensator, compatible
with the limiting filtration, such that the limiting process solves the corresponding SPDE.

The main ingredient is a noise bridge. Together with a shifted dyadic
approximation of the solution, the increments of the solution path over dyadic
time cells approximate the compensated noise increments; since the time
integral of the drift is continuous in $H$, no regularisation of the drift is
needed. The jumps are then recovered by a measurable registration
procedure, while injectivity of the jump coefficient permits reconstruction of
the corresponding marks through a Borel inverse obtained from the
Lusin--Suslin theorem. The reconstruction is performed simultaneously at all
truncation levels; consistency of the resulting finite jump measures and a
monotone patching argument yield a Poisson random measure with the full
L\'evy measure. Under the injectivity assumption, no enlargement of the
represented probability space is required.
\end{abstract}

\smallskip


\smallskip


\smallskip

\section{Introduction}\label{sec:intro}

Compactness methods are a standard tool for proving the existence of solutions to
nonlinear stochastic partial differential equations. Typically, one constructs a
sequence of approximate solutions, derives uniform a priori estimates, and uses
them to establish tightness of the corresponding laws on a suitable path space.
A Skorokhod-type representation theorem then places a subsequence on a new
probability space where the convergence is almost sure, which is particularly
useful for passing to the limit in nonlinear terms. For Wiener-driven equations
this strategy is by now well established; see, for instance,
\cite{DaPratoZabzcy}.

For equations driven by jump noise, however, an additional difficulty appears.
The limiting process has to be accompanied by a Poisson random measure having
the prescribed intensity and the correct relation to the limiting filtration.
Indeed, a Poisson random measure contains both the jump times and the
corresponding marks, and these data are needed in order to identify the
stochastic integral and to conclude that the limiting process solves the
original equation.

To explain the problem, let $Z$ denote the mark space and suppose that a
sequence of pairs $(u_n,\eta_n)$ is transported, by a Skorokhod-type
representation, to a new probability space on which the corresponding copies
converge almost surely. Such convergence does not by itself justify the passage
to the limit in
\[
 (u,\eta)\longmapsto
 \int_0^\cdot\!\!\int_Z g(u(s),z)\,\widetilde{\eta}(dz,ds).
\]
The stochastic integral is constructed through predictable approximations on
the underlying filtered probability space, and the convergence supplied by the
representation theorem does not provide a continuity principle for this map in
the topologies used here. Consequently, transporting the solution and the
random measure jointly is not sufficient, by itself, to identify the limiting
stochastic integral with the integral driven by the transported limiting
measure. This is particularly relevant in the present variational setting,
where
\[
 u\in D([0,T];H)\cap L^\alpha(0,T;V),
\]
so that the path is c\`adl\`ag only in the weaker space $H$, while the drift
naturally takes values in $V^\ast$.

Our approach avoids this continuity problem by reconstructing the driving
Poisson random measure from the approximating solution paths themselves. We
consider a Gelfand triple
\[
 V\hookrightarrow H\equiv H^\ast\hookrightarrow V^\ast
\]
and approximating equations of the form
\[
 du_n(t)
 =
 A_n(t,u_n(t))\,dt
 +
 \int_Z g(u_n(t),z)\,\widetilde{\eta}(dz,dt),
 \qquad
 u_n(0)=u_0^n,
\]
where the approximations are defined on one filtered probability space and are
driven by the same Poisson random measure. The essential idea is to recover the
martingale increments from quantities depending only on the solution path.

At a dyadic level $m$, the
solution is replaced by its shifted Haar approximation $u_{n,m}:=\mathrm{Proj}_m u_n.$ The shift by one dyadic interval makes this approximation predictable. On a
dyadic interval $I_k^m$, the increment $\ell_{k,n,m}$ of the solution over $I_k^m$, reduced by the drift increment $\int_{I_k^m}A_n(s,u_n(s))\,ds$, equals the compensated-noise increment
\[
 \int_{I_k^m}\!\!\int_Z
 g(u_n(s),z)\,\widetilde{\eta}(dz,ds).
\]
This relation is the \emph{noise bridge} underlying the reconstruction. 
The approximation and dyadic
parameters have to be chosen in the appropriate order. Importantly, the proof does not
require a convergence rate that is uniform in the approximation index. The
relevant good-event probabilities converge to one, and this convergence in
probability is sufficient for the subsequent representation argument.

The next step is to recover the individual jumps and their marks. We decompose
the mark space into layers
\[
 \Delta Z_j
 =
 \bigl\{z\in Z:\varepsilon_{j+1}<|z|_Z\leq\varepsilon_j\bigr\},
 \qquad \varepsilon_j\downarrow0,
\]
and associate with each layer the admissible jump set $g(v,\Delta Z_j)\subset H.$ A reconstructed increment is registered whenever it lies within a prescribed
tolerance of one of these sets. A measurable nearest-layer rule resolves
possible ambiguities. Under the injectivity and separation assumptions imposed
on $g$, the corresponding mark can then be recovered through the Borel inverse
of $(v,z)\longmapsto (v,g(v,z)),$ provided by the Lusin--Suslin theorem. Hence the reconstruction is measurable
and, under the injectivity assumption, requires no enlargement of the
probability space.

For every fixed truncation level $J$, the reconstruction is exact on events
whose probabilities tend to one. More precisely, the atoms of the reconstructed
measure can be matched with those of the corresponding truncated path-jump
measure, with the matching error tending to zero. This yields convergence, in
the Prokhorov metric, of the reconstructed finite point measures and therefore
convergence of their laws to the law of a Poisson random measure with the
prescribed truncated intensity. In particular, tightness of the reconstructed
measures is obtained from convergence of their laws rather than merely from a
bound on the number of atoms, which would be insufficient when the mark space
is infinite-dimensional.

\paragraph{Main result.}
We consider the limiting equation
\[
 du(t)
 =
 A(t,u(t))\,dt
 +
 \int_Z g(u(t),z)\,\widetilde{\eta}(dz,dt),
 \qquad
 u(0)=u_0,
\]
on the Gelfand triple $V\hookrightarrow H\equiv H^\ast\hookrightarrow V^\ast,$ where $V$ is reflexive and separable and $H$ is a separable Hilbert space. The
drift $A$ satisfies the hemicontinuity, local monotonicity, coercivity and
growth assumptions of Hypothesis~\ref{hyp}, while the jump coefficient $g$
satisfies the structural assumptions of Hypotheses~\ref{ass:GF} and the
corresponding separation condition.

Starting from approximating solutions $(u_n)$ satisfying the uniform a priori
bounds and having tight laws on
\[
 D([0,T];H)\cap L^\alpha(0,T;V),
\]
Theorem~\ref{Theoremmain} constructs a filtered probability space $\mathfrak A^\ast
 =
 \bigl(\Omega^\ast,\mathcal F^\ast,
       (\mathcal F_t^\ast)_{t\in[0,T]},\PP^\ast\bigr)$ carrying a process $u_\ast$ with the limiting law and a time-homogeneous
Poisson random measure $\eta_\ast$ with intensity $\nu$. The measure
$\eta_\ast$ has compensator $\nu(dz)\,dt$, is compatible with the limiting
filtration, and the process $u_\ast$ solves the limiting SPDE driven by
$\widetilde{\eta}_\ast$. Thus the theorem provides the jump-noise counterpart
of the usual compactness step for stochastic evolution equations in the
variational framework.

Several features of the construction are important. First, the Lévy measure may
have infinite mass near the distinguished point $0^\ast$. Hence the
reconstruction is first performed after truncating the small jumps. Treating
each truncation level by a separate representation theorem would lead to
different probability spaces or different subsequences and would destroy the
consistency needed to recover one global random measure. We therefore represent
the solution jointly with \emph{all} truncation levels on a single countable
product space; see Proposition~\ref{prop:JS-pairs-allJ}. Tightness of the
product follows from tightness of its coordinates, and the limiting truncated
measures are shown to be consistent on their common layers; see
Lemma~\ref{lem:J-consistency}. They can consequently be patched pathwise, and
the monotone construction of Proposition~\ref{prop:full-levy-patch} yields a
Poisson random measure with the full Lévy intensity.

Second, identifying the law of the reconstructed measure is not sufficient:
the resulting noise must also have the correct nonanticipative relation to the
limiting solution. We formulate this through compatibility of the solution and
the random measure. The relevant compatibility identities hold for the
approximating pair, are transferred to the represented variables by equality of
laws, and pass to the limit. This yields independence of future Poisson
increments from the limiting past and hence the required compatibility with
the filtration.

A further point concerns the value of the solution immediately before a jump.
The recovery of a mark uses $g\bigl(u(s-),z\bigr),$ whereas a c\`adl\`ag path is a priori only $H$-valued at a fixed time. The
required $V$-regularity at the random jump times is supplied by the
probabilistic structure: the compensation formula shows that, almost surely,
the Poisson random measure has no atom at a time at which the relevant pre-jump
state fails to belong to $V$. Such a conclusion would not follow from the
pathwise regularity alone.

The assumptions of Theorem~\ref{Theoremmain} are chosen to make the
reconstruction mechanism transparent rather than to be optimal. In particular,
the present argument uses injectivity of the mark map and bounded support of
the Lévy measure, while several possible relaxations and extensions are
discussed in Subsection~\ref{sec:open-problems}.

\paragraph{Relation to existing work.}
The pathwise construction of the jump measure of a finite-dimensional
c\`adl\`ag process is classical; see
\cite{IkedaWatanabe1981,JacodShiryaev2003,conttankov}. The correspondence
between martingale problems and weak solutions of stochastic equations with
jumps is also well understood in finite dimensions; see \cite{kurtz2011}.
Representations of integer-valued random measures with prescribed compensators
as images of Poisson random measures on enlarged probability spaces go back to
Grigelionis~\cite{Grigelionis1971} and
El Karoui--Lepeltier~\cite{ElKarouiLepeltier1977}.

Kallianpur and Xiong~\cite{KallianpurXiong1995} proved a representation theorem
for purely discontinuous martingales with values in the dual of a nuclear
space. Brze\'zniak, Kosmala, Motyl and
Razafimandimby~\cite{BKMR} subsequently obtained a corresponding result for
martingales taking values in a separable Banach space. In their application to
the stochastic Navier--Stokes equations, the construction is performed in an
auxiliary space $U^\ast$, with a compact embedding $H\hookrightarrow U^\ast,$ and the noise coefficient is controlled directly in the $H$-norm.

The present work differs in that the reconstruction is carried out directly in
the variational framework
\[
 V\hookrightarrow H\hookrightarrow V^\ast,
\]
where the path is c\`adl\`ag only in $H$ and the drift takes values in
$V^\ast$. We require the jump coefficient to act from $V$ into $H$ and assume
bounds of the form
\[
 \int_Z |g(v,z)|_H^p\,\nu(dz)
 \leq C_g\bigl(1+|v|_V^p\bigr),
 \qquad
 \int_Z |g(v_1,z)-g(v_2,z)|_H^p\,\nu(dz)
 \leq L_g |v_1-v_2|_V^p.
\]
Thus $g$ need not extend to a Lipschitz map on $H$. This permits, for example,
coefficients of transport type
\[
 g(v,z)=\bigl(c+|v|_V\bigr)\,z,\qquad c>0,
 \]
which naturally map a stronger spatial space $V$ into $H$ but need not define
bounded maps from $H$ to itself. However, the factor in front has to be bounded away from zero.

The main additional ingredient is the reconstruction mechanism described above:
the shifted dyadic approximation and the increments of the path over the dyadic cells create a bridge
from the solution path to the martingale increments; the layerwise registration
recovers the jumps; the Lusin--Suslin inverse recovers the marks; and the
simultaneous treatment of all truncation levels produces a single Poisson
random measure with the full Lévy intensity and the correct filtration
compatibility. Under the injectivity assumption, this construction is carried
out without adding auxiliary randomness to the represented probability space.

The variational theory for locally monotone stochastic equations underlying the
present setting originates in the work of Liu--R\"ockner
\cite{LiuRockner2010,LiuRockner2015} and Brze\'zniak--Liu--Zhu
\cite{BrzezniakLiuZhu2011}. The stochastic integration theory in Banach spaces
of martingale type $2$ used here is developed in
\cite{Brz+Haus_2009}. A L\'evy--Ciesielski construction related in spirit to
the present dyadic approximation, though arising in a different setting, is
given in \cite{HausenblasKumarToelle2026}.

\paragraph{Organisation of the paper.}
Section~\ref{Sec-2} introduces the notation and probabilistic framework,
including integer-valued random measures, stochastic integration with respect
to compensated Poisson random measures, and the notion of compatibility.
Section~\ref{sec:representation} introduces the variational setting, states the
assumptions on the drift, jump coefficient and approximations, and presents the
main representation theorem together with possible extensions.
The proof of Theorem~\ref{Theoremmain} is carried out in Section~4 and is
organised into Steps~(I)--(XIV): the first steps construct the pathwise
reconstruction, the middle steps perform the simultaneous representation and
identify the reconstructed random measures, aMarie4
nd the final steps establish
compatibility and pass to the limiting SPDE.

The appendices collect the auxiliary results used in the proof.
Appendix~\ref{app:levy-linde} develops the topology, compactness and tightness
theory for the relevant random measures, including the countable product of
truncated jump measures. Appendix~\ref{haarproj} studies the shifted Haar
projection used in the predictable approximation.
 Finally,
Appendix~\ref{app:Skorohod-space} records the required facts about the
Skorokhod space and the measurability of the jump data of c\`adl\`ag paths.

\section{Notation, Function Spaces, and Stochastic Preliminaries}\label{sec:notation}
\label{Sec-2}
\label{sec:preliminaries}

This section collects the notation, measure-theoretic framework, and
stochastic analysis tools used throughout the paper. We first introduce the
measurable structures associated with integer-valued measures, then define the
function spaces used for stochastic processes and stochastic integrands. We
subsequently recall time-homogeneous Poisson random measures, their compensated
versions, and stochastic integration with respect to compensated Poisson random
measures. Finally, we introduce the notion of compatibility used to formulate
nonanticipativity in a law-invariant form.

\subsection{General notation and spaces of integer-valued measures}

We begin by fixing some notation used throughout the paper.

\begin{notation}\label{notationref}
The symbol $\RR$ denotes the real numbers, $\RR^+:=\{x\in\RR:x>0\}, \ \RR_0^+:=\RR^+\cup\{0\}.$ By $\NN$ we denote the set of natural numbers including $0$, and $\bar\NN:=\NN\cup\{\infty\}$.
\end{notation}

\begin{notation}\label{notationref2}
Let $(\CF_t)_{t\in\INT}$ be a filtration on a measurable space
$(\Omega,\CF)$ and let $\theta$ be a measure on $(\Omega,\CF)$.
We denote by $\CF_t^\theta$ the augmentation of $\CF_t$ by the
$\theta$-null sets in $\CF$.
\end{notation}

We call a measurable space $(Z,\MscrZ)$ a \emph{Polish measurable space} if
there exists a complete separable metric $\varrho$ on $Z$ such that
\[
\MscrZ=\mathscr B(Z,\varrho),
\]
where $\mathscr B(Z,\varrho)$ denotes the corresponding Borel
$\sigma$-field.

\begin{notation}\label{notationref4}
The set of all finite nonnegative measures on $(Z,\MscrZ)$ is denoted by
$M_+(Z)$, while $\mathcal P(Z)$ denotes the set of probability measures on
$(Z,\MscrZ)$.

Let $\{Z_n\}_{n\in\NN}\subset\MscrZ$ satisfy $Z_n\uparrow Z.$ We denote by $M_{\bar\NN}(\{Z_n\})$ the family of all $\bar\NN$-valued
measures $\theta$ on $(Z,\MscrZ)$ such that
\[
\theta(Z_n)<\infty,
\qquad n\in\NN.
\]
The $\sigma$-field on $M_{\bar\NN}(\{Z_n\})$ generated by the evaluation
maps
\[
i_B:
M_{\bar\NN}(\{Z_n\})
\longrightarrow
\bar\NN,
\qquad
i_B(\mu):=\mu(B),
\qquad B\in\MscrZ,
\]
is denoted by $\CM_{\bar\NN}(\{Z_n\})$.
\end{notation}

The preceding space will serve as a state space for integer-valued random
measures. The following result provides the topological structure needed in
the subsequent compactness arguments. Further information is given in appendix \ref{app:levy-linde}.

\begin{lemma}\label{measure_lemma}
Let $(Z,\MscrZ)$ be a Polish measurable space and let
$\{Z_n\}_{n\in\NN}\subset\MscrZ$ satisfy $Z_n\uparrow Z$.
Then $M_{\bar\NN}(\{Z_n\})$ can be endowed with a Polish topology whose
Borel $\sigma$-field coincides with
$\CM_{\bar\NN}(\{Z_n\})$.
\end{lemma}

Throughout the paper, $\mathfrak A=(\Omega,\CF,\BF,\PP), \ \BF=(\CF_t)_{t\in[0,T]},$ denotes a filtered probability space satisfying the usual conditions; that
is, $(\Omega,\CF,\PP)$ is complete and $(\CF_t)_{t\in[0,T]}$ is complete
and right-continuous.

\subsection{Function spaces for stochastic integration}

Let $Y$ be a separable Banach space. For $q\in[1,\infty)$, we denote by
$L^q(Z,\nu;Y)$ the space of equivalence classes of measurable mappings
\[
\xi:(Z,\MscrZ)\longrightarrow(Y,\mathcal B(Y))
\]
such that
\[
|\xi|_{L^q(Z,\nu;Y)}
:=
\left(
\int_Z |\xi(z)|_Y^q\,\nu(dz)
\right)^{1/q}
<\infty.
\]

Similarly, $L^q(\Omega;Y)$ and $L^q(\Omega_T;Y)$ have their usual
meanings, where
\[
\Omega_T:=[0,T]\times\Omega
\]
is equipped with the product $\sigma$-field
$\mathscr B([0,T])\otimes\CF$.
By $L^0(\Omega;Y)$ we denote the space of measurable mappings from
$(\Omega,\CF)$ to $Y$.

Let $\mathcal N(0,T;Y)$ denote the space of equivalence classes of
progressively measurable processes
\[
\xi:[0,T]\times\Omega\longrightarrow Y.
\]
For $q\in(1,\infty)$, we set
\begin{align}
\mathcal N^q(0,T;Y)
&:=
\left\{
\xi\in\mathcal N(0,T;Y):
\int_0^T|\xi(t)|_Y^q\,dt<\infty
\quad\PP\text{-a.s.}
\right\},
\label{def:Nq}
\\
\mathcal M^q(0,T;Y)
&:=
\left\{
\xi\in\mathcal N(0,T;Y):
\EE\int_0^T|\xi(t)|_Y^q\,dt<\infty
\right\}.
\label{def:Mq}
\end{align}

Let $\mathcal N_{\rm step}(0,T;Y)$ denote the space of all
$\xi\in\mathcal N(0,T;Y)$ for which there exists a partition
\[
0=t_0<t_1<\cdots<t_n=T
\]
such that $\xi$ is constant on each interval $(t_{k-1},t_k]$ and its value
there is $\CF_{t_{k-1}}$-measurable. We set
\[
\mathcal M_{\rm step}^q(0,T;Y)
:=
\mathcal M^q(0,T;Y)
\cap
\mathcal N_{\rm step}(0,T;Y).
\]

The space $\mathcal M^q(0,T;Y)$ is a closed subspace of
\[
L^q([0,T]\times\Omega;Y)
\cong
L^q\bigl(0,T;L^q(\Omega;Y)\bigr),
\]
and $\mathcal M_{\rm step}^q(0,T;Y)$ is dense in
$\mathcal M^q(0,T;Y)$.

\subsection{Stochastic preliminaries}

We next recall the probabilistic framework for Poisson random measures and
their compensated versions. Since the reconstruction argument depends
essentially on the relation between the noise and the underlying filtration,
this relation is included explicitly in the definition.

\paragraph{Poisson random measures}

\newcommand{\NNnull}{{\NN_0}}
Let $(Z,\MscrZ)$ be a Polish measurable space, let $\nu$ be a
$\sigma$-finite measure on $(Z,\MscrZ)$, and let
$\{Z_n\}_{n\in\NN}\subset\MscrZ$ satisfy
\[
Z_n\uparrow Z,
\qquad
\nu(Z_n)<\infty,
\qquad n\in\NN.
\]

\begin{definition}\label{def-Prm}
{\rm(see \cite[Definition~I.8.1]{ikeda1981})}
A \emph{time-homogeneous Poisson random measure} on $(Z,\MscrZ)$ with
intensity measure $\nu$, over the filtered probability space
$\mathfrak A$, is a measurable mapping
\[
\eta:
(\Omega,\CF)
\longrightarrow
\Bigl(
M_{\NNnull }(\{Z_n\times(0,T]\}),
\CM_{\NNnull }(\{Z_n\times(0,T]\})
\Bigr)
\]
such that the following properties hold.

\begin{enumerate}[label=\upshape(\roman*)]

\item
For every $B\in\MscrZ\otimes\mathcal B((0,T])$ satisfying $(\nu\otimes\Leb)(B)<\infty$, the random variable $\eta(B)$ is
Poisson distributed with parameter $(\nu\otimes\Leb)(B)$.

\item
The random measure $\eta$ is independently scattered: if
$B_1,\ldots,B_m\in\MscrZ\otimes\mathcal B((0,T])$ are pairwise disjoint
and have finite $(\nu\otimes\Leb)$-measure, then $\eta(B_1),\ldots,\eta(B_m)$ are mutually independent.

\item
For every $U\in\MscrZ$ with $\nu(U)<\infty$, the counting process
\[
N(t,U):=\eta(U\times(0,t]),
\qquad t\in[0,T],
\]
is $(\CF_t)$-adapted. Moreover, for every $0\le s<t\le T$,
\[
N(t,U)-N(s,U)
=
\eta(U\times(s,t])
\]
is independent of $\CF_s$ and is Poisson distributed with parameter
$(t-s)\nu(U)$.

\end{enumerate}
\end{definition}

\begin{remark}\label{compensator}
The intensity measure is uniquely determined by the law of $\eta$. Indeed,
for every $A\in\MscrZ$ with $\nu(A)<\infty$ and every $t\in(0,T]$,
\[
\EE\bigl[\eta(A\times(0,t])\bigr]
=
t\,\nu(A).
\]
The corresponding compensator of $\eta$ on $Z\times(0,T]$ is $\nu(dz)\,dt.$
\end{remark}

The preceding definition is formulated directly in terms of the random
measure. The following equivalent characterisation in terms of the
associated counting processes will be useful later.

\begin{lemma}\label{eqdefprm}
Let $\eta:
\Omega
\longrightarrow
M_{\NN_0}(\{Z_n\times(0,T]\})$ be measurable and define
\[
N(t,U):=\eta(U\times(0,t]),
\qquad
t\in[0,T],\quad U\in\MscrZ.
\]
Then $\eta$ is a time-homogeneous Poisson random measure with intensity
measure $\nu$ over $\mathfrak A$ if and only if the following properties
hold.

\begin{enumerate}[label=\upshape(\alph*)]

\item
For every $U\in\MscrZ$ with $\nu(U)<\infty$ and every $t\in[0,T]$,
$N(t,U)$ is Poisson distributed with parameter $t\nu(U)$. Moreover, for
every finite family of pairwise disjoint sets
$U_1,\ldots,U_m\in\MscrZ$ of finite $\nu$-measure,
\[
N(t,U_1),\ldots,N(t,U_m)
\]
are mutually independent.

\item
For every $t\in[0,T]$, the mapping $N(t,\cdot):
\Omega
\longrightarrow
M_{\NNnull }(\{Z_n\})$ is $\CF_t$-measurable.

\item
For every $0\le s<t\le T$, the random measure
\[
N(t,\cdot)-N(s,\cdot)
=
\eta(\,\cdot\,\times(s,t])
\]
is independent of $\CF_s$ and has the law of a Poisson random measure on
$Z$ with intensity measure $(t-s)\nu$.

\end{enumerate}
\end{lemma}

\begin{proof}
The implication from Definition~\ref{def-Prm} to
\emph{(a)}--\emph{(c)} follows immediately by restricting $\eta$ to
measurable rectangles.

Conversely, the stated properties yield the required Poisson laws and
independence first for finite unions of measurable rectangles of the form
$U\times(s,t]$. Since these rectangles generate
$\MscrZ\otimes\mathcal B((0,T])$, a standard monotone-class argument gives
properties \emph{(i)} and \emph{(ii)} of Definition~\ref{def-Prm}, while
\emph{(b)}--\emph{(c)} give the required adaptedness and independence of
future increments. We refer to
\cite[Definition~I.8.1 and Proposition~I.8.2]{ikeda1981}; see also
\cite[Chapter~3]{applebaum2009}.
\end{proof}

\paragraph{Compensated Poisson random measures}

Let $\eta$ be a time-homogeneous Poisson random measure on $(Z,\MscrZ)$
with intensity measure $\nu$. For every $B\in\MscrZ\otimes\mathcal B((0,T])$ satisfying $(\nu\otimes\Leb)(B)<\infty$, the corresponding compensated
Poisson random measure is defined by
\[
\tilde\eta(B)
:=
\eta(B)-(\nu\otimes\Leb)(B).
\]
In particular, for a measurable rectangle $U\times I$ satisfying
$\nu(U)\Leb(I)<\infty$,
\[
\tilde\eta(U\times I)
=
\eta(U\times I)-\nu(U)\Leb(I).
\]

For every $\sigma$-finite measure $\nu$ on $(Z,\MscrZ)$, a
time-homogeneous Poisson random measure with intensity measure $\nu$ exists
on a suitable filtered probability space; see, for example,
\cite{ikeda1981,applebaum2009}.

\subsection{Stochastic integration with respect to compensated Poisson random measures}
\label{subsec-Stochastic preliminaries}

We briefly recall the construction of the stochastic integral following
\cite{Brz+Haus_2009}. Let $p\in(1,2]$ and let $E$ be a Banach space of
martingale type $p$.

For $\xi\in
\mathcal M_{\rm step}^p
\bigl(0,T;L^p(Z,\nu;E)\bigr),$ the stochastic integral with respect to $\tilde\eta$ is defined by
\begin{equation}
\widetilde I(\xi)
=
\sum_{j=1}^n
\int_Z
\xi(t_j,z)\,
\tilde\eta(dz,(t_{j-1},t_j]).
\label{eqn-int-def}
\end{equation}

It follows from \cite{Brz+Haus_2009} that $\widetilde I$ is a bounded
linear map from
\[
\mathcal M_{\rm step}^p
\bigl(0,T;L^p(Z,\nu;E)\bigr)
\]
into $L^p(\Omega;E)$. More precisely, there exists a constant $C>0$,
depending only on $p$ and $E$, such that
\begin{equation}
\EE\bigl|\widetilde I(\xi)\bigr|_E^p
\le
C\,
\EE\int_0^T\int_Z
|\xi(t,z)|_E^p\,\nu(dz)\,dt.
\label{I1}
\end{equation}
Moreover, $\EE[\widetilde I(\xi)]=0.$ Since$ \mathcal M_{\rm step}^p
\bigl(0,T;L^p(Z,\nu;E)\bigr)$ is dense in $\mathcal M^p
\bigl(0,T;L^p(Z,\nu;E)\bigr),$ the stochastic integral extends uniquely to the latter space.

\begin{theorem}\label{th-ineq}
Let $p\in(1,2]$ and let $E$ be a Banach space of martingale type $p$.
Then there exists a unique bounded linear operator
\[
I:
\mathcal M^p(0,T;L^p(Z,\nu;E))
\longrightarrow
L^p(\Omega,\CF;E)
\]
such that
\[
I(\xi)=\widetilde I(\xi),
\qquad
\xi\in
\mathcal M_{\rm step}^p(0,T;L^p(Z,\nu;E)).
\]

For every
\[
\xi\in
\mathcal M^p(0,T;L^p(Z,\nu;E)),
\]
the process $\left(I(1_{(0,t]}\xi)\right)_{t\in[0,T]}$ is an $E$-valued $p$-integrable martingale and admits a c\`adl\`ag
modification. Moreover, there exists a constant $C>0$, depending only on
$p$ and $E$, such that
\begin{equation}
\EE
\sup_{t\in[0,T]}
\left|
I(1_{(0,t]}\xi)
\right|_E^p
\le
C\,
\EE\int_0^T\int_Z
|\xi(s,z)|_E^p\,\nu(dz)\,ds.
\label{I4}
\end{equation}
\end{theorem}

As usual, we write
\[
\int_0^t\int_Z
\xi(s,z)\,\tilde\eta(dz,ds)
:=
I(1_{(0,t]}\xi),
\qquad
t\in[0,T].
\]

\subsection{Compatibility and nonanticipativity}

A further issue arises when a Skorokhod-type representation theorem is
used. On the original probability space, nonanticipativity is expressed
through the filtration: the solution is adapted and future increments of the
driving Poisson random measure are independent of the past. A representation
theorem, however, transports laws rather than the original filtration.
Therefore, we need a formulation of nonanticipativity that can be expressed
through the joint law of the solution and the noise. This is provided by
Kurtz's notion of compatibility; see \cite{kurtz1998,kurtz2011}.

\begin{definition}[Compatibility of a pair]\label{jcomp}
Let $X$ be a process with paths in
\[
\DD([0,T];H)\cap L^\alpha(0,T;V),
\]
and let $\mu$ be an integer-valued random measure on $Z\times(0,T]$,
both defined on the same probability space.

The pair $(X,\mu)$ is called \emph{compatible} if, for every
$0\le s<t\le T$, every $p,q,r\in\NN\setminus\{0\}$, every $s_1,\ldots,s_p,s'_1,\ldots,s'_q\in[0,s],$ every $D_1,\ldots,D_q,C_1,\ldots,C_r\in\MscrZ$ with $\nu(D_l)<\infty, \ \nu(C_l)<\infty,$ and every pair of bounded Borel functions
\[
h:H^p\times\NN^q\longrightarrow\RR,
\qquad
\phi:\NN^r\longrightarrow\RR,
\]
one has
\begin{align}
&\EE\Bigl[
h\bigl(
X(s_1),\ldots,X(s_p),
\mu(D_1\times(0,s'_1]),\ldots,
\mu(D_q\times(0,s'_q])
\bigr)
\nonumber\\
&\hspace{2.5cm}\times
\phi\bigl(
\mu(C_1\times(s,t]),\ldots,
\mu(C_r\times(s,t])
\bigr)
\Bigr]
\label{eq:compat-def}
\\
&\quad=
\EE\Bigl[
h\bigl(
X(s_1),\ldots,X(s_p),
\mu(D_1\times(0,s'_1]),\ldots,
\mu(D_q\times(0,s'_q])
\bigr)
\Bigr]
\nonumber\\
&\hspace{2.5cm}\times
\EE\Bigl[
\phi\bigl(
\mu(C_1\times(s,t]),\ldots,
\mu(C_r\times(s,t])
\bigr)
\Bigr].
\nonumber
\end{align}
Thus, the past of the pair $(X,\mu)$ up to time $s$ carries no additional
information about the future increments of $\mu$.
\end{definition}

\begin{remark}\label{def:independence}
One can also replace $X$ by a filtration $(\CF_t)_{t\in[0,T]}$. Here, $((\CF_t)_{t\in[0,T]},\mu)$ are compatible means that the increments of the process $\mu$ are independent from $(\CF_t)_{t\in[0,T]}$.
\end{remark}

\begin{remark}\label{rem:compat-equiv}
Define
\[
\mathscr F_s^{X,\mu}
:=
\sigma\Bigl(
X(r),
\mu(D\times(0,r]):
0\le r\le s,\;
D\in\MscrZ,\;
\nu(D)<\infty
\Bigr).
\]
By a monotone-class argument, Definition~\ref{jcomp} is equivalent to the
conditional-expectation identity
\[
\EE\!\left[
\phi\bigl(
\mu(C_1\times(s,t]),\ldots,
\mu(C_r\times(s,t])
\bigr)
\,\middle|\,
\mathscr F_s^{X,\mu}
\right]
=
\EE\!\left[
\phi\bigl(
\mu(C_1\times(s,t]),\ldots,
\mu(C_r\times(s,t])
\bigr)
\right]
\]
for every bounded Borel $\phi$ and every finite family
$C_1,\ldots,C_r\in\MscrZ$ of finite $\nu$-measure.
In particular, on the original probability space, if $u_n$ is
$(\CF_t)_{t\in[0,T]}$-adapted and $\eta_n$ is a time-homogeneous Poisson random measure
over $\mathfrak{A}$ with filtration $(\CF_t)_{t\in[0,T]}$, then the pair $(u_n,\eta_n)$ is compatible.
\end{remark}

The essential feature of Definition~\ref{jcomp} is that compatibility
depends only on the joint law of the pair.

\begin{remark}[Transfer mechanism]\label{rem:compat-pathspace}
Suppose that $(X_1,\mu_1)$ and $(X_2,\mu_2)$ are defined on possibly
different probability spaces and satisfy
\[
\Law(X_1,\mu_1)
=
\Law(X_2,\mu_2).
\]
Then $(X_1,\mu_1)$ is compatible if and only if $(X_2,\mu_2)$ is
compatible. Indeed, both sides of \eqref{eq:compat-def} are expectations of
bounded Borel functionals of the pair.

This observation is used in the reconstruction argument as follows. On the
original probability space, compatibility is available for
$(u_n,\eta_n)$. The reconstructed measures are compared with $\eta_n$ on
the good events of the reconstruction procedure, and the corresponding
compatibility identities are transferred to the represented pairs through
equality of laws. Their passage to the limiting pair is established in
Step~(XII), where the required convergence of the corresponding bounded
Borel functionals is proved. Thus, compatibility of the limiting pair is
not assumed, but is obtained as part of the reconstruction argument.
\end{remark}

\begin{remark}\label{rem:compat-dyadic}
The compatibility condition can equivalently be tested against bounded
Borel functionals of finitely many dyadic coefficients of the solution,
since these coefficients are Borel functions of the path. This observation
is useful in the reconstruction procedure, where the shifted dyadic
approximations are used to encode the past of the solution. No additional
compatibility assumption on the dyadic coefficients is required.
\end{remark}

		\section{A Representation Theorem for Poisson Random Measures}\label{sec:representation}

Before going to the infinite-dimensional setting, let us recap the finite-dimensional setting.
		The construction of the jump measure associated with a c\`adl\`ag path in
		\(\RR^d\) is classical and can be carried out pathwise; see, for example,
		Tankov~\cite{conttankov}, Applebaum~\cite{zbMATH05553180}, and
		Ikeda--Watanabe~\cite{ikeda1981}. Let $X\in \mathbb D([0,T];\RR^d).$ For \(s\in(0,T]\), define the jump of \(X\) at time \(s\) by
		\[
		\Delta X_s:=X_s-X_{s-}.
		\]
		Since every c\`adl\`ag path possesses at most countably many jumps and only
		finitely many jumps larger than a prescribed threshold \(\varepsilon>0\) on
		a compact interval, one may define the truncated jump measure
		\[
		\mu^{X,\varepsilon}
		:=
		\sum_{0<s\le T}
		\mathbf 1_{\{\|\Delta X_s\|>\varepsilon\}}
		\,\delta_{(s,\Delta X_s)}.
		\]
		This is a finite point measure on $(0,T]\times\bigl(\RR^d\setminus B_\varepsilon(0)\bigr).$ Passing to the limit \(\varepsilon\downarrow0\), one obtains the jump
		measure
		\[
		\mu^X
		=
		\sum_{0<s\le T:\,\Delta X_s\neq0}
		\delta_{(s,\Delta X_s)},
		\]
		which is an integer-valued random measure on $(0,T]\times(\RR^d\setminus\{0\}).$ Thus, \(\mu^X\) records both the jump times and the corresponding jump
		amplitudes. If \(X\) is a L\'evy process, then \(\mu^X\) is in fact a
		Poisson random measure whose compensator is given by $dt\otimes \nu(dz),$ where \(\nu\) denotes the L\'evy measure.

		More generally, finite-dimensional jump processes are often constructed via
		martingale problems. If the limiting process solves the martingale problem
		associated with the integro-differential operator
		\[
		Lf(x)
		=
		b(x)\cdot \nabla f(x)
		+
		\int_Z
		\Bigl[
		f(x+g(x,z))
		-f(x)
		-\nabla f(x)\cdot g(x,z)
		\Bigr]
		\,\nu(dz),
		\]
		then, under the usual measurability and integrability assumptions, one may
		realise the solution on a suitable filtered probability space together with
		a Poisson random measure \(\eta\) having intensity measure $\nu(dz)\,dt,$ such that
		\[
		X_t
		=
		X_0
		+
		\int_0^t b(X_s)\,ds
		+
		\int_0^t\int_Z
		g(X_{s-},z)\,
		\widetilde{\eta}(dz,ds).
		\]
		This correspondence between martingale solutions and weak solutions of
		stochastic equations with jumps is well understood; see
		Kurtz~\cite{Kurtz2010}.

		The situation becomes substantially more delicate for stochastic partial
		differential equations. In contrast to the finite-dimensional case, the
		solution typically exhibits different spatial and temporal regularity
		properties. More precisely, one often has
		\[
		u\in L^\alpha(0,T;V)
		\]
		for a more regular space \(V\), while the pathwise time regularity is
		available only in a larger Hilbert space \(H\),
		\[
		u\in \mathbb D([0,T];H).
		\]
		Thus the solution is integrable in time with values in the stronger space
		\(V\), but possesses c\`adl\`ag trajectories only in the weaker space
		\(H\).

		Now, let us turn our attention to the infinite dimensional setting.
		The natural framework for the representation theorem developed in this
		section is provided by a Gelfand triple
		\[
		V
		\hookrightarrow
		H
		\equiv H^\ast
		\hookrightarrow
		V^\ast,
		\]
		where \(V\) is a reflexive and separable Banach space and \(H\) is a separable Hilbert
		space. We identify \(H\) with its dual \(H^\ast\). Let
		\((Z,\MscrZ )\) be a standard measurable space equipped with a
		$\sigma$-finite L\'evy measure \(\nu\). We further assume that both
		\(V\) and \(H\) are Banach spaces of martingale type \(2\) and possess
		the UMD property.

		The theory of maximal monotone operators on Hilbert and Banach spaces is
		classical; see Br\'ezis~\cite{Brezis1973} and
		Barbu~\cite{Barbu1976,Barbu2010}. The extension to operators
		\[
		A:V\to V^\ast
		\]
		defined on Gelfand triples, together with the coercivity,
		monotonicity, hemicontinuity, and growth assumptions employed throughout
		this work, is treated systematically in Lions~\cite{Lions1969} and
		Showalter~\cite{Showalter1997}. The stochastic framework with locally
		monotone coefficients is developed in
		Liu--R\"ockner~\cite{LiuRockner2010,LiuRockner2015} and
		Brzez\'niak--Liu--Zhu~\cite{BrzezniakLiuZhu2011}.

		We are interested in stochastic evolution equations of the form
		\begin{equation*}
			\left\{
			\begin{aligned}
					du(t)
			&=
			A\bigl(t,u(t)\bigr)\,dt
			+
			\int_Z g(u(t),z)\,\widetilde\eta(dz,dt),
			\\
			u(0)
			&=0.
		\end{aligned}
			\right.
		\end{equation*}In order to avoid the technical diffucilties, we are assuming our initial data to be $0$.

		The assumptions imposed on the coefficients \(A\) and  \(g\)
		are stated in the following subsections.

	\paragraph{\bf Structural assumptions on the drift operator:}
		To be more precise, we assume that the operator $A$ satisfies the hemicontinuity, local monotonicity, coercivity, and growth conditions stated below.

		\begin{hypo}\label{hyp}

			Let $V \hookrightarrow H \equiv H^\ast \hookrightarrow V^\ast$ be a Gelfand triple, where $V \hookrightarrow H $ is continuously and densely, $H$ and $V$ are separable and reflexiv,
			and let $A:[0,T]\times V\to V^\ast$ be a {jointly} measurable, i.e.\ $\mathscr B([0,T])\otimes\mathscr B(V)/\mathscr B(V^\ast)$-measurable, (possibly nonlinear) operator. Assume that there exist constants
			\[
			\alpha\ge 2,\qquad
			\beta\ge 0,\qquad
			\theta>0,\qquad
			C_0\ge0,
			\]
	a function $f\in L^1(0,T;\mathbb{R}_+),$ and a measurable, locally bounded function $\rho:V\to [0,\infty),$ such that the following conditions hold.

			\begin{enumerate}

				\item[\rm(H1)] \emph{Hemicontinuity.}
				For a.e.\ $t\in[0,T]$ and all \(w_1,w_2,w_3\in V\), the mapping $\lambda
				\longmapsto
				\bigl\langle A(t,w_1+\lambda w_2),\,w_3\bigr\rangle$ is continuous on \(\mathbb{R}\).

				\item[\rm(H2)] \emph{Local monotonicity.} For  a.e.\ $t\in[0,T]$ {and all \(w_1,w_2\in V\),}
				\[
				2\bigl\langle A(t,w_1)-A(t,w_2),\,w_1-w_2\bigr\rangle \le
				\bigl(f(t)+\rho(w_2)\bigr)
				|w_1-w_2|_H^2 ,
				\]
				where, for some constant $C>0$, the function $\rho$ satisfies
				\begin{align*}
					\rho(w) \leq C(1+|w|_V^\alpha)(1+|w|_H^\beta), \ \text{ for } \ w \in V.
				\end{align*}
				\item[\rm(H3)] \emph{Coercivity.}  For a.e.\ $t\in[0,T]$ {and every \(w\in V\),}
				\[
				2\langle A(t,w),w\rangle
				\le
				C_0|w|_H^2
				-\theta |w|_V^\alpha
				+f(t).
				\]

				\item[\rm(H4)] \emph{Growth condition.}  For a.e.\ $t\in[0,T]$ {and every \(w\in V\),}
				\[
				|A(t,w)|_{V^\ast}^{\frac{\alpha}{\alpha-1}}
				\le
				\bigl(f(t)+C_0|w|_V^\alpha\bigr)
				\bigl(1+|w|_H^\beta\bigr).
				\]

			\end{enumerate}

		\end{hypo}

		\paragraph{Assumptions on the jump coefficient.}
In the next section we list the properties of the jump coefficient.

		\begin{hypo}\label{ass:GF}
			Let $Z$ be a complete separable metric space with metric $d_Z$, with a
			distinguished point $0^\ast\in Z$, and write $|z|_Z:={d_Z}(z,0^\ast)$, and let
			$\mathscr{Z}$ be the associated Borel
			$\sigma$-algebra given by $\mathscr{Z} := \mathscr{B}(Z)$.
			Let $\nu$ be a L\'evy measure with bounded support and no mass at $0^\ast$. {To} be more precise, we assume that  $\nu(\{0^\ast\})=0$ and  $\supp\nu\subset\{|z|_Z\le R_0\}$.
			 Then, we assume that
			$g : V \times Z \to H$ is jointly measurable.
 Assume that there
			exist constants $C_{g}, L_{g} > 0$ such that, for all
			$v, v_{1}, v_{2} \in V$, and for $p= 2$
	\DEQSZ\label{eq:g-second-moment1}
			\int_{Z} |g(v,z)|_{H}^{p}\,\nu(dz)
			\le C_{g}\bigl( 1 + |v|_{V}^{p} \bigr),
		\EEQSZ
				and
		\DEQSZ\label{eq:g-second-moment2}
			\int_{Z} |g(v_{1},z) - g(v_{2},z)|_{H}^{p}\,\nu(dz)
			\le L_{g}\,|v_{1} - v_{2}|_{V}^{p} .
		\EEQSZ
				In addition, we assume the following \emph{weak nondegeneracy} property, i.e.,
			that there is a continuous, strictly positive
			function $\rho_g:[0,\infty)\to(0,\infty)$ with
			\DEQSZ\label{it:nondeg}
			(\ast\ast) \qquad
	|g(u,z)|_{H} \;\ge\; {\rho_g}\bigl(|u|_V\bigr)\,d_Z(z,0^\ast),
			\qquad u\in V,\ z\in Z.
					\EEQSZ
In addition, we assume the corresponding upper bound: there is $C_g>0$ with
\DEQSZ\label{it:upper}
(\ast\!\ast\!\ast)\qquad
|g(v,z)|_H\;\le\;C_g\bigl(1+|v|_V\bigr)\,d_Z(0^\ast,z),
\qquad v\in V,\ z\in Z .
\EEQSZ

			\smallskip

			For every $\epsilon > 0$ and $ R' > 0$:
			\begin{enumerate}[label=(\arabic*)]
				\item\label{it:inj} For every fixed $v \in V$, the map
				$z \mapsto g(v,z)$ is injective on
				$Z \setminus \{ d_Z(z,0^\ast) \le \epsilon \}$.
				\item\label{it:cont} \emph{(Continuity in $Z$.)}  For every
				$v \in V$, $z \mapsto g(v,z)$ is continuous on
				$Z \setminus \{  d_Z(z,0^\ast)\le \epsilon \}$.
				\item\label{it:lip} \emph{(Locally uniform Lipschitz continuity.)}
				There is $\Lt = \Lt(\epsilon, R') > 0$ such that
				\[
				\sup_{\substack{z \in Z \setminus \{{|z|_Z} \le \epsilon\} \\
						|z|_{Z} \le R'}}
				\bigl| g(v_{1},z) - g(v_{2},z) \bigr|_{H}
				\;\le\; \Lt(\epsilon, R')\,|v_{1} - v_{2}|_{V},
				\qquad v_{1}, v_{2} \in V .
				\]
			\end{enumerate}
		\end{hypo}

\begin{remark}
	\label{rem:GF-consequences}
	Let us shortly explain the consequences of  \eqref{it:nondeg} and \eqref{it:upper}:
	\begin{enumerate}[label=\upshape(\roman*)]
		\item For $R>0$ and $\epsilon>0$ we set
		\DEQSZ\label{eq:cRe}
		c(R,\epsilon)\;:=\;\epsilon\,\min_{0\le r\le R}\rho_g(r)\;>\;0,
		\EEQSZ
		the minimum being attained and strictly positive since $\rho_g$ is
		continuous and strictly positive on the compact interval $[0,R]$. By
		\eqref{it:nondeg}, $|g(v,z)|_H\ge c(R,\epsilon)$ for all $v\in V$ with
		$|v|_V\le R$ and all $z\in Z$ with $|z|_Z\ge\epsilon$. Comparing
		\eqref{it:nondeg} and \eqref{it:upper} at $v=0$ gives
		$\rho_g(0)\le C_g$.
		\item Since $\supp\nu\subset\{|z|_Z\le R_0\}$ and $\nu$ is a L\'evy
		measure, $\int_Z|z|_Z^2\,\nu(dz)<\infty$. Hence \eqref{it:upper}
		implies \eqref{eq:g-second-moment1} with $p=2$, with the constant
		$2C_g^2\int_Z|z|_Z^2\,\nu(dz)$ in place of $C_g$.
		Enlarging $C_g$ if necessary, we use one and the same constant $C_g$
		in \eqref{eq:g-second-moment1} and \eqref{it:upper}.
		\item \emph{(Uniform small-jump tail.)} For $p=2$ and $\varepsilon>0$ set
		\DEQSZ\label{eq:XV-tail}
		\kappa_p(\varepsilon)
		:=2^{p-1}C_g^p\int_{\{0<|z|_Z\le\varepsilon\}}\bigl(d_Z(0^\ast,z)\bigr)^p\,\nu(dz).
		\EEQSZ
		Then, by \eqref{it:upper},
		\[
		\int_{\{0<|z|_Z\le\varepsilon\}}|g(v,z)|^p_H\,\nu(dz)
		\le\kappa_p(\varepsilon)\bigl(1+|v|^p_V\bigr),\qquad v\in V,
		\]
		and $\kappa_p(\varepsilon)\downarrow0$ as $\varepsilon\downarrow0$ by
		monotone convergence, since $\nu(\{0^\ast\})=0$. The decay is uniform in
		$v$; it is used for the tiny jumps in Step~(VIII) and for the passage
	$J\to\infty$ in Step~(XIV), see \eqref{eq:XV-Mnat-bound}, \eqref{eq:XV-DeltaJ} and Lemma~\ref{lem:XV-smalljumps}.
	\end{enumerate}
\end{remark}

The weak nondegeneracy \eqref{it:nondeg} bounds $|g(v,z)|_H$ from below;
to tell \emph{which} mark produced a given jump we need that $g(v,\cdot)$
separates marks quantitatively, uniformly on $V$-balls.

\begin{hypo}\label{ass:g-sep}
	For all $R>0$ and $\epsilon>0$ there is a nondecreasing function
	$\theta_{R,\epsilon}:(0,\infty)\to(0,\infty)$ such that
	\[
	\bigl|g(v,z)-g(v,z')\bigr|_H\;\ge\;\theta_{R,\epsilon}\bigl(d_Z(z,z')\bigr)
	\qquad\text{for all } v\in V,\ |v|_V\le R,\ z,z'\in Z,\
	|z|_Z,|z'|_Z\ge\epsilon,\ { z\neq z'} .
	\]
	We take $\theta_{R,\epsilon}$ maximal; then $\theta_{R,\epsilon}$ is
	nonincreasing in $R$ and nondecreasing in $\epsilon$.
\end{hypo}

\begin{remark}
	Hypothesis~\ref{ass:g-sep} implies the injectivity of
	Hypothesis~\ref{ass:GF}\,\ref{it:inj} and the continuity of the inverse
	of $g(v,\cdot)$ on $\{|z|_Z\ge\epsilon\}$, uniformly in $|v|_V\le R$.
	For instance, if $g(v,z)=B(v)z$ with $Z\subset H$ and
	$|B(v)h|_H\ge\beta(R)|h|_H$ for $|v|_V\le R$, then
	$\theta_{R,\epsilon}(a)=\beta(R)\,a$.
\end{remark}
		\medskip

		\paragraph{\bf Approximation of the operator $A(s,\cdot)$:}In addition, in the spirit of an abstract space discretisation, respectively, as defined in \cite{cox01,cox02},  let us assume that we are given a family
$\{V_n:n\in\NN\}$ of closed subspaces of $V$ with
\begin{equation}\label{eq:Vn-dense}
	\bigcup_{n\in\NN}V_n\ \text{ dense in } V,
\end{equation}
and a family $\{A_n:n\in\NN\}$ of jointly measurable, i.e.\ $\mathscr B([0,T])\otimes\mathscr B(V_n)/\mathscr B(V^\ast)$-measurable, operators
$A_n:[0,T]\times V_n\to V^\ast$ such that $u^n_0\in V_n$ and $u_n$ is
$V_n$-valued,\footnote{By \eqref{originalun}, every jump of $u_n$ is of the form $\Delta u_n(s)=g(u_n(s-),z)$, where $(s,z)$ is an atom of $\eta$. Hence $u_n$ can be $V_n$-valued only if $g$ maps $V_n\times Z$ into $V_n$; this is assumed in Theorem~\ref{Theoremmain}. The condition is void if $V_n=V$ for all $n$. It holds, for instance, for finite-dimensional spaces $V_1\subset V_2\subset\cdots$ and $g(v,z)=\sigma(v)\,z$, where the marks form a closed set $Z\subset V_1$ (with the metric of $H$ and $0^\ast=0$) and $\sigma:V\to[c,\infty)$, $c>0$, is Lipschitz continuous on $V$.} and such that Hypothesis~\ref{hyp} holds for each $A_n$
\emph{on $V_n$}, uniformly in $n$, in the following sense.

\begin{hypo}[Hypothesis~\ref{hyp}, localised to $V_n$]\label{hyp:An}
	For every $n\in\NN$ the operator $A_n:[0,T]\times V_n\to V^\ast$ is jointly measurable, i.e.\ $\mathscr B([0,T])\otimes\mathscr B(V_n)/\mathscr B(V^\ast)$-measurable.
	There exist constants $\alpha\ge2$, $\beta\ge0$, {$\theta>0$},
	$C_0\ge0$, a function $f\in L^1(0,T;\RR_+)$ and a measurable, locally
	bounded {$\rho:V\to[0,\infty)$} --- all of them \emph{independent of
	$n$} --- such that for every $n\in\NN$ and a.e.\ $t\in[0,T]$:
	\begin{enumerate}
		\item[\rm(H1$_n$)] \emph{Hemicontinuity.} For all
		$w_1,w_2,w_3\in V_n$ the map
		$\lambda\mapsto\bigl\langle A_n(t,w_1+\lambda w_2),w_3\bigr\rangle$
		is continuous on $\RR$.

		\item[\rm(H2$_n$)] \emph{Local monotonicity.} For all
		$w_1,w_2\in V_n$,
		\[
		2\bigl\langle A_n(t,w_1)-A_n(t,w_2),\,w_1-w_2\bigr\rangle
		\le\bigl(f(t)+{\rho}(w_2)\bigr)\,|w_1-w_2|_H^2 .
		\]

		\item[\rm(H3$_n$)] \emph{Coercivity.} For all $w\in V_n$,
		\[
		2\bigl\langle A_n(t,w),w\bigr\rangle
		\le C_0|w|_H^2-{\theta}|w|_V^\alpha+f(t).
		\]

		\item[\rm(H4$_n$)] \emph{Growth condition.} For all $w\in V_n$,
		\[
		\bigl|A_n(t,w)\bigr|_{V^\ast}^{\frac{\alpha}{\alpha-1}}
		\le\bigl(f(t)+C_0|w|_V^\alpha\bigr)\bigl(1+|w|_H^\beta\bigr).
		\]
	\end{enumerate}
	The constants $\alpha,\beta,{\theta},C_0$ and the functions $f,{\rho}$
	are those of Hypothesis~\ref{hyp}.
\end{hypo}

   		\bigskip

\paragraph{The main result.}
Before stating the main result, let us describe in words what it says and
which assumptions enter. The starting point is a sequence of approximating
problems. For every $n\in\NN$ the process $u_n$ solves the equation
\eqref{originalun} on one fixed filtered probability space $\mathfrak A$,
driven by one and the same Poisson random measure $\eta$, but with a drift
$A_n$ in place of the drift $A$ of the limit problem. Such sequences arise,
for instance, from Galerkin or finite element discretisations of the
equation, or from approximations of the nonlinearity in the drift. For them
one can usually derive uniform a~priori bounds and, from these, the tightness
of the laws of $u_n$ on $\DD([0,T];H)\cap L^\alpha(0,T;V)$. Convergence of
the sequence itself on $\mathfrak A$, on the other hand, is in general not
available, unless additional information such as pathwise uniqueness of the
limit problem can be used.

Tightness gives, by a Skorokhod--Jakubowski type argument, a new probability
space on which copies of $u_n$ converge almost surely. The difficulty is that
the noise does not come along with this procedure in a usable form: the
stochastic integral is not a continuous functional of the pair $(u_n,\eta)$,
so it is not clear that the limit process is driven by a Poisson random
measure on the new space at all. Theorem~\ref{Theoremmain} closes this gap.
It states that the limit can be realised on a filtered probability space
which carries, besides the limit process $u_\ast$, a Poisson random measure
$\eta_\ast$ with the same intensity $\nu$, compatible with the filtration,
such that $u_\ast$ solves the limit equation driven by $\widetilde\eta_\ast$.
In other words, the limit of a tight family of approximating solutions is a
martingale solution of the limit problem. The Poisson random measure
$\eta_\ast$ is not obtained abstractly; it is reconstructed from the jumps of
the paths.

		\begin{theorem}\label{Theoremmain}
			Let the Gelfand triplet: $V\hookrightarrow H \hookrightarrow V^\ast$, where \(V\) is a reflexive and separable Banach space and \(H\) is a separable Hilbert
			space. 		Let $Z$ be a complete separable metric space with metric $d_Z$, with a
				distinguished point $0^\ast\in Z$, and write $|z|_Z:={ d_Z}(z,0^\ast)$, and let
			$\mathscr{Z}$ be the associated Borel
			$\sigma$-algebra given by $\mathscr{Z} := \mathscr{B}(Z)$.
	Let \(\nu\) be a
			\(\sigma\)-finite measure on \((Z,\MscrZ )\), {finite on $Z\setminus\{|z|_Z\le \eps\}$ for every $\eps>0$}. Let
			\[
			\eta:\Omega\to
			M_{\NNnull }(\{Z_k\times(0,T]\})
			\]
			be a time-homogeneous Poisson random measure on \(Z\) with intensity measure
			\(\nu\), defined over the filtered probability space $\mathfrak A=(\Omega,\CF,(\CF_t)_{t\ge0},\PP).$
			Let $A:[0,T]\times V\to V^\ast$ satisfy Hypothesis~\ref{hyp}.  Let $g:V\times Z\to H$ be jointly measurable and satisfy Hypotheses~\ref{ass:GF} and~\ref{ass:g-sep}.

In addition, let us assume that we have a family of operators $\{ A_n:n\in\NN\}$  satisfying Hypotheses \ref{hyp:An} such that  the approximating data converge to the
limiting data in the following weak sense: for every $\phi\in V$,
\begin{equation}\label{eq:DG}
	\int_0^T\bigl|\bigl\langle A_n(s,u_n(s))-A(s,u_n(s)),\phi\bigr\rangle\bigr|\,ds
	\longrightarrow0
	\qquad\text{in probability, as }n\to\infty.
\end{equation}%

Let us assume the initial condition to be $0$. 
			Assume that $g(v,z)\in V_n$ for all $n\in\NN$, $v\in V_n$ and $z\in Z$, and that the process $u_n$ is $(\CF_t)$-adapted $V_n$-valued, with c\`adl\`ag paths in $H$  and satisfies, for \(t\in[0,T]\),
\DEQSZ\label{originalun}
			u_n(t)
			=
			\int_0^t A_n(s,u_n(s))\,ds
			+
			\int_0^t\int_Z
			g(u_n(s),z)\,\widetilde\eta(dz,ds),
\EEQSZ
			in  $H$, \(\PP\)-a.s.
			Assume that there exists \(C>0\), independent of \(n\), such that
			\[
			\sup_{n\in\NN}
			\EE\left[
			\sup_{t\in[0,T]} |u_n(t)|_H^2
			+
			\int_0^T |u_n(t) |_V^\alpha\,dt
			\right]
			\le C.
			\]
		If $\alpha=2$, assume the uniform integrability moreover
		\begin{align}\label{eq:UI-V}
			\lim_{R\to\infty}\,\sup_{n\in\NN}\,
			\EE\int_0^T|u_n(t)|_V^2\,\mathds 1_{\{|u_n(t)|_V>R\}}\,dt=0 ;
		\end{align}
		for $\alpha>2$ this follows from the previous bound, since the left-hand
		side is at most \break $R^{2-\alpha}\sup_n\EE\int_0^T|u_n(t)|_V^\alpha\,dt$.
		In addition, assume that $\{ \Law(u_n):n\in\NN\}$ is tight on\footnote{The space $\DD([0,T];H)$ is equipped with the usual Skorokhod $J_1$-topology, see Appendix \ref{app:Skorohod-space}. The space $L^\alpha(0,T;V)$ carries its norm topology, and the intersection carries the join topology of Lemma~\ref{lem:state-space}\,\ref{SS:i}. Step~(XIV) uses the norm topology. If the laws were tight only for the weak topology of $L^\alpha(0,T;V)$, the drift could not be identified as in Step~5(a) there.} $\DD([0,T];H)\cap L^\alpha(0,T;V)$.

			Then, there exists a filtered probability space $\mathfrak A_\ast
			= (\Omega_\ast,\CF_\ast,(\CF_t^\ast)_{t\in[0,T]},\PP_\ast),$ a process
			\[
			u_\ast\in \DD([0,T];H)\cap L^\alpha(0,T;V),
			\]
			and a time-homogeneous Poisson random measure $\eta_\ast$ on \(Z\) with intensity measure \(\nu\), such that:

			\begin{enumerate}[label=\upshape(\roman*)]
				\item \(u_\ast\) is \((\CF_t^\ast)_{t\in[0,T]}\)-progressively measurable;

				\item \(\eta_\ast\) has compensator $\mathscr{Z}
				\times \mathscr{B}([0,T])\ni (U,I)\mapsto \gamma(U,I)=\nu(U)\,\Leb(I)$;

				\item \(\eta_\ast\) is compatible with the filtration
				\((\CF_t^\ast)_{t\in[0,T]}\), that is, for every \(0\le s<t\), $\eta_\ast(U\times(s,t])$
				is independent of the $\sigma$-algebra  \(\CF_s^\ast\) for every \(U\in\MscrZ \) with
				\(\nu(U)<\infty\);

				\item \(u_\ast\) satisfies
				\[
				u_\ast(t)
				=
				\int_0^t A(s,u_\ast(s))\,ds
				+
				\int_0^t\int_Z
			{g(u_\ast(s),z)}\,\widetilde\eta_\ast(dz,ds),
				\]
				for every \(t\in[0,T]\), \(\PP_\ast\)-a.s.
			\end{enumerate}

		\end{theorem}

\subsection{Open problems and possible extensions}
\label{sec:open-problems}

The hypotheses of Theorem~\ref{Theoremmain} were chosen to make the
reconstruction as transparent as possible, not because they are
believed to be sharp. We collect here the extensions that we expect to
hold, together with the point at which the present proof breaks down
in each case, and some questions that appear to be open or not tackled.

\medskip
\emph{(1) Removing the injectivity condition.}
The Hypothesis~\ref{ass:GF} \ref{it:inj}, that is, the injectivity of
$z\mapsto g(v,z)$ outside every neighbourhood of $0\in Z$ is used
once to make the proof simpler.  It makes the jumps
$z$ recoverable from the observed jump amplitude
$\Delta u(s)=g(u(s-),z)$ in $H$, via the Borel inverse provided by the
Lusin--Suslin theorem (Lemma~\ref{lem:DH-meas}). It is a
restrictive assumption and it fails, for instance, whenever $g(v,\cdot)$ is even in $z$, or whenever $\dim H<\infty$ and $Z$ is
higher-dimensional. Without it, the solution path determines only the
\emph{image} measure on $H$, i.e.\ the point measure of jump times and
amplitudes, and the fibres $g(v,\cdot)^{-1}\{x\}$ may contain more
than one jump. We conjecture that the representation theorem
nevertheless persists in the following weaker form: On an
\emph{extension} of $\mathfrak A_1$ carrying an independent sequence
of uniform random variables, one can construct a Poisson random
measure $\eta_1$ with intensity $\nu$ such that the limit equation
holds. The jump size being chosen, at each registered jump, according to
the disintegration of $\nu$ over the fibres of $g(v,\cdot)$. The
program is analogous to the classical representation of
integer-valued random measures with given compensator as images of
Poisson random measures on an extension (Grigelionis \cite[Theorem 1]{Grigelionis1971}, El~Karoui--Lepeltier \cite{ElKarouiLepeltier1977}; see also Ikeda--Watanabe \cite[Lemma 7.1 and Lemma 7.2]{IkedaWatanabe1981}), but two points do not follow from the classical
theory, that is,  the measurability of the fiber disintegration
\emph{jointly} in $(v,x)$, and the compatibility (in the sense of
Definition~\ref{jcomp}) of the added randomness with the pair
filtration through the diagonal limit $m\to\infty$. The advantage of the injectivity condition is that no extension of the probability space is required.
However, the condition can be relaxed, and proving without \ref{it:inj} should be possible similarly to what is done in  \cite{FKEHMH}.

\medskip
\emph{(2) L\'evy measures with unbounded support.}
The standing assumption $\supp\nu\subset\{|z|_Z\le R_0\}$
(Hypothesis~\ref{ass:GF}) enters through the moment estimates and it
guarantees
$$
\displaystyle\int_Z|g(v,z)|^2_H\,\nu(dz)<\infty
$$
from the growth
condition alone and keeps the moduli of the registration algorithm
uniform over the layers. It is not structural: the layer construction
itself only uses $\nu\bigl(\{|z|_Z>\epsilon\}\bigr)<\infty$, which
holds for every L\'evy measure. The natural extension could be done by treating the large jump separately. The large jumps are isolated; between two large jumps, Theorem \ref{Theoremmain}  can be applied.

\medskip
\emph{(3) Mixed and more general noises.}
Two further directions suggest themselves. First, equations driven by
a Wiener process \emph{and} a Poisson random measure: the jump part
should be recoverable by the present registration algorithm, and the
Gaussian part through a quadratic-variation analogue. Here, one would get a {L\'evy--Ciesielski} construction of the Brownian motion, done in another context in \cite{HausenblasKumarToelle2026}, the difficult point is the joint compatibility of the two
reconstructed noises with one pair filtration. Secondly, cylindrical
L\'evy noise without a jump measure of locally finite intensity in
$H$: there the very notion of registering individual jumps fails, and
a reconstruction would have to proceed at the level of the
characteristics rather than the atoms; it is not clear to us
whether a pathwise reconstruction theorem can hold at all in that
generality, and a counterexample would be as interesting as a proof.

\section{The proof of Theorem \ref{Theoremmain}}

			\textbf{Overview of the proof strategy.}
The proof is organised in fourteen steps. \emph{For a fixed index $n$ and a fixed dyadic level $m$}, Steps~(I)--(VII) develop the reconstruction of the jump structure from the solution path. In order, they introduce the shifted dyadic projection (I), the martingale increments (II), the noise bridge (III), the registration algorithm (IV), the random measure of registered jumps (V), a measurable inverse of $g$ (VI), and the reconstructed Poisson random measure on $Z$ (VII).

Step~(VIII) then chooses $n$, $m$, and the tolerance as functions of a single parameter (Lemma~\ref{lem:subsequence}) and proves that the reconstruction is exact on, so-called,
good events (Lemma~\ref{lem:X-good-event}). Step~(IX) applies the Jakubowski--Skorokhod representation theorem, and Step~(X) introduces the path jump measure. Steps~(XI)--(XIII) identify the limits of the random measures, construct the filtration, and prove that the limiting noise is a time-homogeneous Poisson random measure. Finally, Step~(XIV) shows that the limiting process is a solution.

			Following the nomenclature of e.g. Daley and Vere-Jones \cite{DaleyVereJones2003}, we will call a jump in $H$ a jump, and a jump in the underlying metric space $Z$ a mark.

				Now, we begin the proof.
	
					\begin{proof}[Proof of Theorem \ref{Theoremmain}] The proof of this theorem has been divided into several steps. 
				\smallskip
				\noindent

			\begin{steps}

				\item \textbf{Construction of the progressively measurable approximation to $u_n$ by the shifted dyadic projection and of the filtration on the new probability space $\mathfrak{A}$:}

				Due to the definition of the stochastic integral, our choice is to define a progressively measurable sequence for each $u_n$, by taking the shifted Haar Projection.

				Let \(n,m\in\NN\) be fixed and let
				\DEQSZ\label{def_shifted_dyadic}
				\myp_{n,m}(s)
				&=&
				\sum_{k=0}^{{2^m-1}}\xi_{k,n}^{(m)}\,\mathbf 1_{(k2^{-m}T,\;(k+1)2^{-m}T]}(s),
				\qquad s\in[0,T],
				\EEQSZ
				where \(\xi_{k,n}^{(m)}:\Omega\to V\) are \(\mathscr F\)-measurable random variables for $k\ge1$
				given by
				\DEQSZ\label{def_xi}
				\xi_{k,n}^{(m)} &:=&\frac{2^m}{T}\, \int_{t^m_{k-1}}^{t^m_{k}} u_n(s)\, ds,\quad k=1,\dots ,2^m-1,\, m\in\NN .
				\EEQSZ
	For $k=0$, we put
	\DEQSZ\label{init0}
	\xi_{0,n}^{(m)}:=0,
	\EEQSZ  
	For $m\in\NN$ let $\tau_m:=T2^{-m}$, $t^m_k:=k\tau_m$ and
	$I^m_k:=(t^m_k,t^m_{k+1}]$, $k=0,\dots,2^m-1$; we call $I^m_k$ the
	\emph{$k$-th cell} of order $m$.

				Then, for all $n\in\NN$,  \(\myp_{n,m} \xrightarrow[m\to\infty]{}  u_n\) in \(L^1(0,T;V)\) $\PP$-a.s. (Lemma~\ref{lem:Pro-L1} with $E=V$), and for all $n,m\in\NN$, the process \(\myp_{n,m}\) is simple and predictable (Lemma~\ref{lem:Pro-props}).
\pagebreak[2]
						\item \textbf{Construction of the approximate martingale increments.}

			In this step we reconstruct the martingale part of the process.
				Since we only have the processes at our disposal, we approximate the increments of the martingale part
\DEQSZ\label{def:N}
N_n(t)&:=&\int_0^t\int_Z g(u_n(s),z)\,\tilde\eta(dz,ds),\qquad t\in[0,T],
\EEQSZ
by the increments of the process itself; the integral is interpreted as in Theorem~\ref{th-ineq}. 
By \eqref{originalun}, for $0\le t_1\le t_2\le T$,
\DEQSZ\label{mymn}
u_n(t_2)-u_n(t_1)&=&N_n(t_2)-N_n(t_1)+\int_{t_1}^{t_2}A_n\bigl(s,u_n(s)\bigr)\,ds .
\EEQSZ
Although $A_n$ takes values only in $V^\ast$, the drift integral
$D_n(t):=\int_0^tA_n(s,u_n(s))\,ds=u_n(t)-u^n_0-N_n(t)$ is continuous in $H$, so that its
increments over the cells are $H$-valued random variables, see Lemma~\ref{lem:drift-meas} below.

						\item \textbf{Construction of the  Jumps.}

						In this step, we reconstruct, at the dyadic level, the martingale part from the
transported solution processes. Let
\DEQSZ\label{def:lkn}
\ell_{k,n,m} =u_n(t_{k+1}^m)-u_n(t_{k}^m)-\int_{I^m_k}A_n\bigl(s,u_n(s)\bigr)\,ds .
\EEQSZ 

The purpose of the \emph{noise bridge} is to recover the driving noise
from the solution paths alone. A candidate for the
martingale part is therefore the computable quantity
$\ell_{k,n,m}$.

The drift integral in \eqref{def:lkn} is taken over the cell, i.e.\ over $I^m_k=(t^m_k,t^m_{k+1}]$. It is a Bochner integral in $V^\ast$, and it is not clear a~priori that it defines an $H$-valued measurable functional of the path. The next lemma shows that it does, without any further assumption on $A_n$. Here $\mathcal X:=L^\alpha(0,T;V)\cap\DD([0,T];H)$ is the path space of Lemma~\ref{lem:state-space}, and we use the convention $A_n(s,v):=0$ for $v\in H\setminus V_n$.

\begin{lemma}[Measurability of the drift increments]\label{lem:drift-meas}
Let Hypothesis~\ref{hyp:An} hold, where $A_n:[0,T]\times V_n\to V^\ast$ is $\mathscr B([0,T])\otimes\mathscr B(V_n)/\mathscr B(V^\ast)$-measurable. Fix $n,m\in\NN$ and put
\[
G_n:=\Bigl\{x\in\mathcal X:\ \int_0^T\bigl|A_n(s,x(s))\bigr|_{V^\ast}\,ds<\infty\Bigr\},
\qquad
D^m_{k,n}(x):=\int_{I^m_k}A_n\bigl(s,x(s)\bigr)\,ds\in V^\ast,\quad x\in G_n,
\]
for $k=0,\dots,2^m-1$, and
\[
G^H_{n,m}:=\bigl\{x\in G_n:\ D^m_{k,n}(x)\in H\ \text{ for all }k=0,\dots,2^m-1\bigr\}.
\]
Then the following hold.
\begin{enumerate}[label=\upshape(\roman*)]
\item\label{DM:i} $G_n$ and $G^H_{n,m}$ are Borel subsets of $\mathcal X$.
\item\label{DM:ii} For every $k$, the map $\mathcal X\ni x\mapsto\mathds 1_{G^H_{n,m}}(x)\,D^m_{k,n}(x)\in H$ is $\mathscr B(\mathcal X)/\mathscr B(H)$-measurable. Consequently, the \emph{increment functional}
\[
\ell^m_{k,n}(x):=\mathds 1_{G^H_{n,m}}(x)\,\Bigl(x(t^m_{k+1})-x(t^m_k)-D^m_{k,n}(x)\Bigr),\qquad x\in\mathcal X,
\]
is a Borel map from $\mathcal X$ into $H$.
\item\label{DM:iii} $\PP\bigl(u_n\in G^H_{n,m}\bigr)=1$, and $\PP$-a.s., for all $k$,
\[
D^m_{k,n}(u_n)=u_n(t^m_{k+1})-u_n(t^m_k)-\bigl(N_n(t^m_{k+1})-N_n(t^m_k)\bigr),
\qquad
\ell_{k,n,m}=\ell^m_{k,n}(u_n)=N_n(t^m_{k+1})-N_n(t^m_k).
\]
In particular, $\ell_{k,n,m}$ is an $H$-valued, $\CF_{t^m_{k+1}}$-measurable random variable.
\end{enumerate}
\end{lemma}

\begin{proof}
\emph{Step 1 (The evaluation map).} For fixed $s$, the map $x\mapsto x(s)\in H$ is Borel on $\DD([0,T];H)$, see \cite[Theorem~12.5]{Billingsley1999}, hence on $\mathcal X$. For fixed $x$, the map $s\mapsto x(s)$ is right-continuous. Hence $(s,x)\mapsto x(s)$ is the pointwise limit of the maps $(s,x)\mapsto x\bigl((\lceil 2^ms/T\rceil T2^{-m})\wedge T\bigr)$, $j\to\infty$, each of which is $\mathscr B([0,T])\otimes\mathscr B(\mathcal X)$-measurable, and is therefore jointly measurable.

\emph{Step 2 (From $H$ to $V_n$).} The embedding $V\hookrightarrow H$ is continuous and injective between Polish spaces. By the Lusin--Suslin theorem \cite[Thm.~15.1]{Kechris1995}, $V\in\mathscr B(H)$ and $\mathscr B(V)=\mathscr B(H)\cap V$. Since $V_n$ is closed in $V$, also $V_n\in\mathscr B(H)$ and $\mathscr B(V_n)=\mathscr B(H)\cap V_n$. By Step~1, the set $\{(s,x):x(s)\in V_n\}$ is measurable, and on it $(s,x)\mapsto(s,x(s))$ is measurable into $[0,T]\times V_n$. Together with the joint measurability of $A_n$ and the convention $A_n(s,v)=0$ for $v\notin V_n$, the map
\[
[0,T]\times\mathcal X\ni(s,x)\longmapsto a(s,x):=A_n\bigl(s,x(s)\bigr)\in V^\ast
\]
is $\mathscr B([0,T])\otimes\mathscr B(\mathcal X)/\mathscr B(V^\ast)$-measurable.

\emph{Step 3 (The integral is Borel into $V^\ast$).} Since $V$ is reflexive and separable, $V^\ast$ is separable. Hence $(s,x)\mapsto|a(s,x)|_{V^\ast}$ is jointly measurable, and $G_n$ is a Borel set by Tonelli's theorem. For $x\in G_n$ the function $a(\cdot,x)$ is Bochner integrable in $V^\ast$, so $D^m_{k,n}(x)$ is well defined. Fix a countable set $D_V$ which is dense in $V$. For $\phi\in D_V$, Fubini's theorem shows that
\[
x\longmapsto\mathds 1_{G_n}(x)\,\bigl\langle D^m_{k,n}(x),\phi\bigr\rangle
=\mathds 1_{G_n}(x)\int_{I^m_k}\bigl\langle a(s,x),\phi\bigr\rangle\,ds
\]
is Borel on $\mathcal X$. Since $|y|_{V^\ast}=\sup\{\langle y,\phi\rangle:\phi\in D_V,\ |\phi|_V\le1\}$ and $V^\ast$ is separable, $\mathscr B(V^\ast)$ is generated by the maps $y\mapsto\langle y,\phi\rangle$, $\phi\in D_V$. Hence $x\mapsto\mathds 1_{G_n}(x)D^m_{k,n}(x)$ is Borel from $\mathcal X$ into $V^\ast$.

\emph{Step 4 (From $V^\ast$ to $H$).} For $y\in V^\ast$ put
\[
|y|_H:=\sup\bigl\{\langle y,\phi\rangle:\ \phi\in D_V,\ |\phi|_H\le1\bigr\}\in[0,\infty].
\]
This is a countable supremum of continuous functions on $V^\ast$, hence Borel. The set $\{\phi\in V:|\phi|_H<1\}$ is open in $V$, so $D_V$ is dense in it, and $V$ is dense in $H$. Therefore $|y|_H$ is the norm of $H$ if $y\in H$, and $|y|_H=\infty$ if $y\in V^\ast\setminus H$. Consequently $H=\{y\in V^\ast:|y|_H<\infty\}\in\mathscr B(V^\ast)$. Moreover, for every $h\in H$ the map $y\mapsto|y-h|_H$ is Borel on $V^\ast$, so every ball of $H$ belongs to $\mathscr B(V^\ast)\cap H$; since $H$ is separable, $\mathscr B(H)=\mathscr B(V^\ast)\cap H$. Now \ref{DM:i} and the first part of \ref{DM:ii} follow from Step~3:
\[
G^H_{n,m}=G_n\cap\bigcap_{k=0}^{2^m-1}\bigl\{x:\ \bigl|\mathds 1_{G_n}(x)D^m_{k,n}(x)\bigr|_H<\infty\bigr\}
\]
is Borel, and $x\mapsto\mathds 1_{G^H_{n,m}}(x)D^m_{k,n}(x)$ is Borel into $V^\ast$ with values in $H$, hence Borel into $H$. The second part of \ref{DM:ii} follows, since $x\mapsto x(t^m_{k+1})-x(t^m_k)$ is Borel into $H$ by Step~1.

\emph{Step 5 (Along the solution).} $\PP$-a.s., $u_n\in L^\alpha(0,T;V)$ and $\sup_{t\le T}|u_n(t)|_H<\infty$. By (H4$_n$), $A_n(\cdot,u_n(\cdot))\break \in L^{\alpha/(\alpha-1)}(0,T;V^\ast)\subset L^1(0,T;V^\ast)$, so $u_n\in G_n$, $\PP$-a.s. By \eqref{mymn}, which is an identity in $V^\ast$, $\PP$-a.s.\ for all $k$,
\[
D^m_{k,n}(u_n)=u_n(t^m_{k+1})-u_n(t^m_k)-\bigl(N_n(t^m_{k+1})-N_n(t^m_k)\bigr).
\]
The right-hand side belongs to $H$, since $u_n$ and $N_n$ are $H$-valued. Hence $u_n\in G^H_{n,m}$, $\PP$-a.s., and the formula for $\ell_{k,n,m}$ follows from \eqref{def:lkn}. Finally, $N_n(t^m_{k+1})-N_n(t^m_k)$ is $\CF_{t^m_{k+1}}$-measurable, and the filtration is complete. This proves \ref{DM:iii}.
\end{proof}

\medskip

				\item \textbf{Construction of the jumps:}
				We have constructed an approximation of the difference $\ell_{k,m,n}$. However, we still have to construct the jumps on the mark space themselves.

				Here, we want to introduce the following notation.
				\noindent For a parameter $\epsilon$ we define $Z_\epsilon:=\{z\in Z: |z|\le \epsilon\}$ and $Z_\epsilon^c:=Z\setminus Z_{\epsilon}$.
				Fix a strictly decreasing sequence $\{\ep_j>0:j\in\NN_0\}$ with $\epsilon_0> R_0$ and $\epsilon_j\downarrow 0$, so that $Z^c_{\ep_j}\subset Z^c_{\epsilon_{j+1}}$, and such that
				$$
				\nu\bigl(Z^c_{\ep_{j+1}}\setminus Z^c_{\epsilon_{j}}\bigr)<\infty
				\quad\mbox{and}\quad
				\nu\bigl(\{z\in Z: |z|_Z=\epsilon_j\}\bigr)=0, \qquad j\in\NN_0.
				$$
				Next, let us define for $j\in\NN_0$,
				\DEQSZ\label{eq:layers}
				\diff Z_j:= Z^c_{\epsilon_{j+1}}\setminus  Z^c_{\epsilon_{j}}
				=\{z\in Z:\epsilon_{j+1}<|z|_Z\le \epsilon_j\}.
				\EEQSZ

				\begin{remark}\label{rem:support}
	By Hypothesis~\ref{ass:GF}, $\supp\nu\subset\{|z|_Z\le R_0\}$.
	Since $\epsilon_0 > R_0$, every layer is bounded,
	\[
	\overline{\Delta Z_j}
	\subset
	\bigl\{ z \in Z : \epsilon_{j+1} \le |z|_Z \le \epsilon_j \bigr\},
	\]
	and the layers cover all large jumps:
	$$
	\nu\Bigl( Z^c_{\epsilon_{J+1}}\setminus \bigcup_{\jt=0}^{J}\Delta Z_\jt\Bigr)
	=\nu\bigl(\{|z|_Z>\epsilon_0\}\bigr)=0 .
	$$
	In particular, $\displaystyle\nu\bigl(\bigcup_{\jt=0}^{J}\Delta Z_\jt\bigr)=\nu(Z^c_{\epsilon_{J+1}})<\infty$.
\end{remark}

					Observe, due to the construction of the operator $\Pro_m$, the dyadic coefficients $\xi^{(m)}_{k,n}$ of \eqref{def_xi} are $V$-valued (hence in particular $H$-valued) for all $k=0,1,2,\cdots,2^m-1$.
					Let us define for given $n, m \in \N$, $k \in \{0,\dots,2^{m}-1\}$ and
					$j \in \{0,\dots,J\}$ a
					(random) \emph{admissible jump set} by
					\begin{equation}\label{eq:DH}
						\Delta H^{k,m,n}_{j}
						:= g\bigl(\xi^{(m)}_{k,n}, \Delta Z_{j}\bigr)
						= \bigl\{ g\bigl(\xi^{(m)}_{k,n}, z\bigr) : z \in \Delta Z_{j} \bigr\}
						\;\subset\; H .
					\end{equation}

					By Hypothesis~\ref{ass:GF} and $\xi^{(m)}_{k,n} \in V$ a.s.,
					\begin{equation}\label{eq:g-int}
						\int_{Z} \bigl|g\bigl(\xi^{(m)}_{k,n}, z\bigr)\bigr|_{H}^{p}
						\,\nu( dz)
						\;\le\; C_{g}\bigl(1 + |\xi^{(m)}_{k,n}|_{V}^{p}\bigr) < \infty
						\qquad \PP\text{-a.s.},
					\end{equation}
					and consequently
					\[
					\EE\bigg[ \int_{Z} \bigl(1 \wedge
					|g(\xi^{(m)}_{k,n},z)|_{H}^{p}\bigr)\,\nu(d z) \bigg]< \infty,
					\qquad
					\sup_{m\in\N} \EE\bigg[ \int_{0}^{T}\!\!\int_{Z}
					\bigl(1 \wedge |g( u_{n,m}(s),z)|_{H}^{p}\bigr)
					\,\nu(d z) ds\bigg] < \infty .
					\]

						The next lemma settles the measurability questions attached to the
					random sets \eqref{eq:DH}.  It also provides the measurable
					inversion of  $g$, which will be used later.
					First we analyse the measurability of the admissible jump sets.
					\begin{lemma}
						\label{lem:DH-meas}
						Let Hypothesis~\ref{ass:GF} hold and fix
						$J \in \N$.  For $j \in \{0,\dots,J\}$ let us consider the graph map
						\[
						G_{j} : V \times \Delta Z_{j} \longrightarrow V \times H,
						\qquad
						G_{j}(v,z) := \bigl( v,\, g(v,z) \bigr).
						\]
						Then the following hold:
						\begin{enumerate}[label=\textup{(\roman*)}]
							\item\label{L:i} $G_{j}$ is injective and Borel  measurable; its range
							$Ran(G_{j}) := G_{j}\bigl(V\times\Delta Z_{j}\bigr)
							=\{(v,g(v,z)):\ v\in V,\ z\in\Delta Z_{j}\}$ is a Borel subset of $V \times H$, and
							$G_{j}^{-1} : Ran(G_{j}) \to V \times \Delta Z_{j}$ is also Borel measurable.
							\item\label{L:ii} For every $v \in V$, the set
							$\Delta H_{j}(v) := g(v, \Delta Z_{j})
							= \{ x \in H : (v,x) \in Ran(G_{j})  \}$
							is a Borel subset of $H$, and the sets
							$\Delta H_{0}(v), \dots, \Delta H_{J}(v)$ are pairwise disjoint.
							In particular, for every $\omega \in \Omega$, the admissible
							jump sets $\Delta H^{k,m,n}_{j}(\omega)
							= \Delta H_{j}\bigl(\xi^{(m)}_{k,n}(\omega)\bigr)$ are Borel and
							pairwise disjoint in $j$.
							\item\label{L:iii}
							For every
							$V$-valued random variable $\xi$ and every $H$-valued random
							variable $Y$ on $(\Omega,\cF)$, the event
							\[
							\bigl\{ Y \in \Delta H_{j}(\xi) \bigr\}
							= (\xi, Y)^{-1}(Ran(G_{j}) )
							\;\in\; \sigma(\xi, Y) .
							\]
							\item\label{L:iv}
							The map
							\[
							V \times H \ni (v,x) \longmapsto
							d_{j}(v,x) := \dist_{H}\bigl( x,\, \Delta H_{j}(v) \bigr)
							= \inf_{z \in \Delta Z_{j}} \bigl| x - g(v,z) \bigr|_{H}
							\]
							is Borel measurable on $V \times H$.  Consequently, for $\xi$, $Y$
							as in \ref{L:iii}, the random variables
							$d_{j}(\xi, Y)$, $\min_{0 \le j \le J} d_{j}(\xi,Y)$, and, for every
							$\delta > 0$, the registration events
							$\{ d_{j}(\xi,Y) \le \delta \}$ are $\sigma(\xi,Y)$-measurable.
							\item\label{L:v}
							Let 	$\pi_{Z}$ denote the projection onto the $Z$-coordinate. On the event
							$\{ Y \in \Delta H_{j}(\xi) \}$, the map
							\[
							\zhat_{j} := \pi_{Z} \circ G_{j}^{-1}(\xi, Y),
							\]
							is a
							well-defined, $\sigma(\xi,Y)$-measurable $\Delta Z_{j}$-valued
							random variable satisfying $g(\xi, \zhat_{j}) = Y$ on this event;
							it is the unique mark in $\Delta Z_{j}$ with this property.
						\end{enumerate}
					\end{lemma}

					\begin{proof}
						\ref{L:i}  We know that $V$ and $H$ are separable Banach spaces, hence both spaces are Polish.
						The space $Z$ is also Polish and $\Delta Z_{j} \in \mathscr Z$ {by definition}, since
					the set 	$\supp\nu$ is closed and the annulus is a Borel subset.  Hence, we know that
						$V \times \Delta Z_{j}$ and $V \times H$ are standard Borel spaces.
						The map $G_{j}$ is Borel because $g$ is jointly measurable
						(Hypothesis~\ref{ass:GF}) and the identity on $V$ is continuous.
					Let us analyse what happens with the injectivity.
					If $G_{j}(v,z) = G_{j}(v',z')$, the first coordinate
						gives $v = v'$, and then $g(v,z) = g(v,z')$ with
						$z, z' \in \Delta Z_{j} \subset
						Z \setminus \{|z| \le \epsilon_{J+1}\}$. It follows,  by
						Hypothesis~\ref{ass:GF} \ref{it:inj} applied with
						$\epsilon = \epsilon_{J+1}$, $z = z'$.  The Lusin--Suslin theorem
						\cite[Thm.~15.1]{Kechris1995} states that an injective Borel map
						between standard Borel spaces has Borel range and Borel inverse.
						This gives all three assertions in (i).

					\smallskip

						\ref{L:ii} From the fact that $\Delta H_{j}(v)$ is the $v$-section of the Borel set
						$Ran(G_{j}) $, and sections of Borel sets are Borel sets, it follows that
						$\Delta H_{j}(v)$ is a Borel set.	Next, we will show that these sets are disjoint.
						If	$x \in \Delta H_{i}(v) \cap \Delta H_{j}(v)$ with $i \neq j$, then
						$x = g(v, z_{i}) = g(v, z_{j})$ with $z_{i} \in \Delta Z_{i}$,
						$z_{j} \in \Delta Z_{j}$. To be more precise,  both elements  lie in
						$Z \setminus \{|z| \le \epsilon_{J+1}\}$, so
						Hypothesis~\ref{ass:GF}-\ref{it:inj} forces $z_{i} = z_{j}$,
						contradicting $\Delta Z_{i} \cap \Delta Z_{j} = \emptyset$.  The
						statement for the random sets follows by substituting
						$v = \xi^{(m)}_{k,n}(\omega)$.

											\smallskip

						\ref{L:iii} This point follows  from \ref{L:i}. To be more precise,  one has to take into account that $Ran(G_{j}) $ is a Borel set  and
						$(\xi, Y)$ is a $\sigma(\xi,Y)/\cB(V \times H)$-measurable map.

											\smallskip

						\ref{L:iv} Since $Z$ is separable, it follows that  $\Delta Z_{j}$ is separable. Fix a
						countable dense subset $D_{j} \subset \Delta Z_{j}$.  We claim
						\begin{equation}\label{eq:dense}
							d_{j}(v,x)
							= \inf_{z \in D_{j}} \bigl| x - g(v,z) \bigr|_{H}
							\qquad \text{for all } (v,x) \in V \times H .
						\end{equation}
						Indeed,  the direction ``$\le$'' is trivial. For the direction  ``$\ge$'', let
						$z \in \Delta Z_{j}$ and $z_{\ell} \in D_{j}$ with
						$z_{\ell} \to z$ in $Z$.  Since
						$\Delta Z_{j} \subset Z \setminus \{|z| \le \epsilon_{j+1}\}$ and
						$g(v,\cdot)$ is continuous there by
						Hypothesis~\ref{ass:GF}-\ref{it:cont},
						 it follows that
						$g(v,z_{\ell}) \to g(v,z)$ in $H$. Hence
						$|x - g(v,z_{\ell})|_{H} \to |x - g(v,z)|_{H}$, which proves
						\eqref{eq:dense}.  For each fixed $z$, the map
						$(v,x) \mapsto |x - g(v,z)|_{H}$ is Borel: $v \mapsto g(v,z)$ is
						Borel (a section of the jointly measurable $g$), the norm distance
						is jointly continuous in $(x, g)$, and a Carath\'eodory function
						(measurable in $v$, continuous in $x$) on a product of a measurable
						space with a separable metric space is jointly measurable.  A
						countable infimum of Borel functions is Borel, so \eqref{eq:dense}
						shows $d_{j}$ is a Borel function on $V \times H$.  The remaining assertions
						follow by composing $d_j$ with the measurable map $(\xi, Y)$ and taking
						finite minima.

	\smallskip

						\ref{L:v} By \ref{L:i}, we know that  $\pi_{Z} \circ G_{j}^{-1}$ is Borel measurable on
						$Ran(G_{j}) $. Composing the mapping  with $(\xi,Y)$ on the   event $\{(\xi,Y) \in Ran(G_{j}) \}$, which is $\sigma(\xi,Y)$-measurable
						by \ref{L:iii}, gives measurability.
						The identity $g(\xi,\zhat_{j}) = Y$ holds because
						$G_{j}(\xi, \zhat_{j}) = (\xi, Y)$ by definition of the inverse;
						uniqueness of the jump is Hypothesis~\ref{ass:GF}-\ref{it:inj}.
					\end{proof}

					 The
					layer $Z_j$, $j\in\NN$, of Lemma~\ref{lem:DH-meas} are defined for elements in $V$. To be more precise, the
					admissible sets $\Delta H_\jt(v)$, the graph $Ran(G_{j}) \subset V\times H$
					and the inverse $\zhat_\jt$ all presuppose a $V$-valued first
					argument. The \emph{registration algorithm} later on evaluates them at
					the dyadic coefficients $\xi^{(m)}_{k,n}$, which are $V$-valued by
					construction.

					\begin{lemma}\label{lem:DH-meas-H}
						Let Hypothesis~\ref{ass:GF} hold and fix $J\in\NN$ and
						$\jt\in\{0,\dots,J\}$.
						\begin{enumerate}[label=\upshape(\roman*)]
							\item\label{HH:i} Every closed ball
							$B^V_k:=\{x\in V:|x|_V\le k\}$ is a closed subset of $H$;
							consequently $V=\bigcup_{k\in\NN}B^V_k\in\mathscr B(H)$, and the
							extended norm $|\cdot|_V:H\to[0,\infty]$, set to $+\infty$ on
							$H\setminus V$, is $\mathscr B(H)$-measurable.

							\item\label{HH:ii} The range $Ran(G_{j}) $ of
							Lemma~\ref{lem:DH-meas} \ref{L:i} is a Borel subset of
							$H\times H$, and
							$G_\jt^{-1}:Ran(G_{j})\to V\times\Delta Z_\jt$ is Borel measurable as
							a map defined on this subset of $H\times H$.

							\item\label{HH:iii} Adopt the convention
							$\Delta H_\jt(v):=\emptyset$ for $v\in H\setminus V$;
							equivalently, replace every membership test
							$Y\in\Delta H_\jt(\xi)$ by $(\xi,Y)\in Ran(G_{j}) $. Then the
							statements \ref{L:ii},
							 \ref{L:iii} and \ref{L:v} of
							Lemma~\ref{lem:DH-meas} remain valid verbatim for arbitrary
							\emph{$H$-valued} random variables $\xi$: the event
							$\{(\xi,Y)\in Ran(G_{j}) \}$ belongs to $\sigma(\xi,Y)$; the sets
							$\Delta H_0(v),\dots,\Delta H_J(v)$ are pairwise disjoint Borel
							subsets of $H$ for every $v\in H$; and on
							$\{(\xi,Y)\in Ran(G_{j}) \}$ the jump
							$\zhat_\jt:=\pi_Z\circ G_\jt^{-1}(\xi,Y)$ is a
							$\sigma(\xi,Y)$-measurable $\Delta Z_\jt$-valued random variable
							with $g(\xi,\zhat_\jt)=Y$; on this event, we have automatically
							$\xi\in V$.
						\end{enumerate}
					\end{lemma}

					\begin{proof}
						\ref{HH:i} Let $x_m\in B^V_k$ with $x_m\to x$ in $H$. Since $V$ is
reflexive, the $V$-bounded sequence $(x_m)$ has a subsequence
$x_{m'}\rightharpoonup y$ weakly in $V$, with
$|y|_V\le\liminf_{m'}|x_{m'}|_V\le k$ by weak lower semicontinuity
of the norm. Since the embedding $V\hookrightarrow H$ is continuous,
$x_{m'}\rightharpoonup y$ weakly in $H$ as well, and the weak limit
is unique, so that $y=x\in B^V_k$. Hence $B^V_k$ is closed in $H$,
$V\in\mathscr B(H)$, and $\{|\cdot|_V\le k\}=B^V_k$ gives the
measurability of the extended norm.

						\ref{HH:ii} By \ref{HH:i}, we know that  $V\times\Delta Z_\jt$ is a Borel subset
of the Polish space $H\times Z$. Hence $V\times\Delta Z_\jt$ is a standard Borel space, and
$G_\jt$ is an injective Borel map from $V\times\Delta Z_\jt$ into the Polish space
$H\times H$ (injectivity exactly as in
Lemma~\ref{lem:DH-meas} \ref{L:i}); the Lusin--Suslin theorem
\cite[Thm.~15.1]{Kechris1995} yields that the range is Borel in
$H\times H$ and that the inverse is Borel.

						\ref{HH:iii} The proof is split into four steps.

						\emph{Step 1 (The set identity).} For $v\in V$ the set
						$\Delta H_\jt(v)$ is, by Lemma~\ref{lem:DH-meas} \ref{L:ii}, the
						$v$-section $\{x\in H:(v,x)\in Ran(G_\jt)\}$ of $R_\jt$. For
						$v\in H\setminus V$ the convention sets
						$\Delta H_\jt(v)=\emptyset$, and the $v$-section of $Ran(G_j)$ is
						empty as well, since every point of
						$Ran(G_j)\subseteq V\times H$ has first
						coordinate in $V$. Hence, for \emph{every} $v\in H$,
						$\Delta H_\jt(v)=\{x\in H:(v,x)\in Ran(G_j)\}$, and consequently,
						$\omega$-by-$\omega$,
						$\{Y\in\Delta H_\jt(\xi)\}=\{(\xi,Y)\in Ran(G_j)\}$ for an arbitrary
						$H$-valued random variable $\xi$; the convention is the unique
						extension of the layer sets to $H$ for which this identity holds.

						\emph{Step 2 (Extension of \ref{L:iii}).} The map
						$\omega\mapsto(\xi(\omega),Y(\omega))$ is
						$\sigma(\xi,Y)/\mathscr B(H\times H)$-measurable, since $H$ is
						separable, so that
						$\mathscr B(H\times H)=\mathscr B(H)\otimes\mathscr B(H)$. By
						\ref{HH:ii}, $Ran(G_j)\in\mathscr B(H\times H)$.
						Taking the preimage of the Borel set
						$Ran(G_j)$ under the measurable map $(\xi,Y)$ gives
						\[
						\{Y\in\Delta H_\jt(\xi)\}=(\xi,Y)^{-1}(Ran(G_j))\in\sigma(\xi,Y).
						\]

						\emph{Step 3 (Extension of \ref{L:ii}).} For $v\in V$ the sets
						$\Delta H_0(v),\dots,\Delta H_J(v)$ are Borel and pairwise
						disjoint by Lemma~\ref{lem:DH-meas} \ref{L:ii}, the disjointness
						resting on the injectivity of $g(v,\cdot)$ on
						$\{|z|_Z>\epsilon_{J+1}\}$, Hypothesis~\ref{ass:GF}\,\ref{it:inj}.
						For $v\in H\setminus V$ all $J+1$ sets equal $\emptyset$; the
						empty set is Borel and disjoint from every set, so the assertions
						hold trivially.

						\emph{Step 4 (Extension of \ref{L:v}; on the event, automatically
						$\xi\in V$).} Put $E:=\{(\xi,Y)\in Ran(G_j)\}$, so that
						$E\in\sigma(\xi,Y)$ by Step~2. On $E$ the random variable
						$\zhat_\jt:=\pi_Z\circ G_\jt^{-1}(\xi,Y)$ is the composition of
						the $\sigma(\xi,Y)$-measurable map $(\xi,Y)$, the map
						$G_\jt^{-1}$ is Borel on $Ran(G_j)$ \emph{as a subset of
						$H\times H$} by \ref{HH:ii}  and the continuous projection
						$\pi_Z$; hence it is a $\sigma(\xi,Y)$-measurable
						$\Delta Z_\jt$-valued random variable on $E$. Moreover, on $E$,
						\[
						(\xi,Y)=G_\jt\bigl(w,\zhat_\jt\bigr),
						\qquad
						(w,\zhat_\jt):=G_\jt^{-1}(\xi,Y)\in V\times\Delta Z_\jt .
						\]
						Since the first coordinate of $G_\jt$ is the identity, comparing
						first coordinates gives $w=\xi$; in particular $\xi\in V$ on $E$,
						with no separate hypothesis membership of $(\xi,Y)$ in the
						range $Ran(G_j)\subseteq V\times H$ forces the first coordinate into
						$V$. Comparing second coordinates gives $g(\xi,\zhat_\jt)=Y$.
						Uniqueness of the jump in $\Delta Z_\jt$ follows from
						Hypothesis~\ref{ass:GF} \ref{it:inj}, exactly as in
						Lemma~\ref{lem:DH-meas} \ref{L:v}.
					\end{proof}
						Since the recovered jump term  carries some error terms,
it lies in general only \emph{near}, but not \emph{in}, one of the admissible jump sets $\Delta H^{k,m,n}_j$. However, one has to select a set $\Delta H^{k,m,n}_j$ to which it will belong. Here, due to the error,
						we  relax the membership test and admit a tolerance: An increment is registered as a jump as soon as its distance to some layer is at most $\delta$. If it lies  $\delta$-close to two adjacent layers simultaneously, we assign the membership to the layer with the minimal distance.

					\begin{definition}
						\label{def:selector}
						Let $\delta>0$ be a tolerance parameter, let $\xi$ be a $V$-valued and
						$Y$ an $H$-valued random variable on $(\Omega,\mathcal F)$, and
						recall the distance functionals
						\[
						d_j(\xi,Y)=\dist_H\bigl(Y,\Delta H_j(\xi)\bigr)
						=\inf_{z\in\Delta Z_j}\bigl|Y-g(\xi,z)\bigr|_H,
						\qquad j=0,\ldots,J,
						\]
						of Lemma~\ref{lem:DH-meas} \ref{L:iv}.
						\begin{enumerate}
							\item[\textup{(a)}] The \emph{registration event} at tolerance
							$\delta$ is
							\[
							\mathrm R^{\delta}(\xi,Y):=\Bigl\{\min_{0\le j\le J} d_j(\xi,Y)\le\delta\Bigr\}.
							\]
							On $\mathrm R^{\delta}(\xi,Y)$ we say that the increment $Y$ \emph{passes the
								tolerance test}, i.e.\ it is registered as a jump.
							\item[\textup{(b)}] The \emph{nearest-layer selector} is the
							$\{0,\ldots,J\}\cup\{\partial\}$-valued random variable
							\[
							\hat\jmath(\xi,Y):=
							\begin{cases}
								\min\Bigl\{ j\in\{0,\ldots,J\}\;:\;
								d_j(\xi,Y)=\min\limits_{0\le i\le J} d_i(\xi,Y)\Bigr\},
								& \text{on } \mathrm R^{\delta}(\xi,Y),\\[6pt]
								\partial, & \text{on } \Omega\setminus \mathrm R^{\delta}(\xi,Y),
							\end{cases}
							\]
							where $\partial$ is a so-called cemetery symbol (no registration).  The minimum
							over the argmin set is a deterministic rule for deciding the outcome when two (or more) sets are equally scored, so that the decision is still well defined; any other
							measurable rule may be used instead.
						\end{enumerate}
						In the setting of Step~(V) we apply this with
						$\xi=\xi^{(m)}_{k,n}$ and $Y=\xhat_{k,n,m}$ from
						\eqref{def_xhat}, and write
						$\mathrm r^{\delta}_{k,m,n}:=\mathrm R^{\delta}(\xi^{(m)}_{k,n},\xhat_{k,n,m})$
						and
						$\hat\jmath_{k,m,n}:=\hat\jmath(\xi^{(m)}_{k,n},\xhat_{k,n,m})$.
					\end{definition}
					In the next Lemma we prove the measurability of the selector rule.

					\begin{lemma}
						\label{lem:selector-meas}
						Let Hypothesis~\ref{ass:GF} hold and let $\xi$, $Y$ be as in
						Definition~\ref{def:selector}. Then, for every $\delta>0$:
						\begin{enumerate}
							\item[\textup{(i)}] $\mathrm R^{\delta}(\xi,Y)\in\sigma(\xi,Y)$, and
							$\{\hat\jmath(\xi,Y)=j\}\in\sigma(\xi,Y)$ for every
							$j\in\{0,\ldots,J\}\cup\{\partial\}$; in particular
							$\hat\jmath(\xi,Y)$ is a $\sigma(\xi,Y)$-measurable random variable.
							\item[\textup{(ii)}] If $\xi=\xi^{(m)}_{k,n}$ and
							$Y=\xhat_{k,n,m}$, then $\mathrm r^{\delta}_{k,m,n}$ and
							$\hat\jmath_{k,m,n}$ are $\mathscr F_{t^m_{k+1}}$-measurable.
							\item[\textup{(iii)}] On the exact-membership event
							$\{Y\in\Delta H_j(\xi)\}$ one has $d_j(\xi,Y)=0$, hence
							$\mathrm R^{\delta}(\xi,Y)$ occurs and, if in addition
							$d_i(\xi,Y)>0$ for all $i\neq j$, then $\hat\jmath(\xi,Y)=j$; on
							$\{\hat\jmath=j\}\cap\{Y\in\Delta H_j(\xi)\}$ the recovered jump
							$\hat z_j=\pi_Z\circ G_j^{-1}(\xi,Y)$ of
							Lemma~\ref{lem:DH-meas} \ref{L:v} is well defined.
						\end{enumerate}
					\end{lemma}

					\begin{proof}
						(i) By Lemma~\ref{lem:DH-meas} \ref{L:iv}, each $d_j(\xi,Y)$,
						$j=0,\ldots,J$, is $\sigma(\xi,Y)$-measurable; hence so is
						$d^{(J)}:=\min_{0\le i\le J}d_i(\xi,Y)$, and
						$\mathrm R^{\delta}(\xi,Y)=\{d^{(J)}\le\delta\}\in\sigma(\xi,Y)$. For
						$j\in\{0,\ldots,J\}$,
						\[
						\{\hat\jmath(\xi,Y)=j\}
						=\mathrm R ^{\delta}(\xi,Y)\cap\{d_j(\xi,Y)=d^{(J)}\}\cap
						\bigcap_{i<j}\{d_i(\xi,Y)>d^{(J)}\},
						\]
						a finite intersection of $\sigma(\xi,Y)$-measurable events, and
						$\{\hat\jmath=\partial\}=\Omega \setminus \mathrm R^\delta(\xi,Y)$. Since
						$\{0,\ldots,J\}\cup\{\partial\}$ is finite, measurability of the
						level sets is measurability of $\hat\jmath$.

						(ii) $\xi^{(m)}_{k,n}$ is $\mathscr F_{t^m_k}$-measurable and
						$\xhat_{k,n,m}$ is $\mathscr F_{t^m_{k+1}}$-measurable, so
						$\sigma(\xi^{(m)}_{k,n},\xhat_{k,n,m})\subset
						\mathscr F_{t^m_{k+1}}$, and (ii) follows from (i).

						(iii) The first two assertions are immediate from the definitions and
						the last one is Lemma~\ref{lem:DH-meas} \ref{L:v}.
					\end{proof}

					\subsubsection*{Algorithm: extraction of the jumps on each time interval in $H$}
					\label{alg:jump-extraction}

					Fix $m,n\in\NN$, the number of layers $J\in\NN$, a tolerance $\delta_r>0$ and a
					cut-off $R_r>0$.  Later on, in Step~(VIII) all parameters are fixed as
functions of one parameter $r$: $n=n_r$, $m=m_r$,
and $\delta_r$, $R_r$ (see
 Lemma~\ref{lem:subsequence}, and
					Remark~\ref{rem:dictionary}).
The dependence on the common parameter will be chosen in such a way, that we will get convergence of the here constructed random measure to a Poisson random measure with L\'evy measure $\nu$.
In particular, we assume here for $r\ge r_0(J)$ that they satisfy
					\begin{equation}\label{eq:delta-window}
						2\delta_r\;<\;c(R_r,\epsilon_{J+1}),
						\qquad \delta_r\downarrow 0,\quad R_r\uparrow\infty.
					\end{equation}
			Under this assumption, we know 		by Lemma~\ref{lem:subsequence}\,\ref{DS:iii}, p.\ \pageref{lem:subsequence}, where $c(R,\epsilon)$ is
					the constant \eqref{eq:cRe} of the weak nondegeneracy
					\eqref{it:nondeg}, that  every admissible jump $g(v,z)$ with $|v|_V\le R_r$
					and $z\in\Delta Z_0\cup\dots\cup\Delta Z_J$ has $H$-norm at least
					$c(R_r,\epsilon_{J+1})$. Hence, $\delta_r$ is sufficiently small.

Now, let us present the exact algorithm. The
					following {steps} are carried out independently for every dyadic
					time interval $I^m_k=(t^m_k,t^m_{k+1}]$, $k=0,\ldots,2^m-1$.

					\medskip
					\noindent\textbf{Input:} At each time interval $I^m_k$,
					the frozen state $\xi^{(m)}_{k,n}$
					(which is $\mathscr F_{t^m_k}$-measurable) and the two
values $u_n(t^m_k)$, $u_n(t^m_{k+1})$ and the drift increment $D^m_{k,n}(u_n)=\int_{I^m_k}A_n(s,u_n(s))\,ds$ of Lemma~\ref{lem:drift-meas} are given.

					\medskip
					\noindent\textbf{Output:} At each time interval $I^m_k$, the algorithm says either \emph{``no jump in $I^m_k$''}, or a
					triple $(\hat\jmath,\Upsilon_{k,m,\jty},\hat x)$ consisting of a layer index
					$\hat\jmath\in\{0,\ldots,J\}$, a jump time
					$\Upsilon_{k,m,\jty}\in I^m_k$ and a jump size $\hat x\in H$ will be constructed.

					\begin{enumerate}
						\item[\textbf{(A1)}] \textbf{Increment.} Compute (compare \eqref{def:lkn})
\DEQSZ\label{mart_incr}
\ell:=\ell_{k,n,m}
=u_n(t^m_{k+1})-u_n(t^m_k)-{\int_{I^m_k}}A_n\bigl(s,u_n(s)\bigr)\,ds ,
\EEQSZ
By Lemma~\ref{lem:drift-meas}\,\ref{DM:iii}, $\ell=N_n(t^m_{k+1})-N_n(t^m_k)$, $\PP$-a.s., and $\ell$ is an $H$-valued random variable.

						\item[\textbf{(A2)}] \textbf{Decompensation.} Add the compensator
						back:
					\DEQSZ\label{def:jump}
						\hat x:=\hat x_{k,n,m}=\ell+\lbar^{(J)}_{k,m,n},
						&\qquad&
						\lbar^{(J)}_{k,m,n}
						=\tau_m\sum_{\jt=0}^{J}\int_{\Delta Z_\jt}
						g\bigl(\xi^{(m)}_{k,n},z\bigr)\,\nu(dz),
					\EEQSZ

						\item[\textbf{(A3)}] \textbf{Distances to the layers.} For
						$\jt=0,\ldots,J$ compute
						\[
						d_\jt:=d_\jt\bigl(\xi^{(m)}_{k,n},\hat x\bigr)
						=\inf_{z\in\Delta Z_\jt}
						\bigl|\hat x-g\bigl(\xi^{(m)}_{k,n},z\bigr)\bigr|_H
						=\inf_{z\in D_\jt}
						\bigl|\hat x-g\bigl(\xi^{(m)}_{k,n},z\bigr)\bigr|_H,
						\]
						where $D_\jt\subset\Delta Z_\jt$ is a countable dense subset, cf.\
						\eqref{eq:dense}; set $d^{(J)}:=\min_{0\le\jt\le J}d_\jt$.

						\item[\textbf{(A4)}] \textbf{Registration (thresholding with a
							$V$-cut-off).} One decides whether there is a jump in $I^m_k$: the
						recovered increment must lie within the tolerance of one of the
						admissible jump sets, and the frozen state must lie in the $V$-ball
						of radius $R_r$, on which the lower bound \eqref{it:nondeg} is
						uniform.
						\begin{itemize}
							\item If $\bigl|\xi^{(m)}_{k,n}\bigr|_V>R_r$ or $d^{(J)}>\delta_r$:
							output \emph{``no jump in $I^m_k$''} and stop.
							\item If $\bigl|\xi^{(m)}_{k,n}\bigr|_V\le R_r$ and
							$d^{(J)}\le\delta_r$: a jump is registered; continue.
						\end{itemize}
						Since $\xi^{(m)}_{k,n}$ is $\mathscr F_{t^m_k}$-measurable, the
						cut-off does not affect the measurability of the registration
						event. By \eqref{eq:delta-window}, a cell without a jump of mark in
						$Z^c_{\epsilon_{J+1}}$ and with small errors is not registered, while
						a true jump is registered as soon as the errors are below
						$\delta_r/2$; the precise statement is
						Lemma~\ref{lem:X-good-event} in Step~(VIII).

						\item[\textbf{(A5)}] \textbf{Layer selection (nearest layer,
							deterministic rule).} Set
						\[
						\hat\jmath:=\min\bigl\{\jt\in\{0,\ldots,J\}: d_\jt=d^{(J)}\bigr\}.
						\]
Thus, if several layers are equally close to \(\hat x\), we will select the layer with the
smallest index.
So, the rule for which layer is selected is \emph{deterministic}.

						\item[\textbf{(A6)}] \textbf{Recording the $H$-valued atom.} Assign
						the jump time $\Upsilon_{k,m.J}:=\tfrac12(t^m_k+t^m_{k+1})$. Here, we chose the midpoint and record the atom
						$\delta_{(\hat x,\,\Upsilon_{k,m,\jty})}$ in the layer-$\hat\jmath$ time-interval
						measure $\hat N^{k,n,m,\hat\jmath}_{\up}$ of
						\eqref{eq:Nhat-delta}.

\item[\textbf{(A7)}] \textbf{Recovering the jumps  in $Z$.}
\begin{itemize}
	\item If $\hat x\in\Delta H_{\hat\jmath}(\xi^{(m)}_{k,n})$
	\emph{exactly}, the jump  is recovered deterministically and uniquely by
	the measurable inversion of Lemma~\ref{lem:DH-meas} \ref{L:v}:
	$$\hat z:=\pi_Z\circ G_{\hat\jmath}^{-1}\bigl(\xi^{(m)}_{k,n},\hat x\bigr)
	\in\Delta Z_{\hat\jmath}, $$
	and $g(\xi^{(m)}_{k,n},\hat z)=\hat x$.
	\item In general (A4) guarantees only
	$d_{\hat\jmath}(\xi^{(m)}_{k,n},\hat x)\le\delta_r$, which does not imply that
	$(\xi^{(m)}_{k,n},\hat x)\in Ran(G_{\hat\jmath})$, so neither the inversion
	of \ref{L:v} given by $\zhat_\jt:=\pi_Z\circ G_\jt^{-1}(\xi,Y)$
	nor  \eqref{eq:f-explicit}, which returns the default $z_0$ off the layer applies to $\hat x$ itself. The jump is
	recovered by the \emph{tolerance-adapted inverse}. Let us choose a measurable
	$\hat z\in\Delta Z_{\hat\jmath}$ with
	\[
	\bigl|g(\xi^{(m)}_{k,n},\hat z)-\hat x\bigr|_H
	\;\le\;d_{\hat\jmath}(\xi^{(m)}_{k,n},\hat x)+2^{-m},
	\]
	e.g.\ the first point of a fixed countable dense subset of
	$\Delta Z_{\hat\jmath}$ achieving this bound (which is measurable, since $Z$ is
	separable).
\end{itemize}

					\end{enumerate}

					\begin{remark}\label{rem:single-atom}%
					For a given $J$, 	the algorithm constructs at most one jump per time interval and layer, in
						accordance with the definition \eqref{eq:Nhat-delta}. This choice is made without loss of generality. By a Borel--Cantelli argument, $\PP$-a.s.\ for all
						sufficiently large $m$ every time interval $I^m_k$ contains at most one jump
						of size in $\bigcup_{\jt\le J}\Delta Z_\jt$, so that for $m\ge m_0(\omega)$, only one jump occurs.
					\end{remark}

					By applying  Lemma~\ref{lem:DH-meas} with $\xi = \xi^{(m)}_{k,n}$, which is
					$\cF_{t^{m}_{k}}$-measurable, and
					$Y = \xhat_{k,n,m} = \ell_{k,n,m}
					+ \bar\ell_{k,m,n}$, which is
					$\cF_{t^{m}_{k+1}}$-measurable, the following result can be proven.
					\begin{corollary}%
						\label{cor:adapted}
						For every
						$j \in \{0,\dots,J\}$ and $\delta > 0$ the following items are satisfied:
						\begin{enumerate}[label=\textup{(\alph*)}]
							\item the membership event
							$\{\xhat_{k,n,m} \in \Delta H^{k,m,n}_{j}\}$, the distances
							$d_{j}(\xi^{(m)}_{k,n}, \xhat_{k,n,m})$, the registration
							event $\{ \min_{j \le J} d_{j} \le \delta \}$ and the nearest-layer
							selector are all $\cF_{t^{m}_{k+1}}$-measurable;
							\item on the membership event, the recovered jump
							$\zhat_{j} = \pi_{Z} \circ G_{j}^{-1}(\xi^{(m)}_{k,n},
							\xhat_{k,n,m})$ is an
							$\cF_{t^{m}_{k+1}}$-measurable $\Delta Z_{j}$-valued random
							variable with
							$g(\xi^{(m)}_{k,n}, \zhat_{j}) = \xhat_{k,n,m}$.
						\end{enumerate}
					\end{corollary}

					\begin{proof}
					Since	$\sigma(\xi^{(m)}_{k,n}, \xhat_{k,n,m}) \subset
						\cF_{t^{m}_{k+1}}$, everything follows from
						Lemma~\ref{lem:DH-meas}\ref{L:iii} \ref{L:v}; the measurability of the registration event and of the nearest-layer selector is Lemma~\ref{lem:selector-meas}\,(ii).
					\end{proof}

					\begin{remark}\label{rem:shift}
						It follows from Corollary \ref{cor:adapted}, that
						the random measures restricted to a time interval $[t_k^m,t_{k+1}^m)$ of Step~(VI) built from these
						ingredients are $\cF^{1,n}_{t^{m}_{k+1}}$-measurable, and the
						associated counting processes are adapted to the one-time-interval shifted
				filtration $\bigl(\cF_{(t+\tau_{m}) \wedge T}\bigr)_{t}$, $\tau_m=2^{-m}T$.
					\end{remark}

					\begin{remark}\label{rem:interfaces}
						By part \ref{L:ii} the sets $\Delta H_{j}(v)$ are disjoint, but
						their \emph{closures} may intersect along the images of the spheres
						$\{|z|_{Z} = \epsilon_{j}\}$.  Since the levels were chosen with
						$\nu(\{|z|_{Z} = \epsilon_{j}\}) = 0$.
					\end{remark}

					\item   \textbf{Definition of a random measure on $\mathfrak{A}$:}
					\noindent

		Recall from Corollary~\ref{cor:adapted} the recovered increment is given by
						\DEQSZ\label{def_xhat}
						\xhat_{k,n,m}:=\ell_{k,n,m}+\lbar^{(J)}_{k,m,n},
						\qquad
						\lbar^{(J)}_{k,m,n}:=\tau_m\sum_{\jt=0}^{J}\int_{\Delta Z_\jt} g\bigl(\xi^{(m)}_{k,n},z\bigr)\,\nu(dz),
						\EEQSZ
						where $\ell_{k,n,m}$ is the increment defined in \eqref{mart_incr}.
Now, $\xhat_{k,n,m}$ approximates the actual jump size of the noise in the time interval $I^m_k$, up to some  error terms.
					Next, let us define on $\mathfrak{A}$ the random measure
					$$	\hat{N}_{\mathbf{p}}^{k,n,m,j}
					:\mathscr{B}(H)\times  \mathscr{B}([0,T])\to \mathbb N_0
					$$
					for $B\in\mathscr B(H)$ and $ I\in\mathscr B([0,T])$ by
					\begin{equation}\label{eq:Nhat-delta}
						\hat{N}_{\mathbf{p}}^{k,n,m,j}	(B\times I)
						\;:=\;
						{\chi^{\jt,m}_k} \cdot \mathds{1}_{ \{  \Upsilon_{k,m,\jty}\in I\}}   \cdot   \mathds{1}_{ \{ \xhat_{k,n,m}\in B\}}
		.
					\end{equation}
 Here $\chi^{\jt,m}_k$ is the \emph{registration indicator} of the cell $I^m_k$ for the layer $\jt$. It implements the three rules of (A4)--(A5), the cut-off, the tolerance and the nearest-layer selection:
\begin{equation}\label{eq:chi}
\chi^{\jt,m}_k:=\mathds 1_{\{|\xi^{(m)}_{k,n}|_V\le R_r\}}\cdot\mathds 1_{\{d^{(J)}(\xi^{(m)}_{k,n},\xhat_{k,n,m})\le\delta_r\}}\cdot\mathds 1_{\{\hat\jmath(\xi^{(m)}_{k,n},\xhat_{k,n,m})=\jt\}}\in\{0,1\},
\end{equation}
with $d^{(J)}$ and $\hat\jmath$ as in (A3) and (A5). 
					To show that $\Omega \ni \omega\mapsto \hat{N}_{\mathbf{p}}^{k,n,m,j}$ is a random measure, we have to show that, for every fixed $B\in\mathscr B(H)$ and $I\in\mathscr B([0,T])$, the mapping $\Omega\ni\omega\mapsto \hat{N}_{\mathbf{p}}^{k,n,m,j}(\omega, B\times I)\in\mathbb N_0$ is measurable
					on $(\Omega,\CF)$. Note, for a fixed $\omega$, $\hat{N}_{\mathbf{p}}^{k,n,m,j}(\omega,\cdot)$ is an $\mathbb N_0$-valued measure by construction.
					\begin{proposition}%
						\label{prop:Nhat-meas}
						Under the assumption that $g\colon V\times Z\to H$ is
						$(\mathscr{B}(V)\otimes\mathscr{B}(Z),\mathscr{B}(H))$-measurable
						and $ \xi_{k,n}^{(m)}\colon\Omega\to V$ is $\CF$-measurable,
						the map
						\[
				\Omega \ni 		\omega
						\;\longmapsto\;
						\hat{N}_{\mathbf{p}}^{k,n,m,j}(\omega,\,B\times I)
						\]
						is $\CF$-measurable for every fixed
						$B\in\mathscr{B}(H)$ and $I\in\mathscr{B}(I_{k,m})$. Moreover,
						$\hat{N}_{\mathbf{p}}^{k,n,m,j}$ is an integer-valued
						random measure on $(H\times I_{k,m},\,\mathscr{B}(H)\otimes
						\mathscr{B}(I_{k,m}))$ in the sense of Kallenberg
						\textup{\cite[Ch.~1]{Kallenberg2017}}.
					\end{proposition}

						\begin{proof}
							By \eqref{def_xi}, \eqref{mart_incr}, and \eqref{def_xhat}, the random variable $\xhat_{k,n,m}$ is built from the values $u_n(t^m_k)$, $u_n(t^m_{k+1})$, the drift increment $D^m_{k,n}(u_n)$,  and the dyadic coefficient $\xi^{(m)}_{k,n}$.  Since $g$ is jointly measurable and the $\xi^{(m)}_{k,n}$ are
							$\mathscr{F}$-measurable, the map
							$\Omega\ni \omega\mapsto \xhat_{k,n,m}(\omega)$ is
							$\mathscr{F}/\mathscr{B}(H)$-measurable.
 By Lemma~\ref{lem:selector-meas}\,(i), applied with $\xi=\xi^{(m)}_{k,n}$ and $Y=\xhat_{k,n,m}$, the registration event $\{d^{(J)}\le\delta_r\}$ and the level sets of the selector $\hat\jmath$ belong to $\sigma(\xi^{(m)}_{k,n},\xhat_{k,n,m})\subset\CF $, and so does the cut-off event $\{|\xi^{(m)}_{k,n}|_V\le R_r\}$. Hence $\chi^{\jt,m}_k$ is $\CF $-measurable.
							For fixed $B\in\mathscr{B}(H)$, the map
							$\Omega\ni \omega\mapsto\mathds{1}_B(\xhat_{k,n,m}(\omega))$ is therefore
							$\mathscr{F}$-measurable, and the factor
							$\mathds{1}_I(\Upsilon_{k,m,\jty})$ is deterministic, so the product in \eqref{eq:Nhat-delta}
							is
							$\mathscr{F}$-measurable.  Since, for each $\omega$,
							$\hat{N}_{\mathbf{p}}^{k,n,m,j}(\omega,\,\cdot)$ is
							either the Dirac mass $\delta_{(\xhat_{k,n,m},\Upsilon_{k,m,\jty})}$ or the zero
							measure, it is an integer-valued measure,
							hence a random measure in the sense of
							\cite[Def.~1.1]{Kallenberg2017}.
						\end{proof}

						These random measures are extended to the whole time interval in the following way:
						For $B\in \mathscr{B}(H)$ and $I\in\mathscr{B}([0,T])$, we put
						\DEQSZ\label{defrandom}
					\hat{N}_{\mathbf{p}}^{n,m,j}	(B\times I)
						:= \sum_{k=0}^{2^m-1} \hat{N}_{\mathbf{p}}^{k,n,m,j}\lk( B
						\times \lk(I\cap I_{k,m}\rk)\rk).
						\EEQSZ
						Finally, we add up all the jumps over $\Delta Z_\jt$, $\jt =0,1,2,\ldots, J$. To be more precise, for
						$B\in \mathscr{B}(H)$ and $I\in\mathscr{B}([0,T])$  we set
						\DEQSZ\label{defrandomJ}
						\hat{N}_{\mathbf{p}}^{n,m,J}	(B\times I)
						:= \sum_{\jt=0}^{J} 	\hat{N}_{\mathbf{p}}^{n,m,j}	(B\times I).
						\EEQSZ
						Since $\sum_{\jt=0}^J\chi^{\jt,m}_k\le1$ by \eqref{eq:chi}, each time interval $I^m_k$ carries at most one atom of $\hat N^{n,m,J}_{\mathbf p}$, and $	\hat{N}_{\mathbf{p}}^{n,m,J}	$ is again an $\mathbb N_0$-valued random measure on $\mathscr B(H)\otimes\mathscr B([0,T])$.

				Next, let us introduce a filtration on $\MA$ generated by \(\myp_{n}\) as an element of $\mathcal X=L^\alpha(0,T;V)\cap\mathbb \DD([0,T];H)$, and a discrete filtration generated by the dyadic coefficients and the values of $u_n$ at the dyadic times,
						given by  $(\mathscr F_{t}^{n})_{t\in[0,T]}$ and $( \mathscr F_{t}^{n,m})_{t\in[0,T]}$, where for $t\in[0,T]$ we define
\begin{align}\label{def_dis_n_filtration}
							\mathscr F_{t}^{n}
					&	:=
						\sigma\bigl(\{\myp_{n}(s):0\le s\le t\}\bigr)^{\PP\text{-comp}}\quad \mbox{and}\\ \quad \mathscr F_{t}^{n,m}
						&:=
						\sigma\bigl(\{\xi^{(m)}_{j,n},\ u_n(t^m_j),\ {D^m_{j-1,n}(u_n)}:\ 0\le j\le 2^m,\ t^m_j\le t\}\bigr)^{\PP\text{-comp}}.
					\end{align}
							\begin{remark}
							In fact, 		$(\mathscr F_{t}^{ n,m})_{t\in[0,T]}$ is a discrete filtration and for all $t\in[0,T]$,
							$$\mathscr F^{n,m}_t=\mathscr F^{ n,m}_{k^m(t)},
							$$
							where  $k^m(t):=\sup \{ s\in \pi_T^m: s\le t\}$.
						\end{remark}
						The discrete filtration $\mathscr F^{n,m}$ is generated by the dyadic
						coefficients and by the values of $u_n$ at the dyadic times that have
						already passed, and by the drift increments $D^m_{j-1,n}(u_n)=\int_{I^m_{j-1}}A_n(s,u_n(s))\,ds$; these are $H$-valued and $\mathscr F^n_{t^m_j}$-measurable. Since $\xi^{(m)}_{j,n}$ is the average of $u_n$ over
						$I^m_{j-1}$, we have
						$\mathscr F^{n,m}_t\subset\mathscr F^n_t\subset\mathscr F_t$. In
						particular:
						\begin{enumerate}[label=\upshape(\roman*)]
							\item $\xi^{(m)}_{k,n}$ is $\mathscr F^{n,m}_{t^m_k}$-measurable, so
							$u_{n,m}$ is $(\mathscr F^{n,m}_t)$-predictable;
							\item the computable increment $\ell_{k,n,m}=u_n(t^m_{k+1})-u_n(t^m_k)-D^m_{k,n}(u_n)$ of \eqref{mart_incr}, and hence $\xhat_{k,n,m}$, are
							$\mathscr F^{n,m}_{t^m_{k+1}}$-measurable;
							\item for $U\in\mathscr Z$ with $\nu(U)<\infty$, the increment
							$\eta(U\times I^m_k)$ is independent of $\mathscr F^{n,m}_{t^m_k}$,
							since $\mathscr F^{n,m}_{t^m_k}\subset\mathscr F_{t^m_k}$.
						\end{enumerate}

			The next proposition records the measurability of the reconstructed random measures with respect to the discrete filtration.

						\begin{proposition}\label{prop:measurable}
							For any $m,n\in\NN$,
							\begin{enumerate}
								\item[\textup{(i)}] for all
								$B \in \cB(H)$ and $I \in \cB([0,T])$ with $I\subset I_{k,m}$, the random variable $\hat{N}_{\mathbf{p}}^{n,m,j}(B \times I)$ is
								{ $\mathscr F^{n,m}_{t^m_{k+1}}$}-measurable; for general $I\in\cB([0,T])$, $\hat{N}_{\mathbf{p}}^{n,m,j}(B\times I)$ and $\hat{N}_{\mathbf{p}}^{n,m,J}(B\times I)$ are { $\mathscr F^{n,m}_T$}-measurable;
								\item[\textup{(ii)}] $\hat{N}_{\mathbf{p}}^{n,m,Jx}$ is an
								integer-valued random measure in the sense of
								\cite[Ch.~1]{Kallenberg2017};
							\end{enumerate}
						\end{proposition}
						\begin{proof}
				Observe that
							$$\ell_{k,n,m}=u_n(t^m_{k+1})-u_n(t^m_k)-D^m_{k,n}(u_n)
							$$
is built from $u_n(t^m_k)$, $u_n(t^m_{k+1})$ and the drift increment $D^m_{k,n}(u_n)$, and is
							therefore $\mathscr F^{n,m}_{t^m_{k+1}}$-measurable by
							\eqref{def_dis_n_filtration}. It follows  by
							\eqref{def_xi} and \eqref{eq:g-int} that  $\lbar^{(J)}_{k,m,n}$ is
							$\sigma(\xi^{(m)}_{k,n})\subset\mathscr F^{n,m}_{t^m_k}$-measurable. Hence $\xhat_{k,n,m}$ is
							$\mathscr F^{n,m}_{t^m_{k+1}}$-measurable, and by
							Lemma~\ref{lem:DH-meas}, items \ref{L:iii}-\ref{L:iv} so are $d_{\jt}$, $d^{(J)}$, the
							registration event and the selector; the cut-off event $\{|\xi^{(m)}_{k,n}|_V\le R_r\}$ is even $\mathscr F^{n,m}_{t^m_k}$-measurable, so $\chi^{\jt,m}_k$ is $\mathscr F^{n,m}_{t^m_{k+1}}$-measurable.  Since
							$\ind_{I}(\sigma^{m}_{k})$ is deterministic, (i) follows.  For
							(ii), the mass assigned to $B \times (0,t]$ involves only time intervals with
							{$\sigma^{m}_{k} \le t$, hence $t^{m}_{k+1} \le t + \tau_{m}$,}  with $\tau_m={2^{-m}}T$,  and a
							finite sum of $\NN_0$-valued random measures with measurable masses is
							an integer-valued random measure.
						\end{proof}

\renewcommand{\n}{n}

						\item
						\textbf{Construction of the inverse of $g$.}

						Throughout this step, $m,n\in\NN$, $J\in\NN$ and the layers $\Delta Z_\jt$, $\jt\in\{0,\dots,J\}$, of \eqref{eq:layers} are fixed.

						\medskip
						Recall from Step~(I) (cf.\ \eqref{def_xi}) that the shifted dyadic interpolant $u_{\n,m}:=\mathrm{Proj}_m u_{n}$ is the $V$-valued, piecewise constant process
						\[
u_{\n,m}(s)=\xi^{(m)}_{k,\n},\qquad s\in I^m_k=(t^m_k,t^m_{k+1}],\quad k=0,\dots,2^m-1,
\]
 where
\[
{\xi^{(m)}_{0,\n}=0,\qquad\xi^{(m)}_{k,\n}=\frac{2^m}{T}\int_{t^m_{k-1}}^{t^m_k}u_\n(r)\,dr,\quad k\ge1,}
						\]
						and that $\xi^{(m)}_{k,\n}$ is $\mathscr F^{ \n,m}_{t^m_k}$-measurable. Let us  define the \emph{frozen jump field}
						\[
						\Theta^{m,\n}(s,z):= g\bigl(u_{\n,m}(s),z\bigr),\qquad (s,z)\in[0,T]\times Z .
						\]

						Since $g:V\times Z\to H$ is jointly measurable (Hypothesis~\ref{ass:GF}) and $(s,\omega)\mapsto u_{\n,m}(s,\omega)$ is $\mathscr B([0,T])\otimes\mathscr F$-measurable (piecewise constant in time with $\mathscr F$-measurable coefficients), the map
						\[
						\Theta^{m,\n}:[0,T]\times Z\times\Omega \to H
						\]
						is jointly measurable. Moreover, on each time interval $I^m_k$ the field depends on $\omega$ only through the frozen coefficient $\xi^{(m)}_{k,\n}$, which is $\mathscr F^{ \n,m}_{t^m_k}$-measurable. Hence, for every fixed $z\in Z$, the process $(\omega,s)\mapsto\Theta^{m,\n}(s,z)$ is $\mathscr P^{\n}_m$-measurable, i.e.\ predictable with respect to the filtration  $(\mathscr F^{ \n,m}_t)_{t\in[0,T]}$.

						\medskip
						For $v\in V$ and $\jt\in\{0,\dots,J\}$ define the finite measure\footnote{Here, for a function $f$ $f_{\#}\mu:=\mu\circ f^{-1}$ denotes the image (pushforward) measure, see \cite{Federer1969}.}
						\DEQSZ\label{eq:mu-layer}
						\mu_\jt(v;B):=\int_{\Delta Z_\jt}\mathds 1_B\bigl(\underbrace{g(v,z)}_{\in H}\bigr)\,\nu(dz)
						=\Bigl({g(v,\cdot)}_{\#}\bigl(\nu|_{\Delta Z_\jt}\bigr)\Bigr)(B),
						\qquad B\in\mathscr B(H).
						\EEQSZ
						{ Since $g:V\times Z\to H$ is jointly measurable (Hypothesis~\ref{ass:GF}) and $\nu(\Delta Z_\jt)<\infty$, Tonelli's theorem shows that $V\ni v\mapsto\mu_\jt(v;B)$ is $\mathscr B(V)$-measurable for every $B\in\mathscr B(H)$. Hence $\mu_\jt$ is a finite Borel kernel from $(V,\mathscr B(V))$ to $(H,\mathscr B(H))$ with total mass $\mu_\jt(v;H)=\nu(\Delta Z_\jt)$, concentrated on the admissible jump set $\Delta H_\jt(v)=g(v,\Delta Z_\jt)$ of Lemma~\ref{lem:DH-meas} \ref{L:ii}.}
						We now substitute the frozen state. For $s\in I^m_k$, $k=0,\dots,2^m-1$, and $x\in H$ set
						\[
						\zhat^{m,\n,\jt}(s,x)(\omega):=
						\begin{cases}
							\zhat_\jt\bigl(\xi^{(m)}_{k,\n}(\omega),x\bigr), & x\in\Delta H_\jt\bigl(\xi^{(m)}_{k,\n}(\omega)\bigr),\\[2pt]
							z_0, & \text{otherwise.}
						\end{cases}
						\]
						By Lemma~\ref{lem:DH-meas} \ref{L:i}, \ref{L:iii} and \ref{L:v} (Lusin--Suslin) and the $\mathscr F^{ \n,m}_{t^m_k}$-measurability of $\xi^{(m)}_{k,\n}$, the map $(\omega,s,x)\mapsto\zhat^{m,\n,\jt}(s,x)(\omega)$ is $\mathscr P^{\n}_m\otimes\mathscr B(H)$-measurable.  In particular, for $s\in I^m_k$ it is $\mathscr F^{ \n,m}_{t^m_k}\otimes\mathscr B(H)$-measurable. This predictability is the nonanticipativity needed later on. {Recovered increments $\xhat_{k,n,m}$ are in general only $\delta_r$-close to the layer; outside the membership set,}  the tolerance-adapted inverse of step (A7) is used instead; its measurability is part of (A7).

						\medskip
						Since $\zhat_\jt(v,g(v,z))=z$ for every $z\in\Delta Z_\jt$ (Lemma~\ref{lem:DH-meas} \ref{L:v}), pushing the layer image kernel \eqref{eq:mu-layer} back to $Z$ under the measurable inverse returns the restricted L\'evy measure: for every $v\in V$ and every $\jt\in\{0,\dots,J\}$,
						\DEQSZ\label{eq:Q-total-mass}
						\bigl(\zhat_\jt(v,\cdot)\bigr)_{\#}\,\mu_\jt(v;\cdot)&=&\nu|_{\Delta Z_\jt}.
						\EEQSZ
						Applied with $v=\xi^{(m)}_{k,\n}(\omega)$,  \eqref{eq:Q-total-mass} says that, at the frozen state, the inverse map carries the layer image of $\nu$ back to $\nu|_{\Delta Z_\jt}$.

						\begin{remark}\label{rem:kallenberg-disint}
							If the injectivity assumption of Hypothesis~\ref{ass:GF}, item~\ref{it:inj}, is dropped, a conditional inverse still exists, but only through a  disintegration (done by Ikeda and Watanabe \cite{IkedaWatanabe1989}): for each fixed $v\in V$ one disintegrates the finite measure
							\[
							\rho_\jt(v;\cdot):=\bigl(\nu\big|
_{\Delta Z_\jt}\bigr)\circ\bigl(z\mapsto(g(v,z),z)\bigr)^{-1}
							\qquad\text{on }H\times Z
							\]
							over its first marginal $\mu_\jt(v;\cdot)$. Since $H$ and $Z$ are Polish, hence standard Borel, the disintegration theorem for kernels on standard Borel spaces {(see \cite[Lemma~3.22]{kallenberg2021foundations})} yields a probability kernel $Q_\jt(v;x,dz)$, unique for $\mu_\jt(v;\cdot)$-a.e.\ $x$, satisfying the  identity 
$$\int_{\Delta Z_\jt}\varphi\bigl(g(v,z),z\bigr)\,\nu(dz)=\int_H\int_Z\varphi(x,z)\,Q_\jt(v;x,dz)\,\mu_\jt(v;dx).
$$
 A jointly measurable version in $(v,x)$ then has to be produced by a measurable-selection argument. Under item~\ref{it:inj} all of this is unnecessary: the measurable inverse $\zhat_\jt$ of Lemma~\ref{lem:DH-meas} \ref{L:v} provides the explicit (Dirac) inverse, jointly Borel in $(v,x)$, which is unique.
						\end{remark}

						\newcommand{\mykappa}{\kappa}

Under Hypothesis~\ref{ass:GF}, item~\ref{it:inj}, the generalised
inverse is  chosen \emph{explicitly}:
with the measurable inverse $\zhat_\jt$ of Lemma~\ref{lem:DH-meas} \ref{L:v} we take, for $t\in I^m_k$,
\DEQSZ\label{eq:f-explicit}
f^{n,m,j }(t,x,\omega)=
\begin{cases}
	\zhat_j\bigl(\xi^{(m)}_{k,n}(\omega),\,x\bigr), & x\in\Delta H^{k,m,n}_j(\omega),\\[2pt]
	0^\ast, & \text{otherwise.}
\end{cases}
\EEQSZ

						\item \textbf{Explicit construction of the Poisson random
							measure on $Z$.}
						In this step we construct an approximation of a Poisson random
						measure on $(Z,\mathscr{Z})$ with L\'evy measure $\nu$, by
						inverting the map $g$ via its explicit measurable inverse (Lemma~\ref{lem:DH-meas} \ref{L:v}).

	\medskip

						\noindent\textbf{The registration data at level $m$.}
						By the single-atom convention, each time interval $I^m_k$
carries at most one registered jump. 
We recall the registration indicator $\chi^{\jt,m}_k:\Omega\to\{0,1\}$ of $I^m_k$ from \eqref{eq:chi}:
\[
\chi^{\jt,m}_k=\mathds 1_{\{|\xi^{(m)}_{k,n}|_V\le R_r\}}\cdot\mathds 1_{\mbox{\rm r}^{\delta_r}_{k,m,n}}\cdot\mathds 1_{\{\hat\jmath_{k,m,n}=\jt\}}\in\{0,1\}.
\]
Here
\DEQSZ\label{registerevent}
\mbox{\rm r}^{\delta_r}_{k,m,n}
:=\mbox{\rm r}^{\delta_r}\bigl(\xi^{(m)}_{k,n},\xhat_{k,n,m}\bigr)
=\bigl\{\min_{0\le \jt\le J}
d_\jt\bigl(\xi^{(m)}_{k,n},\xhat_{k,n,m}\bigr)\le\delta_r\bigr\}
\EEQSZ
denotes the registration event of Definition~\ref{def:selector}: the recovered increment $\xhat_{k,n,m}$ of the interval $I^m_k$ lies within $H$-distance $\delta_r$ of at least one admissible jump set $\Delta H_\jt\bigl(\xi^{(m)}_{k,n}\bigr)$. The cut-off discards the cells whose frozen state lies outside the $V$-ball of radius $R_r$, on which the lower bound \eqref{it:nondeg} is uniform, and the nearest-layer selector $\hat\jmath_{k,m,n}$ attributes the jump to the layer realising the minimal distance.

The corresponding considerations for $u_n$ are  given in
						Definition~\ref{def:selector}\footnote{Here $\hat\jmath_{k,m,n}
:=\hat\jmath\bigl(\xi^{(m)}_{k,n},\xhat_{k,n,m}\bigr)$
denotes the nearest-layer selector of Definition~\ref{def:selector}: on the registration event ${\mbox{\rm r}^{\delta_r}_{k,m,n}}$ it equals
the smallest index $\jt\in\{0,\dots,J\}$ realising
$\min_{0\le j\le J}d_j\bigl(\xi^{(m)}_{k,n},\xhat_{k,n,m}\bigr)$,
i.e.\ it gives back the registered mark of the interval $I^m_k$ to the layer
whose admissible jump set $\Delta H_\jt\bigl(\xi^{(m)}_{k,n}\bigr)$ is
closest in $H$ to the recovered increment.
};
the associated jump time
						$\Upsilon_{k,m,\jty}=\tfrac12(t^m_k+t^m_{k+1})$ is the \emph{deterministic}
						midpoint, and the recovered jump size is $\xhat_{k,n,m}$. By Corollary~\ref{cor:adapted} (Lemma~\ref{lem:selector-meas}
						applied with $\xi=\xi^{(m)}_{k,n}$ and $Y=\xhat_{k,n,m}$), the
						variables $\chi^{\jt,m}_k$ and $\xhat_{k,n,m}$ are
						$\mathscr F^{n,m}_{t^m_{k+1}}$-measurable, i.e.\ the registration of the
						interval $I^m_k$ is decided at its {right} endpoint. The number of
						registered layer-$\jt$ jumps at level $m$ over $\mathfrak{A}$ is
					\DEQSZ\label{defrandom2}
						\iota^{\jt,m}_{n}
						:=\hat N^{n,m,\jt}_{\up}\bigl(H\times[0,T]\bigr)
						=\sum_{k=0}^{2^m-1}\chi^{\jt,m}_k\;\le\;2^m .
				\EEQSZ

						\medskip

						In the next step, we construct an approximation of a Poisson random measure on $(Z, \MscrZ )$ having L\'evy measure $\nu$.
					The map $f^{n,m,\jt}$ defined in \eqref{eq:f-explicit} is jointly measurable in $(t,x,\omega)$ and predictable in $(t,\omega)$ by the inverse-map paragraph of Step~(VI).

						\medskip
					Let us now construct the  random measure on $Z\times[0,T]$.
					The $Z$-valued measure
						$$\hat{N}^{n,m,\jmath}_{\uq}
						:\mathscr{Z}\times\mathscr{B}([0,T])\to\mathbb{N}_0$$
						is defined  from the registration data of the
						algorithm: its atoms are the recovered jumps in $Z$ of step (A7),
						placed at the deterministic times $\Upsilon_{k,m,\jty}$.
						For all $C\in\mathscr{Z}$ and $I\in\mathscr{B}([0,T])$ we set
		\begin{equation}\label{zuvor}
							\hat{N}^{n,m,\jmath}_{\uq} (C\times I)
							:=\sum_{k=0}^{2^m-1}\chi^{\jmath,m}_k\,
							\mathds{1}_C\bigl(\hat z_{k,m,\jmath}\bigr)\,
							\mathds{1}_I\bigl(\Upsilon_{k,m,\jty}\bigr),
						\end{equation}
						where $\chi^{\jmath,m}_k$ is the registration indicator of
						\eqref{defrandom2} and $\hat z_{k,m,\jmath}\in\Delta Z_\jmath$
						is the mark in $Z$ recovered in step (A7). In particular, it is  the exact inverse
						$\zhat_\jmath\bigl(\xi^{(m)}_{k,n},\xhat_{k,n,m}\bigr)$
						in case of the exact-membership event, or otherwise, it is the tolerance-adapted inverse.
						 Each summand is measurable by
						Corollary~\ref{cor:adapted} and the measurable selection in
						(A7), so \eqref{zuvor} is an $\NN_0$-valued random measure.
						On the exact-membership events, \eqref{zuvor} coincides with the
						pushforward of $	\hat{N}^{n,m,\jmath}_{\uq}$ under the map
						$f^{n,m,\jmath}$ of \eqref{eq:f-explicit}; the direct
						definition is used because the recovered increments are in
						general only {$\delta_r$}-close to the layer, where the
						pushforward under $f^{n,m,\jmath}$ would return the
						default value $z_0$.

						\medskip
						\noindent Let us take the sum  over $\jmath=0,\ldots,J$.
						Since the sets $\Delta Z_\jmath$ are disjoint by construction,
						the sets $\Delta H^{k,m,n}_\jmath
						= \{g(\xi^{(m)}_{k,n},z):z\in\Delta Z_\jmath\}$ are
						also disjoint for different $\jmath$, and the sum
					\begin{equation}\label{nochmals_zuvor}
							\hat N^{n,m,J}_{\uq}(C\times I)
							:=\sum_{\jmath=0}^{J}
							\hat{N}^{n,m,\jmath}_{\uq}(C\times I),
							\quad\text{with the layer measures of \eqref{zuvor},}
						\end{equation}
						is again a well-defined $\mathbb{N}_0$-valued random measure on
						$\mathscr{Z}\times\mathscr{B}([0,T])$ over $\mathfrak{A}$,
						counting jumps in $Z$ originating from each annular region
						$\Delta Z_\jmath$, $\jmath=0,\ldots,J$.

\bigskip

\item\label{step:diagonal}\textbf{Fixing the scaling: $n$, $m$
	and $\delta$ as functions of one parameter $\seque$.}

Up to now the index $n$ of the approximating equation, the
dyadic level $m$ and the
tolerance $\delta$ of the registration algorithm were free. In this step
all of them are fixed as functions of \emph{one} parameter, and we show
that, outside events of small probability, the registration algorithm of
Step~(IV) recovers exactly the jumps of the driving noise with marks in
$Z^c_{\epsilon_{J+1}}$, for every truncation level $J$.

\medskip

\medskip
\noindent\textbf{(a) The exact decomposition of the error.}
Let us first analyse the exact decomposition of the error.
Fix $J\in\NN$ and parameters $n,m\in\NN$ and
$\varepsilon_\flat\in(0,\epsilon_{J+1}]$. For brevity we write $\xhat_k:=\xhat_{k,n,m}$.
By \eqref{mart_incr},
\begin{equation}\label{eq:ell-new}
	\ell_{k,n,m}
	=u_n(t^m_{k+1})-u_n(t^m_k)-\int_{{I^m_k}}A_n\bigl(s,u_n(s)\bigr)\,ds .
\end{equation}
We split the mark space into the \emph{tiny} marks
$Z_{\varepsilon_\flat}=\{|z|_Z\le\varepsilon_\flat\}$, the \emph{band}
$B:=\{\varepsilon_\flat<|z|_Z\le\epsilon_{J+1}\}$ and the \emph{large}
marks $Z^c_{\epsilon_{J+1}}$, and define, for $0\le k<2^m$,
the error of the small jumps $T_k^{m,n}$ given by
\begin{align}\label{eq:errors_tiny}
T_k^{m,n}&:=\int_{{I^m_k}}\!\int_{Z_{\varepsilon_\flat}}
	g\bigl(u_n(s),z\bigr)\,\widetilde\eta(dz,ds)
\end{align}
the error of the compensators on $I_k^m$, i.e.  
$$
\int_{I^m_k}\int_{Z^c_{\varepsilon_\flat}}g(u_n(s),z)\,\nu(dz)\,ds
$$
and the calculated compensator $	\lbar^{(J)}_{k,m,n}$ given in \eqref{def_xhat} (compare also \eqref{def:jump})
given by
\begin{align}
	\label{eq:errors_comp}
	C_k^{m,n}&:=-\int_{{I^m_k}}\!\int_{B}g\bigl(u_n(s),z\bigr)\,\nu(dz)\,ds
	-\int_{{I^m_k}}\!\int_{Z^c_{\epsilon_{J+1}}}
	\bigl[g\bigl(u_n(s),z\bigr)-g(\xi^{(m)}_{k,n},z)\bigr]\,\nu(dz)\,ds,
\end{align}
and, for an atom $(s,z)$ of $\eta$ with $s\in {I^m_k}$, the
	error arising by averaging the coefficients is given by
\begin{align}\label{eq:errors_freezing}
	F(s,z)&:=g(u_n(s-),z)-g(\xi^{(m)}_{k,n},z).
\end{align}
Here we write $s-$: by Step~0 of Lemma~\ref{lem:jump-identity}, $\PP$-a.s.\
$u_n(s-)\in V$ at every atom $(s,z)$ of $\eta$. The left limits and the choice
of the integrand are explained in the next remark.

\begin{remark}\label{rem:left-limits}
	The stochastic integral with respect to $\widetilde\eta$ is defined on
	$\mathcal M^p(0,T;L^p(Z,\nu;H))$ (Theorem~\ref{th-ineq}), i.e.\ for
	equivalence classes of integrands modulo $\PP\otimes\Leb\otimes\nu$. We
	therefore write all compensated integrals, such as $T_k^{n,m}$, $S_\flat$ and the
	noise term of \eqref{originalun}, with $g(u_n(s),z)$. The left limit enters
	only where the integrand is \emph{evaluated at individual atoms}: in sums
	over the atoms of $\eta$, in the averaging error $F(s,z)$, and in the jump
	identity $\Delta u_n(s)=g(u_n(s-),z)$. There one uses the predictable
	representative $g(u_n(s-),z)\mathds1_{\{u_n(s-)\in V\}}$ of $g(u_n(s),z)$,
	which belongs to the same class since $u_n(s-)=u_n(s)$ off the at most
	countable set of jump times (Step~1 of Lemma~\ref{lem:jump-identity}). The
	path $u_n$ is c\`adl\`ag in $H$ but not in $V$; by Step~0 of
	Lemma~\ref{lem:jump-identity} and
	Lemma~\ref{lem:pathwise-atoms}\,\ref{PA:null}, $\PP$-a.s.\ $u_n(s-)\in V$
	at every atom $(s,z)$ of $\eta$, including the infinitely many small ones,
	so the indicator equals $1$ there. At an atom $(s,z)$ the value $u_n(s)$
	already contains the jump $g(u_n(s-),z)$, which is why the left limit, and
	not $u_n(s)$, is evaluated there.
\end{remark}

In the next Lemma  we tackle the exact decomposition of the recovered increment.

\begin{lemma}
	\label{lem:exact-decomp}
	$\PP$-a.s., for every $0\le k<2^m$,
	\begin{align}
		\xhat_{k,n,m}&=\sum_{\substack{(s,z)\in\operatorname{supp}\eta\\ s\in {I^m_k},\ |z|_Z>\varepsilon_\flat}}
		\bigl[g(\xi^{(m)}_{k,n},z)+F(s,z)\bigr]
		+T_k^{m,n}+C_k^{m,n} .\label{eq:xhat-exact}
	\end{align}
\end{lemma}
\begin{proof}

By \eqref{eq:ell-new} and the increment identity \eqref{mymn}, $\ell_{k,n,m}=\int_{I^m_k}\int_Zg(u_n(s),z)\,\widetilde\eta(dz,ds)$, $\PP$-a.s. On
	$Z^c_{\varepsilon_\flat}$ the intensity is finite, so there the
	compensated integral is the $\eta$-integral minus the compensator. This
	splitting evaluates the integrand at the atoms and therefore uses the
	predictable representative $\tilde\xi(s,z)=g(u_n(s-),z)\mathds1_{\{u_n(s-)\in V\}}$
	of $g(u_n(s),z)$ (Lemma~\ref{lem:pathwise-atoms}\,\ref{PA:pathwise} and
	Step~1 of Lemma~\ref{lem:jump-identity}); the indicator equals $1$ at the
	atoms by Step~0 there. Since $u_n(s-)=u_n(s)$ for a.e.\ $s$, the
	compensator may be written with $u_n(s)$, and
	\[
	\int_{{I^m_k}}\!\int_Z g(u_n(s),z)\,\widetilde\eta(dz,ds)
	=T_k^{m,n}+\sum_{\substack{s\in {I^m_k}\\ |z|_Z>\varepsilon_\flat}}g(u_n(s-),z)
	-\int_{{I^m_k}}\!\int_{Z^c_{\varepsilon_\flat}}g(u_n(s),z)\,\nu(dz)\,ds .
	\]
The sum can be written as a finite sum, since there exists only a finite number of jumps with $|z|_Z>\varepsilon_\flat$.
	Adding $\lbar^{(J)}_{k,m,n}=\int_{{I^m_k}}\int_{Z^c_{\epsilon_{J+1}}}g(\xi^{(m)}_{k,n},z)\,\nu(dz)\,ds$
	(Remark~\ref{rem:support}) and writing
	$g(u_n(s-),z)=g(\xi^{(m)}_{k,n},z)+F(s,z)$ gives \eqref{eq:xhat-exact}.
\end{proof}

\medskip
\noindent\textbf{(b) The total error.}
For fixed parameters $n,m\in\NN$, $J\in\NN$ and
$\varepsilon_\flat\in(0,\epsilon_{J+1}]$ we collect all errors in one random variable,
\begin{align}\label{eq:e-total}
	e_J&:=2\sup_{t\in[0,T]}\bigl|S_\flat(t)\bigr|_H+\max_{k}|C_k^{m,n}|_H
	+\sum_{\substack{(s,z)\in\operatorname{supp}\eta\\ |z|_Z>\varepsilon_\flat}}
	|F(s,z)|_H,
\end{align}
with
\begin{align}
	S_\flat(t)&:=\int_0^t\!\!\int_{Z_{\varepsilon_\flat}}g(u_n(s),z)\,\widetilde\eta(dz,ds).
\end{align}
Here we took into account that  $|T_k^{m,n}|_H\le2\sup_t|S_\flat(t)|_H$ and, by \eqref{eq:xhat-exact},
in every cell containing at most one atom with mark in
$Z^c_{\varepsilon_\flat}$,
\begin{equation}\label{eq:xhat-close}
	\Bigl|\xhat_{k,n,m}-\mathbf 1_{\{\text{atom }(s,z)\text{ in }I^m_k\text{ with }|z|_Z>\varepsilon_\flat\}}\,g(\xi^{(m)}_{k,n},z)\Bigr|_H\le e_J .
\end{equation}

\medskip
\noindent\textbf{(c) Bounds for fixed parameters.} We now bound each error
term of the decomposition, for fixed parameters.

Throughout we use the constants
\[
C_\alpha:=\sup_{n\in\NN}\EE\int_0^T|u_n(s)|_V^\alpha\,ds,
\qquad
C_2:=\sup_{n\in\NN}\EE\int_0^T\bigl(1+|u_n(s)|_V^2\bigr)\,ds,
\qquad
C_1:=\sup_{n\in\NN}\EE\int_0^T|u_n(s)|_V\,ds,
\]
which are finite by the a~priori estimate and $\alpha\ge2$. Next, for
$\epsilon>0$, the masses $\nu(Z^c_\epsilon)<\infty$ and
$\nu_J:=\nu(Z^c_{\epsilon_{J+1}})$ are also finite. For fixed $n$ we abbreviate the
projection error
\[
\Theta_n(m):=\EE\bigl\|u_n-\myp_{n,m}\bigr\|_{L^1(0,T;V)} .
\]

By the Cauchy--Schwarz inequality in $\nu$ and
\eqref{eq:g-second-moment2} with $p=2$, for every $\epsilon>0$,
\begin{equation}\label{eq:g-L1nu}
	\int_{Z^c_{\epsilon}}\bigl|g(v_1,z)-g(v_2,z)\bigr|_H\,\nu(dz)
	\le\bigl(\nu(Z^c_\epsilon)L_g\bigr)^{1/2}\,|v_1-v_2|_V,\qquad v_1,v_2\in V.
\end{equation}

First, let us calculate an error bound for fixed parameters.
\begin{lemma}\label{lem:error-bounds}
	Let $J\in\NN$, $n,m\in\NN$ with $\tau_m\le1$ and
	$\varepsilon_\flat\in(0,\epsilon_{J+1}]$, and put
	$\nu_\flat:=\nu(Z^c_{\varepsilon_\flat})$. Then
	\begin{align}\label{eq:Ee-fixed}
		\EE\,e_J\;\le\;&4\bigl(\kappa_2(\varepsilon_\flat)C_2\bigr)^{1/2}
		+\bigl(\nu_\flat C_gC_2\tau_m\bigr)^{1/2}
		+2\bigl(\nu_\flat L_g\bigr)^{1/2}\Theta_n(m).
	\end{align}
\end{lemma}

\begin{proof}
	\emph{The small jumps:}  The accumulated compensated sum of the small jumps $S_\flat$ is a square integrable $H$-valued
	martingale. By Doob's inequality, the It\^o isometry and
	\eqref{eq:XV-tail} with $p=2$, we know that
	$$\EE\sup_t|S_\flat(t)|_H^2\le4\,\EE\int_0^T\int_{Z_{\varepsilon_\flat}}
	|g(u_n(s),z)|_H^2\,\nu(dz)\,ds\le4\kappa_2(\varepsilon_\flat)C_2.
	$$
	Hence,
	$\EE\,2\sup_t|S_\flat(t)|_H\le4(\kappa_2(\varepsilon_\flat)C_2)^{1/2}$.
\medskip

	\emph{The Compensators.} Since $\nu(B)\le\nu_\flat$, the Cauchy--Schwarz
	inequality in $\nu$ and \eqref{eq:g-second-moment1} give
	$$\int_B|g(u_n(s),z)|_H\,\nu(dz)\le(\nu_\flat C_g)^{1/2}(1+|u_n(s)|_V^2)^{1/2},
$$
	and H\"older's inequality on $I^{m}_k$ gives
	$$\max_k\int_{I^m_k}(1+|u_n|_V^2)^{1/2}ds\le\tau_m^{1/2}
	\bigl(\int_0^T(1+|u_n|_V^2)ds\bigr)^{1/2}.
$$
Since
	$\myp_{n,m}=\xi^{(m)}_{k,n}$ on $I^{n,m}_k$, the second part of $C^{n,m}_k$ is
	bounded, uniformly in $k$, by \eqref{eq:g-L1nu} with $\epsilon=\epsilon_{J+1}$. Hence
	$\EE\max_k|C^{n,m}_k|_H\le(\nu_\flat C_gC_2\tau_m)^{1/2}
	+(\nu_\flat L_g)^{1/2}\Theta_n(m)$, since $\nu_J\le\nu_\flat$.

\medskip

	\emph{Averaged coefficients:} The integrand
	$$(s,z)\mapsto|g(u_n(s-),z)-g(\myp_{n,m}(s),z)|_H\mathbf 1_{Z^c_{\varepsilon_\flat}}(z)
	$$
	is nonnegative and predictable ($\myp_{n,m}$ is left-continuous and
	adapted, and $\myp_{n,m}(s)=\xi^{(m)}_{k,n}$ for $s\in I^m_k$). The compensation
	formula, $u_n(s-)=u_n(s)$ for a.e.\ $s$, and \eqref{eq:g-L1nu} with
	$\epsilon=\varepsilon_\flat$ give
	$$\EE\sum|F(s,z)|_H\le(\nu_\flat L_g)^{1/2}\Theta_n(m).
	$$

	Summing up the three bounds gives \eqref{eq:Ee-fixed}.
\end{proof}

\begin{remark}\label{rem:projection-fixed-n}
	For every fixed $n$, $\Theta_n(m)\to0$ as $m\to\infty$. Indeed,
	$\myp_{n,m}\to u_n$ in $L^1(0,T;V)$ $\PP$-a.s.\ by Step~(I) and Lemma~\ref{lem:Pro-L1}, and
	$\|u_n-\myp_{n,m}\|_{L^1(0,T;V)}\le2\|u_n\|_{L^1(0,T;V)}+\tau_m|u^m_0|_V$,
	where $\tau_m|u^m_0|_V{\to0}$ by
	\eqref{init0}. Hence $\Theta_n(m)\to0$ by dominated convergence.
\end{remark}

\medskip
\noindent\textbf{(d) Choice of the parameters.}
As announced at the beginning, in the next step we pass to the limit along a
single parameter $r$, and all errors have to converge to zero as
$r\to\infty$. The errors depend on several parameters, among them $n$
and $m$, and these cannot be chosen independently: some error
terms tend to zero in $m$ only for fixed $n$. For instance, the
projection error $\Theta_n(m)$ tends to zero as $m\to\infty$ only for fixed
$n$ (Remark~\ref{rem:projection-fixed-n}), so $m$ has to be chosen after $n$;
Therefore the order in which the parameters are chosen matters, and each
parameter is chosen depending on the ones chosen before it. A possible choice
is given in the next lemma.

To make all parameters functions of the single index $r$, we fix a sequence $(\seque)_{r\in\NN}\subset(0,1]$ with
\begin{equation}\label{eq:seque}
	\seque\downarrow0\quad (r\to\infty),
	\qquad
	\sum_{r\in\NN}\seque^{1/2}<\infty
	\qquad(\text{e.g. }\seque=2^{-r}).
\end{equation}
The summability in \eqref{eq:seque} makes it possible to apply the Borel--Cantelli lemma later and is the only place we need the exact rate.

By the tightness assumption of Theorem~\ref{Theoremmain} and Prokhorov's theorem, there are an infinite set $\NNtilde \subset\NN$ and a Borel probability measure $\rho$ on $\mathcal X$ such that $\Law(u_n)\to\rho$ weakly on $\mathcal X$ as $n\to\infty$, $n\in\NNtilde $. We fix $\NNtilde $ and $\rho$; the indices $n_r$ below are taken from $\NNtilde $.

\begin{lemma}\label{lem:subsequence}
	There exist sequences of parameters, chosen in the order
	\[
	\seque\;\longrightarrow\;\bigl(R_r,a_r,\delta_r,\varepsilon_{\flat,r}\bigr)
	\;\longrightarrow\;n_r\;\longrightarrow\;
	m_r ,
	\]
	such that the following hold.
	\begin{enumerate}[label=\upshape(P\,\roman*)]
		\item\label{DS:i} $n_r\uparrow\infty$, $m_r\uparrow\infty$,
		$R_r\uparrow\infty$, and $\delta_r,a_r,\varepsilon_{\flat,r}\downarrow0$;
		moreover $\delta_r\le\seque^{1/2}$ and $2^{-m_r}\le\delta_r$.
		\item\label{DS:ii} With $n=n_r$, $m=m_r$ and
		$\nu_{\flat,r}:=\nu(Z^c_{\varepsilon_{\flat,r}})$,
		\begin{align*}
			&\kappa_2(\varepsilon_{\flat,r})\le\seque\,\delta_r^2,
			\qquad
			\bigl(\nu_{\flat,r}C_gC_2\,\tau_{m_r}\bigr)^{1/2}\le\seque\,\delta_r,
			\qquad
			2\bigl(\nu_{\flat,r}L_g\bigr)^{1/2}\,\Theta_{n_r}(m_r)\le\seque\,\delta_r,
			\\
&
			T\,\tau_{m_r}\,\nu_{\flat,r}^2+\tau_{m_r}\nu_{\flat,r}\le2\seque .
		\end{align*}
		\item\label{DS:iii} For every $J\in\NN$ there is $r_0(J)\ge J$ such
		that for all $r\ge r_0(J)$:
		\begin{align*}
			c\bigl(R_r,\epsilon_{J+1}\bigr)&>2\delta_r,
			\qquad
			\theta_{R_r,\epsilon''_{r,J}}(a_r)\ge3\delta_r,
			\qquad
			\varepsilon_{\flat,r}\le\epsilon''_{r,J},
			\\
			a_r&<\tfrac12\min_{0\le j\le J}\bigl(\epsilon_j-\epsilon_{j+1}\bigr),
		\end{align*}
		where $c(R,\epsilon)$ is the constant \eqref{eq:cRe} and
		\[
		\epsilon''_{r,J}:=\frac{c(R_r,\epsilon_{J+1})}{4C_g(1+R_r)}\;\le\;\frac{\epsilon_{J+1}}4 .
		\]
	\end{enumerate}
\end{lemma}

\begin{proof}
	Put $\underline\rho(R):=\min_{0\le s\le R}\rho_g(s)>0$; it is
	nonincreasing. Comparing \eqref{it:nondeg} with \eqref{it:upper} gives
	$\underline\rho(R)\le\rho_g(0)\le C_g$, hence
	$\epsilon''_{r,J}\le\epsilon_{J+1}/4$. Choose, for given $r$,
	\begin{enumerate}[label=(\arabic*)]
		\item $R_r:=\max\{R\in\{1,\dots,r\}:\ \underline\rho(R)\ge\seque^{1/4}\}$,
		and $R_r:=1$ if this set is empty. For fixed $R$ one has
		$\underline\rho(R)\ge\seque^{1/4}$ for large $r$, hence
		$R_r\to\infty$.
		\item $a_r:=\seque^{1/4}$.
		\item $\delta_r:=\min\bigl\{\seque^{1/2},\ \tfrac13\min_{J\le r}
		\theta_{R_r,\epsilon''_{r,J}}(a_r)\bigr\}>0$ (a minimum of finitely
		many positive numbers). Since $\epsilon''_{r,J}$ is nonincreasing in
		$J$, $\varepsilon''_r:=\epsilon''_{r,r}=\min_{J\le r}\epsilon''_{r,J}$.
		\item $\varepsilon_{\flat,r}\in(0,\varepsilon''_r]$ with
		$\kappa_2(\varepsilon_{\flat,r})\le\seque\delta_r^2$; possible since
		$\kappa_2(\varepsilon)\downarrow0$ as $\varepsilon\downarrow0$, see
		\eqref{eq:XV-tail}.
		\item $n_r>n_{r-1}$ arbitrary, with $n_r\in\NNtilde $.
		
		\item $m_r>m_{r-1}$ so large that $\tau_{m_r}\le1$, $2^{-m_r}\le\delta_r$
		and the remaining inequalities in \ref{DS:ii} hold. This is possible
		since $\tau_m\to0$ and, $n_r$ being fixed,
$\Theta_{n_r}(m)\to0$ by Remark~\ref{rem:projection-fixed-n}.
	\end{enumerate}
	\ref{DS:i} and \ref{DS:ii} hold by construction (the sequences in
	\ref{DS:i} can be made monotone by passing to running minima/maxima).
	For \ref{DS:iii} let $r_0(J)\ge J$ be such that, for $r\ge r_0(J)$,
	$\seque^{1/4}<\min\{\epsilon_{J+1}/2,\ \underline\rho(1)\}$ and
	$a_r<\tfrac12\min_{j\le J}(\epsilon_j-\epsilon_{j+1})$. Then the set in
	(1) is nonempty, so
	$c(R_r,\epsilon_{J+1})=\epsilon_{J+1}\underline\rho(R_r)
	\ge\epsilon_{J+1}\seque^{1/4}>2\seque^{1/2}\ge2\delta_r$. Since $J\le r$,
	(3) gives $\theta_{R_r,\epsilon''_{r,J}}(a_r)\ge3\delta_r$, and (4)
	gives $\varepsilon_{\flat,r}\le\varepsilon''_r\le\epsilon''_{r,J}$.
\end{proof}

\begin{corollary}[Error bound at level $r$]\label{cor:error-r}
	For $J\in\NN$ and $r\ge r_0(J)$ let $e_{r,J}$ denote the total error $e_J$
	of \eqref{eq:e-total} for the parameters $n=n_r$,
	$m=m_r$ and $\varepsilon_\flat=\varepsilon_{\flat,r}$ of
	Lemma~\ref{lem:subsequence}. Then
	\begin{equation}\label{eq:Ee-bound}
		\EE\,e_{r,J}\le C_\ast\,\seque^{1/2}\,\delta_r,
		\qquad C_\ast:=4C_2^{1/2}+{2} .
	\end{equation}
	The constant $C_\ast$ does not depend on $J$.
\end{corollary}

\begin{proof}
	By \ref{DS:iii}, $\varepsilon_{\flat,r}\le\epsilon''_{r,J}\le\epsilon_{J+1}$,
	so Lemma~\ref{lem:error-bounds} applies. By \ref{DS:ii}, the first term on
	the right-hand side of \eqref{eq:Ee-fixed} is at most
	$4C_2^{1/2}\seque^{1/2}\delta_r$, and each of the {two} remaining terms is at
	most $\seque\delta_r$. Since $\seque\le\seque^{1/2}$, this gives
	\[
	\EE\,e_{r,J}\le4C_2^{1/2}\seque^{1/2}\delta_r+{2}\seque\delta_r
	\le C_\ast\,\seque^{1/2}\,\delta_r .
	\]
\end{proof}

\begin{remark}\label{rem:dictionary}
	From now on the level index is $r$. At level $r$ the solution is
	$u_{n_r}$, the shifted dyadic projection is $\myp_{n_r,m_r}$, and the algorithm of
	Step~(IV) runs with tolerance $\delta_r$ and cut-off $R_r$. Its output
	at truncation level $J$, constructed in Steps~(V)--(VII), is written
	$\Phi^J_r(u_{n_r})$, the random
	measure on $Z^c_{\epsilon_{J+1}}\times[0,T]$ with one atom
	$(\zhat_k,\Upsilon_{k,m_r,\jty})$ for every registered cell $k$.
\end{remark}

\medskip
\noindent\textbf{(e) The good event.}\label{godevent}
The good event $\mathcal G_{r,J}$ is the event on which everything goes
right for the algorithm to recover exactly the jumps with marks in
$Z^c_{\epsilon_{J+1}}$ (layers $0,\dots,J$) at level $r$. It is the
intersection of four events, each of which rules out a specific way in
which the reconstruction could fail.

\begin{remark}
Recall that, by (A4) in Step~(IV), a jump is registered in $I^\mr _k$ only if
$|\xi^{(\mr )}_{k,\nr }|_V\le R_r$ and $d^{(J)}\bigl(\xi^{(\mr )}_{k,n},\xhat_{k,n,\mr }\bigr)\le\delta_r$;
the cut-off $R_r$ is needed because the lower bound \eqref{it:nondeg}
depends on $|v|_V$.
\end{remark}

For $J\in\NN$ and $r\ge r_0(J)$ (with the notation of (d)) let
$\mathcal G_{r,J}:=\bigcap_{i=1}^4\mathcal G^{(i)}_{r,J}$, where
\begin{enumerate}[label=\upshape(G\arabic*)]
	\item $\mathcal G^{(1)}_{r,J}$: every cell $I_k$ contains at most one
	atom of $\eta$ with mark in $Z^c_{{\varepsilon_{\flat,r}}}$, and $I_0$
	contains no atom with mark in $Z^c_{\epsilon_{J+1}}$;
	\item $\mathcal G^{(2)}_{r,J}$: no atom $(s,z)$ of $\eta$ satisfies
	$\bigl||z|_Z-\epsilon_j\bigr|\le a_r$ for some $j\in\{1,\dots,J+1\}$;
	\item $\mathcal G^{(3)}_{r,J}$: $|\myp_{n_r,\mr }(s)|_V\le {R_r}$ for every atom
	$(s,z)$ of $\eta$ with $s>\tau_\mr $ and $|z|_Z>\epsilon_{J+1}$;
	\item $\mathcal G^{(4)}_{r,J}$: $e_{r,J}\le{\delta_r}/2$.
\end{enumerate}
$\mathcal G^{(1)}_{r,J}$ excludes two jumps in the same cell: their sum
is in general not of the form $g(\xi^{(\mr )}_{k,n},z)$, so neither of them could be
recognised. The first cell is excluded for large jumps, since its frozen
state is the initial datum $0$ 
and not a cell average.
$\mathcal G^{(2)}_{r,J}$ keeps the marks away from the boundaries of the
layers: a mark just inside a layer produces almost the same jump as one
just outside, so the algorithm could pick the wrong layer or register a
band jump as a large one. With the margin $a_r$,
Hypothesis~\ref{ass:g-sep} provides the gap
$\theta_{R_r,\epsilon''_{r,J}}(a_r)\ge3\delta_r$.
$\mathcal G^{(3)}_{r,J}$ ensures that the $V$-cut-off
of (A4) in Step~(IV) does not discard a true jump; the cut-off itself is
needed because the lower bound \eqref{it:nondeg} on $g$ depends on
$|v|_V$. Finally, $\mathcal G^{(4)}_{r,J}$ says that the errors
collected in $e_{r,J}$ (the tiny-jump martingale, the compensators, and the averaging errors) are too small to push the
recovered increment $\xhat_k$ into, or out of, the tolerance window
around the admissible jump sets.

The next Lemma gives an estimate of the probability of the good events.

\begin{lemma}\label{lem:X-good-event}
	Let $J\in\NN$ and $r\ge r_0(J)$. Then
\DEQSZ\label{eq:X-good-prob}
\lqq{ 		\PP\bigl(\mathcal G^c_{r,J}\bigr)\le c_{r,J}}
\\
&		:=&2\seque+2C_\ast\seque^{1/2}
		+T\,\nu\Bigl(\bigcup_{j=1}^{J+1}\bigl\{\bigl||z|_Z-\epsilon_j\bigr|\le a_r\bigr\}\Bigr)
		+\nu_JC_\alpha R_r^{-\alpha}
		\;\xrightarrow[r\to\infty]{}\;0 ,\notag
\EEQSZ
	and on $\mathcal G_{r,J}$ the following hold.
	\begin{enumerate}[label=\upshape(\alph*)]
		\item\label{X:good-a} \emph{(No false registration.)} Every registered
		cell contains exactly one atom of $\eta$ with mark in
		$Z^c_{\epsilon_{J+1}}$.
		\item\label{X:good-b} \emph{(Correct registration.)} Every atom $(s,z)$
		of $\eta$ with $z\in\Delta Z_j$, $j\le J$, lies in a cell $I_k$,
		$k\ge1$, containing no other atom with mark in
		$Z^c_{{\varepsilon_{\flat,r}}}$; this cell is registered in the correct
		layer $\hat\jmath=j$, and the recovered mark satisfies
		\[
		\bigl|g(\xi^{(\mr )}_{k,\nr },\zhat_k)-g(\xi^{(\mr )}_{k,\nr },z)\bigr|_H\le\delta_r+2^{-m_r},
		\qquad d_Z(\zhat_k,z)\le a_r .
		\]
		\item\label{X:good-c} \emph{(Time error.)} The recorded time satisfies
		$|\Upsilon_{k,m_r,\jty}-s|\le\tau_{m_r}/2$.
	\end{enumerate}
	In particular, the atoms of $\Phi^J_r(u_{n_r})$ and of
	$\eta|_{Z^c_{\epsilon_{J+1}}\times[0,T]}=\mathcal J^J(u_{n_r})$
	(Lemma~\ref{lem:jump-identity}) are in one-to-one correspondence, and
	for the Prokhorov metric $d_{\mathrm v}$ built on
	$\rho((z,s),(z',s'))=d_Z(z,z')+|s-s'|$,
	\begin{equation}\label{eq:X-Dm-bound}
		d_{\mathrm v}\bigl(\Phi^J_r(u_{n_r}),\mathcal J^J(u_{n_r})\bigr)
		\le b_r:=a_r+\tau_{m_r}/2
		\qquad\text{on }\mathcal G_{r,J} .
	\end{equation}
\end{lemma}

\noindent\emph{Idea of the proof.} On $\mathcal G_{r,J}$ every cell
$I_k$ with $|\xi^{(\mr )}_{k,\nr }|_V\le R_r$ is in one of four situations, each handled
by a margin:
\begin{itemize}
	\item \emph{only tiny jumps:} $|\xhat_k|_H\le\delta_r/2$, while every
	admissible jump has size at least $c(R_r,\epsilon_{J+1})>2\delta_r$;
	the cell is not registered;
	\item \emph{a band jump with $|z|_Z\le\epsilon''_{r,J}$:} by
	\eqref{it:upper} the jump is at most $\frac14c(R_r,\epsilon_{J+1})$ in
	$H$; not registered;
	\item \emph{a band jump with $\epsilon''_{r,J}<|z|_Z\le\epsilon_{J+1}$:}
	by $\mathcal G^{(2)}_{r,J}$ its mark lies at least $a_r$ below
	$\epsilon_{J+1}$, so by Hypothesis~\ref{ass:g-sep} its image is at
	distance at least $\theta(a_r)-\delta_r/2\ge\frac52\delta_r$ from all
	layers; not registered;
	\item \emph{a large jump in layer $j$:} it is registered (distance at
	most $\delta_r/2$), in the correct layer (all other layers are at
	distance at least $\frac52\delta_r$), and the recovered mark lies within
	$a_r$ of the true one.
\end{itemize}
This gives a one-to-one correspondence between the true and the
registered jumps, with mark error at most $a_r$ and time error at most
$\tau_{m_r}/2$, hence \eqref{eq:X-Dm-bound}.

\begin{proof}
	\emph{Probabilities.} The counts
	$\eta(Z^c_{{\varepsilon_\flat,r}}\times I_k)$ are Poisson with parameter
	${\nu_{\flat,r}}\tau_\mr $, so
	$\PP((\mathcal G^{(1)})^c)\le2^\mr ({\nu_{\flat,r}}\tau_\mr )^2+\nu_J\tau_\mr 
	\le T\tau_\mr {\nu_{\flat,r}}^2+\tau_\mr {\nu_{\flat,r}}\le2\seque$ by \ref{DS:ii}. Next,
	$\PP((\mathcal G^{(2)})^c)$ is at most the expected number of atoms in
	the indicated set, i.e.\ $T\nu(Z_J)$; this tends to $0$ as
	$a_r\downarrow0$, since the sets decrease to
	$\bigcup_{j\le J+1}\{|z|_Z=\epsilon_j\}$, which is $\nu$-null, and have
	finite $\nu$-measure: by \ref{DS:iii}, $a_r=\seque^{1/4}<\epsilon_{J+1}/2$ for
	$r\ge r_0(J)$, so they are contained in $\{|z|_Z\ge\epsilon_{J+1}/2\}$. For $\mathcal G^{(3)}$, the
	integrand $\mathbf 1\{|\myp_{n,\mr }(s)|_V>R\}$ is predictable, so by the
	compensation formula, Chebyshev's inequality and Jensen's inequality on
	each cell ($|\xi^{(\mr )}_{k,\nr }|_V^\alpha\le\tau_\mr ^{-1}\int_{I_{k-1}}|u_n|_V^\alpha$,
	$k\ge1$),
	\[
	\PP\bigl((\mathcal G^{(3)})^c\bigr)\le
	\nu_J\,\EE\int_{\tau_\mr }^T\mathbf 1\{|\myp_{n,\mr }(s)|_V>R\}\,ds
	\le\nu_JR^{-\alpha}\,\EE\int_0^T|u_n(s)|_V^\alpha\,ds\le\nu_JC_\alpha R^{-\alpha}.
	\]
	Finally $\PP((\mathcal G^{(4)})^c)\le2\EE e_{r,J}/{\delta_r}\le2C_\ast\seque^{1/2}$
	by Corollary~\ref{cor:error-r}. The limit in \eqref{eq:X-good-prob}
	follows from $\seque\to0$, $a_r\to0$ and $R_r\to\infty$.

	\emph{The registration.} Work on $\mathcal G_{r,J}$. By
	$\mathcal G^{(1)}$, every cell contains at most one atom with mark in
	$Z^c_{{\varepsilon_\flat,r}}$, so \eqref{eq:xhat-close} applies, with
	$e_{r,J}\le{\delta_r}/2$ by $\mathcal G^{(4)}$. Recall from \eqref{eq:cRe},
	\eqref{it:nondeg} that $|g(v,z')|_H\ge c(R,\epsilon_{J+1})$ for
	$|v|_V\le R$ and $|z'|_Z>\epsilon_{J+1}$, i.e.\ for all $z'$ in the
	layers. By the cut-off in (A4) only cells with $|\xi^{(\mr )}_{k,n}|_V\le R$ can be
	registered; consider such a cell and $z'\in\Delta Z_{j'}$, $j'\le J$.

	\emph{Case 1: no atom with mark in $Z^c_{{\varepsilon_\flat,r}}$.} Then
	$|\xhat_k|_H\le{\delta_r}/2$, and
	$|\xhat_k-g(\xi^{(\mr )}_{k,n},z')|_H\ge c(R,\epsilon_{J+1})-{\delta_r}/2>\tfrac32{\delta_r}$
	by \ref{DS:iii}. So $d^{(J)}>{\delta_r}$: the cell is not registered.

	\emph{Case 2: one atom $(s,z)$ with ${\varepsilon_\flat,r}<|z|_Z\le\epsilon''_{r,J}$.}
	By \eqref{it:upper}, $|g(\xi^{(\mr )}_{k,\nr },z)|_H\le C_g(1+R)\epsilon''_{r,J}
	=\tfrac14c(R,\epsilon_{J+1})$, hence
	$|\xhat_k-g(\xi^{(\mr )}_{k,n},z')|_H\ge\tfrac34c(R,\epsilon_{J+1})-\tfrac12{\delta_r}>{\delta_r}$.
	Not registered.

	\emph{Case 3: one atom $(s,z)$ with $\epsilon''_{r,J}<|z|_Z\le\epsilon_{J+1}$.}
	By $\mathcal G^{(2)}$, $|z|_Z<\epsilon_{J+1}-a_r$, so
	$d_Z(z,z')\ge|z'|_Z-|z|_Z>a_r$. Both marks have norm
	$\ge\epsilon''_{r,J}$, so Hypothesis~\ref{ass:g-sep} and \ref{DS:iii} give
	$|g(\xi^{(\mr )}_{k,n},z)-g(\xi^{(\mr )}_{k,n},z')|_H\ge\theta_{R,\epsilon''_{r,J}}(a_r)\ge3{\delta_r}$,
	hence $|\xhat_k-g(\xi^{(\mr )}_{k,n},z')|_H\ge\tfrac52{\delta_r}>{\delta_r}$. Not registered.

	Together, Cases~1--3 prove \ref{X:good-a}.

	\emph{Case 4: one atom $(s,z)$ with $z\in\Delta Z_j$, $j\le J$.} By
	$\mathcal G^{(1)}$, $k\ge1$; by $\mathcal G^{(3)}$, $|\xi^{(\mr )}_{k,n}|_V\le R$. By
	\eqref{eq:xhat-close}, $d_j\le|\xhat_k-g(\xi^{(\mr )}_{k,n},z)|_H\le{\delta_r}/2$, so the
	cell is registered. By $\mathcal G^{(2)}$ and \ref{DS:iii},
	$\epsilon_{j+1}+a_r<|z|_Z<\epsilon_j-a_r$ (only the lower bound if
	$j=0$), so $d_Z(z,z')\ge a_r$ for every $z'\in\Delta Z_{j'}$, $j'\ne j$,
	and Hypothesis~\ref{ass:g-sep} gives
	$d_{j'}\ge3{\delta_r}-{\delta_r}/2>{\delta_r}/2\ge d_j$. Hence $\hat\jmath=j$. By
	(A7), $|g(\xi^{(\mr )}_{k,n},\zhat_k)-\xhat_k|_H\le d_j+2^{-\mr }$, so
	$|g(\xi^{(\mr )}_{k,n},\zhat_k)-g(\xi^{(\mr )}_{k,n},z)|_H\le{\delta_r}+2^{-\mr }\le2{\delta_r}<3{\delta_r}$. If
	$d_Z(\zhat_k,z)>a_r$, Hypothesis~\ref{ass:g-sep} (both marks have norm
	$>\epsilon_{J+1}\ge\epsilon''_{r,J}$) would give a difference
	$\ge3{\delta_r}$, a contradiction. This proves \ref{X:good-b}, and
	\ref{X:good-c} holds since $\Upsilon_{k,\mr,\jt }$ is the midpoint of $I_k$.

	\emph{The distance.} By \ref{X:good-a}--\ref{X:good-c}, the atoms of
	$\Phi^J_r(u_{n_r})$ and of $\eta|_{Z^c_{\epsilon_{J+1}}\times[0,T]}$
	can be matched one to one with $\rho$-distance at most
	$a_r+\tau_\mr /2$, and the atom-matching estimate for the Prokhorov
	metric gives \eqref{eq:X-Dm-bound}.
\end{proof}

\medskip
\noindent\textbf{(f) Summability.}
By Lemma~\ref{lem:X-good-event}, $\PP(\mathcal G^c_{r,J})\le c_{r,J}\to0$ for
every $J$. Step~(IX) and the later steps use only this convergence in
probability; the almost sure convergence used there is produced on
$\mathfrak A^1$ by the Jakubowski--Skorokhod theorem. The following remark
records what a Borel--Cantelli argument would give.

\begin{remark}\label{rem:tightness-not-enough}
	The parts $2\seque+2C_\ast\seque^{1/2}$ of $c_{r,J}$ are summable by
	\eqref{eq:seque}; the two remaining parts tend to zero but carry no
	rate, since they depend on $\nu$ near the spheres $\{|z|_Z=\epsilon_j\}$
	and on the decay of $\rho_g$. This is the same dichotomy as elsewhere
	in the paper: convergence comes for free, a rate does not. For the
	Borel--Cantelli lemma one passes to a subsequence $(r_l)$ with
	$c_{r_l,J}\le2^{-l}$ for all $J\le l$, which exists since
	$c_{r,J}\to0$ for every $J$. The arguments of Step~(IX) and later only
	use $\PP(\mathcal G^c_{r,J})\to0$, i.e.\ convergence in probability.
\end{remark}

\bigskip
\item \textbf{The application of the
	Jakubowski--Skorokhod representation theorem.}

The aim of this step is twofold. First, we establish tightness of the laws of
\[
\Bigl(u_{n_r},(\Phi^J_r(u_{n_r}))_{J\in\NN}\Bigr),\qquad r\in\NN,
\]
on the corresponding product space; here $r$ is the level index of
Lemma~\ref{lem:subsequence} and $\Phi^J_r$ the reconstruction at level $r$,
see Remark~\ref{rem:dictionary}.
Second, we apply the Jakubowski--Skorokhod representation theorem to this tight family and obtain, on a new probability space
$\mathfrak A^1=(\Omega^1,\mathscr F^1,\PP^1)$, an almost surely convergent sequence of copies together with its limit
$(u^1_\ast,(\eta^{1,J})_{J\in\NN})$. We use the following form of the representation theorem; see \cite[Thm.~2]{Jakubowski1997}.

\begin{theorem}
	\label{thm:jakubowski}

Let $(\mathcal T,\tau)$ be a topological space on which there exists a
countable family $\{F_l:\mathcal T\to[0,1]\}_{l\in\NN}$ of
$\tau$-continuous functions separating the points of $\mathcal T$. Let
$(\chi_m)_{m\in\NN}$ be a sequence of $\mathcal T$-valued random elements,
defined on a probability space $(\Omega,\mathscr F,\PP)$ and measurable
with respect to the $\sigma$-field $\sigma(F_l:l\in\NN)$, whose laws are
uniformly tight: for every $\varepsilon>0$ there is a $\tau$-compact
$K_\varepsilon\subset\mathcal T$ with
$\PP(\chi_m\in K_\varepsilon)\ge1-\varepsilon$ for all $m\in\NN$. Then
there exist a subsequence $(m_i)_{i\in\NN}$, a probability space
$(\Omega^1,\mathscr F^1,\PP^1)$ --- which may be taken to be
$\bigl([0,1],\mathscr B([0,1]),\Leb\bigr)$ --- and $\mathcal T$-valued
random elements $\chi^1_i$, $i\in\NN$, and $\chi^1_\ast$ on it, measurable
with respect to $\sigma(F_l:l\in\NN)$, such that
\begin{itemize}
\item  $	\Law_{\PP^1}\bigl(\chi^1_i\bigr)=\Law_{\PP}\bigl(\chi_{m_i}\bigr)$ on $\sigma(F_l:l\in\NN)$
for every $i\in\NN$,
\item $	\chi^1_i(\omega)\xrightarrow[i\to\infty]{\tau}\chi^1_\ast(\omega)$
for {every} $\omega\in\Omega^1$.

\end{itemize}
\end{theorem}

\medskip

Let
$\mathcal X=L^\alpha(0,T;V)\cap\DD([0,T];H)$ is equipped with the join of
the norm topology of $L^\alpha(0,T;V)$ and the Skorokhod $J_1$ topology of
$\DD([0,T];H)$, so convergence in $\mathcal X$ means convergence in both components. For fixed $J$, let $\MMN(\ZZJ)$ denote the space of finite integer-valued Borel measures on
$\ZZJ:=Z^c_J\times[0,T]$, where $Z^c_J:=\bigcup_{j\le J}\Delta Z_j=Z^c_{\epsilon_{J+1}}\cap\{|z|_Z\le\epsilon_0\}$, equipped with the ordinary weak topology of finite measures. Since $\ZZJ$ is Polish, $\MMN(\ZZJ)$ is Polish as well.

To represent \emph{all levels simultaneously}
we put next
\[
\mathcal T_\infty:=\mathcal X\times\prod_{J\in\NN}\MMN\bigl(\ZZJ \bigr)
\]
and equip $\mathcal T_\infty$ with the product topology. We record the
state-space properties used in the representation argument.

\begin{lemma}\label{lem:state-space}
Let $Z$ be a complete separable metric space
(Hypothesis~\ref{ass:GF}). Then:
\begin{enumerate}[label=\upshape(\roman*)]
\item\label{SS:i} Let $\mathcal X=L^\alpha(0,T;V)\cap\DD([0,T];H)$ be the set of
all $u\in\DD([0,T];H)$ with $u(t)\in V$ for a.e.\ $t$ and
$\int_0^T|u(t)|_V^\alpha\,dt<\infty$, equipped with the \emph{join topology},
i.e.\ the coarsest topology for which both maps
$\mathcal X\to L^\alpha(0,T;V)$, $u\mapsto u$ (taken as an a.e.-class), and
$\mathcal X\to\DD([0,T];H)$, $u\mapsto u$, are continuous. Then $\mathcal X$ is
a Polish space, and
\[
d_{\mathcal X}(u,w):=\|u-w\|_{L^\alpha(0,T;V)}+d^\circ(u,w)
\]
is a complete metric on $\mathcal X$ inducing the join topology, where
$d^\circ$ is Billingsley's metric \eqref{eq:billingsley-metric} on
$\DD([0,T];H)$.

\item\label{SS:ii} $\MMN(\ZZJ)$, the space of finite
$\NN_0$-valued Borel measures on $\ZZJ$ equipped with the weak topology, is Polish. It is a closed subset of the space of finite Borel measures on $\ZZJ$; the total-mass functional $\mu\mapsto\mu(\ZZJ)$ is continuous, and the evaluation maps $\mu\mapsto\mu(B)$, $B\in\mathscr B(\ZZJ)$, are Borel measurable.

\item\label{SS:iib} The countable product space
\[
\prod_{J\in\NN}\MMN\bigl(\ZZJ\bigr)
\]
equipped with the product topology is Polish.

\item\label{SS:iii} $\mathcal T_\infty$ is Polish with the product topology. Moreover, it admits a countable family of $[0,1]$-valued continuous functions separating its points, and this family generates the Borel $\sigma$-field $\mathscr B(\mathcal T_\infty)$.
\end{enumerate}
\end{lemma}

\begin{proof}

\ref{SS:i} Polishness of the two ambient spaces does not by itself give
Polishness of their intersection, so we use the diagonal embedding
\[
\iota:\mathcal X\longrightarrow\mathcal P:=L^\alpha(0,T;V)\times\DD([0,T];H),
\qquad \iota(u):=(u,u).
\]
\emph{The ambient space.} $V$ is a separable Banach space, so
$L^\alpha(0,T;V)$ is a separable Banach space. $H$ is a separable Hilbert
space, so $(\DD([0,T];H),d^\circ)$ is complete and separable
(Appendix~\ref{app:Skorohod-space}; \cite[\S12]{Billingsley1999},
\cite[Ch.~3, \S5]{Ethier1986}). Hence $\mathcal P$ with the product
topology is Polish, and $\|v-v'\|_{L^\alpha(0,T;V)}+d^\circ(w,w')$ is a
complete metric for it.

\emph{$\iota$ is a homeomorphism onto its image.} The map $\iota$ is
injective, since its second component is the identity on paths. By
definition, the join topology is the initial topology of the two
components of $\iota$, i.e.\ the topology induced on $\mathcal X$ by $\iota$
and the product topology. Hence $\iota$ is a homeomorphism from
$\mathcal X$ onto
\[
\mathcal G:=\iota(\mathcal X)
=\bigl\{(v,w)\in\mathcal P:\ w(t)=v(t)\ \text{for a.e.\ }t\in[0,T]\bigr\},
\]
and $d_{\mathcal X}(u,w)$ is the distance of $\iota(u)$ and $\iota(w)$ in
$\mathcal P$. The identification of $\mathcal G$ uses that a path
$w\in\DD([0,T];H)$ belongs to $\mathcal X$ if and only if its a.e.-class
lies in $L^\alpha(0,T;V)$, where $V\hookrightarrow H$.

\emph{$\mathcal G$ is closed in $\mathcal P$.} Let $(v_k,w_k)\in\mathcal G$
converge in $\mathcal P$ to $(v,w)$. Since $V\hookrightarrow H$
continuously and $T<\infty$,
\[
\|v_k-v\|_{L^1(0,T;H)}\le c_{V\hookrightarrow H}\,T^{1-\frac1\alpha}\,
\|v_k-v\|_{L^\alpha(0,T;V)}\longrightarrow0 .
\]
Since $w_k\to w$ in the $J_1$ topology, Lemma~\ref{lem:J1-pointwise}\,(a)
gives $w_k(t)\to w(t)$ in $H$ for every continuity point $t$ of $w$, hence
for a.e.\ $t\in[0,T]$. Moreover, $J_1$-convergent sequences are bounded in
the supremum norm (see Lemma~\ref{lem:J1-pointwise}\,(d)): with the time changes $\lambda_k$ of
Appendix~\ref{app:Skorohod-space},
$\sup_t|w_k(t)|_H=\sup_s|w_k(\lambda_k(s))|_H\le\sup_s|w(s)|_H+\varepsilon_k$
with $\varepsilon_k\to0$, and $\sup_s|w(s)|_H<\infty$ for a c\`adl\`ag path.
By dominated convergence, $w_k\to w$ in $L^1(0,T;H)$. Since $w_k=v_k$
a.e.,
\[
\|w-v\|_{L^1(0,T;H)}\le\|w-w_k\|_{L^1(0,T;H)}+\|v_k-v\|_{L^1(0,T;H)}
\longrightarrow0,
\]
so $w=v$ a.e., i.e.\ $(v,w)\in\mathcal G$.

\emph{Conclusion.} A closed subset of a Polish space is Polish, and the
restriction of a complete metric to a closed subset is complete. Hence
$\mathcal G$, and with it $\mathcal X\cong\mathcal G$, is Polish, and
$d_{\mathcal X}$ is a complete metric inducing the join topology. In
particular, $u_k\to u$ in $\mathcal X$ if and only if $u_k\to u$ in
$L^\alpha(0,T;V)$ and in $\DD([0,T];H)$, and the Borel $\sigma$-field of
$\mathcal X$ is the trace of $\mathscr B(L^\alpha(0,T;V))\otimes\mathscr B(\DD([0,T];H))$
under $\iota$.

\ref{SS:ii} The space of finite Borel measures on a Polish space, with the
weak topology, is Polish (Prokhorov's metric; see
\cite[App.~A2.5]{DaleyVereJones2003}). Closedness of
$\MMN(\ZZJ )$: A weak limit of finite $\NN_0$-valued measures
on a complete separable metric space is again $\NN_0$-valued, see
\cite[App.~A2.6, Prop.~A2.6.III]{DaleyVereJones2003}. The total mass is $\mu\mapsto\mu(\mathds 1)$ with
$\mathds 1\in C_b(\ZZJ )$, hence weakly continuous. Borel measurability of
the evaluations $\mu\mapsto\mu(B)$ is
\cite[App.~A2.6]{DaleyVereJones2003}; cf.\ also
\cite[Ch.~4]{Kallenberg2017}.

\ref{SS:iib} Since each $\MMN(\ZZJ)$ is Polish, the countable product
$\mathbb M:=\prod_{J\in\NN}\MMN(\ZZJ)$, equipped with the product topology, is Polish; see, for example, Kechris~\cite[\S3.A]{Kechris1995}.

\ref{SS:iii} Since $\mathcal T_\infty=\mathcal X\times\mathbb M$, assertion~\ref{SS:iii} follows from~\ref{SS:i}--\ref{SS:iib} and the standard countable separating-family construction for Polish spaces.
\end{proof}

By Lemma~\ref{lem:state-space}\,\ref{SS:iii}, $\mathcal T_\infty$ is Polish,
and the union over the factors of countable separating families, composed
with the projections, is a countable family of continuous functions
separating the points of $\mathcal T_\infty$ and generating
$\mathscr B(\mathcal T_\infty)$.

Now, we apply  Theorem~\ref{thm:jakubowski} to the sequence of tuples
\DEQSZ\label{tuples}
\Bigl(u_{n_r},\ \bigl(\Phi^J_r(u_{n_r})\bigr)_{J\in\NN}
\Bigr),\qquad r\in\NN,
\EEQSZ
on $\mathfrak A$. Once tightness has been established, Theorem~\ref{thm:jakubowski} yields a single probability space $\mathfrak A^1$ and a single subsequence
$\{r_i:i\in\NN\}$ serving all truncation levels simultaneously.
All levels $J\in\NN$ are represented on one common probability space and are
consistent, which will be shown in Lemma~\ref{lem:J-consistency}.
The passage $J\to\infty$ is
\emph{pathwise monotone} and handled in Proposition
\ref{prop:full-levy-patch}, whose Poisson property follows set by set from the
stability of the Poisson family under independent sums and monotone limits.
So the
Poisson random measure $\eta^1$ with the full L\'evy measure $\nu$ is
obtained from the coherent family $(\eta^{1,J})_{J\in\NN}$ by a pathwise
monotone limit.

\medskip

In the first part (A1), we tackle the tightness.
The representation step is
formulated precisely in (A2) below.  The state space is prepared in
Lemma~\ref{lem:state-space},  and Theorem~\ref{thm:jakubowski} is applied to the tuples \eqref{tuples} in
Proposition~\ref{prop:JS-pairs-allJ}.

\medskip

Throughout, write
$\pi_J:=\mathrm{PRM}\bigl(\nu\rest_{Z^c_J}\otimes\Leb\bigr)$ for the law of a Poisson
random measure with intensity $\nu\rest_{Z^c_J}\otimes\Leb$ on $\MMN(\ZZJ)$;
as in Step~(VIII), $\nu_J:=\nu(Z^c_J)$ denotes its total mass per unit time. For more details see Appendix \ref{app:levy-linde}.

For the sake of simplicity, every limit
$r\to\infty$ in the sequel is understood along $(r_i)$. We write
$\mathfrak A^1=(\Omega^1,\mathscr F^1,\PP^1)$ for the represented probability space and
$u^1_i$, $N^{1,J}_i$, $u^1_\ast$ and $\eta^{1,J}$ for the random elements given by Proposition~\ref{prop:JS-pairs-allJ}. To keep the representation simple, we show first tightness and measurability for $\mathcal T:=\mathcal X\times \MMN(\ZZJ)$ and extend the results in a second step to $\mathcal T_\infty$.

\medskip

\noindent\emph{(A1) Verification of the Assumptions of Theorem \ref{thm:jakubowski}:}
For the convenience of the reader we recall the product space here; the construction of the product space and of the metric inducing its topology is given in Subsection~\ref{app:product-space}, Definition~\ref{def:product-vague}.
Let us consider firstly the sequence
\[
\Bigl((\Phi^J_r(u_{n_r}))_{J\in\NN}\Bigr)_{r\in\NN}
\qquad\text{with values in}\qquad
\mathbb M:=\prod_{J\in\NN} \MMN\bigl(\ZZJ\bigr).
\]
\begin{steps_a1}
\item\emph{Separability.} Being a countable product of Polish spaces,
$\mathbb M$ is Polish, hence separable, and
thus carries a countable family of $[0,1]$-valued continuous functions
separating its points, as required in Theorem~\ref{thm:jakubowski}.

\item\emph{Measurability.} By Proposition~\ref{prop:Nhat-meas}
(page~\pageref{prop:Nhat-meas}) and the measurable mark-recovery map, the evaluations of
$\Phi^J_r(u_{n_r})$ on Borel subsets of $\ZZJ$ are $\mathcal F$-measurable. For the weak topology on finite measures, the Borel $\sigma$-field is generated by the evaluation maps. Hence
\[
\Phi^J_r(u_{n_r}):\Omega\longrightarrow\MMN(\ZZJ)
\]
is Borel measurable. Since the Borel $\sigma$-field of a countable product of Polish spaces is the product $\sigma$-field, coordinatewise measurability gives measurability of
$\bigl(\Phi^J_r(u_{n_r})\bigr)_{J\in\NN}:\Omega\to\mathbb M$.

\item \emph{Tightness of the tuples.}
By Proposition~\ref{prop:tight-product}\,\ref{it:tp-joint}$\Leftrightarrow$\ref{it:tp-coord}
it suffices to verify tightness coordinatewise, i.e.\ that for every
$J\in\NN$ the family $\bigl(\Phi^J_r(u_{n_r})\bigr)_{r\in\NN}$ is tight in
$\MMN(\ZZJ)$. We show the stronger statement
$\Law\bigl(\Phi^J_r(u_{n_r})\bigr)\Rightarrow\pi_J$ below, by comparison with
the path jump measure; a weakly convergent sequence of laws on a Polish space
is tight. A first-moment bound on the total mass $\Phi^J_r(u_{n_r})(\ZZJ)$
alone would not suffice when $Z$ is infinite-dimensional, since it does not
control the location of the atoms; see
Remark~\ref{rem:ktw-counterexample}.

\end{steps_a1}

Fix $J\in\NN$ and let $r\ge r_0(J)$, with the parameters of Lemma~\ref{lem:subsequence}.
\par

We work on the space of finite measures on the Polish space
$\ZZJ$ with its ordinary weak topology, metrized by a Prokhorov metric $d_{\mathrm v}$. The only property used below is the standard atom-matching estimate: if
$\mu=\sum_{i=1}^{K}\delta_{p_i}$ and
$\mu'=\sum_{i=1}^{K}\delta_{q_i}$ with
$\rho(p_i,q_i)\le\varepsilon$ for every $i$, then
$d_{\mathrm v}(\mu,\mu')\le\varepsilon$.
Consider, on $\mathfrak A$, the following sequence
\[
\Bigl(u_{n_r},\ \Phi^J_r(u_{n_r})\Bigr)_{r\in\NN}
\qquad\text{with values in}\qquad
\mathcal X\times\MMN(\ZZJ).
\]
Lemma~\ref{lem:X-good-event} compares the reconstructed measure with the truncated jump measure. Indeed, by \eqref{eq:X-Dm-bound} and \eqref{eq:X-good-prob},
\[
\PP\Bigl(d_{\mathrm v}\bigl(\Phi^J_r(u_{n_r}),\mathcal J^J(u_{n_r})\bigr)>b_r\Bigr)
\le c_{r,J},
\qquad b_r\to0,\ c_{r,J}\to0,
\]
while, by Lemma~\ref{lem:jump-identity},
$\Law\bigl(\mathcal J^J(u_{n_r})\bigr)=\pi_J$ is \emph{one and the same law for every} $r$. Hence, for every bounded Lipschitz function $h$ on
$\MMN(\ZZJ)$,
\[
\bigl|\EE\,h\bigl(\Phi^J_r(u_{n_r})\bigr)-\pi_J(h)\bigr|
\le \operatorname{Lip}(h)\,b_r+2\|h\|_\infty\,c_{r,J}\longrightarrow0.
\]
Since bounded Lipschitz functions determine weak convergence on a Polish space,
$\Law\bigl(\Phi^J_r(u_{n_r})\bigr)\Rightarrow\pi_J$. Thus, the family is tight in
$\MMN(\ZZJ)$ by Prokhorov's theorem; the finitely many levels $r<r_0(J)$
do not matter, since a single law on a Polish space is tight.

We now consider the full tuple \eqref{tuples},
\[
\Bigl(u_{n_r},\bigl(\Phi^J_r(u_{n_r})\bigr)_{J\in\NN}\Bigr),\qquad r\in\NN.
\]
The first coordinate is tight by the hypothesis of Theorem~\ref{Theoremmain}. For every fixed $J$, the preceding comparison argument gives tightness of
$\{\Phi^J_r(u_{n_r})\}_r$. Given $\varepsilon>0$, choose a compact
$K_0\subset\mathcal X$ with
$\inf_r\PP(u_{n_r}\in K_0)\ge1-\varepsilon/2$, and for every $J\in\NN$ choose a compact
$K_J\subset\MMN(\ZZJ)$ such that
\[
\inf_r\PP\bigl(\Phi^J_r(u_{n_r})\in K_J\bigr)
\ge1-\varepsilon 2^{-J-2}.
\]
Then $K:=K_0\times\prod_{J\in\NN}K_J$ is compact in the product topology, and the union bound gives
\[
\inf_r\PP\Bigl(
\bigl(u_{n_r},(\Phi^J_r(u_{n_r}))_J\bigr)\in K
\Bigr)\ge1-\varepsilon.
\]
Hence the laws of the tuples are tight on $\mathcal T_\infty$.

\medskip
\noindent\emph{(A2) The representation step.}

Theorem~\ref{thm:jakubowski}, applied to the tuples \eqref{tuples}, gives the
following.

\begin{proposition}\label{prop:JS-pairs-allJ}
Let $n_r,m_r,\delta_r,R_r$ be the parameters of
Lemma~\ref{lem:subsequence}, and let $\Phi^J_r$ be the reconstruction
functional at level $r$ (Remark~\ref{rem:dictionary}), so that
\[
\chi_r:=\Bigl(u_{n_r},(\Phi^J_r(u_{n_r}))_{J\in\NN}\Bigr),
\qquad r\in\NN,
\]
is a sequence of $\mathcal T_\infty$-valued Borel random elements on
$\mathfrak A$.

Then there exists a
subsequence $(r_i)_{i\in\NN}$, a probability space
$\mathfrak A^1=(\Omega^1,\mathscr F^1,\PP^1)$ which
may be taken to be $([0,1],\mathscr B([0,1]),\Leb)$ and Borel random elements
\[
\bigl(u^1_i,(N^{1,J}_i)_{J\in\NN}\bigr),\ i\in\NN,
\qquad\text{and}\qquad
\bigl(u^1_\ast,(\eta^{1,J})_{J\in\NN}\bigr)
\]
with values in $\mathcal T_\infty$, such that:
\begin{enumerate}[label=\upshape(\alph*)]
\item \label{eq:JS-laws} \emph{Equality of laws holds:} for every $i\in\NN$,
$$\Law_{\PP^1}\bigl(u^1_i,(N^{1,J}_i)_J\bigr)
=\Law_{\PP}\bigl(u_{n_{r_i}},(\Phi^J_{r_i}(u_{n_{r_i}}))_J\bigr)\quad \mbox{on}\quad
\mathscr B(\mathcal T_\infty);
$$
\item \label{JS:as} \emph{Convergence $\PP^1$-a.s. is given:}  We have $\PP^1$-a.s. as $i\to\infty$
\begin{equation}\label{eq:JS-as}
\left\{
\begin{aligned}
u^1_i(\omega)\longrightarrow u^1_\ast(\omega)
\quad\text{in }L^\alpha(0,T;V)	\text{ {and}  }\text{ in }	\,	\DD([0,T];H)\ (J_1),
\\
\mbox{ for all }J\in\NN\quad
N^{1,J}_i(\omega)\longrightarrow\eta^{1,J}(\omega)
\quad\text{weakly in }\MMN(\ZZJ ).
\end{aligned}
\right.\end{equation}

\item\label{JS:rho} \emph{Identification of the first marginal:}
It holds	$\Law_{\PP^1}(u^1_\ast)=\rho$, where $\rho$ is the limit law of
$\Law(u_n)$ along $\NNtilde $.

\item\label{JS:point} \emph{The second marginal is a finite
integer-valued random measure:}  For all $J\in \N$, $\PP^1$-a.s., $\eta^{1,J}$ is a finite $\NN_0$-valued measure, and
$\eta^{1,J}(\ZZJ)=\lim_{i\to\infty}N^{1,J}_i(\ZZJ)<\infty$.
In addition,  for every
$B\in\mathscr B(\ZZJ  )$ the mass $\eta^{1,J}(B)$ is a
$\mathscr F^1$-measurable random variable and on a set $\Omega_1$ with $\PP^1(\Omega_1)=1$ (see
\ref{JS:as}), $N^{1,J}_i(A)\to\eta^{1,J}(A)$ for every $A\in\mathscr B(\ZZJ)$
with $\eta^{1,J}(\partial A)=0$.
\end{enumerate}
\end{proposition}

\begin{proof}[Proof of Proposition \ref{prop:JS-pairs-allJ}:]
Throughout, $\pi_J=\mathrm{PRM}\bigl(\nu\rest_{Z^c_J}\otimes\Leb\bigr)$ is the
law of a Poisson random measure with intensity $\nu\rest_{Z^c_J}\otimes\Leb$ on
$\MMN(\ZZJ)$.
In (A1) we verified the assumptions of Theorem~\ref{thm:jakubowski}: $\mathcal T_\infty$ is Polish and the sequence of tuples defined in \eqref{tuples} is measurable and tight.

In particular, the sequence $(\chi_r)_{r\in\NN}$ satisfies the hypotheses of
Theorem~\ref{thm:jakubowski} on $(\mathcal T_\infty,\{F_l\})$. It yields a
subsequence $(r_i)_{i\in\NN}$, a probability space
$\mathfrak A^1=(\Omega^1,\mathscr F^1,\PP^1)$ which may be taken to be
$([0,1],\mathscr B([0,1]),\Leb)$ and $\mathcal T_\infty$-valued random
elements
$\chi^1_i=:\bigl(u^1_i,(N^{1,J}_i)_{J}\bigr)$ and
$\chi^1_\ast=:\bigl(u^1_\ast,(\eta^{1,J})_{J}\bigr)$
with $\Law_{\PP^1}(\chi^1_i)=\Law_{\PP}(\chi_{r_i})$ on
$\mathscr B(\mathcal T_\infty)$ for every $i$ which is assertion (a)
and $\chi^1_i(\omega)\to\chi^1_\ast(\omega)$ for every
$\omega\in\Omega^1$ in the product topology. Convergence in the product
topology of $\mathcal T_\infty$ is exactly componentwise convergence, that is,
$u^1_i\to u^1_\ast$ in $\mathcal X$ (that is, in $L^\alpha(0,T;V)$
\emph{and} in the $J_1$ topology of $\DD([0,T];H)$, by the definition of
the join topology) and, for every $J$, $N^{1,J}_i\to\eta^{1,J}$
weakly in $\MMN(\ZZJ)$. This is
assertion (b).

\emph{Identification of the limit components.}
\emph{First marginal.} By (a), $\Law(u^1_i)=\Law(u_{n_{r_i}})$. Since
$n_{r_i}\in\NNtilde $ and $n_{r_i}\to\infty$, and since $\Law(u_n)\to\rho$ weakly on $\mathcal X$ along $\NNtilde $, $\Law(u^1_i)\to\rho$. On the other
hand $u^1_i\to u^1_\ast$ everywhere, so dominated convergence gives
$\EE^1\big[F(u^1_i)\big]\to\EE^1\big[F(u^1_\ast)\big]$ for every $F\in C_b(\mathcal X)$.
Weak limits on the metrisable space $\mathcal X$ are unique, whence
$\Law(u^1_\ast)=\rho$.

\emph{Point-measure property.} Fix $J$. By construction,
$\Phi^J_r(u_{n_r})$ takes values in the set
$\MMN(\ZZJ)$ of finite $\NN_0$-valued measures, which is
closed in the space of finite Borel measures on $\ZZJ$
(Lemma~\ref{lem:state-space} \ref{SS:ii}). By (a),
$\PP^1\bigl(N^{1,J}_i\in\MMN(\ZZJ)\bigr)=1$ for every $i$.
On the intersection of these full sets the everywhere limit
$\eta^{1,J}(\omega)$ lies in the closed set $\MMN(\ZZJ)$,
the continuity of the total-mass functional gives
$\eta^{1,J}(\ZZJ)=\lim_iN^{1,J}_i(\ZZJ)<\infty$, and the Borel evaluations
$\mu\mapsto\mu(B)$ make every mass $\eta^{1,J}(B)$
$\mathscr F^1$-measurable.
\end{proof}

In the next lemma we show that the family of reconstructed random measures is
consistent, i.e.\ that for $J\le J'$ the measure at level $J$ is the
restriction of the measure at level $J'$ to $Z^c_J\times[0,T]$, $\PP^1$-almost surely.

\begin{lemma}[{Consistency in $J$}]\label{lem:J-consistency}
For all $J'\ge J$, every $\mathcal{D} \in\mathscr Z$ with
$\mathcal{D} \subset\bigcup_{j\le J}\Delta Z_j$ and every $I\in\mathscr B([0,T])$,
\[
\eta^{1,J'}(\mathcal{D}\times I)\;=\;\eta^{1,J}(\mathcal{D}\times I)
\qquad\PP^1\text{-a.s.},
\]
the reconstructed measures at different truncation levels agree on the
common layers.
\end{lemma}

\begin{proof}
Fix $J<J'$ and write $\Gamma_J:=Z^c_J\times[0,T]$ with
$Z^c_J=\bigcup_{j\le J}\Delta Z_j$. For $r\ge\max\{r_0(J),r_0(J')\}$, on
$\mathcal G_{r,J}\cap\mathcal G_{r,J'}$, Lemma~\ref{lem:X-good-event}
applies at both levels. A jump with mark in $\Delta Z_j$, $j\le J$, is
registered at both levels and in the same layer $j$; a jump with mark in
$\Delta Z_{J+1}\cup\dots\cup\Delta Z_{J'}$ is registered only at level $J'$,
in a layer larger than $J$; and no other cell is registered. The recovered
marks at the two levels need not coincide, since the decompensation
$\lbar^{(J)}_{k,m,n}$ in $\xhat_k$ depends on $J$, but by
\eqref{eq:X-Dm-bound} both are within $a_r$ of the true mark. Hence the
atoms of $\Phi^{J'}_r(u_{n_r})\rest_{\Gamma_J}$ and of
$\Phi^J_r(u_{n_r})$ can be matched one to one with $\rho$-distance at most
$2b_r$, i.e.
\[
d_{\mathrm v}\Bigl(\Phi^{J'}_r(u_{n_r})\rest_{\Gamma_J},\,\Phi^J_r(u_{n_r})\Bigr)\le2b_r
\qquad\text{on }\mathcal G_{r,J}\cap\mathcal G_{r,J'} .
\]
The restriction map is Borel, so by
Proposition~\ref{prop:JS-pairs-allJ}\,\ref{eq:JS-laws},
$\PP^1\bigl(d_{\mathrm v}(N^{1,J'}_i\rest_{\Gamma_J},N^{1,J}_i)>2b_{r_i}\bigr)
\le c_{r_i,J}+c_{r_i,J'}\to0$. By the comparison argument of (A1),
$\Law(N^{1,J'}_i)=\Law(\Phi^{J'}_{r_i}(u_{n_{r_i}}))\Rightarrow\pi_{J'}$, and
$N^{1,J'}_i\to\eta^{1,J'}$ $\PP^1$-a.s., so $\Law(\eta^{1,J'})=\pi_{J'}$; since $\epsilon_{J+1}\in\mathcal R(\nu)$, the
limit $\eta^{1,J'}$ has $\PP^1$-a.s.\ no atom on
$\{|z|_Z=\epsilon_{J+1}\}\times[0,T]$, so the restriction to $\Gamma_J$ is
continuous at $\eta^{1,J'}$ and
$N^{1,J'}_i\rest_{\Gamma_J}\to\eta^{1,J'}\rest_{\Gamma_J}$ $\PP^1$-a.s.\
by \ref{JS:as}. Together with $N^{1,J}_i\to\eta^{1,J}$ and $b_{r_i}\to0$
this gives $d_{\mathrm v}(\eta^{1,J'}\rest_{\Gamma_J},\eta^{1,J})=0$
$\PP^1$-a.s., which is the claim.
\end{proof}

\medskip

\item
\noindent\textbf{The path jump measure and the limit object.}
In this step we analyse the limit object.
Let us define
for an arbitrary  path $x\in\DD([0,T];H)$ put $\Delta x(s):=x(s)-x(s-)$ and define, for
$\jt\in\{0,\dots,J\}$,
\[
\mathcal J_\jt(x)
:=\sum_{0<s\le T\atop \Delta x(s)\not=0}
\mathds 1_{\big\{\{(x(s-),\,\Delta x(s))\in Ran(G_j)\}\big\}}\,
\delta_{\bigl(\zhat_\jt(x(s-),\Delta x(s)),\;s\bigr)},
\qquad
\mathcal J^J(x):=\sum_{\jt=0}^{J}\mathcal J_\jt(x).
\]
Observe, for $x\in\DD([0,T];H)$, $ \mathcal J_\jt(x)$
is a point measure on $Z^c_J\times[0,T]$. Here $Ran(G_j)$ is the graph set
of Lemma~\ref{lem:DH-meas},
 \emph{viewed as a Borel subset of
$H\times H$} (Lemma~\ref{lem:DH-meas-H} \ref{HH:ii}), and\footnote{$\zhat$ is defined in Lemma \ref{lem:DH-meas}.}
$\zhat_\jt=\pi_Z\circ G^{-1}_\jt$
is the measurable inverse, defined on
$Ran(G_j)$. Here, equivalently, the membership test is written as

$\Delta x(s)\in\Delta H_\jt(x(s-))$ under the convention
$\Delta H_\jt(v)=\emptyset$ for $v\in H\setminus V$ of
Lemma~\ref{lem:DH-meas-H}  \ref{HH:iii}. With this reading, the summand is well
defined for \emph{every} path, and
$\mathcal J^J:\DD([0,T];H)\to\mathbb M\bigl(Z^c_J\times[0,T]\bigr)$ is
Borel measurable. A c\`adl\`ag path has at most countably many jumps; its
successive jump times and the pairs $(x(s-),\Delta x(s))$ are Borel functionals
on $\DD([0,T];H)$ (Lemma~\ref{lem:jump-measurable} with $E=H$), and the membership tests and the inverse are Borel in
$H\times H$ by Lemma~\ref{lem:DH-meas-H}. The sets $Ran(G_0),\dots,Ran(G_J)$ have
pairwise disjoint $v$-sections, so each jump is registered in at most one
layer.

\begin{remark}\label{rem:J-mass}
For each $\jt\le J$, the total mass
	$x\mapsto\mathcal J_\jt(x)\bigl(Z^c_J\times[0,T]\bigr)\in[0,\infty]$, i.e.\ the number of
	jumps of $x$ registered in the layer $\jt$, is Borel: it is a countable sum of Borel
	indicators over a measurable enumeration of the jumps (Lemma~\ref{lem:jump-measurable}).
	On the Borel set where it is infinite we set $\mathcal J_\jt(x):=0$. Since this
	convention is applied to each layer separately, $\mathcal J_\jt$ does not depend on $J$,
	and $\mathcal J^J(x)=\sum_{\jt\le J}\mathcal J_\jt(x)$ is a finite measure for every path
	$x$, so that $\mathcal J^J$ takes values in $\mathbb M\bigl(Z^c_J\times[0,T]\bigr)$. For
	the solutions $u_n$ the convention is $\PP$-a.s.\ never active: by
	Lemma~\ref{lem:jump-identity}, the total mass of $\mathcal J^J(u_n)$ is Poisson
	distributed with parameter $\nu(Z^c_J)\,T$.
\end{remark}

\medskip

As next we state a lemma whose result are classical.
We nevertheless state and prove the lemma here, for the following
reason.
In our construction we approximate the underlying Poisson random measure by a point process. On the other side,
is the integrand of the stochastic integral  against a Poisson random measure defined  as an equivalence class: Two integrands agreeing $\PP\otimes \Leb\otimes \n$-a.e.~are identified, because the stochastic integral  cannot distinguish them.
The reconstruction, however, must
	\emph{evaluate} integrands at individual atoms $(\tau,z)$, and a single point
	is itself a null set.  For an equivalence class, the value there is undefined.
	So, one may
	fix an everywhere-defined representative, show (via the compensation formula) that the atoms a.s.~avoid the null set where any two representatives disagree, so the value at each atom is well defined independently of the choice, and show that the jump of the integral process at the atom equals that value.
The Lemma is used in this pointwise form in
	Step~0 of Lemma~\ref{lem:jump-identity}, in the decomposition \eqref{eq:xhat-exact}, and in Lemma~\ref{lem:X-good-event}.

\begin{lemma}%
\label{lem:pathwise-atoms}
	Let $\eta$ be a time-homogeneous Poisson random measure on $(Z,\mathscr Z )$ with
	intensity $\nu$ over the filtered (right continuous filtration)  probability space $\mathfrak A_0=(\Omega_0,\CF_0,(\CF_t^0)_{t\in[0,T]},\PP_0)$, let $U\in\mathscr Z $ with $\nu(U)<\infty$, and
	let $p\in(1,2]$ be as in Theorem~\ref{th-ineq}. Let $\mathscr{P}$ be the predictable $\sigma$-algebra.
	\begin{enumerate}[label=\upshape(\roman*)]
		\item\label{PA:null} If $N\subset\Omega\times(0,T]\times Z$ is
		 $\mathscr{P}\otimes \mathscr Z $-measurable with
		$(\PP\otimes\Leb\otimes\nu)(N)=0$, then $\PP$-a.s.\ no atom $(\tau,z)$
		of $\eta$ with $\tau\le T$ lies in $N$.
		\item\label{PA:distinct} $\PP$-a.s., the atoms of $\eta$ in
		$Z\times(0,T]$ have pairwise distinct time coordinates.
		\item\label{PA:version} If $\xi,\xi'$ are predictable $H$-valued
		processes on $\Omega\times(0,T]\times Z$ with $\xi=\xi'$
		$\PP\otimes\Leb\otimes\nu$-a.e., then $\PP$-a.s.\
		$\xi(\tau,z)=\xi'(\tau,z)$ at every atom $(\tau,z)$ with $z\in U$.
		\item\label{PA:pathwise} If $\xi\in\mathcal M^p(0,T;L^p(Z,\nu;H))$
		admits an everywhere-defined predictable representative (again denoted
		$\xi$), then, $\PP$-a.s., for every $t\in[0,T]$,
		\begin{equation}\label{eq:pathwise-integral}
			\int_0^t\!\!\int_U\xi(s,z)\,\widetilde\eta(dz,ds)
			=\sum_{(\tau_i,z_i):\,z_i\in U,\ \tau_i\le t}\xi(\tau_i,z_i)
			\;-\;\int_0^t\!\!\int_U\xi(s,z)\,\nu(dz)\,ds ,
		\end{equation}
		the sum being $\PP$-a.s.\ finite and the compensator integral absolutely
		convergent. In particular, the c\`adl\`ag modification of the left-hand
		side jumps  only at the atom times $\tau_i$ with $z_i\in U$, and its jump at $\tau_i$ is
$\xi(\tau_i,z_i)$, which may vanish; by \ref{PA:version}, the right-hand side does
		not depend on the choice of the predictable representative.
	\end{enumerate}
\end{lemma}

\begin{proof}
Throughout, $\mathscr P$ denotes the predictable $\sigma$-field on
$\Omega\times(0,T]$; it is generated by the sets $F\times(a,b]$ with
$0\le a<b\le T$ and $F\in\CF_a$. For $U'\in\mathscr Z$ we write
$\Omega\times(0,T]\times U'$ for the corresponding subset of
$\Omega\times(0,T]\times Z$, and $(Z^c_{\epsilon_j})_{j\in\NN}$ for the
exhausting sets $Z^c_{\epsilon_j}\uparrow Z\setminus\{0^\ast\}$ with
$\nu(Z^c_{\epsilon_j})<\infty$.

\ref{PA:null} We prove the \emph{compensation formula}: for every
$\mathscr P\otimes\mathscr Z$-measurable $f:\Omega\times(0,T]\times Z\to[0,\infty]$,
\begin{equation}\label{eq:PA-comp}
	\EE\int_{(0,T]}\!\int_Z f(s,z)\,\eta(dz,ds)
	=\EE\int_0^T\!\!\int_Z f(s,z)\,\nu(dz)\,ds .
\end{equation}
Fix $U'\in\mathscr Z$ with $\nu(U')<\infty$, and let $\mathcal R$ be the class
of sets $F\times(a,b]\times B$ with $0\le a<b\le T$, $F\in\CF_a$ and
$B\in\mathscr Z$, $B\subset U'$. It is a $\pi$-system, it contains
$\Omega\times(0,T]\times U'$, and it generates
$(\mathscr P\otimes\mathscr Z)\cap(\Omega\times(0,T]\times U')$. For
$A=F\times(a,b]\times B\in\mathcal R$, Definition~\ref{def-Prm}\,(iii) gives
that $\eta(B\times(a,b])$ is independent of $\CF_a$ and has mean
$(b-a)\nu(B)$, hence
\[
\EE\bigl[\mathds1_F\,\eta(B\times(a,b])\bigr]=\PP(F)\,(b-a)\,\nu(B)
=\EE\int_0^T\!\!\int_Z\mathds1_A\,\nu(dz)\,ds .
\]
Both sides of \eqref{eq:PA-comp}, evaluated at $f=\mathds1_A$, define finite
measures on $(\mathscr P\otimes\mathscr Z)\cap(\Omega\times(0,T]\times U')$,
of total mass $T\nu(U')$. (The left-hand side is measurable in $\omega$, since
$\eta$ is a random measure and $\mathds1_A$ is jointly measurable.) They agree
on $\mathcal R$, so by Dynkin's $\pi$--$\lambda$ theorem they agree on the
generated $\sigma$-field. By linearity and monotone convergence,
\eqref{eq:PA-comp} holds for every nonnegative $f$ vanishing outside
$\Omega\times(0,T]\times U'$. Taking $U'=Z^c_{\epsilon_j}$ and letting
$j\to\infty$ gives \eqref{eq:PA-comp} for $f$ vanishing on
$\Omega\times(0,T]\times\{0^\ast\}$, by monotone convergence. Finally,
$\eta(\{0^\ast\}\times(0,T])$ is Poisson distributed with parameter
$T\nu(\{0^\ast\})=0$ by Definition~\ref{def-Prm}\,(i), so this set carries no
atom, $\PP$-a.s., and both sides of \eqref{eq:PA-comp} ignore it. Applying
\eqref{eq:PA-comp} to $f=\mathds1_N$ gives $\EE\,\eta(N)=0$, i.e.\
$\eta(N)=0$, $\PP$-a.s., which is \ref{PA:null}.

\ref{PA:distinct} Fix $j\in\NN$, put $U'=Z^c_{\epsilon_j}$, and let
$I^m_k=(kT2^{-m},(k+1)T2^{-m}]$. If two atoms in $U'\times(0,T]$ have the same
time coordinate, then for every $m$ some cell $I^m_k$ contains at least two
atoms with mark in $U'$. Since $\eta(U'\times I^m_k)$ is Poisson distributed
with parameter $\mu_m:=\nu(U')T2^{-m}$ and $\PP(\mathrm{Poi}(\mu)\ge2)\le\mu^2$,
\[
\PP\bigl(\text{two atoms in }U'\times(0,T]\text{ share a time}\bigr)
\le\sum_{k=0}^{2^m-1}\mu_m^2=T^2\nu(U')^2\,2^{-m}\xrightarrow[m\to\infty]{}0 .
\]
Taking the union over $j\in\NN$, and recalling from the proof of
\ref{PA:null} that there is no atom with mark $0^\ast$, proves \ref{PA:distinct}.

\ref{PA:version} Apply \ref{PA:null} to the
$\mathscr P\otimes\mathscr Z$-measurable null set
$N=\{\xi\ne\xi'\}\cap(\Omega\times(0,T]\times U)$.

\ref{PA:pathwise} Write $\xi$ for the given everywhere-defined predictable
representative, i.e.\ a $\mathscr P\otimes\mathscr Z$-measurable map
$\xi:\Omega\times(0,T]\times Z\to H$ with
$\EE\int_0^T\int_Z|\xi(s,z)|_H^p\,\nu(dz)\,ds<\infty$. The left-hand side of
\eqref{eq:pathwise-integral} denotes the c\`adl\`ag modification, given by
Theorem~\ref{th-ineq}, of $t\mapsto I\bigl(\mathds1_{(0,t]\times U}\,\xi\bigr)$.

\emph{Step 1 (The right-hand side).} Since $(\PP\otimes\Leb\otimes\nu)$
restricted to $\Omega\times(0,T]\times U$ has total mass $T\nu(U)<\infty$,
H\"older's inequality gives, with $p'=p/(p-1)$,
\begin{equation}\label{eq:PA-L1}
	\EE\int_0^T\!\!\int_U|\xi|_H\,d\nu\,ds
	\le\bigl(T\nu(U)\bigr)^{1/p'}
	\Bigl(\EE\int_0^T\!\!\int_U|\xi|_H^p\,d\nu\,ds\Bigr)^{1/p}<\infty .
\end{equation}
By \eqref{eq:PA-comp} applied to the nonnegative predictable function
$|\xi|_H\mathds1_U$, also $\EE\int_{(0,T]}\int_U|\xi|_H\,d\eta<\infty$.
Moreover $\eta(U\times(0,T])$ is Poisson distributed with finite parameter.
Hence, $\PP$-a.s., $\eta$ has finitely many atoms $(\tau_i,z_i)$ in
$U\times(0,T]$, the compensator integral converges absolutely, and
\[
J_t(\xi):=\sum_{i:\ \tau_i\le t}\xi(\tau_i,z_i)
-\int_0^t\!\!\int_U\xi(s,z)\,\nu(dz)\,ds,\qquad t\in[0,T],
\]
is well defined for all $t$. Since
$|J_t(\xi)|_H\le\int_{(0,T]}\int_U|\xi|_H\,d\eta+\int_0^T\int_U|\xi|_H\,d\nu\,ds$
for all $t$, \eqref{eq:PA-comp} yields
\begin{equation}\label{eq:PA-J-bound}
	\EE\sup_{t\le T}\bigl|J_t(\xi)\bigr|_H\le2\,\EE\int_0^T\!\!\int_U|\xi|_H\,d\nu\,ds .
\end{equation}
The map $\xi\mapsto J(\xi)$ is linear. The compensator term is continuous in
$t$, and the atom sum is a right-continuous step function which, by
\ref{PA:distinct}, jumps exactly at the times $\tau_i$, with jump
$\xi(\tau_i,z_i)$. If $\xi'$ is another everywhere-defined predictable
representative, then \eqref{eq:PA-J-bound} applied to $\xi-\xi'$ gives
$\sup_t|J_t(\xi)-J_t(\xi')|_H=0$, $\PP$-a.s.; this is the independence of the
representative.

\emph{Step 2 (Elementary integrands).} Let $\mathcal E_U$ be the linear span of
the functions $\zeta=\mathds1_F\,\mathds1_{(a,b]}(s)\,\mathds1_B(z)\,e$ with
$0\le a<b\le T$, $F\in\CF_a$, $B\in\mathscr Z$, $B\subset U$, and $e\in H$.
Each such $\zeta$ is a step process in
$\mathcal M^p_{\rm step}(0,T;L^p(Z,\nu;H))$, and $\mathds1_{(0,t]}\zeta$ is of
the same form with $(a,b]$ replaced by $(a\wedge t,b\wedge t]$, the indicator
$\mathds1_F$ being $\CF_{a\wedge t}$-measurable whenever the interval is
nonempty. By \eqref{eqn-int-def},
\[
I\bigl(\mathds1_{(0,t]\times U}\,\zeta\bigr)
=\mathds1_F\,e\;\widetilde\eta\bigl(B\times(a\wedge t,b\wedge t]\bigr)
=\mathds1_F\,e\,\Bigl(\eta\bigl(B\times(a\wedge t,b\wedge t]\bigr)
-\nu(B)\,(b\wedge t-a\wedge t)\Bigr)=J_t(\zeta)
\]
for all $t$, $\PP$-a.s. The process on the left is therefore already
c\`adl\`ag, and by linearity \eqref{eq:pathwise-integral} holds for every
$\zeta\in\mathcal E_U$, for all $t$ simultaneously, $\PP$-a.s.

\emph{Step 3 (Density).} We claim that $\mathcal E_U$ is dense in
$L^p_U:=L^p\bigl(\Omega\times(0,T]\times U,\ (\mathscr P\otimes\mathscr Z)\cap(\Omega\times(0,T]\times U),\ \PP\otimes\Leb\otimes\nu;\ H\bigr)$.
Since $H$ is separable, simple functions $\sum_ke_k\mathds1_{A_k}$ with
$e_k\in H$ and $A_k$ predictable are dense in $L^p_U$, so it suffices to
approximate $\mathds1_A$ for every predictable $A\subset\Omega\times(0,T]\times U$.
Let $\mathcal D$ be the class of such $A$ for which $\mathds1_A$ lies in the
$L^p_U$-closure of the span of $\{\mathds1_R:R\in\mathcal R\}$, with
$\mathcal R$ as in the proof of \ref{PA:null} for $U'=U$. Then
$\mathcal R\subset\mathcal D$; $\mathcal D$ is closed under complements in
$\Omega\times(0,T]\times U$, since
$\mathds1_{A^c}=\mathds1_{\Omega\times(0,T]\times U}-\mathds1_A$; and
$\mathcal D$ is closed under increasing unions, by dominated convergence, the
measure being finite. By Dynkin's $\pi$--$\lambda$ theorem, $\mathcal D$
contains $\sigma(\mathcal R)$, which is the whole predictable $\sigma$-field on
$\Omega\times(0,T]\times U$.

\emph{Step 4 (Passage to the limit).} By Step~3 choose $\zeta_k\in\mathcal E_U$
with $\zeta_k\to\xi\mathds1_U$ in $L^p_U$, hence also in
$L^1(\Omega\times(0,T]\times U;H)$ by \eqref{eq:PA-L1}. By the maximal
inequality \eqref{I4} of Theorem~\ref{th-ineq},
\[
\EE\sup_{t\le T}\Bigl|I\bigl(\mathds1_{(0,t]\times U}\,\zeta_k\bigr)
-I\bigl(\mathds1_{(0,t]\times U}\,\xi\bigr)\Bigr|_H^p
\le C\,\EE\int_0^T\!\!\int_U|\zeta_k-\xi|_H^p\,d\nu\,ds\longrightarrow0,
\]
where both integral processes are the c\`adl\`ag modifications. By
\eqref{eq:PA-J-bound} applied to $\zeta_k-\xi$,
$\EE\sup_t|J_t(\zeta_k)-J_t(\xi)|_H\to0$. Passing to a subsequence, both
convergences hold uniformly in $t$, $\PP$-a.s. Since
$I(\mathds1_{(0,t]\times U}\zeta_k)=J_t(\zeta_k)$ for all $t$ by Step~2, the
limits coincide: $\PP$-a.s., $I(\mathds1_{(0,t]\times U}\xi)=J_t(\xi)$ for
every $t\in[0,T]$, which is \eqref{eq:pathwise-integral}. The finiteness of
the sum, the absolute convergence of the compensator integral, the jumps,
and the independence of the representative were shown in Step~1.
\end{proof}

\begin{remark}%
	\label{rem:pathwise-atoms-classical}
Lemma~\ref{lem:pathwise-atoms} is classical and can be found already in the literature.
However, we state it here, since it is important for our setting. What is important is that they are independent of the
dimension of the jump state $Z$, the arguments are
		purely measure-theoretic, requiring of $(Z,\mathscr Z)$ only that it be a
		measurable space carrying the $\sigma$-finite intensity $\nu$, so they apply
		verbatim when $Z$ is an infinite-dimensional Banach or general Polish space.
		Assertion \ref{PA:null} is the compensation formula, in its unconditional
		Poisson form the Mecke equation, valid for $\sigma$-finite intensities on
		arbitrary measurable state spaces; see \cite[Ch.~II, \S3]{ikeda1981} and
		\cite[Thm.~4.1]{LastPenrose2018}. Assertion \ref{PA:distinct} is built into the
		very formalism of \cite[Ch.~I, \S\S8--9]{ikeda1981}, where a point process
		\emph{is} a $Z$-valued point function on a countable time set one jump per
		time and the Poisson point process with $\sigma$-finite characteristic
		measure is constructed in that class (Thm.~I.9.1), for the jump measures of
		L\'evy processes cf.\ also \cite[Thm.~19.2]{sato1999}. Assertion
		\ref{PA:version} is an immediate corollary of \ref{PA:null}. Assertion
		\ref{PA:pathwise} is the ``integral of the first kind'' over a finite-activity
		window and is tackled by \cite{Ruediger2004} (comparison of the pathwise `natural' integral with the strong-$p$ integral), \cite{AlbeverioRuediger2005,RuedigerZiglio2006} and \cite{HausenblasBanach2010} for integrands with values in separable, respectively martingale-type $p$, Banach spaces, and  for the
		real and Hilbert-space-valued cases, \cite[Ch.~6]{PeszatZabczyk2007}.
\end{remark}

Let us remind, $\eta$ is a time-homogeneous Poisson random measure on $(Z,\mathscr Z )$ with
intensity measure $\nu$, defined over the filtered probability space
$\mathfrak A=(\Omega,\CF,(\CF_t)_{t\ge0},\PP)$ and compatible with
$(\CF_t)$ in the sense of Definition~\ref{def-Prm} having  the law $\mathrm{PRM}(\nu\otimes\Leb)$.
\begin{lemma}%
	\label{lem:jump-identity}
	Assume Hypotheses~\ref{hyp}, \ref{hyp:An} and \ref{ass:GF}, and let $\eta$
	be a time-homogeneous Poisson random measure with intensity $\nu$,
	compatible with $(\CF_t)$ in the sense of Definition~\ref{def-Prm}. Let
	$n\in\NN$ and let $u_n$ be an $(\CF_t)$-adapted solution of
	\eqref{originalun} on $\mathfrak A$, driven by $\eta$, which is c\`adl\`ag in
	$H$ and satisfies $u_n\in L^\alpha(0,T;V)$, $\PP$-a.s., and
	$\EE\int_0^T|u_n(s)|_V^\alpha\,ds<\infty$. Let
	$\mathcal J^J$ be the path jump measure defined above. Each layer measure
	$\mathcal J_\jt$ depends on $\jt$ and $g$ only, not on $J$, which is applied to each layer
	separately). Consequently, for
	$J'\ge J$,
	$\mathcal J^{J'}(x)\restriction_{Z^c_J\times[0,T]}
	=\mathcal J^{J}(x)$ for every $x\in\DD([0,T];H)$.

	Assume $\ep_0\ge R_0$, so that
	$\bigcup_{\jt\le J}\Delta Z_\jt=\{\epsilon_{J+1}<|z|_Z\le\epsilon_0\}$ carries
	all jumps of size $>\epsilon_{J+1}$ (Remark~\ref{rem:support}). Then,
	$\PP$-a.s.,
	\begin{itemize}
		\item $	u_n(\tau-)\in V$ at \emph{every} atom $(\tau,z)$ of $\eta$  {with} $\tau\le T$,
		\item and $$
		\mathcal J^J(u_n)
		=\sum_{\jt=0}^{J}\;
		\int_{(0,T]}\int_{\Delta Z_\jt}\delta_{(z,s)}\;\eta(dz,ds),$$
		i.e.\ $\mathcal J^J(u_n)$ is exactly the restriction of $\eta$ to
		$\bigl(\bigcup_{\jt\le J}\Delta Z_\jt\bigr)\times(0,T]$.
	\end{itemize}
	 In particular,
	$\mathcal J^J(u_n)$ is on $\mathfrak A$, a Poisson random measure on
	$Z^c_J\times[0,T]$ with intensity
	$\nu\rest_{Z^c_J}\otimes\Leb$, i.e.\ $\Law(\mathcal J^J(u_n))=\pi_J$.
\end{lemma}

\begin{proof}
\emph{Step 0 (The pre-jump state lies in $V$ at every atom).}
Fix $n\in\NN$. The process $s\mapsto u_n(s-)$ on $\Omega\times(0,T]$ is
$H$-valued, left-continuous and adapted, hence predictable. Note that the
right-continuous process $u_n(\cdot)$ is progressively measurable but
\emph{not} predictable, so it could not be used in its place below. By
Lemma~\ref{lem:DH-meas-H} \ref{HH:i} we have $V\in\mathscr B(H)$, so that the set
\[
\mathcal D_n:=\bigl\{(\omega,s)\in\Omega\times(0,T]:\ u_n(s-,\omega)\notin V\bigr\}
\]
is a predictable set. It is $(\PP\otimes\Leb)$-null, since, for $\PP$-a.e.\
$\omega$ the path $u_n(\cdot,\omega)$ is c\`adl\`ag in $H$, so
$u_n(s-,\omega)=u_n(s,\omega)$ off the at most countable set of its jump
times, while $u_n(s,\omega)\in V$ for $\Leb$-a.e.\ $s$: membership in
$\mathcal X=\DD([0,T];H)\cap L^\alpha(0,T;V)$ is membership of the
everywhere-defined c\`adl\`ag function itself, i.e.\ it \emph{means}
$u_n(s,\omega)\in V$ for a.e.\ $s\in[0,T]$ together with
$\int_0^T|u_n(s,\omega)|^\alpha_V\,ds<\infty$. Hence
$\Leb\bigl((\mathcal D_n)_\omega\bigr)=0$ for $\PP$-a.e.\ $\omega$, and the
Tonelli Theorem gives $(\PP\otimes\Leb)(\mathcal D_n)=0$.

Since $\nu(\{0^\ast\})=0$ and $Z^c_{\epsilon_j}\uparrow Z\setminus\{0^\ast\}$ with
$\nu(Z^c_{\epsilon_j})<\infty$ (Remark~\ref{rem:support}), the intensity $\nu$ is
$\sigma$-finite on $Z$, so $\mathcal D_n\times Z$ is
$\mathscr P\otimes\mathscr Z$-measurable and, as the product of a
$(\PP\otimes\Leb)$-null set with $Z$, it is a
$(\PP\otimes\Leb\otimes\nu)$-null set, although $\nu(Z)=\infty$ is allowed.
Since $u_n$ is $(\CF_t)$-adapted and $\eta$ is compatible with $(\CF_t)$,
Lemma~\ref{lem:pathwise-atoms}\,\ref{PA:null} applies: $\PP$-a.s.\ \emph{no} atom
$(\tau,z)$ of $\eta$ with $\tau\le T$ lies in $\mathcal D_n\times Z$, i.e.
\[
u_n(\tau-)\in V\qquad\text{at every atom }(\tau,z)\text{ of }\eta,\ \tau\le T,
\qquad\PP\text{-a.s.}
\]

	\emph{Step 1 (Identification of the jumps).}
	Let $\tilde\xi(s,z):=g\bigl(u_n(s-),z\bigr)\mathds1_{\{u_n(s-)\in V\}}$.
	\emph{Measurability:} the embedding $V\hookrightarrow H$ is continuous and
	injective between Polish spaces, so by the Lusin--Souslin theorem
	\cite[Thm.~15.1]{Kechris1995} every Borel subset of $V$ is a Borel subset of
	$H$. Hence the $V$-valued map
	$(\omega,s)\mapsto u_n(s-)\mathds1_{\{u_n(s-)\in V\}}$ (with an arbitrary
	fixed value in $V$ off $\{u_n(s-)\in V\}$) is $\mathscr P$-measurable, and,
	$g$ being $\mathscr B(V)\otimes\mathscr Z$-measurable, $\tilde\xi$ is
	$\mathscr P\otimes\mathscr Z$-measurable and defined everywhere.
	\emph{Version:} since $u_n(s-)=u_n(s)$ off the at most countable set of jump
	times and $u_n(s)\in V$ for a.e.\ $s$, $\tilde\xi$ coincides
	$\PP\otimes\Leb\otimes\nu$-a.e.\ with the integrand $g(u_n(s),z)$ of
	\eqref{originalun}. Both therefore define the same stochastic integral, and
	we work with the predictable representative $\tilde\xi$. As both sides of
	\eqref{originalun} are right-continuous in $V^\ast$, the identity holds,
	$\PP$-a.s., for all $t\in[0,T]$ simultaneously, with the c\`adl\`ag
	modification $M$ of the stochastic integral. Split the noise term of
	\eqref{originalun} into the layer part over
	$U:=\bigcup_{\jt\le J}\Delta Z_\jt$, $\nu(U)<\infty$, and the small-jump part
	over $Z_{\epsilon_{J+1}}$. By \eqref{eq:g-second-moment1}, the moment bound
	in the hypotheses and $\alpha\ge2$,
	\[
	\EE\int_0^T\!\!\int_Z|\tilde\xi(s,z)|_H^2\,\nu(dz)\,ds
	\le C_g\,\EE\int_0^T\bigl(1+|u_n(s)|_V^2\bigr)\,ds<\infty ,
	\]
	so $\tilde\xi\in\mathcal M^2(0,T;L^2(Z,\nu;H))$, and
	Lemma~\ref{lem:pathwise-atoms} applies with $p=2$. By
	Lemma~\ref{lem:pathwise-atoms} \ref{PA:pathwise}, the layer part is the
	pathwise sum \eqref{eq:pathwise-integral}, it jumps exactly at the layer
	atoms, with jump $\tilde\xi(\tau,z)$.

	\emph{The small-jump part.} Let
	$M_\flat(t):=\int_0^t\!\int_{Z_{\epsilon_{J+1}}}\tilde\xi\,d\widetilde\eta$
	(c\`adl\`ag modification). This integral has infinite activity and is
	defined only as an $L^2$-limit, so
	Lemma~\ref{lem:pathwise-atoms}\,\ref{PA:pathwise} does not apply to it
	directly; we truncate. For $k>J$ let
	$M^{(k)}(t):=\int_0^t\!\int_{\{\epsilon_{k+1}<|z|_Z\le\epsilon_{J+1}\}}\tilde\xi\,d\widetilde\eta$.
	The window has finite $\nu$-measure, so by
	Lemma~\ref{lem:pathwise-atoms}\,\ref{PA:pathwise} the process $M^{(k)}$
	jumps exactly at the atoms $(s,z)$ with $\epsilon_{k+1}<|z|_Z\le\epsilon_{J+1}$,
	with jump $\tilde\xi(s,z)$. By Doob's inequality, the It\^o isometry and
	\eqref{eq:g-second-moment1},
	\[
	\EE\sup_{t\le T}\bigl|M_\flat(t)-M^{(k)}(t)\bigr|_H^2
	\le4\,\EE\int_0^T\!\!\int_{\{0<|z|_Z\le\epsilon_{k+1}\}}
	\bigl|g(u_n(s),z)\bigr|_H^2\,\nu(dz)\,ds\xrightarrow[k\to\infty]{}0,
	\]
	by dominated convergence, since $\nu(\{0^\ast\})=0$ and the inner integral
	is dominated by $C_g(1+|u_n(s)|_V^2)$, which is $\PP\otimes\Leb$-integrable
	by the moment bound in the hypotheses and $\alpha\ge2$. Along a subsequence,
	$\sup_t|M_\flat(t)-M^{(k)}(t)|_H\to0$ $\PP$-a.s., and uniform convergence of
	c\`adl\`ag paths implies convergence of all jumps. Together with
	Lemma~\ref{lem:pathwise-atoms}\,\ref{PA:distinct} and the absence of atoms
	with mark $0^\ast$ (Lemma~\ref{lem:pathwise-atoms}\,\ref{PA:null} applied to
	$\Omega\times(0,T]\times\{0^\ast\}$), this gives, $\PP$-a.s.:
	$\Delta M_\flat(s)=\tilde\xi(s,z')$ if $(s,z')$ is an atom of $\eta$ with
	$0<|z'|_Z\le\epsilon_{J+1}$, and $\Delta M_\flat(s)=0$ at every other $s$.
	In particular, $M_\flat$ is continuous at every layer-atom time, since by
	Lemma~\ref{lem:pathwise-atoms}\,\ref{PA:distinct} no small atom shares its
	time with a layer atom.

	\emph{The drift part.}
	The drift part is continuous in $V^\ast$: by (H4$_n$) of Hypothesis~\ref{hyp:An}
	together with $u_n\in L^\alpha(0,T;V)$ and $\sup_{t\le T}|u_n(t)|_H<\infty$
	$\PP$-a.s., one has
	$A_n(\cdot,u_n(\cdot))\in L^{\alpha/(\alpha-1)}(0,T;V^\ast)\subset L^1(0,T;V^\ast)$
	$\PP$-a.s., so $t\mapsto\int_0^t A_n(s,u_n(s))\,ds$ is absolutely continuous as a
	$V^\ast$-valued function and has no jump in $V^\ast$. (It need not be
	$H$-valued at all, which is why the argument is carried out in $V^\ast$.)
	Denoting by $\iota:H\hookrightarrow V^\ast$ the embedding of the Gelfand
	triple and by $M$ the stochastic integral term, both $u_n$ and $M$ are
	c\`adl\`ag in $H$, hence c\`adl\`ag in $V^\ast$ since $\iota$ is continuous,
	and $\iota(\Delta x(t))=\Delta(\iota x)(t)$ for either of them. Taking jumps
	at time $t$ in $V^\ast$ in \eqref{originalun} therefore gives
	$\iota(\Delta u_n(t))=\iota(\Delta M(t))$. Both $\Delta u_n(t)$ and
	$\Delta M(t)$ lie in $H$, and $\iota$ is injective because $V$ is dense in
	$H$ (Hypothesis~\ref{hyp}); consequently
	$\Delta u_n(t)=\Delta M(t)$ in $H$ for every $t$.
		Consequently, by Step~0 (so that the indicator in $\tilde\xi$ equals $1$ at
	every atom), $\PP$-a.s., $\Delta u_n(\tau)=g(u_n(\tau-),z)$ at every layer
	atom $(\tau,z)$; $\Delta u_n(s)=g(u_n(s-),z')$ at every small atom $(s,z')$;
	and $\Delta u_n(s)=0$ at every time $s$ that is not an atom time.

	\emph{Step 2 (Every layer atom is registered in the correct
	layer with the correct mark).}
	Let $(\tau,z)$ be an atom with $z\in\Delta Z_\jt$, $\jt\le J$. By Step~0,
	$u_n(\tau-)\in V$, and by the weak nondegeneracy \eqref{it:nondeg},
	\[
	\bigl|\Delta u_n(\tau)\bigr|_H
	=\bigl|g\bigl(u_n(\tau-),z\bigr)\bigr|_H
	\;\ge\;\rho_g\bigl(|u_n(\tau-)|_V\bigr)\,|z|_Z
	\;>\;0 ,
	\]
	since $\rho_g>0$ and $|z|_Z>\epsilon_{J+1}>0$; only the positivity of
	$\rho_g$ at the finite value $|u_n(\tau-)|_V$ is used, which Step~0
	provides, and no uniform lower bound is needed. So the path jumps at $\tau$.
	Without \eqref{it:nondeg} a vanishing jump would make $\mathcal J^J$ miss an
	atom that the restriction of $\eta$ counts, and the identity would fail. Moreover
	$(u_n(\tau-),\Delta u_n(\tau))=G_\jt(u_n(\tau-),z)\in Ran(G_j)$, so the
	membership test of (X.1) is true for the layer $\jt$, and the recovered mark is
	$\zhat_\jt\bigl(u_n(\tau-),\Delta u_n(\tau)\bigr)=z$ by the uniqueness of the
	mark, Lemma~\ref{lem:DH-meas} \ref{L:v} (in the $H\times H$ reading of
	Lemma~\ref{lem:DH-meas-H} \ref{HH:iii}). By the pairwise disjointness of the
	layers (Lemma~\ref{lem:DH-meas} \ref{L:ii}), the atom is registered exactly
	once.

	\emph{Step 3 (Nothing else is registered).}
	Let $s$ be a jump time of $u_n$ that is not a layer-atom time. By Step~1,
	$s$ is a small-jump atom time, and
	$\Delta u_n(s)=g(u_n(s-),z')$ with $0<|z'|_Z\le\epsilon_{J+1}$ and
	$u_n(s-)\in V$ by Step~0 (atoms with mark $z'=0^\ast$ are excluded,
	$\PP$-a.s., by Lemma~\ref{lem:pathwise-atoms}\,\ref{PA:null}, since
	$\nu(\{0^\ast\})=0$). Registration
	of this jump in a layer $\jt\le J$ would mean
	$g(u_n(s-),z')=g(u_n(s-),z)$ for some $z\in\Delta Z_\jt$, since $z\ne z'$ and
	both are nonzero, this contradicts the injectivity of $g(u_n(s-),\cdot)$ on
	$Z\setminus\{0^\ast \}$, which follows from
	Hypothesis~\ref{ass:GF}\,\ref{it:inj} applied with any
	$\epsilon<\min(|z|_Z,|z'|_Z)$. Hence $\mathcal J^J(u_n)$ has no atoms beyond
	those of Step~2, and summing over the layer atoms gives the identity. The
	exceptional sets used above (Step~0, Lemma~\ref{lem:pathwise-atoms}\,\ref{PA:distinct},
	the subsequence in the small-jump part, the atoms with mark $0^\ast$, and the
	sets where \eqref{originalun} or the pathwise identity
	\eqref{eq:pathwise-integral} fail) are countably many $\PP$-null sets, so
	all statements hold simultaneously, $\PP$-a.s. The last assertion follows
	from Definition~\ref{def-Prm}, since the restriction of a Poisson random
	measure to a measurable set is a Poisson random measure with the restricted
	intensity, and $\nu(\cup_{\jt\le J}\Delta Z_\jt)=\nu(Z^c_J)<\infty$.
\end{proof}

\medskip
To summarise, at each fixed level $J$ the limit object $\eta^{1,J}$ of
Proposition~\ref{prop:JS-pairs-allJ} \emph{is} a Poisson random measure.

\begin{lemma}\label{lem:etaJ-PRM}
	For every $J\in\NN$, the limit $\eta^{1,J}$ of
	Proposition~\ref{prop:JS-pairs-allJ} is a Poisson random measure on
	$\ZZJ$ with intensity $\nu\rest_{Z^c_J}\otimes\Leb$, that is,
	$\Law(\eta^{1,J})=\pi_J$, in the sense of
	Definition~\ref{def-Prm}\,(i)--(ii). The filtration condition (iii) is
	established in Step~(XII).
\end{lemma}

\begin{proof}
	By the pathwise identity $\mathcal J^J(u_{n_r})=\eta\rest_{\ZZJ}$ of
	Lemma~\ref{lem:jump-identity}, $\Law(\mathcal J^J(u_{n_r}))=\pi_J$ for
	every $r$. By Lemma~\ref{lem:X-good-event},
	$d_{\mathrm v}(\Phi^J_r(u_{n_r}),\mathcal J^J(u_{n_r}))\le b_r$ on
	$\mathcal G_{r,J}$ and $\PP(\mathcal G^c_{r,J})\le c_{r,J}$, where
	$b_r,c_{r,J}\to0$. Hence, for every bounded Lipschitz
	$f:\MMN(\ZZJ)\to\RR$,
	\[
	\Bigl|\EE f\bigl(\Phi^J_r(u_{n_r})\bigr)-\int f\,d\pi_J\Bigr|
	\le\mathrm{Lip}(f)\,b_r+2\|f\|_\infty\,c_{r,J}\longrightarrow0,
	\]
	so $\Law(\Phi^J_r(u_{n_r}))\Rightarrow\pi_J$; this is the comparison
	argument of part (A1) of Step~(IX). By
	Proposition~\ref{prop:JS-pairs-allJ}\,\ref{JS:as}, $\PP^1$-a.s.\
	$N^{1,J}_i\to\eta^{1,J}$ along the subsequence $(r_i)$, hence in
	distribution, and $\Law(N^{1,J}_i)=\Law(\Phi^J_{r_i}(u_{n_{r_i}}))$ by
	\eqref{eq:JS-laws}. Limits in distribution being unique on the Polish
	space $\MMN(\ZZJ)$, $\Law(\eta^{1,J})=\pi_J$. Finally, by
	Proposition~\ref{prop:JS-pairs-allJ}\,\ref{JS:point}, $\eta^{1,J}$ is
	$\PP^1$-a.s.\ a finite integer-valued measure whose masses are random
	variables, so ``Poisson random measure'' has its full meaning of
	Definition~\ref{def-Prm}\,(i)--(ii).
\end{proof}

	\medskip
\item  Identification of the limit of the random measure.

Proposition~\ref{prop:full-levy-patch} patches the
family $\{\eta^{1,J}\}_{J\in\NN}$ into a single random measure $\eta^1$ on
$Z\times(0,T]$, the limit being monotone and consistent by
Lemma~\ref{lem:J-consistency}, and $\eta^1$ satisfies
Definition~\ref{def-Prm}\,(i)--(ii) with the \emph{full} intensity $\nu$.

\begin{proposition}[{Patching: the full L\'evy measure}]\label{prop:full-levy-patch}
	Recall that $\nu(\{0^\ast\})=0$ by Hypothesis~\ref{ass:GF}. Set $Z_n:=\{d_Z(z,0^\ast)>\epsilon_{n+1}\}\cup\{0^\ast\}$,
	$n\in\NN$, so that $Z_n\uparrow Z$ and $\nu(Z_n)<\infty$, and define, for
	$\mathcal{D}\in\mathscr Z$ and $I\in\mathscr B([0,T])$,
	\begin{align*}
		\eta^1(\mathcal{D} \times I)
		:=\lim_{J\to\infty}\eta^{1,J}\Bigl(\bigl(\mathcal{D}\cap
		\textstyle\bigcup_{j\le J}\Delta Z_j\bigr)\times I\Bigr)
		\;=\sup\limits_{J\in\NN}\,
		\eta^{1,J}\Bigl(\bigl(\mathcal{D}\cap
		\textstyle\bigcup_{j\le J}\Delta Z_j\bigr)\times I\Bigr)
		\;\in\;\NNnull ,
	\end{align*}
	the limit being monotone (and consistent) by
	Lemma~\ref{lem:J-consistency}. More precisely, the proof of
	Lemma~\ref{lem:J-consistency} gives
	$d_{\mathrm v}\bigl(\eta^{1,J'}\rest_{Z^c_J\times[0,T]},\eta^{1,J}\bigr)=0$,
	i.e.\ $\eta^{1,J'}\rest_{Z^c_J\times[0,T]}=\eta^{1,J}$ as measures, and since
	there are only countably many pairs $J\le J'$, this holds for all of them
	outside one $\PP^1$-null set. Then $\eta^1$ is a measurable mapping
	$\Omega^1\to\mathbb M_{\NNnull }(\{Z_n\times(0,T]\})$ in the sense of
	Notation~\ref{notationref4}, and it satisfies conditions (i) and (ii) of
	Definition~\ref{def-Prm} with intensity measure $\nu$, that is,  for every
	$\mathcal{D}\in\mathscr Z$ with $\nu(\mathcal{D})<\infty$ and every $I$, the mass
	$\eta^1(\mathcal{D}\times I)$ is Poisson distributed with parameter
	$\nu(\mathcal{D})\,\Leb(I)$, and $\eta^1$ is independently scattered.
\end{proposition}

\begin{proof}%
	For $\mathcal{D} \subset\bigcup_{j\le J}\Delta Z_j$ the assertions are those of
	$\eta^{1,J}$ (Lemma~\ref{lem:etaJ-PRM}). For general $\mathcal{D}$
	with $\nu(\mathcal{D})<\infty$, since $\supp\nu\subset\{|z|_Z\le
	R_0\}\subset\{|z|_Z\le\epsilon_0\}$ (Remark~\ref{rem:support}) and
	$\nu(\{0^\ast\})=0$, the sets $\mathcal{D}\cap\Delta Z_j$, $j\in\NN_0$, cover $\mathcal{D}$ up
	to a $\nu$-null set. The masses
	$\eta^1((\mathcal{D}\cap\Delta Z_j)\times I)=\eta^{1,J}((\mathcal{D}\cap\Delta Z_j)\times I)$
	(any $J\ge j$) are Poisson with parameters $\nu(\mathcal{D}\cap\Delta Z_j)\Leb(I)$,
	and they are independent for distinct $j$: for finitely many $j$ this is
	the independent scattering of $\eta^{1,J}$ with $J$ large. Hence their sum
	$\eta^1(\mathcal{D}\times I)$ is, by monotone convergence of the partial sums and
	the stability of the Poisson family under independent sums and monotone
	limits, Poisson with parameter
	$\sum_j\nu(\mathcal{D}\cap\Delta Z_j)\Leb(I)=\nu(\mathcal{D})\Leb(I)<\infty$. In particular
	the mass is a.s.\ finite and $\mathscr F^1$-measurable as an a.s.\ limit
	of measurable maps. For the independent scattering, let
	$B_1,\dots,B_r\in\mathscr Z\otimes\mathscr B([0,T])$ be disjoint. For each
	finite $J$, the truncated masses
	$\eta^{1,J}\bigl(B_l\cap(Z^c_J\times[0,T])\bigr)$, $l\le r$, are independent
	by (ii) for $\eta^{1,J}$, and an a.s.\ monotone limit of vectors with
	independent components has independent components (characteristic
	functions and a $\pi$-system argument).
\end{proof}

\medskip

					\item

							\textbf{Construction of the filtration for $u^1_\ast$ and the random measure on $\mathfrak{A}^1$.}
Let us introduce the filtration on  $\mathfrak A^1$:

\begin{definition}%
	\label{def:pair-filtration}
	Let $(u^1_\ast,\eta^1)$ be as in Proposition~\ref{prop:JS-pairs-allJ} and
	Proposition~\ref{prop:full-levy-patch}, and let
	$\mathscr Z_0\subset\mathscr Z$ be a countable $\pi$-system generating
	$\mathscr Z$ and consisting of sets $\mathcal{D}$ with $\nu(\mathcal{D})<\infty$ and
	$\nu(\partial \mathcal{D})=0$ (e.g.\ finite intersections of sets
	$B\cap\{|z|_Z>\epsilon_j\}$, $j\in\NN$, where $B$ runs through the balls of a
	countable base whose radii avoid the at most countably many exceptional
	values; the cut at $\epsilon_j$ is needed, since balls around $0^\ast$ have
	infinite $\nu$-mass).
	For $t\in[0,T)$ put
\DEQS
	\mathscr F^1_t
	&:=&\bigcap_{s\in(t,T]}
	\sigma\Bigl(u^1_\ast(r),\ \eta^1\bigl(\mathcal{D}\times(0,r]\bigr)\;:\;
	r\le s,\ \mathcal{D}\in\mathscr Z_0\Bigr)\vee\mathcal N_{\PP^1},
\\
\mbox{and}\qquad 	\mathscr F^1_T &:=&\sigma\Bigl(u^1_\ast(r),\eta^1(\mathcal{D}\times(0,r]):r\le T,
	\mathcal{D}\in\mathscr Z_0\Bigr)\vee\mathcal N_{\PP^1}.
\EEQS
	Then $(\mathscr F^1_t)_{t\in[0,T]}$ is right-continuous and complete
	\emph{by construction} and $u^1_\ast$ is $(\mathscr F^1_t)$-adapted with
	c\`adl\`ag paths in $H$, hence progressively measurable and
	$t\mapsto\eta^1(\mathcal{D}\times(0,t])$ is $(\mathscr F^1_t)$-adapted for every
	$\mathcal{D}\in\mathscr Z$, by a monotone class argument from $\mathscr Z_0$.
	Adaptedness of \emph{both} coordinates holds by definition; nothing has
	to be proven at this step.
\end{definition}

\begin{lemma}[{Compatibility, functional form}]\label{lem:compat-transfer}
	{%
	For all $0\le s<t\le T$, $p,q,r\in\NN$,
	$s_1,\dots,s_p,\break s'_1,\dots,s'_q\in[0,s]$, and
	$\mathcal{D}_1,\dots,\mathcal{D}_q,\mathcal{C}_1,\dots,\mathcal{C}_r\in\mathscr Z$ with
	$\nu(\mathcal{D}_l),\nu(\mathcal{C}_l)<\infty$, and all bounded Borel functions
	$h:H^p\times\NN_0^q\to\RR$ and $g:\NN_0^r\to\RR$,
\DEQSZ \label{eq:compat-functional}
\lefteqn{ 		\EE^1\Bigl[
		h\bigl(u^1_\ast(s_1),\dots,u^1_\ast(s_p),
		\eta^1(\mathcal{D}_1\times(0,s'_1]),\dots,\eta^1(\mathcal{D}_q\times(0,s'_q])\bigr)\,}
		\\\nonumber
		&&\qquad \qquad \cdot
		g\bigl(\eta^1(\mathcal{C}_1\times(s,t]),\dots,\eta^1(\mathcal{C}_r\times(s,t])\bigr)
		\Bigr]
		\\\nonumber
	\qquad 	&&
		=\EE^1\bigl[	h\bigl(u^1_\ast(s_1),\dots,u^1_\ast(s_p),
		\eta^1(\mathcal{D}_1\times(0,s'_1]),\dots,\eta^1(\mathcal{D}_q\times(0,s'_q])\bigr)\,\bigr]\;
		\\\nonumber
		&&\quad \times \EE^1\bigl[	g\bigl(\eta^1(\mathcal{C}_1\times(s,t]),\dots,\eta^1(\mathcal{C}_r\times(s,t])\bigr)\bigr].
\EEQSZ
	Consequently, for every $U\in\mathscr Z$ with $\nu(U)<\infty$ and every
	$0\le s<t\le T$, the increment $\eta^1(U\times(s,t])$ is independent of
	$\mathscr F^1_s$, that is, condition (iii) of Definition~\ref{def-Prm} holds
	for $\eta^1$ over $(\mathfrak A^1,(\mathscr F^1_t))$, and, together with
	Proposition~\ref{prop:full-levy-patch}, $\eta^1$ is a time-homogeneous
	Poisson random measure with intensity measure $\nu$ over
	$(\mathfrak A^1,(\mathscr F^1_t))$, compatible with the filtration in the
	sense of conclusion (iii) of Theorem~\ref{Theoremmain}.}
\end{lemma}

\begin{proof}
	{%
	It suffices to prove \eqref{eq:compat-functional} with $\eta^1$ replaced
	by $\eta^{1,J}$ and all sets $\mathcal{C}_l,\mathcal{D}_l\subset\bigcup_{j\le J}\Delta Z_j$,
	for every fixed $J$.  The general case follows by monotone convergence
	($\eta^1$-masses are increasing limits of $\eta^{1,J}$-masses, Lemma~%
	\ref{lem:J-consistency} and independence is preserved under a.s.\
	limits, by characteristic functions). Fix $J$ and proceed in four steps.

	\emph{Step~1 (The identity on $\mathfrak A$, exactly).} For every $n$,
	$u_n$ is $\BF$-adapted (it solves \eqref{originalun}) and $\eta$ is
	compatible with $\BF$, so the $h$-functional of the pair $(u_n,\eta)$ is
	$\CF_s$-measurable, while the vector of increments
	$(\eta(\mathcal{C}_l\times(s,t]))_{l\le r}$ is independent of $\CF_s$ by
	Definition~\ref{def-Prm} (ii) and (iii) (jointly over the finitely many
	sets, cf.\ Lemma~\ref{eqdefprm}). Hence the pair $(u_n,\eta)$ satisfies
	the analogue of \eqref{eq:compat-functional} on $\mathfrak A$, with $\EE$ and
	$\CF_s$ in place of $\EE^1$ and $\mathscr F^1_s$. The law of that vector is an explicit
	Poisson functional of $\nu\otimes\Leb$, so $\EE[g(\cdot)]$ is a constant
	$c_g$ independent of $n$. By Lemma~\ref{lem:jump-identity}, $\eta$ may be
	replaced by $\mathcal J^J(u_n)$ on the truncated sets, so the identity
	holds exactly for the pair $(u_n,\mathcal J^J(u_n))$.

	\emph{Step~2 (From the marked noise to the reconstruction, up to
	$o(1)$).} Restrict to $h$ bounded and continuous, $g$ bounded, and to
	sets in the continuity $\pi$-system $\mathscr Z_0$. Let $\mathcal A$ be
	the finite family of the sets $\mathcal{D}_l\times(0,s'_l]$ and
	$\mathcal{C}_l\times(s,t]$, and, for $\delta>0$,
	\[
	\Gamma_\delta:=\bigcup_{A\in\mathcal A}\bigl\{x\in\ZZJ:
	\mathrm{dist}_\rho(x,A)\le\delta,\ \mathrm{dist}_\rho(x,\ZZJ\setminus A)\le\delta\bigr\},
	\]
	with the metric $\rho$ of \eqref{eq:X-Dm-bound}. The counts
	$\mu(A)$ are integer-valued and \emph{not} Lipschitz in
	$d_{\mathrm v}$: an atom close to $\partial A$ changes the count under an
	arbitrarily small shift. They are compared as follows. By
	Lemma~\ref{lem:X-good-event}, on $\mathcal G_{r,J}$ the atoms of
	$\Phi^J_r(u_{n_r})$ and of $\mathcal J^J(u_{n_r})$ correspond one-to-one,
	with time shift $\le\tau_{m_r}/2$ and mark error $\le a_r$, so partners
	have $\rho$-distance $\le b_r$. If an atom $x$ of $\mathcal J^J(u_{n_r})$
	and its partner $y$ were separated by some $A\in\mathcal A$, then
	$x\in\Gamma_{b_r}$. Hence, on $\mathcal G_{r,J}\cap\{\mathcal
	J^J(u_{n_r})(\Gamma_{b_r})=0\}$, all counts entering $h$ and $g$ coincide
	for the two measures, and, since $\Law(\mathcal J^J(u_{n_r}))=\pi_J$ for
	every $r$,
	\[
	\Bigl|\EE\bigl[(h\,g)\bigl(u_{n_r},\Phi^J_r(u_{n_r})\bigr)\bigr]
	-\EE\bigl[(h\,g)\bigl(u_{n_r},\mathcal J^J(u_{n_r})\bigr)\bigr]\Bigr|
	\;\le\;2\|h\|_\infty\|g\|_\infty\bigl(c_{r,J}+\beta_J(b_r)\bigr)
	\;\longrightarrow\;0,
	\]
	where $\beta_J(\delta):=\pi_J\bigl(\{\mu:\mu(\Gamma_\delta)\ge1\}\bigr)$
	does not depend on $r$. The sets $\Gamma_\delta$ are closed and decrease,
	as $\delta\downarrow0$, to $\bigcup_{A\in\mathcal A}\partial A\subset
	\bigl(\bigcup_l\partial\mathcal{D}_l\cup\bigcup_l\partial\mathcal{C}_l\bigr)\times[0,T]
	\,\cup\,Z^c_J\times\{0,s,t,s'_1,\dots,s'_q\}$, which is
	$\nu\otimes\Leb$-null. Since the measures are integer-valued,
	$\beta_J(\delta)\downarrow\pi_J(\{\mu:\mu(\bigcup_A\partial A)\ge1\})=0$.
	The same bound holds for the factors separately.

	\emph{Step~3 (Transfer and limit $i\to\infty$).} Let
	$\mathbb T_\ast:=\{\tau\in[0,T]:\PP^1(u^1_\ast(\tau)\neq u^1_\ast(\tau-))=0\}$.
	Since $u^1_\ast$ is c\`adl\`ag in $H$, $[0,T]\setminus\mathbb T_\ast$ is at
	most countable \cite[Lemma~3.7.7]{Ethier1986}. In this step assume, in
	addition, $s_1,\dots,s_p\in\mathbb T_\ast\cup\{0,T\}$. By equality of laws
	(Proposition~\ref{prop:JS-pairs-allJ}\,\ref{eq:JS-laws}), the expectations
	for $(u^1_i,N^{1,J}_i)$ equal those for
	$(u_{n_{r_i}},\Phi^J_{r_i}(u_{n_{r_i}}))$. Combining Steps~1 and 2,
	\[
	\EE^1\bigl[(h\,g)(u^1_i,N^{1,J}_i)\bigr]
	=\EE^1\bigl[h(u^1_i,N^{1,J}_i)\bigr]\,c_g+o(1),\qquad i\to\infty .
	\]
	By Proposition~\ref{prop:JS-pairs-allJ}\,\ref{JS:as}, $\PP^1$-a.s.\
	$u^1_i\to u^1_\ast$ in $\DD([0,T];H)$ and $N^{1,J}_i\to\eta^{1,J}$. The
	evaluation $x\mapsto x(\tau)$ is continuous on $\DD([0,T];H)$ at every
	path continuous at $\tau$, and at every path if $\tau\in\{0,T\}$
	\cite[Theorem~12.5]{Billingsley1999}; hence $u^1_i(s_k)\to u^1_\ast(s_k)$
	a.s.\ for $k\le p$. The counts $\mu\mapsto\mu(A)$, $A\in\mathcal A$, are
	continuous at $\eta^{1,J}$ a.s., since by Lemma~\ref{lem:etaJ-PRM}
	$\eta^{1,J}$ has a.s.\ no atoms on the $\nu\otimes\Leb$-null set
	$\bigcup_A\partial A$. Dominated convergence gives
	\eqref{eq:compat-functional} for the restricted class, with
	$\EE^1[g]=c_g$.

	\emph{Step~4 (Extension).} \emph{(a) All times.} Let $s_k\notin\mathbb
	T_\ast\cup\{0,T\}$ for some $k$. Choose $s^{(\ell)}_k\in\mathbb T_\ast$ with
	$s^{(\ell)}_k\downarrow s_k$ (and $s^{(\ell)}_k:=s_k$ if $s_k\in\mathbb
	T_\ast\cup\{0,T\}$), and put $s^{(\ell)}:=\max\{s,s^{(\ell)}_1,\dots,s^{(\ell)}_p\}$,
	so that $s^{(\ell)}\downarrow s$ and $s^{(\ell)}<t$ for large $\ell$. Step~3
	applies with $s^{(\ell)}_k$ and $(s^{(\ell)},t]$. As $\ell\to\infty$,
	$u^1_\ast(s^{(\ell)}_k)\to u^1_\ast(s_k)$ by right-continuity, and
	$\eta^{1,J}(\mathcal C_l\times(s^{(\ell)},t])\to\eta^{1,J}(\mathcal C_l\times(s,t])$,
	the counts being integer-valued and eventually constant. Dominated
	convergence gives \eqref{eq:compat-functional} for all times.
	\emph{(b) All functions and sets.} Bounded continuous $h$, bounded $g$
	and the $\pi$-system $\mathscr Z_0$ determine the joint law of the
	finitely many observables. A $\pi$--$\lambda$ argument extends
	\eqref{eq:compat-functional} to all bounded Borel $h,g$ and, letting
	$J\to\infty$ as explained at the beginning of the proof, to $\eta^1$ and
	all sets with finite $\nu$-mass. By the monotone class theorem, it then
	extends to all bounded
	$\widetilde{\mathscr F}^1_s$-measurable first factors $Y$, where
	$\widetilde{\mathscr F}^1_s$ is the unaugmented $\sigma$-field generated
	by the pair up to time $s$. \emph{(c) Augmentation.} Let $Y$ be bounded
	and $\mathscr F^1_s$-measurable. For every $s'\in(s,t)$, $Y$ coincides
	a.s.\ with a bounded $\widetilde{\mathscr F}^1_{s'}$-measurable random
	variable, so by (b) with $s'$ in place of $s$,
	$\EE^1[Y\,g(\eta^1(\mathcal C_l\times(s',t]))_l]=\EE^1[Y]\,
	\EE^1[g(\eta^1(\mathcal C_l\times(s',t]))_l]$. Letting $s'\downarrow s$,
	the counts converge a.s.\ (they are integer-valued and eventually
	constant), and dominated convergence gives the identity for $s$.}
\end{proof}

Finally, let us define  the predictable $\sigma$-algebra $\mathscr{P}^{1}$  on $\Omega^1 \times [0,T]$ by
\DEQSZ\label{def:predfil}
\mathscr{P}^1_{\ast}:=\sigma\left(\left\{ F^1_{s}\times (s,t]:F^1_{s}\in\mathscr{F}^{1}_{s},\ 0\le s< t\le T\, \rk\}\cup \lk\{ F_0\times \{0\}: F_0\in\mathscr{F}_0^1 \right\}\right) .
\EEQSZ

							\newcommand{\Nbar}{\bar{N}}
							\newcommand{\nueps}{\nu_{\varepsilon}}
							\newcommand{\aslim}{\xrightarrow{\;\mathrm{a.s.}\;}}
							\newcommand{\dto}{\xrightarrow{\;d\;}}
							\newcommand{\wto}{\xrightarrow{\;w\;}}
							\newcommand{\lto}{\longrightarrow}

							\item

\textbf{{The limit noise is a time-homogeneous Poisson random
		measure whose intensity measure is the given L\'evy measure $\nu$.}}

The next proposition summarises what we have proven in Steps~(I)--(XII), where we have
 shown that the limit noise is a Poisson random measure
	with L\'evy measure $\nu$.

\begin{proposition}\label{prop:XIII-levy}
Let
$$\bigl(u^1_\ast,(\eta^{1,J})_{J\in\NN}\bigr)
$$
 be the limit
		given in  Proposition~\ref{prop:JS-pairs-allJ}, and let $\eta^1$ be
		the patched measure of Proposition~\ref{prop:full-levy-patch}. Then, there exists a probability space $\mathfrak A^\ast$ such that
		$\eta^1$ is a time-homogeneous Poisson random measure on
		$(Z,\mathscr Z)$ over $\mathfrak A^\ast$ with intensity measure $\nu$ and
		compensator
		$$
\gamma:\mathscr{Z}\times\mathscr{B}([0,T])\ni(B,I) \mapsto \nu(B)\,\Leb(I).
$$
 Moreover, $\eta^1$ is compatible to the filtration $(	\mathscr F^1_t)_{t\in[0,T]}$  In particular the
		compensated Poisson random measure $\tilde\eta^1$ defines a L\'evy process with L\'evy
		measure $\nu$.
\end{proposition}

\begin{proof}
	The proof consists of four assertions.

		\emph{(XIII.1) Existence and a.s.~convergence:} The existence of the limit pair and the
		$\PP^1$-a.s.\ convergences $u^1_i\to u^1_\ast$ and
		$N^{1,J}_i\to\eta^{1,J}$ (for every $J$) are part of the
		\emph{conclusion} of Proposition~\ref{prop:JS-pairs-allJ}.

		\emph{(XIII.2) Identification of the law at fixed $J$.}
		By Lemma~\ref{lem:etaJ-PRM} (end of Step~(X)), for every $J$ the limit
		$\eta^{1,J}$ is a Poisson random measure with intensity
		$\nu\rest_{Z^c_J}\otimes\Leb$, i.e.\ $\Law(\eta^{1,J})=\pi_J$.

		\emph{(XIII.3) Consistency in $J$ and passage to the full L\'evy measure.} The restrictions are a.s.\ consistent for $J'\ge J$
		(Lemma~\ref{lem:J-consistency}), and the monotone patching of
		Proposition~\ref{prop:full-levy-patch} produces $\eta^1$ with intensity
		$\nu$ on $Z\setminus\{0^\ast\}$; the $\varepsilon\to0$ diagonal and Linde's
		weak-convergence criterion for infinitely divisible laws are not used.

		\emph{(XIII.4) Compatibility with the filtration.}
		Adaptedness of both coordinates holds by construction of the pair
		filtration (Definition~\ref{def:pair-filtration}). The independence of
		the future increments $\eta^1(\mathcal{C}\times(s,t])$, $s<t$, from
		$\mathscr F^1_s$ is the functional compatibility identity. It holds on
		$\mathfrak A$ for $(u_n,\eta)$, carries over to
		$(u_{n_r},\Phi^J_r(u_{n_r}))$ up to $o(1)$ by
		Lemma~\ref{lem:X-good-event}, is transferred to $(u^1_i,N^{1,J}_i)$ by
		equality of laws, and passes to the a.s.\ limit
		(Lemma~\ref{lem:compat-transfer}). This is
		Definition~\ref{def-Prm} (iii).

		Combining (XIII.1)--(XIII.4) proves the proposition.
\end{proof}

\begin{corollary}\label{cor:usual-conditions}
		Set $\mathfrak A^\ast:=(\Omega^1,\mathscr F^1,(\mathscr F^1_t)_{t\in[0,T]},
		\PP^1)$, $u_\ast:=u^1_\ast$ and $\eta_\ast:=\eta^1$. Then
		$(\mathscr F^1_t)$ satisfies the usual conditions,
		$u_\ast\in\DD([0,T];H)\cap L^\alpha(0,T;V)$ $\PP^1$-a.s.\ with
		$\Law(u_\ast)=\rho$, and $u_\ast$ is
		$(\mathscr F^1_t)$-progressively measurable and by
		Proposition~\ref{prop:XIII-levy}, $\eta_\ast$ is a time-homogeneous
		Poisson random measure on $(Z,\mathscr Z)$ with intensity measure
		$\nu$ and compensator $\nu(dz)\,dt$, compatible with
		$(\mathscr F^1_t)$. This establishes conclusions (i), (ii) and (iii) of
		Theorem~\ref{Theoremmain}; conclusion (iv), the equation,
		remains and is the object of the final step.
\end{corollary}

						\item

\noindent\textbf{The limit process $u^1_\ast$ is a
							solution to the SPDE.}

We show that $u^1_\ast$ solves the equation of Theorem~\ref{Theoremmain} on $\mathfrak A^\ast$ (Corollary~\ref{cor:usual-conditions}). The proof is divided into six steps. In Steps~1 and~2 we transfer what we need from $\mathfrak A$ to $\mathfrak A^1$, and in Steps~3--5 we pass to the limit.

In Step~1 we replace for a fixed $r\in\NN$ the good event $\mathcal G_{r,J}$ (page~\pageref{godevent}, Lemma~\ref{lem:X-good-event}), on which the registration algorithm recovers $\eta\rest_{\ZZJ}$ up to registration errors, by an event that depends on the path alone: up to a null set, $\mathcal G_{r,J}$ is contained in $\{u_{n_r}\in\Gamma_{r,J}\}$, where $\Gamma_{r,J}$ is a Borel set of paths in $\mathcal X$. We write the error for a fixed $r\in\NN$  as a residual $\rho^J_r$, which is also a Borel functional of the path. Hence, for any fixed $r\in\NN$, the good sets can be transferred from $\mathfrak A$ to $\mathfrak A^1$ by equality of laws, and we estimate the residual on the transferred good sets in expectation (Lemma~\ref{lem:XV-transfer}). 
In Step~2 we transfer the a~priori bounds, which we have on $\MA$ to $\mathfrak A^\ast$ and derive the corresponding bounds for $u^1_\ast$.
In Step~3 we fix $J\in\NN$ and treat the martingale part: as $i\to\infty$, the stochastic integral of $\langle g(\bar u^1_i(s),z),\phi\rangle$ with respect to the compensated measure $\widetilde N^{1,J}_i$ converges in probability to the stochastic integral of $\langle g(u^1_\ast(s-),z),\phi\rangle$ with respect to $\widetilde\eta^{1,J}$ (Lemma~\ref{lem:XV-mart-fixedJ}). 
Step~4 prepares the limit $J\to\infty$: the stochastic integrals over $Z^c_J$ with respect to $\widetilde\eta^{1,J}$ converge in $\mathcal M^2(0,T;H)$ to the integral over $Z\setminus\{0\}$ with respect to $\widetilde\eta^1$ (Lemma~\ref{lem:XV-smalljumps}). 
In Step~5\,(a) we show that, for every $\phi\in V$,
\[
\int_0^t\langle A_{n_{r_i}}(s,u^1_i(s)),\phi\rangle\,ds\longrightarrow\int_0^t\langle A(s,u^1_\ast(s)),\phi\rangle\,ds\qquad\text{as }i\to\infty,
\]
uniformly in $t\in[0,T]$, in $\PP^1$-probability. In Step~5\,(b) we combine Steps~1--3 and~5\,(a), let first $i\to\infty$ and then, by Step~4, $J\to\infty$, and obtain the limit equation \eqref{eq:XV-chi-eq}. Step~6 completes the proof of conclusion~(iv) of Theorem~\ref{Theoremmain}. Now, we are ready to start with the proof.
\par

						\emph{Step 1: Transfer of the  identity at a fixed $r\in\NN$.}

On $\mathfrak A$ the reconstructed measure $\Phi^J_r(u_{n_r})$ coincides with the true jumps $\eta\rest_{\ZZJ}$ only on the good event $\mathcal G_{r,J}$ of Lemma~\ref{lem:X-good-event}, and there only up to the registration errors.
 The good event $\mathcal G_{r,J}$ is an event on $\Omega$, see \pageref{godevent}. Up to a null set, it is contained in the event $\{u_{n_r}\in\Gamma_{r,J}\}$, where $\Gamma_{r,J}$ is a Borel set of paths in $\mathcal X=L^\alpha(0,T;V)\cap\DD([0,T];H)$, see {\ref{XV:D-good}} below. Since $\Gamma_{r,J}$ depends only on the path, the event $\{u_{n_r}\in\Gamma_{r,J}\}$ can be transferred from $\MA$ to $\MA^\ast $ by equality of laws. So we write the error as a residual, which is a Borel functional of the path, and we control the residual on the transferred good set.
\smallskip
\noindent  Fix $J\in\NN$ and $r\ge r_0(J)$.
 We define the following \emph{path functionals}:

\begin{enumerate}[label=\upshape(D\arabic*)]
\item\label{XV:D-proj}  the \emph{dyadic coefficients} $\xi^\mr _k:L^\alpha(0,T;V)\ni x\mapsto \xi_k(x)\in
 V$, $k=0,\dots,2^\mr -1$, where  $\xi_k(x)$ is given by
\[
\xi^\mr _0(x):=0,\qquad \xi^\mr _k(x):=\frac1{\tau_\mr }\int_{I^\mr _{k-1}}x(s)\,ds,\quad k\ge1,
\]
and the \emph{shifted dyadic projection} \footnote{See Appendix~\ref{haarproj} with $E=V$ and initial value $u^m_0$.}
$$\bar x_r:=\mathrm{Proj}_\mr  x:=\sum_{k=0}^{2^\mr -1}\xi^\mr_k(x)\,\mathbf 1_{I^\mr _k}.
$$
For $k\ge1$ the map $\xi^\mr_k$ is linear, and H\"older's inequality gives $|\xi_k(x)|_V\le\tau_\mr ^{-1/\alpha}\|x\|_{L^\alpha(0,T;V)}$. Hence the maps $\xi^\mr_k$ and $x\mapsto\bar x_r\in L^\infty(0,T;V)$ are Lipschitz continuous on $L^\alpha(0,T;V)$, and in particular continuous on $\mathcal X$. Moreover, $\xi^\mr_k(u_{n_r})=\xi^{(\mr )}_{k,n_r}$ and $\mathrm{Proj}_\mr  u_{n_r}=\myp_{n_r,m_r}$;
\item the \emph{residual}  $\rho^J_r:\mathcal X\times [0,T]\ni (x,t)\mapsto \rho^J_r(x)(t) \in V^\ast$, given by
\begin{align}\label{eq:XV-rho}
\langle\rho^J_r(x)(t),\phi\rangle&:=\langle x(t),\phi\rangle
 -\int_0^t\langle A_{n_r}(s,x(s)),\phi\rangle\,ds\notag\\
&\quad -\int_0^t\!\!\int_{Z^c_J}\langle g(\bar x_r(s),z),\phi\rangle\,\bigl(\Phi^J_r(x)-\nu\rest_{Z^c_J}\otimes\Leb\bigr)(dz,ds),\quad \forall \phi\in V,
\end{align}
{ where $\bar x_r$ is the projection of \ref{XV:D-proj}, $\Phi^J_r(x)$ is the output of the registration algorithm at fixed  $r$ (Remark~\ref{rem:dictionary}), with one atom $(\zhat_k,\Upsilon_{k,\mr,\jt })$ for every registered cell $I^\mr_k$, $\Upsilon_{k,\mr,jt }$ being the midpoint of $I^\mr _k$, and $\nu\rest_{Z^c_J}$ is the restriction of $\nu$ to $Z^c_J$. Here and below, $A_n(s,v):=0$ for $v\in V\setminus V_n$. This convention is never active along $u_n$, which is $V_n$-valued, and since $\{x:\ x(s)\in V_n\text{ for a.e.\ }s\}$ is closed in $L^\alpha(0,T;V)$, equality of laws shows that it is $\PP^1$-a.s.\ never active along $u^1_i$ either;}
\item the \emph{small-jump part} $M^\Jr _\natural:\mathcal X\times [0,T]\ni (x,t)\mapsto  M^\Jr _\natural(x)(t) \in V^\ast$, given by
\begin{align}\label{eq:XV-Mnat}
\langle M^\Jr _\natural(x)(t),\phi\rangle:=\langle x(t),\phi\rangle
 -\int_0^t\langle A_{n_r}(s,x(s)),\phi\rangle\,ds
 -\langle L^J(x)(t),\phi\rangle, \quad \forall \phi\in V,
\end{align}
where $\mathcal J^J(x)$ is the path jump measure of Step~(X), and $L^J :\mathcal X\times [0,T]\ni (x,t)\mapsto L^J (x)(t) \in V^\ast$ is defined by
\[
\langle L^J (x)(t),\phi\rangle:=\int_0^t\!\!\int_{Z^c_J }\langle g(x(s-),z),\phi\rangle\,\mathcal J^J (x)(dz,ds)
-\int_0^t\!\!\int_{Z^c_J }\langle g(x(s),z),\phi\rangle\,\nu(dz)\,ds
\]
and is the large-jump part of the path;
\item the \emph{error functional}  $E_{r,\Jr }:\mathcal X\ni x\mapsto E_{r,\Jr }(x)\in[0,\infty]$, given by
\begin{align}\label{eq:XV-E}
E_{r,\Jr }(x)&:=\sum_{(s,z)\in\operatorname{supp}\mathcal J^J (x)}\bigl|g(x(s-),z)-g(\bar x_r(s),z)\bigr|_H\notag\\
&\quad+\bigl(\delta_r+2^{-m_r}\bigr)\,\mathcal J ^J (x)\bigl(\ZZJ\bigr)
+\bigl(\nu_J L_g\bigr)^{1/2}\int_0^T|x(s)-\bar x_r(s)|_V\,ds,
\end{align}
where $\nu_J =\nu(Z^c_J )$, as in Step~(VIII);
\item\label{XV:D-good} the \emph{good set on path space} $\Gamma_{r,J}$, the set of all $x\in\mathcal X$ such that\footnote{Compare with the good events on \pageref{godevent}.}
\begin{enumerate}[label=\upshape($\Gamma$\arabic*)]
\item $\mathcal J^J(x)$ has no atom in $I^\mr _0$ and at most one atom in every cell $I^\mr _k$;
\item a cell $I^\mr _k$ carries an atom of $\Phi^J_r(x)$ if and only if it contains an atom of $\mathcal J^J(x)$;
\item if $(s,z)$ is the atom of $\mathcal J^J(x)$ in $I^\mr _k$ and $(\zhat_k,\Upsilon_{k,\mr,jt })$ is the atom of $\Phi^J_r(x)$ in $I^\mr _k$, then 
$$\bigl|g(\xi^\mr_k(x),\zhat_k)-g(\xi^\mr _k(x),z)\bigr|_H\le\delta_r+2^{-m_r};
$$
\end{enumerate}
\item finally, the sets
$O^J_r(t):=\{x\in\mathcal X:\ \mathcal J^J(x)(Z^c_J\times I^\mr _{k})\ge1\}$ for $t\in I^\mr _k$, and $O^J_r(0):=\emptyset$.
\end{enumerate}

The following considerations show that the maps defined above are defined on the path space alone. In particular, they do not depend on the probability space.
At an atom $(s,z)$ of $\mathcal J^J(x)$ one has $x(s-)\in V$ (by the convention $\Delta H_\jt(v)=\emptyset$ for $v\in H\setminus V$ of Step~(X)) and $\Delta x(s)=g(x(s-),z)$. So $L^J(x)$ is well defined for every $x$. The compensator 
integrals are finite, since $\int_{Z^c_J}|g(v,z)|_H\,\nu(dz)\le(\nu_JC_g)^{1/2}(1+|v|_V^2)^{1/2}$ by \eqref{eq:g-second-moment1}, $x\in L^2(0,T;V)$ and $\bar x_r\in L^\infty(0,T;V)$. For fixed $(t,\phi)$ all these functionals are Borel on $\mathcal X$, and $\Gamma_{r,J}$ and $O^J_r(t)$ are Borel sets, for the following reasons. 
The map $x\mapsto\bar x_r$ is continuous by \ref{XV:D-proj}, and the maps $x\mapsto x(t)$ and $x\mapsto\int_0^t\langle A_{n_r}(s,x(s)),\phi\rangle\,ds$ are Borel. 
 For the last map, and for the same map with $A$ in place of $A_{n_r}$ (used in Step~5(a)), note the following. The map $(s,x)\mapsto x(s)\in H$ is $\mathscr B([0,T])\otimes\mathscr B(\mathcal X)$-measurable, since it is right-continuous in $s$ and Borel in $x$. The subspaces $V_n$ and $V$ are Borel subsets of $H$, and their Borel $\sigma$-fields are the traces of $\mathscr B(H)$, by the Lusin--Suslin theorem \cite[Thm.~15.1]{Kechris1995}. Finally, $A_n$ and $A$ are jointly measurable by Hypotheses~\ref{hyp:An} and~\ref{hyp}; we put $A(s,v):=0$ for $v\in H\setminus V$. Hence the integrands are jointly measurable, and the integrals are Borel in $x$ by Fubini's theorem, with the value $0$ where they do not converge absolutely. The maps $\Phi^J_r$ and $\mathcal J^J$ are Borel (Propositions~\ref{prop:Nhat-meas} and~\ref{prop:measurable}, Step~(X)).
For $\Phi^J_r$ this holds on all of $\mathcal X$: by Lemma~\ref{lem:drift-meas} the increment functionals $\ell^{m_r}_{k,n_r}$ are Borel from $\mathcal X$ into $H$, the maps $x\mapsto\xi^{m_r}_k(x)$ are continuous by \ref{XV:D-proj}, and the distances $d_\jt$, the selector $\hat\jmath$ and the tolerance-adapted inverse of (A7) are Borel on $V\times H$ (Lemmas~\ref{lem:DH-meas} and~\ref{lem:selector-meas}). The atoms of a finite point measure can be enumerated measurably (Lemma~\ref{lem:jump-measurable}), and an integral against a finite point measure is a finite sum. The compensator integrals are Borel by Fubini's theorem.

 Hence the good sets, the residual and all other objects above are defined on the path space and do not depend on the underlying probability space.
As a consequence, the following lemma transfers properties of these objects from $\MA$ to $\MA^\ast$.
 In the lemma, $|\cdot|_H$ is also applied to elements of $V^\ast$. To this end, we fix a countable set $D\subset V$ which is dense in $V$, and put
\[
|y|_H:=\sup\bigl\{\langle y,\phi\rangle:\ \phi\in D,\ |\phi|_H\le1\bigr\}\in[0,\infty],\qquad y\in V^\ast;
\]
by the density of $D$ in $V$ and of $V$ in $H$, this is the norm of $H$ if $y\in H$, and $|y|_H=\infty$ if $y\in V^\ast\setminus H$. Since $D$ is countable, $x\mapsto|\rho^J_r(x)(t)|_H$ and $x\mapsto|M^J_\natural(x)(t)|_H$ are Borel on $\mathcal X$ for every $t$.

\begin{lemma}\label{lem:XV-transfer}
Let $J\in\NN$ and $i\in\NN$ with $r_i\ge r_0(J)$. Put $\bar u^1_i:=\mathrm{Proj}_{m_{r_i}}u^1_i$ and $\widetilde N^{1,J}_i:=N^{1,J}_i-\nu\rest_{Z^c_J}\otimes\Leb$.
\begin{enumerate}[label=\upshape(\alph*)]
\item\label{XV-tr:a}  $N^{1,J}_i=\Phi^J_{r_i}(u^1_i)$, $\PP^1$-a.s.
\item\label{XV-tr:b}  we have $\PP^1$-a.s., for all $t\in[0,T]$ and $\phi\in V$,
\begin{align}\label{eq:XV-approx}
\langle u^1_i(t),\phi\rangle&=
\int_0^t\langle A_{n_{r_i}}(s,u^1_i(s)),\phi\rangle\,ds\notag\\
&\quad+\int_0^t\!\!\int_{Z^c_J}\langle g(\bar u^1_i(s),z),\phi\rangle\,\widetilde N^{1,J}_i(dz,ds)+\langle\rho^J_{r_i}(u^1_i)(t),\phi\rangle .
\end{align}
\item\label{XV-tr:c} \emph{The transferred good sets} (compare the good event on page~\pageref{godevent}). For $t\in[0,T]$ let $G^1_{i,J}:=\{u^1_i\in\Gamma_{r_i,J}\}$ and
\[
C_i=C_i(t):=G^1_{i,J}\setminus\{u^1_i\in O^J_{r_i}(t)\}.
\]
In other words,  $C_i$  is  the good set transferred to $\MA^\ast$ without the event that $t\in[0,T]$ lies in a cell of order $m_{r_i}$ containing an atom of $\mathcal J^J(u^1_i)$.\footnote{Remark~\ref{rem:XV-excluded-cell} explains why this event is excluded.} Then
\begin{align}
\PP^1(C_i^c)&\le c_{r_i,J}+\nu_J\tau_{m_{r_i}},\label{eq:XV-Ci}\\
|\rho^J_{r_i}(u^1_i)(t)|_H\,\mathds 1_{C_i}&\le|M^J_\natural(u^1_i)(t)|_H+E_{r_i,J}(u^1_i),\label{eq:XV-rho-bound}
\end{align}
where $c_{r,J}$ is the bound \eqref{eq:X-good-prob}, and
\begin{align}
\EE^1\bigl[|M^\Jr_\natural(u^1_i)(t)|_H^2\bigr]&\le\kappa_2(\epsilon_{J+1})\,C_2,\label{eq:XV-Mnat-bound}\\
\EE^1\bigl[E_{r_i,J}(u^1_i)\bigr]&\le\bar c_{r_i,J}:=\bigl(C_\ast+2T\nu_J+1\bigr)\,\delta_{r_i},\label{eq:XV-E-bound}
\end{align}
where  $C_2$ comes from Step~(VIII)\,(c), $\kappa_2$ comes from \eqref{eq:XV-tail}, and $C_\ast$ comes from \eqref{eq:Ee-bound}. In particular, for every fixed $J$, $\PP^1(C_i^c)\to0$ and $\EE^1[E_{r_i,J}(u^1_i)]\to0$ as $i\to\infty$.
\end{enumerate}
\end{lemma}

\begin{proof}
\ref{XV-tr:a} Since $\MMN(\ZZJ)$ is separable and metrisable, the graph $\mathrm{Gr}:=\{(x,\mu):\ \mu=\Phi^J_{r_i}(x)\}$ of the Borel map $\Phi^J_{r_i}:\mathcal X\to\MMN(\ZZJ)$ is a Borel subset of $\mathcal X\times\MMN(\ZZJ)$. As $(u_{n_{r_i}},\Phi^J_{r_i}(u_{n_{r_i}}))\in\mathrm{Gr}$ everywhere on $\Omega$, Proposition~\ref{prop:JS-pairs-allJ}\,\ref{eq:JS-laws} gives $\PP^1\bigl((u^1_i,N^{1,J}_i)\in\mathrm{Gr}\bigr)=1$.

\ref{XV-tr:b} By the definition \eqref{eq:XV-rho} of the residual and by \ref{XV-tr:a}, \eqref{eq:XV-approx} is an identity, simultaneously for all $(t,\phi)$.
\ref{XV-tr:c} \emph{Part 1: a pathwise estimate.} Let $x\in\Gamma_{r,J}$ and $t\in[0,T]$ with $x\notin O^J_r(t)$. We show
\begin{equation}\label{eq:XV-pathwise}
\bigl|\rho^J_r(x)(t)-M^J_\natural(x)(t)\bigr|_H\le E_{r,J}(x).
\end{equation}
By \eqref{eq:XV-rho} and \eqref{eq:XV-Mnat}, for every $\phi\in V$,
\begin{align*}
\langle\rho^J_r(x)(t)-M^J_\natural(x)(t),\phi\rangle
&=\int_0^t\!\!\int_{Z^c_J}\langle g(x(s-),z),\phi\rangle\,\mathcal J^J(x)(dz,ds)
-\int_0^t\!\!\int_{Z^c_J}\langle g(\bar x_r(s),z),\phi\rangle\,\Phi^J_r(x)(dz,ds)\\
&\quad-\int_0^t\!\!\int_{Z^c_J}\langle g(x(s),z)-g(\bar x_r(s),z),\phi\rangle\,\nu(dz)\,ds ,
\end{align*}
so $\rho^J_r(x)(t)-M^J_\natural(x)(t)\in H$. Let $(s,z)$ be an atom of $\mathcal J^J(x)$, and let $I^m_k$ be the cell containing $s$. By ($\Gamma$1) and ($\Gamma$2), $I^m_k$ carries exactly one atom $(\zhat_k,\Upsilon_{k,m,jt})$ of $\Phi^J_r(x)$, and every atom of $\Phi^J_r(x)$ arises in this way. Since $x\notin O^J_r(t)$, $t\notin I^m_k$; hence $I^m_k\subset[0,t]$ or $I^m_k\subset(t,T]$, and $s\le t$ if and only if $\Upsilon_{k,m,jt}\le t$. Thus the two point-measure integrals run over the same cells. Since $\bar x_r=\xi^\mr_k(x)$ on $I^\mr_k$, ($\Gamma$3) gives
\begin{align*}
&\Bigl|\int_0^t\!\!\int_{Z^c_J} g(x(s-),z)\,\mathcal J^J(x)(dz,ds)
-\int_0^t\!\!\int_{Z^c_J} g(\bar x_r(s),z)\,\Phi^J_r(x)(dz,ds)\Bigr|_H\\
&\le\sum_{(s,z)\in\operatorname{supp}\mathcal J^J(x)}\Bigl(\bigl|g(x(s-),z)-g(\xi^\mr_k(x),z)\bigr|_H+\bigl|g(\xi^\mr_k(x),z)-g(\xi^\mr_k(x),\zhat_k)\bigr|_H\Bigr)\\
&\le\sum_{(s,z)\in\operatorname{supp}\mathcal J^J(x)}\bigl|g(x(s-),z)-g(\bar x_r(s),z)\bigr|_H+\bigl(\delta_r+2^{-m_r}\bigr)\,\mathcal J^J(x)\bigl(\ZZJ\bigr),
\end{align*}
where in the sums $k=k(s)$ is the index of the cell containing $s$. Finally, by \eqref{eq:g-L1nu} with $\epsilon=\epsilon_{J+1}$,
\[
\Bigl|\int_0^t\!\!\int_{Z^c_J}\bigl[g(x(s),z)-g(\bar x_r(s),z)\bigr]\,\nu(dz)\,ds\Bigr|_H\le\bigl(\nu_JL_g\bigr)^{1/2}\int_0^T|x(s)-\bar x_r(s)|_V\,ds .
\]
Adding the two bounds gives \eqref{eq:XV-pathwise}.

\emph{Part 2: the estimates on $\mathfrak A$.} Let $r\ge r_0(J)$. By Lemma~\ref{lem:jump-identity}, $\mathcal J^J(u_\nr)=\eta\rest_{\ZZJ}$, $\PP$-a.s.
\begin{enumerate}[label=\upshape(2.\roman*)]
\item\label{XV-tr:2i} $\PP(u_\nr\notin\Gamma_{r,J})\le c_{r,J}$. Indeed, up to a $\PP$-null set the following holds on $\mathcal G_{r,J}$. (G1) gives ($\Gamma$1), since $Z^c_J\subset Z^c_{{\varepsilon_\flat,r}}$. Lemma~\ref{lem:X-good-event}\,\ref{X:good-a} and~\ref{X:good-b} give ($\Gamma$2): a registered cell contains exactly one atom with mark in $Z^c_J$, every such atom lies in a registered cell, and a registered cell carries exactly one atom of $\Phi^J_r$. Lemma~\ref{lem:X-good-event}\,\ref{X:good-b} gives ($\Gamma$3), since $\xi^\mr_k(u_\nr)=\xi^{(\mr)}_{k,\nr}$. Now use Lemma \ref{lem:X-good-event},  and \eqref{eq:X-good-prob}.
\item\label{XV-tr:2ii} For every $t$, $\PP$-a.s.,
\[
M^\Jr_\natural(u_\nr)(t)=\int_0^t\!\!\int_{\{0<|z|_Z\le\epsilon_{J+1}\}}g(u_\nr(s),z)\,\widetilde\eta(dz,ds),
\
L^J(u_\nr)(t)=\int_0^t\!\!\int_{Z^c_J}g(u_\nr(s),z)\,\widetilde\eta(dz,ds).
\]
Indeed, by \eqref{originalun} the first two terms of \eqref{eq:XV-Mnat} together equal to
$$
\int_0^t\int_Z g(u_\nr(s),z)\,\widetilde\eta(dz,ds).
$$
 On $Z^c_J$ the intensity is finite, so there the compensated integral is the pathwise $\eta$-integral of the predictable representative $g(u_\nr(s-),z)$ minus the compensator (Lemma~\ref{lem:pathwise-atoms}\,\ref{PA:pathwise} and Remark~\ref{rem:left-limits}), and $\eta\rest_{\ZZJ}=\mathcal J^J(u_\nr)$. By the It\^o isometry and Remark~\ref{rem:GF-consequences}\,(iii), we get
\[
\EE\bigl[|M^\Jr_\natural(u_\nr)(t)|_H^2\bigr]=\EE\int_0^t\!\!\int_{\{0<|z|_Z\le\epsilon_{J+1}\}}|g(u_\nr(s),z)|_H^2\,\nu(dz)\,ds\le\kappa_2(\epsilon_{J+1})\,C_2 ,
\]
and, by \eqref{eq:g-second-moment1}, $\EE\bigl[|\langle L^J(u_\nr)(t),\phi\rangle|^2\bigr]\le|\phi|_H^2\,C_g\,C_2$ for $\phi\in V$.
\item\label{XV-tr:2iii} $\EE[E_{r,J}(u_\nr)]\le\bar c_{r,J}$. $\PP$-a.s., the first term of \eqref{eq:XV-E} equals the sum of $|F(s,z)|_H$ over the atoms of $\eta$ with $|z|_Z>\epsilon_{J+1}$, where $F$ is the averaging error \eqref{eq:errors_freezing}. Since $\varepsilon_{\flat,r}\le\epsilon_{J+1}$, this sum is bounded by $e_{r,J}$ (see \eqref{eq:e-total}), so by \eqref{eq:Ee-bound} its expectation is at most $C_\ast\seque^{1/2}\delta_r\le C_\ast\delta_r$. The second term has expectation $(\delta_r+2^{-m_r})\,T\nu_J\le2T\nu_J\delta_r$ by \ref{DS:i}. For the third term, $\nu_J\le\nu_{\flat,r}$ and \ref{DS:ii} give $(\nu_JL_g)^{1/2}\,\Theta_{n_r}(m_r)\le\seque\delta_r/2\le\delta_r$.
\item\label{XV-tr:2iv} For $t\in I^\mr _k$, $\PP(u_\nr\in O^J_r(t))=\PP\bigl(\eta(Z^c_J\times I^\mr _k)\ge1\bigr)=1-e^{-\nu_J\tau_\mr }\le\nu_J\tau_\mr $.
\end{enumerate}

\emph{Part 3: transfer.} All sets and functionals in Parts~1 and Part~2 are Borel on $\mathcal X$, and we know  by Proposition~\ref{prop:JS-pairs-allJ}\,\ref{eq:JS-laws} that
$$
\Law_{\PP^1}(u^1_i)=\Law_{\PP}(u_{n_{r_i}}).
$$
 Hence \ref{XV-tr:2i} and \ref{XV-tr:2iv} give \eqref{eq:XV-Ci}, and \ref{XV-tr:2ii} and \ref{XV-tr:2iii} give \eqref{eq:XV-Mnat-bound} and \eqref{eq:XV-E-bound}. On $C_i$ the path $u^1_i(\omega)$ satisfies the assumptions of Part~1, so \eqref{eq:XV-rho-bound} follows from \eqref{eq:XV-pathwise} and the triangle inequality. The limits follow from $c_{r,J}\to0$ (Lemma~\ref{lem:X-good-event}), $\tau_{m_r}\to0$ and $\delta_r\to0$.
\end{proof}

\begin{remark}\label{rem:XV-excluded-cell}
Let us clarify why the cell of $t$ is excluded.
 The registration algorithm sees only the increments of the path over the dyadic cells, not the time of a jump inside a cell. Therefore $\Phi^J_r(x)$ places the recovered atom at the midpoint $\Upsilon_{k,\mr,jt }$ of its cell $I^\mr _k$, while the corresponding atom $(s,z)$ of $\mathcal J^J(x)$ sits at the true time $s\in I^\mr _k$. Since $\bar x_r=\xi_k(x)$ on $I^\mr _k$, both atoms are evaluated at the same frozen state; only their times differ. If $t\notin I^\mr _k$, then $s\le t$ if and only if $\Upsilon_{k,\mr,jt }\le t$, so in Part~1 of the proof of Lemma~\ref{lem:XV-transfer} the pair contributes only the small errors collected in $E_{r,J}(x)$. If $t\in I^\mr _k$, the two times can lie on different sides of $t$. For $t$ in the shaded interval $[s,\Upsilon_{k,\mr,jt })$, the true jump is already counted, but its reconstruction is not:
\begin{center}
\begin{tikzpicture}[x=1.2cm,y=1cm,>=stealth,
    lab/.style={below,text height=1.6ex,text depth=0.5ex}]
  \fill[gray!25] (2.2,-0.11) rectangle (5,0.11);
  \draw[gray] (-0.8,0) -- (0,0);
  \draw[gray] (10,0) -- (10.8,0);
  \draw[thick] (0,0) -- (10,0);
  \draw[thick] (0,-0.15) -- (0,0.15);
  \draw[thick] (10,-0.15) -- (10,0.15);
  \draw[dashed] (3.6,-0.15) -- (3.6,0.6);
  \fill (2.2,0) circle (2.5pt);
  \draw[fill=white] (5,0) circle (2.5pt);
  \node[lab] at (0,-0.15) {$t^\mr _k$};
  \node[lab] at (2.2,-0.15) {$s$};
  \node[lab] at (3.6,-0.15) {$t$};
  \node[lab] at (5,-0.15) {$\Upsilon_{k,\mr,jt}$};
  \node[lab] at (10,-0.15) {$t^\mr _{k+1}$};
  \node[align=center,font=\small] (true) at (1.0,1.5)
    {true jump $g(x(s-),z)$\\ already counted};
  \draw[->] (true.south) -- (2.12,0.17);
  \node[align=center,font=\small] (rec) at (7.4,1.5)
    {reconstructed jump $g(\xi^\mr _k(x),\zhat_k)$\\ not yet counted};
  \draw[->] (rec.south) -- (5.1,0.17);
  \draw (0,-0.85) -- (0,-0.95) -- (10,-0.95) -- (10,-0.85);
  \node[below,font=\small] at (5,-0.95) {cell $I^\mr_k=(t^\mr_k,t^\mr_{k+1}]$};
\end{tikzpicture}
\end{center}
Then, up to the small errors collected in $E_{r,J}(x)$, the difference $\rho^J_r(x)(t)-M^J_\natural(x)(t)$ equals the whole jump $g(x(s-),z)$. By \eqref{it:nondeg}, $|g(x(s-),z)|_H\ge\rho_g(|x(s-)|_V)\,\epsilon_{J+1}$, which is not small for fixed $J$. The case $\Upsilon_{k,m,jt}\le t<s$ is the same with the roles exchanged. Hence no bound of the form \eqref{eq:XV-pathwise} can hold for $x\in O^J_r(t)$, and this event has to be removed from the good set. This costs little: by Lemma~\ref{lem:jump-identity}, the number of atoms of $\mathcal J^J(u_{n_r})$ in the cell of $t$ is Poisson distributed with mean $\nu_J\tau_{\mr }$, which gives the term $\nu_J\tau_{m_{r_i}}$ in \eqref{eq:XV-Ci}.
\end{remark}

						\emph{Step 2: Uniform bounds on $\MA^\ast$.}
Since $\Law(u^1_i)=\Law(u_{n_{r_i}})$ on $\mathcal X$ (Proposition~\ref{prop:JS-pairs-allJ}\,\ref{eq:JS-laws}) and the mapping
	$$\mathcal X\ni x\mapsto\sup_{t\in[0,T]}|x(t)|_H^2+\int_0^T|x(s)|_V^\alpha\,ds
	$$
 is Borel on $\mathcal X$, the a~priori estimate of Theorem~\ref{Theoremmain} can be transferred: there exists a constant $C$ such that
\begin{equation}\label{eq:XV-apriori}
\sup_{i\in\mathbb N}\;\EE^1\bigg[\sup_{t\in[0,T]}|u^1_{i}(t)|_H^2+\int_0^T|u^1_{i}(s)|_V^\alpha\,ds\bigg]\le C .
\end{equation}
For the projections\footnote{See Appendix \ref{haarproj}.}  $\bar u^1_i=\Pro_{m_{r_i}}u^1_i$, Jensen's inequality on the cells gives
\[
\sup_{0\le t\le T}|\bar u^1_i(t)|_H\le\max\Bigl\{\sup_m|u^m_0|_H,\ \sup_{t\in[0,T]}|u^1_i(t)|_H\Bigr\},
\qquad
\int_{\tau_{m_{r_i}}}^T|\bar u^1_i(s)|_V^\alpha\,ds\le\int_0^T|u^1_i(s)|_V^\alpha\,ds  .
\]
By \ref{prop:JS-pairs-allJ}\,\ref{JS:as}, 
$u^1_i\to u^1_\ast$ in $L^\alpha(0,T;V)$ and in $\DD([0,T];H)$, $\PP^1$-a.s.
	Hence
	$$
	\int_0^T|u^1_i(s)|_V^\alpha\,ds\to\int_0^T|u^1_\ast(s)|_V^\alpha\,ds,\quad \mbox{for}\quad i\to \infty.
	$$
	 Moreover, we have  at every continuity point $t$ of $u^1_\ast$ by Lemma~\ref{lem:J1-pointwise}\,(a)
	 $$
	 |u^1_\ast(t)|_H=\lim_i|u^1_i(t)|_H\le\liminf_i\sup_{s\in[0,T]}|u^1_i(s)|_H
	 $$
These points are dense and $u^1_\ast$ is right-continuous, so
$$
\sup_{t\in[0,T]} |u^1_\ast(t)|_H\le\liminf_i\sup_{t\in[0,T]}|u^1_i(t)|_H.
$$
Fatou's lemma and \eqref{eq:XV-apriori} give
						\begin{equation}\label{eq:XV-limit-apriori}
							\EE^1\bigg[\sup_{t\in[0,T]}|u^1_\ast(t)|_H^2\bigg]\le C\ \ \text{ and }
							\ \ 
							\EE^1\bigg[\int_0^T|u^1_\ast(s)|_V^\alpha\,ds\bigg]\le C .
						\end{equation}

 The process $u^1_\ast$ is $(\mathscr F^1_t)$-progressively measurable by Corollary~\ref{cor:usual-conditions}.
	Hence,
	$$[0,T]\ni s\mapsto A(s,u^1_\ast(s)),
	$$
	is also $(\mathscr F^1_t)$-progressively measurable  by the joint measurability of $A$ (Hypothesis~\ref{hyp}).
	Moreover, $A(s,v):=0$ for $v\in H\setminus V$, see Step~1.
	\par

						\emph{Step 3:  The martingale part at fixed $J$.} For our requirement, we state and prove the following result:

\begin{lemma}\label{lem:XV-mart-fixedJ}
Fix $J\in\NN$, $\phi\in V$ and $t\in[0,T]$. Then, as $i\to\infty$, we have in probability over $\MA^1$
\begin{equation}\label{eq:XV-mart-conv}
\int_0^t\!\!\int_{Z^c_J}\langle g(\bar u^1_i(s),z),\phi\rangle\,\widetilde N^{1,J}_i(dz,ds)
\longrightarrow
\int_0^t\!\!\int_{Z^c_J}\langle g(u^1_\ast(s-),z),\phi\rangle\,\widetilde\eta^{1,J}(dz,ds)
.
\end{equation}
\end{lemma}
\begin{proof}
\newcommand{\QP}{\Pro_{m'}}%
Let $\Lt_J:=\Lt(\epsilon_{J+1},\epsilon_0)$ be the constant of Hypothesis~\ref{ass:GF}\,\ref{it:lip}. Since $Z_J^c=\{\epsilon_{J+1}<|z|_Z\le\epsilon_0\}$,
\begin{equation}\label{eq:XV-lipJ}
\sup_{z\in Z_J^c}|g(v_1,z)-g(v_2,z)|_H\le\Lt_J\,|v_1-v_2|_V,\qquad v_1,v_2\in V .
\end{equation}
As in Remark~\ref{rem:left-limits}, $g(u^1_\ast(s-),z)$ stands for the predictable representative $g(u^1_\ast(s-),z)\mathds 1_{\{u^1_\ast(s-)\in V\}}$.

For $m'\in\NN$ let $\QP$ be the shifted dyadic projection of order $m'$ with initial value $u_0^{m'}$ (Appendix~\ref{haarproj} with $E=V$), and let $\xi_k^{m'}$, $k=0,\dots,2^{m'}-1$, be its coefficients, as in \ref{XV:D-proj} with $m'$ in place of $m_r$. Denote the left-hand side and the right-hand side of \eqref{eq:XV-mart-conv} by $v_i$ and $v$, respectively, and define the same integrals with projected integrands,
\[
\begin{aligned}
	v_i^{m'}&:=\int_0^t\!\!\int_{Z_J^c}\langle g(\QP u_i^1(s),z),\phi\rangle\,\widetilde N_i^{1,J}(dz,ds),
\\
	v^{m'}&:=\int_0^t\!\!\int_{Z_J^c}\langle g(\QP u_\ast^1(s),z),\phi\rangle\,\widetilde\eta^{1,J}(dz,ds).
\end{aligned}
\]
Since $\nu(Z_J^c)<\infty$, all four integrals are pathwise: each is the integral against the point measure minus the integral against the compensator $\nu\rest_{Z_J^c}\otimes\Leb$. For $v_i$ and $v_i^{m'}$ this is the definition of $\widetilde N_i^{1,J}$ (Lemma~\ref{lem:XV-transfer}); note that $N_i^{1,J}=\Phi^J_{r_i}(u^1_i)$ is the reconstructed finite point measure, not a Poisson random measure. For $v$ and $v^{m'}$ it is Lemma~\ref{lem:pathwise-atoms}\,\ref{PA:pathwise}, since $\eta^{1,J}=\eta^1\rest_{\ZZJ}$ has the compensator $\nu\rest_{Z_J^c}\otimes\Leb$ with respect to $(\mathscr F_t^1)$ (Proposition~\ref{prop:full-levy-patch} and Corollary~\ref{cor:usual-conditions}) and the integrands are predictable (see (iii) below).

\emph{Properties of $\QP$.} For $x,y\in L^1(0,T;V)$ and $k\ge1$,
\[
|\xi_k^{m'}(x)-\xi_k^{m'}(y)|_V
\le\tau_{m'}^{-1}\|x-y\|_{L^1(0,T;V)},
\]
while $\xi_0^{m'}=u_0^{m'}$ does not depend on $x$. Lemma~\ref{lem:Pro-props}\,\ref{Pro:jensen} gives
\[
\|\QP x-\QP y\|_{L^1(0,T;V)}\le\|x-y\|_{L^1(0,T;V)},\
\|x-\QP x\|_{L^1(0,T;V)}\le2\|x\|_{L^1(0,T;V)}+\tau_{m'}|u_0^{m'}|_V.
\]
Put $\theta_\ast(m'):=\EE^1\|u_\ast^1-\QP u_\ast^1\|_{L^1(0,T;V)}$. By Lemma~\ref{lem:Pro-L1}, \eqref{init0}, the preceding bound, and \eqref{eq:XV-limit-apriori}, dominated convergence yields $\theta_\ast(m')\to0$ as $m'\to\infty$.

\emph{Strategy.} For every $m'\in\NN$,
\begin{equation}\label{eq:XV-mart-split}
|v_i-v|\le|v_i-v_i^{m'}|+|v_i^{m'}-v^{m'}|+|v^{m'}-v| .
\end{equation}
We fix $m'$ and show: (i) $v_i^{m'}\to v^{m'}$, $\PP^1$-a.s., as $i\to\infty$; (ii) $\limsup_{i\to\infty}\EE^1\bigl[|v_i-v_i^{m'}|\wedge1\bigr]\le C\,\theta_\ast(m')$; (iii) $\EE^1|v^{m'}-v|\le C\,\theta_\ast(m')$, where $C$ depends on $J$ and $\phi$, but not on $i$ and $m'$. Only then we let $m'\to\infty$. Thus the approximation by the simple predictable integrands $g(\QP\,\cdot\,,z)$ enters through the bounds (ii) and (iii), which are uniform in $i$, and the limits $i\to\infty$ and $m'\to\infty$ are not interchanged.

\emph{(i) Convergence for fixed $m'$.} We show that $v_i^{m'}\to v^{m'}$, $\PP^1$-a.s., as $i\to\infty$, treating the compensator parts and the point-measure parts separately. For the compensators, \eqref{eq:g-L1nu} with $\epsilon=\epsilon_{J+1}$ and the contraction property of $\QP$ give
\[
\begin{aligned}
	&\left|\int_0^t\!\int_{Z_J^c}\langle g(\QP u_i^1(s),z)-g(\QP u_\ast^1(s),z),\phi\rangle\,\nu(dz)\,ds\right|\\
	&\qquad\le |\phi|_H(\nu_JL_g)^{1/2}\|u_i^1-u_\ast^1\|_{L^1(0,T;V)}\longrightarrow0
\end{aligned}
\]
by Proposition~\ref{prop:JS-pairs-allJ}\,\ref{JS:as}. For the point-measure parts put, for $(z,s)\in\ZZJ$,
\[
f_i(z,s):=\langle g(\QP u_i^1(s),z),\phi\rangle\,\mathbf 1_{(0,t]}(s),\qquad
f(z,s):=\langle g(\QP u_\ast^1(s),z),\phi\rangle\,\mathbf 1_{(0,t]}(s).
\]
Since $\QP x=\xi_k^{m'}(x)$ on $I_k^{m'}$, the estimate for $\xi_k^{m'}$ and \eqref{eq:XV-lipJ} give
\[
\sup_{\ZZJ}|f_i-f|\le|\phi|_H\,\Lt_J\,\tau_{m'}^{-1}\,\|u_i^1-u_\ast^1\|_{L^1(0,T;V)}\longrightarrow0 .
\]
The function $f$ is bounded on $\ZZJ$ by \eqref{it:upper}, and, by Hypothesis~\ref{ass:GF}\,\ref{it:cont}, it is continuous at every point of $\ZZJ$ whose time coordinate does not lie in the finite set $D:=\{t\}\cup\{t_k^{m'}:k=0,\dots,2^{m'}\}$. By Lemma~\ref{lem:etaJ-PRM}, $\EE^1\eta^{1,J}(Z_J^c\times D)=\nu_J\Leb(D)=0$, so $\eta^{1,J}(Z_J^c\times D)=0$, $\PP^1$-a.s. Since $N_i^{1,J}\to\eta^{1,J}$ weakly (Proposition~\ref{prop:JS-pairs-allJ}\,\ref{JS:as}), the mapping theorem \cite[Theorem~2.7]{Billingsley1999}, applied to the normalised measures, gives $\int f\,dN_i^{1,J}\to\int f\,d\eta^{1,J}$ (if $\eta^{1,J}(\ZZJ)=0$, this is trivial, since $f$ is bounded and $N_i^{1,J}(\ZZJ)\to0$). Hence, by Proposition~\ref{prop:JS-pairs-allJ}\,\ref{JS:point},
\[
\Bigl|\int f_i\,dN_i^{1,J}-\int f\,d\eta^{1,J}\Bigr|
\le\sup_{\ZZJ}|f_i-f|\;N_i^{1,J}(\ZZJ)+\Bigl|\int f\,dN_i^{1,J}-\int f\,d\eta^{1,J}\Bigr|\longrightarrow0 ,
\]
which proves (i). Note that the projection is essential here: for the integrand $g(u^1_\ast(s-),z)$ itself, the discontinuities in time would sit exactly at the atoms of $\eta^{1,J}$, and weak convergence could not be used.

\emph{(ii) Error at the approximating level.}
By \eqref{eq:XV-lipJ} and \eqref{eq:g-L1nu},
$|v_i-v_i^{m'}|\le |\phi|_H F_i$, where
\[
\begin{aligned}
	F_i:={}&\Lt_J\int_0^T\!\int_{Z_J^c}
	|\bar u_i^1(s)-\QP u_i^1(s)|_V\,N_i^{1,J}(dz,ds)\\
	&+(\nu_JL_g)^{1/2}
	\|\bar u_i^1-\QP u_i^1\|_{L^1(0,T;V)}.
\end{aligned}
\]
By Lemma~\ref{lem:XV-transfer}\,\ref{XV-tr:a},
$F_i=F(u_i^1)$ $\PP^1$-a.s., where $F$ is the Borel functional on
$\mathcal X$ obtained by replacing $N_i^{1,J}$, $\bar u_i^1$,
and $u_i^1$ with $\Phi_{r_i}^J(x)$, $\bar x_{r_i}$, and $x$,
respectively. Equality of laws (Proposition~\ref{prop:JS-pairs-allJ}\,\ref{eq:JS-laws}) therefore gives
\[
\EE^1[F_i\wedge1]=\EE[F(u_n)\wedge1],
\qquad n:=n_{r_i}.
\]
Take $i$ so large that $r_i\ge r_0(J)$ and $m:=m_{r_i}\ge m'$. Both
$\bar u_n=\Pro_m u_n$ and $\QP u_n$ are constant on each
cell of order $m$. On $\{u_n\in\Gamma_{r_i,J}\}$, conditions
($\Gamma$1) and ($\Gamma$2) pair the atoms of
$\Phi_{r_i}^J(u_n)$ with those of
$\mathcal J^J(u_n)=\eta\rest_{\ZZJ}$
(Lemma~\ref{lem:jump-identity}), one pair in each occupied cell.
Thus the point-measure integral in $F(u_n)$ is unchanged if
$\Phi_{r_i}^J(u_n)$ is replaced by $\eta\rest_{\ZZJ}$.
The integrand is predictable by
Lemma~\ref{lem:Pro-props}\,\ref{Pro:pred}, so the compensation formula \eqref{eq:PA-comp} yields
\[
\begin{aligned}
	\EE\int_0^T\!\int_{Z_J^c}
	|\bar u_n(s)-\QP u_n(s)|_V\,\eta(dz,ds)
	&=\nu_J\EE\|\bar u_n-\QP u_n\|_{L^1(0,T;V)}\\
	&\le\nu_J\bigl(\Theta_n(m)+\Theta_n(m')\bigr).
\end{aligned}
\]
Together with Part~\ref{XV-tr:2i} of the proof of
Lemma~\ref{lem:XV-transfer}, this implies
\[
\EE^1[F_i\wedge1]
\le c_{r_i,J}
+\bigl(\Lt_J\nu_J+(\nu_JL_g)^{1/2}\bigr)
\bigl(\Theta_{n_{r_i}}(m_{r_i})+
\Theta_{n_{r_i}}(m')\bigr).
\]
As $i\to\infty$, $c_{r_i,J}\to0$, and the term with $\Theta_{n_{r_i}}(m_{r_i})$ vanishes as well: by Part~\ref{XV-tr:2iii} of the proof of Lemma~\ref{lem:XV-transfer}, $(\nu_JL_g)^{1/2}\Theta_{n_r}(m_r)\le\delta_r$ for $r\ge r_0(J)$, and $\Lt_J\nu_J+(\nu_JL_g)^{1/2}=\bigl(\Lt_J(\nu_J/L_g)^{1/2}+1\bigr)(\nu_JL_g)^{1/2}$. For the remaining term, equality
of laws gives
\[
\Theta_{n_{r_i}}(m')
=\EE^1\|u_i^1-\QP u_i^1\|_{L^1(0,T;V)}.
\]
The map $x\mapsto\|x-\QP x\|_{L^1(0,T;V)}$ is continuous on
$L^1(0,T;V)$ and bounded by
$2\|x\|_{L^1(0,T;V)}+\tau_{m'}|u_0^{m'}|_V$.
Hence \eqref{eq:XV-apriori},
Proposition~\ref{prop:JS-pairs-allJ}\,\ref{JS:as}, and
Vitali's theorem imply
$\Theta_{n_{r_i}}(m')\to\theta_\ast(m')$. Altogether,
\[
\limsup_{i\to\infty}\EE^1[F_i\wedge1]\le\bigl(\Lt_J\nu_J+(\nu_JL_g)^{1/2}\bigr)\,\theta_\ast(m').
\]

\emph{(iii) Error at the limit.}
By Corollary~\ref{cor:usual-conditions}, $u_\ast^1$ is
progressively measurable. Consequently $\QP u_\ast^1$ is
predictable by the argument of
Lemma~\ref{lem:Pro-props}\,\ref{Pro:pred}, and so are
$u_\ast^1(s-)$ and $\mathds 1_{\{u^1_\ast(s-)\in V\}}$. The integrand of $v-v^{m'}$ coincides with $\langle g(u_\ast^1(s),z)-g(\QP u_\ast^1(s),z),\phi\rangle$ for a.e.\ $s$ (see the beginning of the proof). Hence the pathwise representation, the compensation formula and \eqref{eq:g-L1nu} give
\[
\begin{aligned}
	\EE^1|v-v^{m'}|
	&\le 2\EE^1\int_0^t\!\int_{Z_J^c}
	\bigl|\langle g(u_\ast^1(s),z)
	-g(\QP u_\ast^1(s),z),\phi\rangle\bigr|\,\nu(dz)\,ds\\
	&\le 2|\phi|_H(\nu_JL_g)^{1/2}\theta_\ast(m').
\end{aligned}
\]

\emph{Conclusion.} For $a,b,c\ge0$ we have $(a+b+c)\wedge1\le a\wedge1+b\wedge1+c\wedge1$, and $(|\phi|_HF_i)\wedge1\le(1+|\phi|_H)(F_i\wedge1)$. Hence \eqref{eq:XV-mart-split}, (ii), (iii), and (i) together with dominated convergence for the middle term give
\[
\limsup_{i\to\infty}\EE^1\bigl[|v_i-v|\wedge1\bigr]\le\bigl(1+|\phi|_H\bigr)\bigl(\Lt_J\nu_J+3(\nu_JL_g)^{1/2}\bigr)\,\theta_\ast(m') .
\]
The left-hand side does not depend on $m'$. Letting $m'\to\infty$ gives $\EE^1[|v_i-v|\wedge1]\to0$, i.e.\ $v_i\to v$ in $\PP^1$-probability, which is \eqref{eq:XV-mart-conv}.
\end{proof}

						\emph{Step 4:  Passing  $J\to\infty$.} Now, our aim is to pass to the limit $J\to \infty$. Here, the  following result is necessary.

\begin{lemma}\label{lem:XV-smalljumps}
Let $\tilde\eta^1$ be the compensated patched measure of Proposition~\ref{prop:full-levy-patch}, and let $Y$ be an $(\mathscr F^1_t)$-adapted process with c\`adl\`ag paths in $H$ such that $\EE^1\int_0^T|Y(s)|_V^2\,ds<\infty$. For $J\in\NN$ and $t\in[0,T]$ put
\[
X^J(t):=\int_0^t\!\!\int_{Z^c_J}g(Y(s-),z)\,\tilde\eta^{1,J}(dz,ds).
\]
Then $(X^J)_{J\in\NN}$ is a Cauchy sequence in $\mathcal M^2(0,T;H)$, and its limit is
\[
\int_0^t\!\!\int_{Z\setminus\{0\}}g(Y(s-),z)\,\tilde\eta^1(dz,ds),\qquad t\in[0,T].
\]
\end{lemma}
\begin{proof}
For $J'\ge J$ the difference $X^{J'}-X^J$ is the stochastic integral of $g(Y(s-),z)$ over the annulus $A_{J,J'}:=\{\epsilon_{J'+1}<|z|_Z\le\epsilon_{J+1}\}$. The integrand is predictable, and $\eta^1$ is a Poisson random measure compatible with $(\mathscr F^1_t)_{t\in[0,T]}$ by Proposition~\ref{prop:full-levy-patch}. So the maximal inequality \eqref{I4} with $p=2$ and $E=H$, whose constant we denote by $c$, gives
\[
\EE^1\Bigl[\sup_{t\in[0,T]}\bigl|X^{J'}(t)-X^J(t)\bigr|_H^2\Bigr]\le c\,\EE^1\int_0^T\!\!\int_{A_{J,J'}}|g(Y(s),z)|_H^2\,\nu(dz)\,ds ,
\]
since $Y(s-)=Y(s)$ for all but countably many $s$. By the uniform small-jump tail of Remark~\ref{rem:GF-consequences}, see \eqref{eq:XV-tail},
\[
\int_{A_{J,J'}}|g(Y(s),z)|_H^2\,\nu(dz)\le\kappa_2(\epsilon_{J+1})\bigl(1+|Y(s)|_V^2\bigr).
\]
Hence
\[
\EE^1\Bigl[\sup_{t\in[0,T]}\bigl|X^{J'}(t)-X^J(t)\bigr|_H^2\Bigr]\le c\,\kappa_2(\epsilon_{J+1})\Bigl(T+\EE^1\int_0^T|Y(s)|_V^2\,ds\Bigr)\xrightarrow[J\to\infty]{}0,
\]
uniformly in $J'\ge J$. Completeness of $\mathcal M^2(0,T;H)$ gives the limit, which coincides with the integral over $Z\setminus\{0\}$ by the construction of the extension of the stochastic integral.
\end{proof}
\par

						\emph{Step 5: Identification of drift.} The proof has been divided into several steps.

\quad\emph{Step 5 (a): Convergence of the drift.}
 Let $\phi\in V$ and put, for $t\in[0,T]$,
\[
d_i(t):=\int_0^t\langle A_{n_{r_i}}(s,u^1_i(s)),\phi\rangle\,ds,\qquad
d(t):=\int_0^t\langle A(s,u^1_\ast(s)),\phi\rangle\,ds .
\]
By (H4$_n$) and (H4), these integrals are finite, $\PP^1$-a.s.; for $d$ see also the bound in (i) below. We show
\begin{equation}\label{eq:XV-drift-conv}
\sup_{t\in[0,T]}|d_i(t)-d(t)|\longrightarrow0\qquad\text{in }\PP^1\text{-probability, as }i\to\infty .
\end{equation}
\par

{\emph{(i) Weak convergence along the paths.} Put $\alpha':=\alpha/(\alpha-1)$, and let $\Omega^1_0$ be the set of all $\omega\in\Omega^1$ such that $u^1_i(\omega)\to u^1_\ast(\omega)$ in $L^\alpha(0,T;V)$ and in $\DD([0,T];H)$. By Proposition~\ref{prop:JS-pairs-allJ}\,\ref{JS:as}, $\PP^1(\Omega^1_0)=1$. Fix $\omega\in\Omega^1_0$ and write $v_i:=u^1_i(\omega,\cdot)$ and $v:=u^1_\ast(\omega,\cdot)$. By Lemma~\ref{lem:J1-pointwise}\,(d),
\[
R:=\sup_{i\in\NN}\sup_{t\in[0,T]}|v_i(t)|_H+\sup_{t\in[0,T]}|v(t)|_H<\infty ,
\]
and by (H4), for every $i\in\NN$,
\[
\int_0^T|A(s,v_i(s))|_{V^\ast}^{\alpha'}\,ds\le(1+R^\beta)\Bigl(\|f\|_{L^1(0,T)}+C_0\sup_{j\in\NN}\int_0^T|v_j(s)|_V^\alpha\,ds\Bigr)=:K<\infty .
\]
The same bound holds for $v$. We show that
\begin{equation}\label{eq:XV-weak-path}
	A(\cdot,v_i)\rightharpoonup A(\cdot,v)\quad\text{weakly in }L^{\alpha'}(0,T;V^\ast),\ \text{as }i\to\infty .
\end{equation}
Since $V$ is reflexive and $1<\alpha'<\infty$, the space $L^{\alpha'}(0,T;V^\ast)$ is reflexive, and its dual is $L^\alpha(0,T;V)$. So every subsequence of $(v_i)_{i\in\NN}$ has a further subsequence $(v_{i_l})_{l\in\NN}$ such that $A(\cdot,v_{i_l})\rightharpoonup y$ weakly in $L^{\alpha'}(0,T;V^\ast)$ for some $y$, and $v_{i_l}(s)\to v(s)$ in $H$ for a.e.\ $s\in[0,T]$. It suffices to show that $y=A(\cdot,v)$.

Let $w\in L^\alpha(0,T;V)\cap L^\infty(0,T;H)$. Then $f+\rho(w)\in L^1(0,T)$ by the growth bound on $\rho$ in (H2), and $A(\cdot,w)\in L^{\alpha'}(0,T;V^\ast)$ by (H4). Integrating (H2) over $[0,T]$ gives
\[
2\int_0^T\bigl\langle A(s,v_{i_l}(s))-A(s,w(s)),v_{i_l}(s)-w(s)\bigr\rangle\,ds\le\int_0^T\bigl(f(s)+\rho(w(s))\bigr)\,|v_{i_l}(s)-w(s)|_H^2\,ds .
\]
On the left, $A(\cdot,v_{i_l})$ converges weakly in $L^{\alpha'}(0,T;V^\ast)$, and $v_{i_l}-w$ converges strongly in $L^\alpha(0,T;V)$. On the right, we use dominated convergence with the bound $|v_{i_l}-w|_H^2\le2\bigl(R^2+\|w\|_{L^\infty(0,T;H)}^2\bigr)$. Letting $l\to\infty$ gives
\[
2\int_0^T\bigl\langle y(s)-A(s,w(s)),v(s)-w(s)\bigr\rangle\,ds\le\int_0^T\bigl(f(s)+\rho(w(s))\bigr)\,|v(s)-w(s)|_H^2\,ds .
\]
Take $w=v-\varepsilon h$ with $h\in L^\infty(0,T;V)$ and $\varepsilon\in(0,1]$, and divide by $\varepsilon$:
\[
2\int_0^T\bigl\langle y(s)-A(s,v(s)-\varepsilon h(s)),h(s)\bigr\rangle\,ds\le\varepsilon\int_0^T\bigl(f(s)+\rho(v(s)-\varepsilon h(s))\bigr)\,|h(s)|_H^2\,ds .
\]
By the growth bound on $\rho$, the integral on the right is bounded uniformly in $\varepsilon\in(0,1]$. On the left, (H1) gives $\langle A(s,v(s)-\varepsilon h(s)),h(s)\rangle\to\langle A(s,v(s)),h(s)\rangle$ as $\varepsilon\downarrow0$, for a.e.\ $s\in[0,T]$. Here (H1) is used at the point $v(s)=u^1_\ast(\omega,s)$, which depends on $s$. This is why the null set of times in (H1) must not depend on $w_1,w_2,w_3$. By (H4), $|\langle A(s,v(s)-\varepsilon h(s)),h(s)\rangle|\le\|h\|_{L^\infty(0,T;V)}\,G(s)^{1/\alpha'}$ for $\varepsilon\in(0,1]$, where
\[
G:=\bigl(f+2^{\alpha-1}C_0\,(|v|_V^\alpha+|h|_V^\alpha)\bigr)\bigl(1+(R+\|h\|_{L^\infty(0,T;H)})^\beta\bigr)\in L^1(0,T).
\]
So dominated convergence, as $\varepsilon\downarrow0$, gives $\int_0^T\langle y(s)-A(s,v(s)),h(s)\rangle\,ds\le0$. Replacing $h$ by $-h$, we get equality for all $h\in L^\infty(0,T;V)$. This set is dense in $L^\alpha(0,T;V)$. Hence $y=A(\cdot,v)$, and \eqref{eq:XV-weak-path} is proved.

\emph{(ii) Uniformity in $t$.} For $t\in[0,T]$, the function $\mathds 1_{[0,t]}\phi$ belongs to $L^\alpha(0,T;V)$. So \eqref{eq:XV-weak-path} gives $\int_0^t\langle A(s,v_i(s)),\phi\rangle\,ds\to\int_0^t\langle A(s,v(s)),\phi\rangle\,ds$ for every $t\in[0,T]$. By H\"older's inequality, for $0\le t'\le t\le T$ and all $i\in\NN$,
\[
\Bigl|\int_{t'}^t\langle A(s,v_i(s)),\phi\rangle\,ds\Bigr|\le|\phi|_V\,K^{1/\alpha'}\,|t-t'|^{1/\alpha},
\]
and the same holds for $v$. An equicontinuous sequence of functions on $[0,T]$ that converges pointwise converges uniformly. Since $\omega\in\Omega^1_0$ was arbitrary,
\begin{equation}\label{eq:XV-drift-A}
\PP^1\lk( 	\sup_{t\in[0,T]}\Bigl|\int_0^t\bigl\langle A(s,u^1_i(s))-A(s,u^1_\ast(s)),\phi\bigr\rangle\,ds\Bigr|\stackrel{i\to\infty}\longrightarrow0\rk)=1.
\end{equation}
\par}

\emph{(iii) Consistency of the drifts.} The map
\[
\mathcal X\ni x\mapsto\int_0^T\bigl|\bigl\langle A_{n_{r_i}}(s,x(s))-A(s,x(s)),\phi\bigr\rangle\bigr|\,ds
\]
is Borel on $\mathcal X$ (see the measurability paragraph of Step~1). Since $\Law_{\PP^1}(u^1_i)=\Law_{\PP}(u_{n_{r_i}})$, we have for every $\ep>0$
\begin{align*}
&\PP^1\lk(\int_0^T\bigl|\bigl\langle A_{n_{r_i}}(s,u^1_i(s))-A(s,u^1_i(s)),\phi\bigr\rangle\bigr|\,ds>\ep\rk)\\
&\quad=\PP\lk(\int_0^T\bigl|\bigl\langle A_{n_{r_i}}(s,u_{n_{r_i}}(s))-A(s,u_{n_{r_i}}(s)),\phi\bigr\rangle\bigr|\,ds>\ep\rk)\xrightarrow[i\to\infty]{}0 ,
\end{align*}
where the limit follows from \eqref{eq:DG} along the subsequence $(n_{r_i})$. Hence the integral tends to $0$ in $\PP^1$-probability.
\par

\emph{(iv) Conclusion.} For every $t\in[0,T]$,
\[
|d_i(t)-d(t)|\le\int_0^T\bigl|\bigl\langle A_{n_{r_i}}(s,u^1_i(s))-A(s,u^1_i(s)),\phi\bigr\rangle\bigr|\,ds+\Bigl|\int_0^t\bigl\langle A(s,u^1_i(s))-A(s,u^1_\ast(s)),\phi\bigr\rangle\,ds\Bigr| .
\]
By (iii), the first term tends to $0$ in probability. By \eqref{eq:XV-drift-A}, the supremum over $t$ of the second term tends to $0$, $\PP^1$-a.s. This proves \eqref{eq:XV-drift-conv}.

\quad\emph{Step 5 (b): (The limit equation).}
 Fix $\phi\in D$ and $J\in\NN$. Let $S\subset[0,T)$ be the at most countable set of times $t$ with $\PP^1(u^1_\ast(t)\neq u^1_\ast(t-))>0$, and let $t\in[0,T]\setminus S$. Put
\[
\Delta^J(t):=\langle u^1_\ast(t),\phi\rangle-\int_0^t\langle A(s,u^1_\ast(s)),\phi\rangle\,ds-\int_0^t\!\!\int_{Z^c_J}\langle g(u^1_\ast(s-),z),\phi\rangle\,\widetilde\eta^{1,J}(dz,ds)=:a-d-n ,
\]
and write \eqref{eq:XV-approx} as $a_i-d_i-n_i=\langle\rho^J_{r_i}(u^1_i)(t),\phi\rangle$. Here $d_i=d_i(t)$ and $d=d(t)$ are as in Step~5(a), and
\[
a_i:=\langle u^1_i(t),\phi\rangle,\qquad n_i:=\int_0^t\!\!\int_{Z^c_J}\langle g(\bar u^1_i(s),z),\phi\rangle\,\widetilde N^{1,J}_i(dz,ds).
\]
We show
\begin{equation}\label{eq:XV-DeltaJ}
\EE^1|\Delta^J(t)|\le|\phi|_H\bigl(\kappa_2(\epsilon_{J+1})\,C_2\bigr)^{1/2}.
\end{equation}

\emph{Convergence of the three terms.} If $t<T$, then, as $t\notin S$, $t$ is $\PP^1$-a.s.\ a continuity point of $u^1_\ast$, so $u^1_i(t)\to u^1_\ast(t)$ in $H$ by the $J_1$-convergence (Lemma~\ref{lem:J1-pointwise}); for $t=T$ this holds as well. Since $u_0^{n_{r_i}}\to u_0$ in $H$, $a_i\to a$, $\PP^1$-a.s. By Step~5(a), $d_i\to d$ in probability, and by Lemma~\ref{lem:XV-mart-fixedJ}, $n_i\to n$ in probability. Let $C_i=C_i(t)$ be the sets of Lemma~\ref{lem:XV-transfer}\,\ref{XV-tr:c}; by \eqref{eq:XV-Ci}, $\mathds 1_{C_i}\to1$ in probability. Hence
\[
\mathds 1_{C_i}\langle\rho^J_{r_i}(u^1_i)(t),\phi\rangle=\mathds 1_{C_i}(a_i-d_i-n_i)\longrightarrow a-d-n=\Delta^J(t)\qquad\text{in }\PP^1\text{-probability.}
\]

\emph{The residual on the good sets.} This is where the good sets enter. By \eqref{eq:XV-rho-bound}, \eqref{eq:XV-Mnat-bound} and \eqref{eq:XV-E-bound},
\[
\EE^1\bigl|\mathds 1_{C_i}\langle\rho^J_{r_i}(u^1_i)(t),\phi\rangle\bigr|\le|\phi|_H\Bigl(\bigl(\kappa_2(\epsilon_{J+1})C_2\bigr)^{1/2}+\bar c_{r_i,J}\Bigr)
\xrightarrow[i\to\infty]{}|\phi|_H\bigl(\kappa_2(\epsilon_{J+1})C_2\bigr)^{1/2}.
\]
Along a subsequence, the convergence in probability above holds $\PP^1$-a.s., and Fatou's lemma gives \eqref{eq:XV-DeltaJ}. No uniform integrability of $a_i$, $d_i$ or $n_i$ is needed.

\emph{The limit $J\to\infty$.} By Lemma~\ref{lem:XV-smalljumps}, applied to  $Y=u^1_\ast$ (note that $\EE^1\int_0^T|u^1_\ast(s)|_V^2\,ds<\infty$ by \eqref{eq:XV-limit-apriori} and $\alpha\ge2$),
\[
\int_0^t\!\!\int_{Z^c_J}\langle g(u^1_\ast(s),z),\phi\rangle\,\widetilde\eta^{1,J}(dz,ds)\longrightarrow\int_0^t\!\!\int_{Z\setminus\{0\}}\langle g(u^1_\ast(s),z),\phi\rangle\,\widetilde\eta^1(dz,ds)
\quad\text{in }L^2(\Omega^1),
\]
and $\kappa_2(\epsilon_{J+1})\to0$. Hence \eqref{eq:XV-DeltaJ} gives, $\PP^1$-a.s.,
\[
\langle u^1_\ast(t),\phi\rangle=
\int_0^t\langle A(s,u^1_\ast(s)),\phi\rangle\,ds+\int_0^t\!\!\int_{Z\setminus\{0^\ast\}}\langle g(u^1_\ast(s),z),\phi\rangle\,\widetilde\eta^1(dz,ds).
\]
This holds for every $\phi\in D$ and every $t\in[0,T]\setminus S$, a set which is dense in $[0,T]$ and contains $T$. Both sides are c\`adl\`ag in $t$ (the stochastic integral has a c\`adl\`ag modification by Theorem~\ref{th-ineq}), so a common null set serves first for a countable dense set of times and then for all $t\in[0,T]$. By the density of $D$ in $V$, we obtain for all $t\in[0,T]$, $\PP^1$-a.s.,

						\begin{equation}\label{eq:XV-chi-eq}
							u^1_\ast(t)
							=
							{\int_0^t A(s,u^1_\ast(s))\,ds}
							+\int_0^t\!\!\int_{Z\setminus\{0^\ast\}}
							g(u^1_\ast(s),z)\,\tilde\eta^1(dz,ds)
						\end{equation}
 as an identity in $V^\ast$.

						\emph{Step 6: Conclusion.} By \eqref{eq:XV-chi-eq}, $u^1_\ast$ solves the equation of Theorem~\ref{Theoremmain} on $\mathfrak A^\ast$ (Corollary~\ref{cor:usual-conditions}). Together with \eqref{eq:XV-limit-apriori}, the regularity $u^1_\ast\in\DD([0,T];H)\cap L^\alpha(0,T;V)$ $\PP^1$-a.s., the progressive measurability of $u^1_\ast$ and the properties of $\eta^1$ from Corollary~\ref{cor:usual-conditions} and Proposition~\ref{prop:full-levy-patch}, this establishes conclusion (iv) of Theorem~\ref{Theoremmain} and completes Step~(XIV).\par

						\end{steps}
					\end{proof}

				\appendix

\section{{Random Measures and Their Convergence}}
\label{app:levy-linde}

\newcommand{\aslim}{\xrightarrow{\;\mathrm{a.s.}\;}}
\newcommand{\dto}{\xrightarrow{\;d\;}}
\newcommand{\wto}{\xrightarrow{\;w\;}}

\newcommand{\Nbar}{\bar{N}}
\newcommand{\nueps}{\nu_{\varepsilon}}
\newcommand{\lto}{\longrightarrow}
\newcommand{\MM}{\mathbf M}
\renewcommand{\MMN}{\mathbf M_{\NN_0}}
\newcommand{\MMloc}{\mathbf M_{\mathrm{loc}}}
\renewcommand{\bNN}{{\NN_0}}
This appendix collects the measure-theoretic results used in the reconstruction argument. We recall the topologies and tightness notions for random measures, with particular attention to non-locally compact state spaces, and construct the product space of truncated jump measures used in Step~(IX).

\medskip
\subsection[Random measures and the vague topology]{{Random measures and the vague topology}}\label{app:convergenceofmeasures}

Let us recall the general notion of a random measure. Let
$(S,\mathscr S)$ be a measurable space and let
$(\Omega,\mathcal F,\mathbb P)$ be a probability space.

\begin{definition}\label{def_rm}
A \emph{random measure} $\xi$ on $(S,\mathscr S)$ is a mapping $\xi:\Omega\times\mathscr S\longrightarrow[0,\infty]$ such that

\begin{enumerate}[label=\upshape(\roman*)]

	\item for every $\omega\in\Omega$, the mapping $B\longmapsto\xi(\omega,B)$ is a measure on $(S,\mathscr S)$;

	\item for every $B\in\mathscr S$, the mapping $\omega\longmapsto\xi(\omega,B)$ is $\mathcal F$-measurable.

\end{enumerate}
\end{definition}

\begin{remark}\label{def_rm_filt}
Let $(\Omega,\mathcal F,(\mathcal F_t)_{t\in[0,T]},\mathbb P)$ be a filtered probability space. A random measure $N$ on
$S\times(0,T]$ is called adapted to $(\mathcal F_t)_{t\in[0,T]}$ if, for every $B\in\mathscr S$ and $t\in[0,T]$, the random variable $N(B\times(0,t])$ is $\mathcal F_t$-measurable.
\end{remark}

For general background on random measures and point processes, we refer to
Kallenberg \cite{Kallenberg2017}, Daley--Vere-Jones
\cite{DaleyVereJones2003}, and Br\'emaud \cite{Bremaud1981}; for marked
Poisson processes, see Last--Penrose \cite{LastPenrose2018}.

\medskip

Why are we introducing a new topology?
The random measures in this paper are
Poisson random measures on $Z\times[0,T]$ with intensity $\nu\otimes\mathrm{Leb}$,
and their approximations. Two features rule out the standard topologies. First,
if $\nu(Z)=\infty$, a Poisson random measure has, in every time interval,
infinitely many atoms near the point $0^\ast$; it is not a finite measure, so
the weak topology of finite measures cannot be used.
Secondly, the mark space $Z$
may be infinite-dimensional and is then not locally compact, so the classical
vague topology, recalled next, degenerates. The $w^\#$-topology of
Subsection~\ref{app:wsharp} removes both obstacles: measures are tested only
against bounded continuous functions that vanish outside a \emph{bounded} set,
and the metric is chosen such that ``bounded'' means ``bounded away from
$0^\ast$'' (Subsection~\ref{par:setting-Sn}).

For point processes one works with integer-valued measures. When the underlying space is locally compact and Polish, the natural topology is the vague topology on locally finite measures; for general Polish spaces we will instead use the $w^\#$-topology on boundedly finite measures. The two settings are stated separately below.

\begin{definition}[Vague topology and random measures]\label{def:vague-rm}
Let $S$ be a locally compact, second countable Hausdorff space (hence Polish),
and let $\widehat{\mathcal S}$ denote the class of relatively compact Borel
subsets of $S$. We write $\MMloc(S)$ for the set of Borel measures $\mu$ on $S$
that are locally finite, i.e.\ $\mu(B)<\infty$ for all $B\in\widehat{\mathcal S}$,
and we equip $\MMloc(S)$ with the \emph{vague topology}, i.e.\ the coarsest
topology making the maps
\[
\MMloc(S)\ni\mu\longmapsto\mu(f):=\int_S f\,d\mu\in\RR,
\qquad f\in C_c(S),
\]
continuous. Then $\MMloc(S)$ is a Polish space and its Borel $\sigma$-field
coincides with $\sigma\bigl(\mu\mapsto\mu(B):B\in\widehat{\mathcal S}\bigr)$; see
\cite[Thm.~4.2]{Kallenberg2017}. A \emph{random measure} on $S$ over a
probability space $(\Omega,\mathcal F,\PP)$ is a measurable map
$\xi:\Omega\to\MMloc(S)$, and $\mathbf M_{\mathrm{loc},\NN_0}(S)\subset\MMloc(S)$
denotes the closed subset of $\NN_0$-valued (integer-valued) locally finite
measures.
\end{definition}

For an infinite-dimensional normed space, local compactness fails and $C_c(S)=\{0\}$; hence the ordinary vague topology is too weak for our purposes. This motivates the $w^\#$-topology introduced below, which is available on arbitrary complete separable metric spaces and agrees with the vague topology in the proper locally compact case. For more information, see Daley--Vere-Jones \cite[App.~A2.6]{DaleyVereJones2003}.

\subsection[The $w^\#$-topology]{{The $w^\#$-topology}}
\label{app:wsharp}

Throughout this paragraph no local compactness is assumed. Following
Daley--Vere-Jones \break \cite[App.~A2.6]{DaleyVereJones2003}, we work with
\emph{boundedly finite} measures, i.e.\ measures that are finite on bounded
sets, and with the $w^\#$-topology generated by the bounded continuous
functions vanishing outside a bounded set. It is convenient to describe the
bounded sets through a localizing sequence.

\smallskip\noindent
\begin{definition}[Localizing sequences and the $w^\#$-topology]
\label{def:wsharp}
Let $E$ be a Polish space. A sequence $(E_n)_{n\in\NN_0}$ of Borel subsets of
$E$ is called \emph{localizing} if $E_0\subset E_1\subset\cdots$,
$\bigcup_nE_n=E$, and there is a complete metric $d_E$ inducing the topology
of $E$ whose bounded sets are exactly the subsets of some $E_n$.
\begin{enumerate}[label=(\alph*)]
	\item $\MM(\{E_n\})$ denotes the set of Borel measures $\mu$ on $E$ with
	$\mu(E_n)<\infty$ for all $n$. $\MMN(\{E_n\})\subset\MM(\{E_n\})$ is the
	subset of those with $\mu(B)\in\NN_0$ for every Borel $B\subset E_n$ and
	every $n$. Note that $\mu(E)=\infty$ is allowed.
	\item $C_b(\{E_n\})$ denotes the bounded continuous $f:E\to\RR$ vanishing
	outside some $E_n$.
	\item The \emph{$w^\#$-topology} on $\MM(\{E_n\})$ is the coarsest topology
	making $\mu\mapsto\int_Ef\,d\mu$ continuous for every $f\in C_b(\{E_n\})$;
	we write $\mu_k\xrightarrow{w^\#}\mu$.
\end{enumerate}
\end{definition}

\smallskip\noindent
For \emph{finite measures} on a Polish space $E$ we write $\MM(E)$ for the
finite Borel measures on $E$, and $\MMN(E)\subset\MM(E)$ for the finite
$\NN_0$-valued ones, i.e.\ the finite point measures. Both carry the weak
topology: $\mu_k\to\mu$ if and only if $\int_Ef\,d\mu_k\to\int_Ef\,d\mu$ for
every $f\in C_b(E)$. Then $\MM(E)$ is Polish, and $\MMN(E)$ is a closed subset,
hence also Polish \cite[App.~A2.5, Prop.~A2.6.III]{DaleyVereJones2003}.

\begin{remark}\label{rem:wsharp-finite-case}
If $E=E_n$ for all $n$, then $\MM(\{E_n\})=\MM(E)$, $C_b(\{E_n\})=C_b(E)$,
and the $w^\#$-topology is the weak topology. The localizing sequence matters
only when the measures of interest have infinite total mass. In our
application, this happens because the number of small jumps accumulates
at~$0^\ast$.
\end{remark}

\begin{proposition}[{Daley--Vere-Jones \cite[App.~A2.6]{DaleyVereJones2003}}]
\label{prop:wsharp-props}
Let $E$ be a Polish space with a localizing sequence $(E_n)_{n\in\NN_0}$ in
the sense of Definition~\ref{def:wsharp}, and let $\widehat{\mathcal E}$
denote the class of Borel sets $B\subset E$ with $B\subset E_n$ for some $n$.
Then the following hold.
\begin{enumerate}[label=\textup{(\roman*)}]
	\item\label{it:ws-equiv}
	For $\mu,\mu_k\in\MM(\{E_n\})$ one has $\mu_k\xrightarrow{w^\#}\mu$ if and
	only if $\mu_k|_B\Rightarrow\mu|_B$ weakly, as finite measures on $E$, for
	every $B\in\widehat{\mathcal E}$ with $\mu(\partial B)=0$. In particular, if
	$\mu(\partial E_n)=0$ for all $n$, then
	\[
	\mu_k\xrightarrow{w^\#}\mu
	\iff
	\mu_k|_{E_n}\Rightarrow\mu|_{E_n}\ \text{weakly for every }n\in\NN_0 .
	\]
	\item\label{it:ws-polish}
	$\MM(\{E_n\})$, equipped with the $w^\#$-topology, is a Polish space, and
	$\MMN(\{E_n\})$ is a closed subset of $\MM(\{E_n\})$, hence also Polish.
	\item\label{it:ws-borel}
	The Borel $\sigma$-field of the $w^\#$-topology on $\MM(\{E_n\})$ coincides
	with $\sigma\bigl(\mu\mapsto\mu(B):\ B\in\widehat{\mathcal E}\bigr)$.
	Consequently, a map $\xi:\Omega\to\MM(\{E_n\})$ is measurable if and only if
	$\omega\mapsto\xi(\omega,B)$ is measurable for every
	$B\in\widehat{\mathcal E}$, i.e.\ if and only if $\xi$ is a random measure
	in the sense of Definition~\ref{def_rm} taking values in $\MM(\{E_n\})$.
\end{enumerate}
\end{proposition}

\begin{proof}
Let $d_E$ be a complete metric as in Definition~\ref{def:wsharp}, whose
bounded sets are exactly the elements of $\widehat{\mathcal E}$. Then
$\MM(\{E_n\})$ is the space of boundedly finite measures on $(E,d_E)$, and
$C_b(\{E_n\})$ is the class of bounded continuous functions vanishing outside
a $d_E$-bounded set. Hence the $w^\#$-topology of Definition~\ref{def:wsharp}
is the $w^\#$-topology of \cite[App.~A2.6]{DaleyVereJones2003} for
$(E,d_E)$, and \ref{it:ws-equiv}--\ref{it:ws-borel} are the statements proved
there. For the second part of~\ref{it:ws-equiv}: ``$\Rightarrow$'' is the first
part with $B=E_n$; for ``$\Leftarrow$'', every $f\in C_b(\{E_n\})$ vanishes
outside some $E_n$, so
$\int f\,d\mu_k=\int f\,d(\mu_k|_{E_n})\to\int f\,d(\mu|_{E_n})=\int f\,d\mu$.
\end{proof}

\begin{remark}[Relation to the vague topology]\label{rem:wsharp-vs-vague}
Let $E$ be a Polish space with a complete metric $d_E$ that is
\emph{proper}, i.e.\ closed bounded subsets of $E$ are compact (for instance
$E=\RR^d$, or $E=Z\times[0,T]$ with $Z\subset\RR^d$ closed), and let $E_n$ be
the closed $d_E$-balls of radius $n$ around a fixed point. Then $E$ is locally
compact, the Borel sets contained in some $E_n$ are exactly the relatively
compact Borel sets, $C_b(\{E_n\})=C_c(E)$, and
\[
\MM(\{E_n\})=\MMloc(E),
\qquad w^\#\text{-topology}=\text{vague topology}.
\]
In this case the two notions may be used interchangeably, and
Proposition~\ref{prop:wsharp-props} reduces to
\cite[Thm.~4.19, Thm.~4.20]{Kallenberg2017}.
	\end{remark}

	\subsection[Point measures on the jump space]{{Point measures on the jump space}}
	\label{par:setting-Sn}
	We now describe the space on which the Poisson random measures of
	the main text live. If $\nu(Z)=\infty$, such a measure has, in every time
	interval, infinitely many atoms near the point $0^\ast$, but only finitely many
	atoms at a distance more than $\varepsilon$ from $0^\ast$, for every
	$\varepsilon>0$. The sets $S_n$ introduced below exclude shrinking
	neighbourhoods of $0^\ast$; on each of them the measure is finite.

	Let $Z$ be a complete separable metric space with metric $d$ and a
	distinguished point $0^\ast\in Z$. Write $|z|_Z:=d(z,0^\ast)$, and let
	$\mathscr Z:=\mathscr B(Z)$ be the Borel $\sigma$-algebra.
Let $\nu$ be a $\sigma$-finite measure on $(Z,\mathscr Z)$ with
$\nu(\{0^\ast\})=0$, with $\nu(\{|z|_Z>\varepsilon\})<\infty$ for every
$\varepsilon>0$, and with bounded support,
$\operatorname{supp}\nu\subset\{|z|_Z\le R_0\}$ (Hypothesis~\ref{ass:GF}).
Thus $\nu$ is finite away from $0^\ast$; its possibly infinite mass sits only
near $0^\ast$, which is why the point $0^\ast$ is removed below. For instance,
these conditions hold if $\nu(\{0^\ast\})=0$, the support is bounded and
$\int_Z(1\wedge|z|_Z^2)\,\nu(dz)<\infty$, since
$r^2\,\nu(\{|z|_Z\ge r\})\le\int_Z(1\wedge|z|_Z^2)\,\nu(dz)$ for $r\in(0,1]$.

Let
\[
\mathcal R(\nu):=\bigl\{r>0:\ \nu(\{|z|_Z=r\})=0\bigr\}
\]
be the set of \emph{$\nu$-continuity levels}. Its complement in $(0,\infty)$ is
at most countable.

Fix $\epsilon_0>R_0$ and choose
$\epsilon_n\in\mathcal R(\nu)$ for every $n\ge1$ so that
$\epsilon_0>\epsilon_1>\epsilon_2>\dots\downarrow0$. Since
$\epsilon_0>R_0$ and $\operatorname{supp}\nu\subset\{|z|_Z\le R_0\}$, one also
has $\nu(\{|z|_Z=\epsilon_0\})=0$. Put
\[
Z^\circ:=\bigl\{z\in Z:\ 0<|z|_Z\le\epsilon_0\bigr\},\qquad
Z_n:=\bigl\{z\in Z:\ \epsilon_{n+1}<|z|_Z\le\epsilon_0\bigr\}
=\bigcup_{j\le n}\Delta Z_j,
\quad n\in\NN_0 .
\]
That is, $Z_n$ is obtained from $Z^\circ$ by removing the closed ball of radius
$\epsilon_{n+1}$ around $0^\ast$; in the notation of Step~(IX), $Z_n=Z^c_n$ and
$Z_n\times[0,T]=\ZZJ$ with $J=n$. Each $Z_n$ is bounded and, as the
intersection of an open and a closed subset of $Z$, a Polish space, and
\[
Z_n\uparrow Z^\circ,\qquad
\nu(Z_n)<\infty,\qquad
\nu(Z\setminus Z^\circ)=0,\qquad
\nu(\partial_{Z^\circ}Z_n)=\nu(\{|z|_Z=\epsilon_{n+1}\})=0 ,
\]
where $\partial_{Z^\circ}$ denotes the boundary relative to $Z^\circ$.
On the time--space level we set
\[
S:=Z^\circ\times[0,T],\qquad
S_n:=Z_n\times[0,T],\quad n\in\NN_0,
\]

so that $S_n\uparrow S$ and $(\nu\otimes\mathrm{Leb})(S_n)=T\,\nu(Z_n)<\infty$.
Random measures on $S$ with possibly infinitely many small jumps take values in the
spaces $\MM(\{S_n\})$ and $\MMN(\{S_n\})$ introduced next.

\medskip

\smallskip\noindent
\emph{Bounded sets and the $w^\#$-topology on $S$.} We call a Borel set
$B\subset S$ \emph{bounded} if $B\subset S_n$ for some $n$, and write
$\widehat{\mathcal S}$ for the class of such sets.

Define
\[
d_S((x,s),(y,t))
:=d(x,y)+\bigl|\,|x|_Z^{-1}-|y|_Z^{-1}\bigr|+|s-t|,
\qquad (x,s),(y,t)\in S .
\]
This is a complete metric inducing the original topology on $S$. Indeed, if
$((x_k,t_k))$ is $d_S$-Cauchy, then $(x_k)$ is $d$-Cauchy and $(t_k)$ is
Cauchy, while $(|x_k|_Z^{-1})$ is Cauchy and hence bounded. Thus
$x_k\to x$ in $Z$, $t_k\to t\in[0,T]$, and the boundedness of
$|x_k|_Z^{-1}$ implies $x\neq0^\ast$; hence $(x,t)\in S$ and the convergence
also holds in $d_S$. The  term $\bigl|\,|x|_Z^{-1}-|y|_Z^{-1}\bigr|$ behaves like the map $r\mapsto \frac 1r$, which sends the interval
$(0,\ep_0]$ to $[\frac 1{\ep_0},\infty)$. If we have $(x_k)\to 0^\ast $ for $k\to\infty$, the problem would be that $0^\ast \not \in Z$ and the $d_S$-distance tends to $\infty$.
So \emph{close} to $0^\ast$ becomes \emph{far away}. Bounded sets are exactly the sets that stay at a positive distance from $0^\ast$, and there $\nu$ is finite.
Hence, $z\mapsto |z|_Z^{-1}$ is continuous on
$Z\setminus\{0^\ast\}$, so $d_S$ induces the original topology on $S$.
In other words, since $|z|_Z\le\epsilon_0$ on $S$,
\[
B\subset S\text{ is }d_S\text{-bounded}
\quad\Longleftrightarrow\quad
\inf_{(z,t)\in B}|z|_Z>0
\quad\Longleftrightarrow\quad
B\subset S_n\text{ for some }n.
\]
Thus $0^\ast$ is sent \emph{to infinity} and $(S_n)_{n\in\NN_0}$ is a localizing
sequence in the sense of Definition~\ref{def:wsharp}.

Accordingly, $\MM(\{S_n\})$ is the set of Borel measures on $S$ that are
finite on every $S_n$, and \break  $\MMN(\{S_n\})\subset\MM(\{S_n\})$ is the subset of
$\NN_0$-valued ones. Both carry the $w^\#$-topology of
Definition~\ref{def:wsharp}: $\mu_k\to\mu$ if and only if
$\int_Sf\,d\mu_k\to\int_Sf\,d\mu$ for every bounded continuous
$f:S\to\RR$ vanishing outside some $S_n$. By
Proposition~\ref{prop:wsharp-props}, $\MMN(\{S_n\})$ is Polish, and its Borel
$\sigma$-algebra is generated by the evaluations $\mu\mapsto\mu(B)$,
$B\in\widehat{\mathcal S}$.

A Poisson random measure $\eta$ on $S$ with intensity $\nu\otimes\mathrm{Leb}$
is a.s.\ an element of $\MMN(\{S_n\})$. Since $\epsilon_{n+1}\in\mathcal R(\nu)$ for every
$n\in\NN_0$, it has a.s.\ no atom on
$\{|z|_Z=\epsilon_{n+1}\}\times[0,T]$ for any $n$.

\subsection[Compactness and tightness]{{Compactness and tightness}}
\label{app:tightness}

We are interested in tightness of measure valued random variables.
We first define tightness and then characterise the (deterministic)
compact sets of $\MM(\{E_n\})$.

\begin{definition}[Tightness of a family of random measures]\label{def:tight-rm}
Let $E$ be a Polish space with a localizing sequence $(E_n)_{n\in\NN_0}$ in
the sense of Definition~\ref{def:wsharp}. A family $\{\xi_m\}_{m\in\NN}$ of
random measures with values in $\MM(\{E_n\})$ (respectively in
$\MMN(\{E_n\})$), defined on probability spaces
$(\Omega_m,\mathcal F_m,\PP_m)$, is called \emph{tight} (in the
$w^\#$-topology) if for every $\varepsilon>0$ there exists a set
$K_\varepsilon\subset\MM(\{E_n\})$ (respectively $K_\varepsilon\subset\MMN(\{E_n\})$),
compact in the $w^\#$-topology, such that
\[
\sup_{m\in\NN}\ \PP_m\bigl(\xi_m\notin K_\varepsilon\bigr)\le\varepsilon .
\]
Since these spaces are Polish
(Proposition~\ref{prop:wsharp-props}\,\ref{it:ws-polish}), this is equivalent
to relative compactness in distribution of $\{\Law(\xi_m)\}_{m\in\NN}$, by
Prokhorov's theorem.
\end{definition}

\begin{lemma}[Relative compactness of deterministic sets in the $w^\#$-topology]
\label{lem:compact-Mn}
Let $E$ be a Polish space with a localizing sequence $(E_n)_{n\in\NN_0}$ in
the sense of Definition~\ref{def:wsharp}. A set
$\mathcal K\subset\MM(\{E_n\})$ is relatively compact in the
$w^\#$-topology if and only if, for every $n\in\NN_0$,
\begin{enumerate}
	\item[(K1)] $\displaystyle\sup_{\mu\in\mathcal K}\mu(E_n)<\infty$, and
	\item[(K2)] for every $\eta>0$ there is a compact set $C\subset E$ with
	$\displaystyle\sup_{\mu\in\mathcal K}\mu(E_n\setminus C)<\eta$.
\end{enumerate}
If $\mathcal K\subset\MMN(\{E_n\})$, then (K2) is equivalent to
\begin{enumerate}
	\item[(K2$'$)] there is a compact set $C_n\subset E$ with
	$\mu(E_n\setminus C_n)=0$ for every $\mu\in\mathcal K$,
\end{enumerate}
i.e.\ on each $E_n$ the total mass stays bounded \emph{and} all atoms lie in one
fixed compact set. If the $E_n$ have compact closure, (K2) holds with
$C=\overline{E_n}$, and relative compactness reduces to~(K1).
\end{lemma}

\begin{proof}
Let $d_E$ be a complete metric as in Definition~\ref{def:wsharp}, and let
$U_k$ be the open $d_E$-balls of radius $k$ around a fixed $x_0\in E$. Each
$E_n$ lies in some $U_k$, and each $U_k$ lies in some $E_n$. Conditions (K1)
and (K2) are inherited by subsets, so they hold for all $E_n$ if and only if
they hold for all $U_k$. In that form, the lemma is
\cite[Prop.~A2.6.IV]{DaleyVereJones2003}. Daley--Vere-Jones use closed balls;
this is equivalent, since the closed ball of radius $k$ lies between $U_k$ and
$U_{k+1}$.

For integer-valued measures, (K2) with $\eta<1$ gives a compact $C_n$ with
$\mu(E_n\setminus C_n)<1$, hence $=0$, for all $\mu\in\mathcal K$. That is
(K2$'$). Conversely, (K2$'$) gives (K2) for every $\eta>0$ with $C=C_n$.
\end{proof}

We now pass from compact sets of measures to tightness of random measures.

\begin{theorem}[Tightness criterion for random measures in the $w^\#$-topology]
\label{thm:kallenberg-tight-wsharp}
Let $E$ be a Polish space with a localizing sequence $(E_n)_{n\in\NN_0}$ in
the sense of Definition~\ref{def:wsharp}, and let $\{\xi_m\}_{m\in\NN}$ be
random measures on $E$, defined on probability spaces
$(\Omega_m,\mathcal F_m,\mathbb P_m)$, with values in $\MM(\{E_n\})$
(respectively in $\MMN(\{E_n\})$). Then the following assertions are
equivalent.
\begin{enumerate}[label=\textup{(\roman*)}]
	\item\label{it:ktw-tight}
	The family $\{\xi_m\}_{m\in\NN}$ is tight in $\MM(\{E_n\})$ (respectively
	in $\MMN(\{E_n\})$) equipped with the $w^\#$-topology, in the sense of
	Definition~\ref{def:tight-rm}.
	\item\label{it:ktw-masses}
	For every $n\in\NN_0$:
	\begin{enumerate}[label=\textup{(\alph*)}]
		\item\label{it:ktw-mass}
		\emph{(control of the masses)} the real random variables $\{\xi_m(E_n)\}_{m\in\NN}$
		are tight in $[0,\infty)$, i.e.
		\[
		\lim_{K\to\infty}\sup_{m\in\NN}\mathbb P_m\bigl(\xi_m(E_n)>K\bigr)=0;
		\]
		\item\label{it:ktw-cc}
		\emph{(compact containment condition)} for all $\varepsilon,\delta>0$ there is a
		compact set $C\subset E$ with
		\[
		\sup_{m\in\NN}\mathbb P_m\bigl(\xi_m(E_n\setminus C)>\delta\bigr)\le\varepsilon .
		\]
	\end{enumerate}
\end{enumerate}
\end{theorem}

\begin{proof}[Proof of Theorem \ref{thm:kallenberg-tight-wsharp}]

	We use the characterisation of relatively compact sets in
	Lemma~\ref{lem:compact-Mn}.

	\emph{\ref{it:ktw-tight}$\Rightarrow$\ref{it:ktw-masses}.} Let
	$\varepsilon>0$, and let $\mathcal K_\varepsilon\subset\MM(\{E_n\})$ be
	compact with $\sup_m\mathbb P_m(\xi_m\notin\mathcal K_\varepsilon)\le\varepsilon$.
	Fix $n\in\NN_0$. By Lemma~\ref{lem:compact-Mn},
	$c_n:=\sup_{\mu\in\mathcal K_\varepsilon}\mu(E_n)<\infty$, and for every
	$\delta>0$ there is a compact set $C\subset E$ with $\mu(E_n\setminus C)<\delta$
	for all $\mu\in\mathcal K_\varepsilon$. Hence, for all $m$,
	\[
	\mathbb P_m\bigl(\xi_m(E_n)>c_n\bigr)\le\mathbb P_m(\xi_m\notin\mathcal K_\varepsilon)\le\varepsilon,
	\qquad
	\mathbb P_m\bigl(\xi_m(E_n\setminus C)>\delta\bigr)\le\mathbb P_m(\xi_m\notin\mathcal K_\varepsilon)\le\varepsilon .
	\]
	Since $\varepsilon>0$ is arbitrary, this is \ref{it:ktw-masses}.

	\emph{\ref{it:ktw-masses}$\Rightarrow$\ref{it:ktw-tight}.} Let
	$\varepsilon>0$. By \ref{it:ktw-masses} we can choose numbers $c_n<\infty$
	and compact sets $C_{n,j}\subset E$, $n\in\NN_0$, $j\in\NN$, with
	\[
	\sup_m\mathbb P_m\bigl(\xi_m(E_n)>c_n\bigr)\le\varepsilon\,2^{-n-2},
	\qquad
	\sup_m\mathbb P_m\bigl(\xi_m(E_n\setminus C_{n,j})>1/j\bigr)\le\varepsilon\,2^{-n-j-2}.
	\]
	Put
	\[
	\mathcal K:=\bigl\{\mu\in\MM(\{E_n\}):\ \mu(E_n)\le c_n\ \text{and}\
	\mu(E_n\setminus C_{n,j})\le 1/j\ \text{for all}\ n\in\NN_0,\ j\in\NN\bigr\}.
	\]
	The set $\mathcal K$ satisfies (K1) and (K2) of Lemma~\ref{lem:compact-Mn}
	(for $\eta>0$ take $j>1/\eta$). Hence $\mathcal K$ is relatively compact, and
	its closure $\overline{\mathcal K}$ is compact. By
	Proposition~\ref{prop:wsharp-props}\,\ref{it:ws-borel}, $\mathcal K$ is a
	Borel set. A union bound gives, for all $m$,
	\[
	\mathbb P_m\bigl(\xi_m\notin\overline{\mathcal K}\bigr)
	\le\mathbb P_m(\xi_m\notin\mathcal K)
	\le\sum_{n\ge0}\varepsilon\,2^{-n-2}+\sum_{n\ge0}\sum_{j\ge1}\varepsilon\,2^{-n-j-2}
	=\frac\varepsilon2+\frac\varepsilon2=\varepsilon .
	\]
	Hence the family is tight. In the integer-valued case we use the compact
	set $\overline{\mathcal K}\cap\MMN(\{E_n\})$ instead; it is compact since
	$\MMN(\{E_n\})$ is closed (Proposition~\ref{prop:wsharp-props}\,\ref{it:ws-polish}),
	and it contains $\xi_m$ whenever $\overline{\mathcal K}$ does.
	\end{proof}

	\begin{remark}
\smallskip\noindent
In the integer-valued case, \ref{it:ktw-cc} with $\delta<1$ reads
$\sup_m\mathbb P_m(\xi_m(E_n\setminus C)>0)\le\varepsilon$, i.e.\ with high
probability all atoms of $\xi_m$ in $E_n$ lie in one fixed compact set.
In particular:
\begin{enumerate}[label=\textup{(\alph*)}]
	\item\label{it:ktw-moment}
	\emph{(First-moment criterion.)} Assertion~\ref{it:ktw-masses} holds as soon
	as, for every $n\in\NN_0$, the mean measures
	$\Lambda_m:=\mathbb E_m\,\xi_m(\cdot)$ satisfy
	\begin{equation}\label{eq:ktw-moment}
		\sup_{m\in\NN}\Lambda_m(E_n)<\infty
		\quad\text{and}\quad
		\forall\varepsilon>0\ \exists\,C\subset E\text{ compact}:\
		\sup_{m\in\NN}\Lambda_m(E_n\setminus C)\le\varepsilon ,
		\tag{$*$}
	\end{equation}
	i.e.\ the restrictions $\{\Lambda_m|_{E_n}\}_m$ are bounded and uniformly
	tight. This follows from Markov's inequality,
	$\mathbb P_m(\xi_m(E_n)>K)\le\Lambda_m(E_n)/K$ and
	$\mathbb P_m(\xi_m(E_n\setminus C)>\delta)\le\Lambda_m(E_n\setminus C)/\delta$.
	\item\label{it:ktw-same-law}
	\emph{(Identical laws.)} If, for every $n$, the law of $\xi_m|_{E_n}$ does
	not depend on $m$ (e.g.\ every $\xi_m$ is a Poisson random measure with the
	same intensity), then \ref{it:ktw-masses} holds: the law of $\xi_m|_{E_n}$ is then a single
	Borel probability measure on the Polish space $\MM(E)$ of finite measures,
	hence tight, and compact subsets of $\MM(E)$ have bounded mass and are
	uniformly tight (Prokhorov's theorem).
	\item\label{it:ktw-single}
	\emph{(Compact localizing sets.)} If every $E_n$ has compact closure, e.g.\
	$E=\{0<|z|\le\epsilon_0\}\times[0,T]\subset\RR^d\times[0,T]$ with
	$E_n=\{\epsilon_{n+1}<|z|\le\epsilon_0\}\times[0,T]$, then
	\ref{it:ktw-cc} holds with $C=\overline{E_n}$, and \ref{it:ktw-masses}
	reduces to the mass condition~\ref{it:ktw-mass}. This is the locally
	compact setting of Kallenberg~\cite[Thm.~4.10]{Kallenberg2017}, cf.\
	Remark~\ref{rem:wsharp-vs-vague}.
\end{enumerate}
Since $\MMN(\{E_n\})$ is a closed subset of $\MM(\{E_n\})$
(Proposition~\ref{prop:wsharp-props}\,\ref{it:ws-polish}), tightness of a
$\MMN(\{E_n\})$-valued family in $\MM(\{E_n\})$ is equivalent to its
tightness in $\MMN(\{E_n\})$.
\end{remark}

\begin{remark}\label{rem:ktw-counterexample}
Without the compact-containment condition~\ref{it:ktw-cc} the theorem fails as
soon as the $E_n$ are not relatively compact, and this is the case in our
setting when $Z$ is infinite-dimensional. Here one can construct a counterexample. Let $Z$ be a separable Hilbert
space, $0^\ast=0$, $(e_m)$ an orthonormal sequence, $E=S$ and $E_n=S_n$, and
fix $r\in(\epsilon_1,\epsilon_0]$. Put $\xi_m:=\delta_{(re_m,T/2)}$. Then
$\xi_m(S_n)\le1$ for all $m,n$, so \ref{it:ktw-mass} holds. The points $re_m$
are $r\sqrt2$-separated, so a compact $C$ contains only finitely many of the
points $(re_m,T/2)$, and $\xi_m(S_0\setminus C)=1$ for all other $m$. Hence
\ref{it:ktw-cc} fails, and the family is not tight. The jumps escape
``sideways'': they stay at fixed distance $r$ from $0^\ast$, never leave
$S_0$, and do not increase in number.
\end{remark}

\begin{remark}\label{rem:tight-summary}
The three levels are:
\begin{itemize}
	\item \emph{Space.} $\MMN(\{S_n\})$ carries the $w^{\#}$-topology, is
	Polish, and its Borel $\sigma$-algebra is generated by the evaluations
	$\mu\mapsto\mu(B)$, $B\in\widehat{\mathcal S}$
	(Proposition~\ref{prop:wsharp-props}).
	\item \emph{Compactness.} A set $\mathcal K\subset\MMN(\{S_n\})$ is
	relatively compact if and only if, for every $n$, the masses $\mu(S_n)$,
	$\mu\in\mathcal K$, are bounded \emph{and} all atoms in $S_n$ of all
	$\mu\in\mathcal K$ lie in one fixed compact set
	(Lemma~\ref{lem:compact-Mn}).
	\item \emph{Tightness.} A family of random measures is tight if and only
	if, for every $n$, the masses $\xi_m(S_n)$ are tight \emph{and}, with
	probability close to one uniformly in $m$, all atoms of $\xi_m$ in $S_n$
	lie in one fixed compact set
	(Theorem~\ref{thm:kallenberg-tight-wsharp}). Bounded masses alone do not
	suffice when $Z$ is infinite-dimensional
	(Remark~\ref{rem:ktw-counterexample}).
\end{itemize}
Thus compactness and tightness on $\MMN(\{S_n\})$ reduce to the scalar
random variables $\xi(S_n)$ and $\xi(S_n\setminus C)$, $C$ compact.
\end{remark}

\subsection[The product space of truncated jump measures]{{The product space of truncated jump measures}}
\label{app:product-space}

In the Jakubowski--Skorokhod step, the approximating solutions are considered jointly with a family of truncated jump measures, where jumps with marks of size $|z|_Z\le\epsilon_{J+1}$ are removed. Since different truncation levels naturally lead to different state spaces, we assemble all levels into a single countable product space.
This is done in the next definition.

Let us go back to the application in Step IX.
 Since every Polish space admits a bounded complete compatible metric,
	we may equip $S_J$ with such a metric and use the constant localizing
	sequence $E_n=S_J$. Then $C_b(\{E_n\})=C_b(S_J)$, and the
	$w^\#$-topology on $\MMN(S_J)$ coincides with the ordinary weak topology
	(Remark~\ref{rem:wsharp-finite-case}).

\begin{definition}[Truncated counting measures and consistent families]
\label{def:product-vague}
Let $Z$ be a Polish space with distinguished point $0^\ast$, and set
$|z|_Z=d_Z(z,0^\ast)$. Let
\[
  \epsilon_0>\epsilon_1>\cdots\downarrow0,
  \qquad
  S_J:=\{z\in Z:\epsilon_{J+1}<|z|_Z\le\epsilon_0\}
         \times[0,T],
  \quad J\in\mathbb N_0.
\]
Then $S_J\subset S_{J+1}$ and
\[
  S:=\bigcup_{J\ge0}S_J
    =\{z\in Z:0<|z|_Z\le\epsilon_0\}\times[0,T].
\]
Each $S_J$, equipped with its subspace topology, is Polish.

\begin{enumerate}
\item
Let $\MMN(S_J)$ denote the space of finite counting
measures on $S_J$, equipped with the topology of weak convergence
of finite measures. This is a Polish space. Choose a complete
compatible metric $d_J$ on $\MMN(S_J)$.

\item

Define
\[
 	\mathbb M:=\prod_{J\in\NN}\MMN(S_J),
\]
equipped with the complete compatible metric
\[
	d_{\mathbb M}\bigl((\mu_J)_J,(\mu'_J)_J\bigr)
		:=\sum_{J=0}^{\infty}2^{-(J+1)}
		\bigl(d_J(\mu_J,\mu'_J)\wedge1\bigr).
\]
Thus convergence in $\mathbb M$ is equivalent to weak
convergence in every coordinate.

\item
The space of consistent families is
\[
  \mathbb M_{\mathrm c}
  :=
  \left\{
    (\mu_J)_{J\ge0}\in\mathbb M:
    \mu_{J'}|_{S_J}=\mu_J
    \text{ whenever }J\le J'
  \right\}.
\]
This is a closed, hence Polish, subspace of $\mathbb M$.

\item
Let $\MMN^{loc}(S,(S_J)_J)$
denote the space of
counting measures $(\mu_J)_J$ on $S$ satisfying $\mu(S_J)<\infty$
for every $J$. Define the assembly map
\[
  \Sigma:\mathbb M_{\mathrm c}
        \longrightarrow\MMN^{loc}(S,(S_J)_J)
\]
by
\[
  \Sigma((\mu_J)_J)(B)
  :=
  \lim_{J\to\infty}\mu_J(B\cap S_J),
  \qquad B\in\mathcal B(S).
\]
Consistency makes this limit monotone. The assembled measure
is the unique counting measure on $S$ satisfying
\[
  \Sigma((\mu_J)_J)|_{S_J}=\mu_J
  \qquad\text{for every }J\ge0.
\]
\end{enumerate}
\end{definition}

\begin{proof}[Closedness of the consistency space]
Fix $J<J'$ and suppose that
\[
  \mu_n\Rightarrow\mu \quad\text{in }\mathcal N_f(S_{J'}),
  \qquad
  \mu_n|_{S_J}\Rightarrow\eta
       \quad\text{in }\mathcal N_f(S_J).
\]
Write $E=S_{J'}$ and $U=S_J$. Since $U$ is relatively open
in $E$, there exist continuous functions
$h_k:E\to[0,1]$ such that
\[
  h_k=0\ \text{on }E\setminus U,
  \qquad h_k\uparrow\mathbf1_U.
\]
For example, for a compatible metric $d$ on $E$, take
\[
  h_k(x)=\min\{1,k\,d(x,E\setminus U)\},
\]
with $h_k\equiv1$ if $U=E$.

For every $f\in C_b(E)$, weak convergence in both spaces gives
\[
  \int_U f h_k\,d\eta
  =\lim_{n\to\infty}\int_U f h_k\,d(\mu_n|_U)
  =\lim_{n\to\infty}\int_E f h_k\,d\mu_n
  =\int_E f h_k\,d\mu.
\]
Letting $k\to\infty$ and using dominated convergence yields
\[
  \int_U f\,d\eta=\int_U f\,d\mu
  \qquad\text{for every }f\in C_b(E).
\]
Regarding $\eta$ as a finite measure on $E$ by extension
by zero, uniqueness of finite Borel measures therefore gives
$\eta=\mu|_U$.

Thus the graph of each restriction map is closed.
Taking the intersection of the corresponding consistency
conditions over all $J<J'$ proves that
$\mathbb M_{\mathrm c}$ is closed in $\mathbb M$.
\end{proof}

\begin{remark}\label{rem:open-closed-truncation}
The sets $S_J$ are the state spaces of Step~(IX), open at the inner sphere,
in accordance with the layers $\Delta Z_j=\{\epsilon_{j+1}<|z|_Z\le\epsilon_j\}$.
Since $\epsilon_{J+1}\in\mathcal R(\nu)$, a Poisson random measure with
intensity $\nu\otimes\mathrm{Leb}$ has a.s.\ no atom on $\partial_J$; this is
what makes the restriction maps continuous at the limits, see
Proposition~\ref{prop:tight-product}\,\ref{it:tp-cons}.
\end{remark}

\begin{proposition}[Tightness in the product space]\label{prop:tight-product}
Let $S_J=Z^c_J\times[0,T]$, $\partial_J$, the product space
$(\mathbb M,d_{\mathbb M})$ and the Borel set $\mathbb M_{\mathrm c}\subset\mathbb M$
of consistent families be as in Definition~\ref{def:product-vague}; recall
that $\mathbb M$ is Polish, with the topology of coordinatewise weak
convergence, and that on each $\MMN(S_J)$ the $w^\#$-topology coincides with
the weak topology. For every $m\in\mathbb N$ let
$\Xi_m=(\Xi_m^J)_{J\in\mathbb N}$ be an $\mathbb M$-valued random element on
$(\Omega_m,\mathcal F_m,\mathbb P_m)$, and let
$\Lambda_m^J:=\mathbb E_m\,\Xi_m^J(\cdot)$ denote the mean measure of
$\Xi_m^J$ on $S_J$.
\begin{enumerate}[label=(\roman*)]
	\item\label{it:tp-joint}
	The family $\{\Xi_m\}_{m\in\mathbb N}$ is tight in $(\mathbb M,d_{\mathbb M})$
	\item[] if and only if
	\item\label{it:tp-coord}
	for every $J\in\mathbb N$ the family $\{\Xi_m^J\}_{m\in\mathbb N}$ is tight
	in $\MMN(S_J)$.
\end{enumerate}
Moreover:
\begin{enumerate}[label=(\roman*),resume]
	\item\label{it:tp-moment}
	\emph{(First-moment criterion.)} \ref{it:tp-coord}, and hence
	\ref{it:tp-joint}, holds as soon as, for every $J\in\mathbb N$, the mean
	measures $\{\Lambda_m^J\}_m$ are bounded and uniformly tight:
	\begin{equation}\label{eq:tp-moment}
		\sup_{m\in\mathbb N}\Lambda_m^J(S_J)<\infty
		\quad\text{and}\quad
		\forall\varepsilon>0\ \exists\,C\subset S_J\ \text{compact}:\
		\sup_{m\in\mathbb N}\Lambda_m^J(S_J\setminus C)\le\varepsilon .
		\tag{$**$}
	\end{equation}
	In particular, \eqref{eq:tp-moment} holds if $\Lambda_m^J$ does not depend
	on $m$, e.g.\ $\Lambda_m^J=(\nu\otimes\mathrm{Leb})|_{S_J}$ for all $m$; and
	\ref{it:tp-coord} holds if the law of $\Xi_m^J$ does not depend on $m$.
	\item\label{it:tp-cons}
	\emph{(Consistency of limits.)} Suppose $\Xi_m\in\mathbb M_{\mathrm c}$
	$\mathbb P_m$-a.s.\ for all $m$, and $\operatorname{Law}(\Xi_{m_k})\to Q$
	weakly on $\mathbb M$ along a subsequence. If
	\[
	Q\bigl(\mu:\ \mu_{J'}(\partial_J)>0\bigr)=0\qquad\text{for all}\ J<J',
	\]
	then $Q(\mathbb M_{\mathrm c})=1$.
\end{enumerate}
\end{proposition}

\begin{proof}
\ref{it:tp-joint}$\Rightarrow$\ref{it:tp-coord}: the coordinate projection
$\pi_J:\mathbb M\to \MMN(S_J)$ is continuous, and continuous
images of tight families are tight.

\ref{it:tp-coord}$\Rightarrow$\ref{it:tp-joint}: let $\varepsilon>0$. By
\ref{it:tp-coord}, choose for each $J$ a compact
$K_J\subset \MMN(S_J)$ with
$\sup_m\mathbb P_m(\Xi_m^J\notin K_J)\le\varepsilon2^{-(J+1)}$. By Tychonoff's
theorem, $K:=\prod_JK_J$ is compact in the product topology, which is the
topology of $d_{\mathbb M}$, and
\[
\sup_m\mathbb P_m(\Xi_m\notin K)
\le\sum_{J=0}^\infty\sup_m\mathbb P_m(\Xi_m^J\notin K_J)\le\varepsilon .
\]

\ref{it:tp-moment}: fix $J$. Markov's inequality gives
$\mathbb P_m(\Xi_m^J(S_J)>K)\le\Lambda_m^J(S_J)/K$ and, as $\Xi_m^J$ is
$\mathbb N_0$-valued,
\[
\mathbb P_m\bigl(\Xi_m^J(S_J\setminus C)>0\bigr)
=\mathbb P_m\bigl(\Xi_m^J(S_J\setminus C)\ge1\bigr)
\le\Lambda_m^J(S_J\setminus C).
\]
Hence \eqref{eq:tp-moment} implies assertion~\ref{it:ktw-masses} of
Theorem~\ref{thm:kallenberg-tight-wsharp}, applied with $E=S_J$ and the
constant localizing sequence $E_n=S_J$ (admissible since $S_J$ is bounded),
so $\{\Xi^J_m\}_m$ is tight.
If
$\Lambda_m^J=\Lambda^J$ does not depend on $m$, then $\Lambda^J$ is a finite
Borel measure on the Polish space $S_J$, hence tight by Ulam's theorem, and
\eqref{eq:tp-moment} holds. If the law of $\Xi_m^J$ does not depend on $m$,
it is a single Borel probability measure on a Polish space, hence tight.

\ref{it:tp-cons}: for $J<J'$ let
$\rho_{J,J'}: \MMN(S_{J'})\to \MMN(S_J)$,
$\mu\mapsto\mu|_{S_J}$. We claim that $\rho_{J,J'}$ is continuous at every
$\mu$ with $\mu(\partial_J)=0$. Indeed, let $\mu_n\to\mu$ weakly on $S_{J'}$
and $f\in C_b(S_J)$, and let $F:=f$ on $S_J$ and $F:=0$ on
$S_{J'}\setminus S_J$. Since $S_J$ is relatively open in $S_{J'}$ and
$S_{J'}\setminus(S_J\cup\partial_J)$ is relatively open as well, $F$ is
bounded and continuous at every point of $S_{J'}\setminus\partial_J$, and
$\partial_J$ is a $\mu$-null set.
{By the extended Portmanteau theorem (equivalently, convergence
	of integrals for bounded Borel functions whose set of discontinuities is
	$\mu$-null),}
$\mu_n|_{S_J}(f)=\mu_n(F)\to\mu(F)=\mu|_{S_J}(f)$.

Let $D:=\{\mu\in\mathbb M:\mu_{J'}(\partial_J)=0\ \forall J<J'\}$, a Borel set
with $Q(D)=1$ by assumption. If $\mu\in D$ lies in the closure
$\overline{\mathbb M_{\mathrm c}}$, pick $\mu^n\in\mathbb M_{\mathrm c}$ with $\mu^n\to\mu$. Then
$\mu^n_J=\rho_{J,J'}(\mu^n_{J'})\to\rho_{J,J'}(\mu_{J'})$ by the claim, while
$\mu^n_J\to\mu_J$, so $\mu_J=\mu_{J'}|_{S_J}$ and $\mu\in\mathbb M_{\mathrm c}$. Thus
$\overline{\mathbb M_{\mathrm c}}\cap D\subset\mathbb M_{\mathrm c}$. By the portmanteau theorem,
$Q(\overline{\mathbb M_{\mathrm c}})\ge\limsup_k\mathbb P_{m_k}(\Xi_{m_k}\in\overline{\mathbb M_{\mathrm c}})=1$.
Hence $Q(\mathbb M_{\mathrm c})\ge Q(\overline{\mathbb M_{\mathrm c}}\cap D)=1$.
\end{proof}

\begin{remark}\label{rem:tp-mass-not-enough}
The bound $\sup_m\Lambda_m^J(S_J)<\infty$ alone does not suffice when $Z$ is
infinite-dimensional: it controls the number of atoms of $\Xi_m^J$ but not
their location. For instance, let $Z$ be a separable Hilbert space, $0^\ast=0$,
$(e_m)$ orthonormal, and let $\Xi_m^J$ be a Poisson random
measure on $S_J$ with intensity ${c}\,\delta_{re_m}\otimes\mathrm{Leb}$, where {$c>0$ and} $\epsilon_{J+1}<r\le\epsilon_0$. Then
$\Lambda_m^J(S_J)={c}T$ for all $m$. Since the points $re_m$ are
$r\sqrt2$-separated, a compact $C$ contains only finitely many of them, and
for all other $m$ one has
$\mathbb P_m(\Xi_m^J(S_J\setminus C)>0)\ge1-e^{-{c}T}$. So the compact-containment condition
\ref{it:ktw-cc} of Theorem~\ref{thm:kallenberg-tight-wsharp} applied with
$E=S_J$, fails, and the family is not tight. If
$Z=\mathbb R^d$, the sets $C_a:=\{\epsilon_{J+1}+a\le|z|\le\epsilon_0\}\times[0,T]$,
$a>0$, are compact subsets of $S_J$, and the second condition in
\eqref{eq:tp-moment} reduces to
$\lim_{a\downarrow0}\sup_m\Lambda_m^J(S_J\setminus C_a)=0$, i.e.\ no mass may
accumulate at the inner sphere $\partial_J$.
\end{remark}

\section{The shifted projection on the Haar basis}
\label{haarproj}
\newcommand{\pro}{P}
Let $\phi:[0,1]\to\RR$ be the Haar wavelet given by
$$
\phi(x):=\bcase 1 &\mbox{ if } 0\le x<\frac{1}{2},
\\-1 &\mbox{ if } \frac{1}{2}\le x <1,
\ecase
$$
Now, define
$\phi_{n,j}(x):=2^n \phi(2^n x-j)$, $x\in[0,1]$, the corresponding multi resolution analysis. Let $X$ be a Banach space. For $f\in L^m(0,T;X)$, $\pro_n f=\sum_{j=0}^{2^n-1}\la f,\phi_{n,j}\ra \, \phi_{n,j}$ its orthogonal projection on
the Haar basis given by  $\{\phi_{n,j}:j=0,\ldots, 2^n-1\}$.

In a second step, we shift the time intervals. In particular, we define for all
$\kappa\in \NN$ and
for $f\in L^m(0,T;X)$    the following shifted  function  $\shiftpro _\zaehler f$ of $f$
$$
\shiftpro _\zaehler f(s)=f(s-\frac{T}{2^n}),\ \forall s\in [0,T]
$$
Let
\DEQSZ\label{projectiondef}
\Pro_\zaehler := \shiftpro _\zaehler\, \pro_\zaehler.
\EEQSZ

\begin{remark}
	It is straightforward to verify that for any $f\in L^m(0,T;X)$, where $X$  is a Banach space, $\Pro_\zaehler f$ can be also written as follows
	\begin{equation*}
		\Pro_\zaehler f(s)=
		\frac{2^n}{T}\sum_{k=0}^{2^n-1}1_{[t_k,t_{k+1}]}(s)\int_{t_k}^{t_{k+1}}f(r)\, dr
	\end{equation*}
	where $t^\zaehler_k=T\frac{k}{2^n}$.
\end{remark}

\paragraph{The shifted dyadic projection with initial value.}
Let $E$ be a separable Banach space, $m\in\NN$, $\tau_m=T2^{-m}$ and
$I^m_k=(t^m_k,t^m_{k+1}]$, $k=0,\dots,2^m-1$. For $v\in L^1(0,T;E)$ and
$v^m_0\in E$ we put
\[
\Pro_mv:=\sum_{k=0}^{2^m-1}\xi_k\,\mathbf 1_{I^m_k},\qquad
\xi_0:=v^m_0,\qquad
\xi_k:=\frac1{\tau_m}\int_{I^m_{k-1}}v(s)\,ds,\quad k\ge1 .
\]
On $(t^m_1,T]$ this is the shifted projection \eqref{projectiondef}; on the
first cell it takes the value $v^m_0$. For $E=V$, $v=u_n$ and $v^m_0=u^m_0$
it is the process $\myp_{n,m}$ of Step~(I), see
\eqref{def_shifted_dyadic}--\eqref{init0}. The proof uses only the following
properties.

\begin{lemma}[Elementary properties]\label{lem:Pro-props}
	Let $v\in L^1(0,T;E)$, $v^m_0\in E$ and $p\ge1$.
	\begin{enumerate}[label=\upshape(\alph*)]
		\item\label{Pro:pred} If $E=V$, $v=u_n$ and $v^m_0=u^m_0$, then
		$\xi_k$ is $\mathscr F_{t^m_k}$-measurable for every $k$, and
		$\Pro_mu_n$ is a simple predictable process.
		\item\label{Pro:jensen} For $k\ge1$,
		$\tau_m|\xi_k|_E^p\le\int_{I^m_{k-1}}|v(s)|_E^p\,ds$; in particular
		$\int_{t^m_1}^T|\Pro_mv(s)|_E^p\,ds\le\int_0^T|v(s)|_E^p\,ds$.
		\item\label{Pro:H} If $E=V$ and $\sup_{s\in[0,T]}|v(s)|_H<\infty$, then
		$|\xi_k|_H\le\sup_{s\in[0,T]}|v(s)|_H$ for $k\ge1$.
	\end{enumerate}
\end{lemma}
\begin{proof}
	(a) For $k\ge1$, $\xi_k$ depends only on $u_n$ on
	$I^m_{k-1}\subset[0,t^m_k]$. Since $u_n$ is progressively measurable,
	$\xi_k$ is $\mathscr F_{t^m_k}$-measurable; $\xi_0=u^m_0$ is deterministic.
	On $I^m_k=(t^m_k,t^m_{k+1}]$ the process $\Pro_mu_n$ equals the
	$\mathscr F_{t^m_k}$-measurable value $\xi_k$. So it is a left-continuous
	step process whose values are known at the left end point of each cell,
	hence simple and predictable.

	(b) By $|\int_{I^m_{k-1}}v\,ds|_E\le\int_{I^m_{k-1}}|v|_E\,ds$ and Jensen's
	inequality for the convex function $x\mapsto x^p$ and the probability
	measure $\tau_m^{-1}ds$ on $I^m_{k-1}$,
	\[
	|\xi_k|_E^p\le\Bigl(\frac1{\tau_m}\int_{I^m_{k-1}}|v(s)|_E\,ds\Bigr)^p
	\le\frac1{\tau_m}\int_{I^m_{k-1}}|v(s)|_E^p\,ds .
	\]
	Multiplying by $\tau_m$ and summing over $k\ge1$ gives the second claim,
	since the cells $I^m_{k-1}$, $k\ge1$, are disjoint subsets of $(0,T]$.

	(c) $|\xi_k|_H\le\tau_m^{-1}\int_{I^m_{k-1}}|v(s)|_H\,ds\le\sup_{s}|v(s)|_H$.
\end{proof}

\begin{lemma}[Convergence of the shifted dyadic projection in $L^1$]
	\label{lem:Pro-L1}
	Let $E$ be a separable Banach space, $v\in L^1(0,T;E)$ and let $(v^m_0)_{m\in\NN}\subset E$ satisfy
	$\tau_m|v^m_0|_{E}\to0$. Then
	\[
	\bigl\|\Pro_mv-v\bigr\|_{L^1(0,T;E)}\xrightarrow[m\to\infty]{}0 .
	\]
	No regularity of $v$ beyond integrability is used; in particular the
	statement holds for every $v\in\DD([0,T];H)$.
\end{lemma}

\begin{proof}
	Extend $v$ by $0$ outside $[0,T]$ and write $S_hv:=v(\cdot-h)$ and
	$E_mv:=\EE[v\mid\mathcal D_m]$, where $\mathcal D_m$ is the
	$\sigma$-algebra on $(0,T]$ generated by $\{I^m_k\}_{k}$; thus
	$(E_mv)(s)=\la v\ra_{I^m_k}$ for $s\in I^m_k$. By construction
	$\Pro_mv=S_{\tau_m}E_mv$ on $(t^m_1,T]$, while
	$\Pro_mv\equiv v^m_0$ on $I^m_0$. Hence
	\[
	\bigl\|\Pro_mv-v\bigr\|_{L^1(0,T;E)}
	\;\le\;
	\underbrace{\bigl\|S_{\tau_m}E_mv-S_{\tau_m}v\bigr\|_{L^1}}_{=:a_m}
	+\underbrace{\bigl\|S_{\tau_m}v-v\bigr\|_{L^1}}_{=:b_m}
	+\underbrace{\tau_m|v^m_0|_{E}+\int_{I^m_0}|v(s)|_{E}\,ds}_{=:c_m}.
	\]
	Now $a_m\le\|E_mv-v\|_{L^1}\to0$, because $\mathcal D_m\uparrow
	\mathscr B((0,T])$ and $E_mv\to v$ in $L^1$ by the martingale
	convergence theorem (for $E$-valued $v$ this follows from the scalar case by approximating $v$ in $L^1(0,T;E)$ by simple functions, since $E_m$ is a contraction on $L^1(0,T;E)$) (equivalently, by the Lebesgue differentiation
	theorem for the dyadic basis); $b_m\to0$ by the continuity of
	translation in $L^1(\RR;E)$; and $c_m\to0$ by hypothesis together with
	the absolute continuity of the integral, $|I^m_0|=\tau_m\to0$.
\end{proof}
\emph{Application to $u_n$.} Fix $n$. Since $u_n(\omega)\in L^1(0,T;V)$
for $\PP$-a.e.\ $\omega$ and
$\tau_m|u^m_0|_V{\to0}$ by
\eqref{init0}, Lemma~\ref{lem:Pro-L1} with $E=V$ gives
$\|\myp_{n,m}-u_n\|_{L^1(0,T;V)}\to0$ as $m\to\infty$, $\PP$-a.s. No rate in
$m$ and no time regularity of $u_n$ is needed. Convergence in
$L^\alpha(0,T;V)$ would in addition require
$\tau_m|u^m_0|_V^\alpha\to0$; it is not used in
the proof.

\begin{remark}\label{rem:Pro-no-rate}
	Lemma~\ref{lem:Pro-L1} is \emph{rate-free}, and necessarily so: the
	modulus of $\|S_hv-v\|_{L^1}$ is not controlled by $\|v\|_{L^1}$, and on
	$\DD([0,T];H)$ there is no modulus valid for all paths. A rate would
	require an increment bound for $v$. The proof does not need one, since
	$n$ is fixed before $m$.
\end{remark}

\section{$L^p$-spaces and the Haar system}\label{app:haar-system}\setcounter{equation}{0} 

In this subsection, we assume that $(Y,|\cdot|)$ is a separable Banach space. We assume for simplicity that $T=1$.
We recall some facts about the approximation properties of the Haar system and its shifted version. For the proofs and other details, we are referring to \cite[Appendix B]{FKEHMH}.

\subsection{Definition of the shifted Haar projection}

For each $n \in \mathbb{N}$, let $\pi_n = \{ s^n_0 = 0 < s^n_1 < \cdots < s^n_{2^n} = 1\}$ denote the dyadic partition of $[0,1]$ given by
\[
s^n_j = j 2^{-n}, \quad j = 0, \dots, 2^n.
\]
Each interval $(s^n_{j-1}, s^n_j]$ is called a \textit{dyadic interval}.

The classical (non-shifted) Haar projection $\mathfrak{h}_n: L^p([0,1];Y) \to L^p([0,1];Y)$ is defined by
\[
\mathfrak{h}_n(x) = \sum_{j=1}^{2^n} \mathbf{1}_{(s^n_{j-1}, s^n_j]} \otimes \iota_{j,n}(x),
\]
where
\[
\iota_{j,n}(x) := \frac{1}{2^{-n}} \int_{s^n_{j-1}}^{s^n_j} x(r) \, dr
\]
denotes the mean value of \(x\) over $(s^n_{j-1}, s^n_j]$.

The \textit{shifted Haar projection} $\mathfrak{h}_n^s : L^p([0,1];Y) \to L^p([0,1];Y)$ is defined by
\[
\mathfrak{h}_n^s(x) = \sum_{j=1}^{2^n-1} \mathbf{1}_{(s^n_j, s^n_{j+1}]} \otimes \iota_{j-1,n}(x),
\]
where $\iota_{0,n}(x) := 0$ for convention.

\begin{proposition}[Convergence rate in fractional Sobolev spaces]\label{prop:haar-conv}
Let $Y$ be a Banach space, $1 < p < \infty$, and let $x \in W^{\alpha,p}([0,1];Y)$ for some $\alpha \in (0,1]$. Then there exists a constant $C>0$, depending only on $p$ and $\alpha$, such that
\[
\| \mathfrak{h}_n^s(x) - x \|_{L^p([0,1];Y)} \leq C 2^{-n\alpha} \|x\|_{W^{\alpha,p}([0,1];Y)},
\]
where
\[
\|x\|_{W^{\alpha,p}}^p = \int_0^1 \int_0^1 \frac{ | x(t) - x(s) |_Y^p }{ |t-s|^{1+\alpha p} } \, ds \, dt
\]
denotes the Gagliardo seminorm.
\end{proposition}

\begin{proof}
	For the proof we are referring to \cite[Lemma B.2]{FKEHMH}
\end{proof}

\begin{proposition}[Uniform convergence on compact sets]\label{prop:haar-uniform}\label{prop:uhaar}
Let $\mathbb{X}$ be a Banach space compactly embedded into $L^p([0,1];Y)$, $1<p<\infty$. Let $\mathcal{X} \subset \mathbb{X}$ be a bounded set. Then for every $\varepsilon>0$, there exists $N \in \mathbb{N}$ such that for all $n \geq N$,
\[
\sup_{x \in \mathcal{X}} \| \mathfrak{h}_n^s(x) - x \|_{L^p([0,1];Y)} \leq \varepsilon.
\]
\end{proposition}

\begin{proof}
Since $\mathbb{X}$ embeds compactly into $L^p([0,1];Y)$ and $\mathcal{X}$ is bounded in $\mathbb{X}$, the closure of $\mathcal{X}$ in $L^p$ is compact.

For each $x \in \mathcal{X}$, we have $\mathfrak{h}_n^s(x) \to x$ strongly in $L^p([0,1];Y)$. 

Since strong convergence on a compact set implies uniform convergence, for any $\varepsilon>0$ there exists $N \in \mathbb{N}$ such that for all $n \geq N$,
\[
\sup_{x \in \mathcal{X}} \| \mathfrak{h}_n^s(x) - x \|_{L^p([0,1];Y)} \leq \varepsilon.
\]
\end{proof}

\section{The Skorohod space}\label{app:Skorohod-space}
For an introduction to the Skorohod space of real-valued paths we refer to
Billingsley \cite{Billingsley1999} and Jacod--Shiryaev \cite{JacodShiryaev2003};
for paths with values in a general Polish space, which is the situation of
this manuscript, to Ethier--Kurtz \cite{Ethier1986}. We state only the results
consumed in the main text.

Let $(E,d_E)$ be a complete separable metric space. $\DD([0,T];E)$ denotes
the space of right continuous functions $x:[0,T]\to E$ with left limits
(c\`adl\`ag paths). Define the admissible time changes of $[0,T]$,
\[
\Lambda_T:=\bigl\{\lambda:[0,T]\to[0,T]:\ \lambda\ \text{continuous, strictly
	increasing},\ \lambda(0)=0,\ \lambda(T)=T\bigr\};
\]
note that every $\lambda\in\Lambda_T$ \emph{fixes both endpoints}.

A sequence $\{x^{(n)}\}_{n\in\NN}\subset\DD([0,T];E)$ converges to
$x\in\DD([0,T];E)$ in the Skorohod ($J_1$) topology iff there are
$\{\lambda_n\}_{n\in\NN}\subset\Lambda_T$ such that
\begin{itemize}
	\item[(i)] $\sup_{0\le s\le T}|\lambda_n(s)-s|\to0$ as $n\to\infty$;
	\item[(ii)] $\sup_{0\le s\le T}
	d_E\bigl(x^{(n)}(\lambda_n(s)),x(s)\bigr)\to0$ as $n\to\infty$.
\end{itemize}

The topology is metrised by \emph{Billingsley's metric}\footnote{Not to be
	confused with the Prokhorov metric, which metrises weak convergence of
	measures and is used, as $d_{\mathrm v}$, on $\mathbb M(\ZZJ )$ in
	Steps~(VIII)--(IX).} with
\[
\|\lambda\|_{\log}:=\sup_{0\le s<t\le T}
\Bigl|\log\frac{\lambda(t)-\lambda(s)}{t-s}\Bigr|,
\qquad
\Lambda_{\log}:=\{\lambda\in\Lambda_T:\|\lambda\|_{\log}<\infty\},
\]
set
\begin{equation}\label{eq:billingsley-metric}
	d^\circ(x,y):=\inf_{\lambda\in\Lambda_{\log}}
	\Bigl(\|\lambda\|_{\log}
	+\sup_{0\le s\le T}d_E\bigl(x(\lambda(s)),y(s)\bigr)\Bigr),
	\qquad x,y\in\DD([0,T];E).
\end{equation}
The metric $d^\circ$ generates the $J_1$ topology, and, $E$ being complete
and separable, $(\DD([0,T];E),d^\circ)$ is a \emph{complete separable}
metric space, hence Polish; see \cite[\S12]{Billingsley1999},
\cite[Ch.~3, \S5]{Ethier1986},
\cite[Lemma~VI.1.30, Cor.~VI.1.32--1.33]{JacodShiryaev2003}. The naive
Skorohod metric, with $\|\lambda-\mathrm{id}\|_\infty$ in place of
$\|\lambda\|_{\log}$, generates the same topology but is \emph{not}
complete; completeness is what Lemma~\ref{lem:state-space}\,\ref{SS:i}
needs.

\begin{lemma}[Pointwise evaluations under $J_1$]\label{lem:J1-pointwise}
	Let $x^{(n)}\to x$ in $\DD([0,T];E)$ with time changes $\lambda_n$ as in
	\textup{(i)--(ii)}. Then:
	\begin{enumerate}[label=\upshape(\alph*)]
		\item $x^{(n)}(t)\to x(t)$ at every $t$ at which $x$ is continuous; the
		complement is at most countable.
		\item $x^{(n)}(0)\to x(0)$ and $x^{(n)}(T)\to x(T)$, \emph{without} any
		continuity assumption, since every $\lambda_n$ fixes the endpoints.
		\item With $y_n:=x^{(n)}\circ\lambda_n$ and
		$\varepsilon_n:=\sup_sd_E(y_n(s),x(s))$, one has
		$\sup_{t}d_E\bigl(y_n(t-),x(t-)\bigr)\le\varepsilon_n$; in a normed
		space $E$ also
		$\sup_t|\Delta y_n(t)-\Delta x(t)|_E\le2\varepsilon_n$, i.e.\ left
		limits and jump sizes converge uniformly along the time changes.
		\item \emph{(Bounded range.)} For every $e\in E$,
		$\sup_{t}d_E\bigl(x^{(n)}(t),e\bigr)\le\sup_{t}d_E\bigl(x(t),e\bigr)+\varepsilon_n$,
		with $\varepsilon_n$ as in~(c), since $\lambda_n$ is a bijection of
		$[0,T]$ and hence $\sup_td_E(x^{(n)}(t),e)=\sup_sd_E(y_n(s),e)$. In
		particular, a $J_1$-convergent sequence is bounded in the supremum norm;
		in a normed space $E$,
		$\sup_t|x^{(n)}(t)|_E\le\sup_t|x(t)|_E+\varepsilon_n$.
	\end{enumerate}
\end{lemma}

\begin{lemma}[Measurable jump structure]\label{lem:jump-measurable}
	For $x\in\DD([0,T];E)$ and $\varepsilon>0$, the number of jumps of size
	$>\varepsilon$, $N_\varepsilon(x):=\#\{s\in(0,T]:
	d_E(x(s),x(s-))>\varepsilon\}$, is finite, and
	$x\mapsto N_\varepsilon(x)$ is Borel on $(\DD([0,T];E),d^\circ)$. Ordering
	these jump times increasingly and letting $\varepsilon\downarrow0$ along a
	fixed sequence yields a measurable enumeration $(s_i(x))_{i\in\NN}$ of all
	jump times, such that the maps
	$x\mapsto s_i(x)\in[0,T]\cup\{\dagger\}$ and
	$x\mapsto\bigl(x(s_i(x)-),\,x(s_i(x))\bigr)\in E\times E$
	are Borel \textup(with a cemetery value $\dagger$ when fewer than $i$ jumps
	exist\textup). Consequently, functionals of the form
	$\sum_i\varphi\bigl(x(s_i-),\Delta x(s_i),s_i\bigr)$ with Borel
	$\varphi\ge0$ are Borel on $\DD([0,T];E)$.
\end{lemma}

\noindent
\textbf{Acknowledgement.} E.H. gratefully acknowledges the support of the Austrian Science Foundation\href{https://www.fwf.ac.at/en/research-radar/10.55776/PAT3076923}{10.55776/PAT3076923}.  A.K. acknowledges  funding by the Austrian Science Fund (FWF)\break  \href{https://www.fwf.ac.at/forschungsradar/10.55776/ESP4373225}{10.55776/ESP4373225}. A.A.P would like to thank Austrian Science Fund (FWF) for supporting a research visit to Montanuniversit\"at Leoben under the ESPRIT Programme \href{https://www.fwf.ac.at/forschungsradar/10.55776/ESP4373225}{10.55776/ESP4373225}.


\noindent	\textbf{Declarations:}

\noindent 	\textbf{Ethical Approval:}   Not applicable.

\noindent  \textbf{Conflict of interest statement: } On behalf of all authors, the corresponding author states that there is no conflict of interest.


\noindent 	\textbf{Authors' contributions: } All authors have contributed equally.

\noindent 	\textbf{Funding:}  Austrian Science Fund (FWF) \href{https://www.fwf.ac.at/en/research-radar/10.55776/PAT3076923}{10.55776/PAT3076923} (E.H.).\\
Austrian Science Fund (FWF) \href{https://www.fwf.ac.at/forschungsradar/10.55776/ESP4373225}{10.55776/ESP4373225} (A.K.).

\noindent 	\textbf{Availability of data and materials: } Not applicable.

		\end{document}